\documentclass[12pt]{article}
\usepackage[leqno]{amsmath}
\usepackage{amssymb,mathrsfs,bm}
\usepackage{hyperref}
\usepackage{mathtools}
\usepackage{fullpage}
\usepackage{epic}
\usepackage{eepic}
\usepackage[english]{babel}
\usepackage[all]{xy}
\usepackage{enumerate}
\usepackage[T1]{fontenc}
\usepackage[latin1]{inputenc}
\usepackage{color}
\usepackage{soul}
\usepackage{mathdots}
\usepackage{mathabx}
\usepackage{dsfont}
\newtheorem{theointro}{Theorem}
\newtheorem{remintro}{Remark}
\newtheorem{theo}{Theorem}[subsection]
\newtheorem{lem}[theo]{Lemma}
\newtheorem{lemdef}[theo]{Lemma-Definition}
\newtheorem{prop}[theo]{Proposition}
\newtheorem{cor}[theo]{Corollary}
\newtheorem{rem}[theo]{Remark}

\def\Q{\mathbb{Q}}
\def\oF{\overline{F}}
\def\R{\mathbb{R}}
\def\C{\mathbb{C}}

\def\ctens{\widehat{\otimes}}

\DeclareMathOperator{\Hom}{Hom}
\DeclareMathOperator{\Tr}{Trace}
\DeclareMathOperator{\Tra}{Tr}

\DeclareMathOperator{\Temp}{Temp}

\DeclareMathOperator{\Ker}{Ker}

\DeclareMathOperator{\Ad}{Ad}
\DeclareMathOperator{\End}{End}
\DeclareMathOperator{\Lie}{Lie}

\DeclareMathOperator{\vol}{vol}

\DeclareMathOperator{\Irr}{Irr}

\DeclareMathOperator{\diag}{diag}

\DeclareMathOperator{\cS}{\mathcal{S}}
\DeclareMathOperator{\cF}{\mathcal{F}}
\DeclareMathOperator{\fS}{\mathfrak{S}}
\DeclareMathOperator{\cA}{\mathcal{A}}
\DeclareMathOperator{\bR}{\mathbb{R}}
\DeclareMathOperator{\bZ}{\mathbb{Z}}
\DeclareMathOperator{\bN}{\mathbb{N}}
\DeclareMathOperator{\bA}{\mathbb{A}}
\DeclareMathOperator{\cH}{\mathcal{H}}
\DeclareMathOperator{\ux}{\underline{x}}
\DeclareMathOperator{\uy}{\underline{y}}
\DeclareMathOperator{\uz}{\underline{z}}
\DeclareMathOperator{\gl}{\mathfrak{gl}}
\DeclareMathOperator{\g}{\mathfrak{g}}
\DeclareMathOperator{\h}{\mathfrak{h}}
\DeclareMathOperator{\n}{\mathfrak{n}}
\DeclareMathOperator{\fb}{\mathfrak{b}}
\DeclareMathOperator{\fg}{\mathfrak{g}}
\DeclarePairedDelimiter\ceil{\lceil}{\rceil}
\DeclarePairedDelimiter\floor{\lfloor}{\rfloor}
\DeclareMathOperator{\cC}{\mathcal{C}}
\DeclareMathOperator{\fc}{\mathfrak{c}}
\DeclareMathOperator{\ind}{ind}
\DeclareMathOperator{\qs}{qs}
\DeclareMathOperator{\stab}{stab}
\DeclareMathOperator{\As}{As}
\DeclareMathOperator{\cU}{\mathcal{U}}
\DeclareMathOperator{\cV}{\mathcal{V}}
\DeclareMathOperator{\cO}{\mathcal{O}}
\DeclareMathOperator{\cZ}{\mathcal{Z}}
\DeclareMathOperator{\cB}{\mathcal{B}}
\DeclareMathOperator{\cW}{\mathcal{W}}
\DeclareMathOperator{\Bil}{Bil}
\DeclareMathOperator{\cont}{cont}
\DeclareMathOperator{\Sym}{Sym}
\DeclareMathOperator{\Ind}{Ind}
\DeclareMathOperator{\ab}{ab}
\DeclareMathOperator{\temp}{temp}

\DeclareMathOperator{\Whitt}{Whitt}
\DeclareMathOperator{\Id}{Id}
\DeclareMathOperator{\Norm}{Norm}
\DeclareMathOperator{\sgn}{sgn}
\DeclareMathOperator{\unit}{unit}
\DeclareMathOperator{\fs}{\mathfrak{s}}
\DeclareMathOperator{\fu}{\mathfrak{u}}
\DeclareMathOperator{\rs}{rs}
\DeclareMathOperator{\trans}{trans}
\DeclareMathOperator{\disc}{disc}
\DeclareMathOperator{\cM}{\mathcal{M}}

\DeclareMathOperator{\GL}{\mathrm{GL}}

\numberwithin{equation}{subsection}
\newcounter{keepeqno}
\newenvironment{num}
 {\setcounter{keepeqno}{\value{equation}}%
  \begin{list}{(\theequation)}{\usecounter{equation}}%
  \setcounter{equation}{\value{keepeqno}}}
 {\end{list}}
 
\newcommand{\eq}[1][r]
   {\ar@<-3pt>@{-}[#1]
    \ar@<-1pt>@{}[#1]|<{}="gauche"
    \ar@<+0pt>@{}[#1]|-{}="milieu"
    \ar@<+1pt>@{}[#1]|>{}="droite"
    \ar@/^2pt/@{-}"gauche";"milieu"
    \ar@/_2pt/@{-}"milieu";"droite"}

\title{Plancherel formula for $\GL_n(F)\backslash \GL_n(E)$ and applications to the Ichino-Ikeda and formal degree conjectures for unitary groups}
\author{Rapha\"el Beuzart-Plessis}

\begin{document}

\maketitle

\begin{abstract}
We establish an explicit Plancherel decomposition for $\GL_n(F)\backslash \GL_n(E)$ where $E/F$ is a quadratic extension of local fields of characteristic zero by making use of a local functional equation for Asai $\gamma$-factors. We also give two applications of this Plancherel formula: first to the global Ichino-Ikeda conjecture for unitary groups by completing a comparison between local relative characters that was left open by W. Zhang \cite{Zh3} and secondly to the Hiraga-Ichino-Ikeda conjecture on formal degrees \cite{HII} in the case of unitary groups.
\end{abstract}

\tableofcontents

\section{Introduction}

Let $E/F$ be a quadratic extension of local fields of characteristic $0$. In this paper, we develop an explicit Plancherel decomposition for the symmetric space $\GL_n(F)\backslash \GL_n(E)$. The proof is, at least superficially, similar to the factorization of global Flicker-Rallis periods (see \cite{Fli2}, \cite[\S 2]{GJR} and \cite[\S 3.2]{Zh3}) and moreover the result is in remarkable agreement with general conjectures of Sakellaridis-Venkatesh \cite{SV} on the spectrum of general spherical varieties.

Let us mention here that Plancherel decompositions of symmetric spaces have already been thoroughly studied in the literature. In particular, most relevant for us is the work of Harinck \cite{Har} giving another explicit Plancherel formula for the case at hand when $E/F=\C/\bR$ (and more generally for symmetric spaces of the form $G(\bR)\backslash G(\C)$) and also the work of Delorme \cite{Del2} (following what was done by Sakellaridis-Venkatesh \cite{SV} for split spherical varieties) on $p$-adic symmetric spaces (which is less explicit since it consider as a black box the ``discrete spectrum'' of certain {\em Levi varieties}). However, the present work is rather orthogonal to those and uses heavily particular features of the pair $(\GL_n(E),\GL_n(F))$ allowing to express the relevant $L^2$-inner products as residues of certain families of local Zeta integrals (similarly, the Flicker-Rallis period is usually studied by considering it as the residue of a global family of Zeta integrals).

We then give two applications of this Plancherel formula: first to the global Ichino-Ikeda conjecture for unitary groups by completing a comparison between local relative characters that was left open by W. Zhang \cite{Zh3} and secondly to the Hiraga-Ichino-Ikeda conjecture on formal degrees for unitary groups \cite{HII}. Besides the Plancherel formula, the main tools to derive these applications are certain local analogs of the Jacquet-Rallis trace formulae completed with a certain comparison of ``relatively unipotent'' orbital integrals. Actually, this relationship between the Plancherel formula for $\GL_n(F)\backslash \GL_n(E)$ and local Jacquet-Rallis trace formulae was already investigated in the case $n=2$ by Ioan Filip in his PhD thesis \cite{Fil}.

In the rest of this introduction, we give more details on the main results.

\subsection{Plancherel formula for $\GL_n(F)\backslash \GL_n(E)$}

We now state the explicit Plancherel formula we obtain for the symmetric space $Y_n=\GL_n(F)\backslash \GL_n(E)$. Let $U(n)$ be a quasi-split unitary group of rank $n$ (with respect to the extension $E/F$) and $\Temp(U(n))/\stab$ be the set of all tempered $L$-packets of $U(n)$. Similarly, we let $\Temp(\GL_n(E))$ be the set of all isomorphism classes of tempered irreducible representations of $\GL_n(E)$. Finally, let $BC_n:\Temp(U(n))/\stab\to \Temp(\GL_n(E))$ be the {\em stable} base-change map if $n$ is odd and the {\em unstable} base-change map if $n$ is even (this last case requires a choice, which however is unimportant for what follows).

Let us emphasize here that the local Langlands correspondence for unitary groups is now fully known by \cite{Mok}, \cite{KMSW} and thus combining it with the local correspondence for $\GL_n$ (\cite{HT},\cite{Hen}, \cite{Sch}), we see that the preceding sentences make perfect sense. Actually, we could have written our Plancherel formula without appealing to the local Langlands correspondence since only the image of $BC_n$ matters and this can be described in a purely representation-theoretic way using local Asai $L$-functions (moreover, as we will see, this is more or less exactly how it shows up in the computations). Nevertheless, we prefer to write everything using this correspondence since we find the resulting formulations more suggestive and moreover this translation will play an important role for the applications (to be described in the second part of this introduction).

We equip $\Temp(U(n))/\stab$ with a natural and canonical measure $d\sigma$ which is locally given by a Haar measure on certain groups of {\em unramified characters} (see Section \ref{Section spectral measures}). A first weak version of our Plancherel formula can be stated as follows.

\begin{theointro}
There exists an isomorphism of unitary representations
$$\displaystyle L^2(Y_n)\simeq \int^\oplus_{\Temp(U(n))/\stab} \cH_\sigma d\sigma$$
where $\sigma\in \Temp(U(n))/\stab\mapsto \cH_\sigma$ is a measurable field of unitary representations of $\GL_n(E)$ with $\cH_\sigma\simeq BC_n(\sigma)$ for every $\sigma\in \Temp(U(n))/\stab$.
\end{theointro}

Note that as a direct corollary of this result, we can describe explicitly the discrete spectrum of $L^2(Y_n)$: the irreducible unitary representations $\pi$ of $\GL_n(E)$ embedding continuously into $L^2(Y_n)$ are precisely the base-change (stable or unstable as before) of discrete series of $U(n)$. Notice that those representations were already shown to embed in $L^2(Y_n)$ by Jerrod Smith \cite{Smith}. Moreover, the above decomposition confirms in the case of the symmetric space $Y_n$ a general conjecture of Sakellaridis-Venkatesh \cite[Conjecture 16.5.1]{SV} on Plancherel decompositions of spherical varieties\footnote{Strictly speaking, the general conjectures of \cite{SV} do not apply to $Y_n$ since in {\it loc. cit.} the group acting $G$ (which here is $R_{E/F}\GL_{n,E}$) was assumed to be {\em split} over $F$. However, in light of another conjecture of Jacquet (see \cite{Pra}) there is a natural way to extrapolate the conjecture of Sakellaridis-Venkatesh to the case at hand by considering the ``$L$-group'' of $Y_n$ to be ${}^L Y_n={}^L U(n)$ and the morphism ${}^L Y_n\to {}^L G$ to be given by base-change (stable or unstable as before).}.

To state our Plancherel formula more precisely, we need to introduce more notation. Let $\cS(\GL_n(E))$ be the {\em Schwartz space} of $\GL_n(E)$ i.e. the space of complex valued functions on $\GL_n(E)$ which are locally constant and compactly supported if $E$ is $p$-adic or which are $C^\infty$ and ``rapidly decreasing with all their derivatives'' in some suitable sense if $E$ is Archimedean (we refer the reader to Section \ref{Section space of functions} for a precise definition; a similar space of functions is defined on the set of $F$-points of any smooth algebraic variety over $F$). For $\pi\in \Temp(\GL_n(E))$, we define a semi-positive Hermitian form $(.,.)_{Y_n,\pi}$ on $\cS(\GL_n(E))$ by
$$\displaystyle (f_1,f_2)_{Y_n,\pi}:=\sum_{W\in \cB(\pi,\psi_n)} \beta(\pi(f_1)W)\overline{\beta(\pi(f_2)W)}.$$
Here, $\psi$ is a non-trivial additive character of $E$ {\em trivial} on $F$, $\psi_n$ the corresponding generic character of the standard maximal unipotent subgroup $N_n(E)$ of $\GL_n(E)$, $\cB(\pi,\psi_n)$ is a suitable orthonormal basis of the Whittaker model $\cW(\pi,\psi_n)$ for the natural scalar product on the associated Kirillov model (i.e. $L^2$-scalar product on $N_n(E)\backslash P_n(E)$ where $P_n$ is the mirabolic subgroup of $\GL_n$) and $\beta$ is the linear form given by
$$\displaystyle \beta(W):=\int_{N_n(F)\backslash P_n(F)} W(p)dp.$$
Using a nontrivial additive character $\psi'$ of $F$, we can define normalized Haar measures on $\GL_n(E)$ and $\GL_n(F)$ (see Section \ref{Section measures}) hence an invariant measure on $Y_n$. To $\sigma\in \Temp(U(n))/\stab$ we associate its adjoint $\gamma$-factor $\gamma(s,\sigma,\Ad,\psi')$ as well as a certain finite group $S_\sigma$ which is just the centralizer of its Langlands parameter when $\sigma$ is a discrete series (for the general case cf.\ \ref{Section Groups of centralizers}). Then, the explicit Plancherel formula we prove for $Y_n$ reads as follows (cf.\ Theorem \ref{theo Plancherel}).

\begin{theointro}\label{theo 2 intro}
For every $\sigma\in \Temp(U(n))/\stab$, the Hermitian form $(.,.)_{Y_n,BC_n(\sigma)}$ factorizes through the natural projection $\cS(\GL_n(E))\to \cS(Y_n)$. Moreover, assuming that all Haar measures have been normalized using the additive character $\psi'$, for every functions $\varphi_1,\varphi_2\in \cS(Y_n)$ we have
$$\displaystyle (\varphi_1,\varphi_2)_{Y_n}=\int_{\Temp(U(n))/stab} (\varphi_1,\varphi_2)_{Y_n,BC_n(\sigma)} \frac{\lvert \gamma^*(0,\sigma,\Ad,\psi')\rvert}{\lvert S_\sigma\rvert} d\sigma$$
where the left-hand side denotes the $L^2$-inner product on $Y_n$, the right-hand side is absolutely convergent and we have set
$$\displaystyle \gamma^*(0,\sigma,\Ad,\psi')=\left( \zeta_F(s)^{n_\sigma}\gamma(s,\sigma,\Ad,\psi')\right)_{s=0},$$
$n_\sigma$ being the order of the zero of $\gamma(s,\sigma,\Ad,\psi')$ at $s=0$.
\end{theointro}

Although we will not explain this in any details, the above explicit version of the Plancherel decomposition of $Y_n$ is pleasantly aligned with certain speculations of Sakellaridis-Venkatesh \cite[\S 17]{SV} on factorization of global automorphic periods. The main reason being that the linear form $\beta$ also appears in the factorization of the global Flicker-Rallis periods (we refer the reader to \cite[\S 2]{GJR} and \cite[\S 3.2]{Zh3} for a precise statement) and, as we will see in the second part of this introduction, the quotient $\frac{\lvert \gamma^*(0,\sigma,\Ad,\psi')\rvert}{\lvert S_\sigma\rvert}$ equals on the nose (again if we normalize measures correctly using the character $\psi'$) the Plancherel densities of unitary groups. In fact, the proof of Theorem \ref{theo 2 intro} is quite similar to the global computations leading to the factorization of Flicker-Rallis periods. Let us explain the main steps. Let $\varphi_1,\varphi_2\in \cS(Y_n)$ and choose (arbitrarily) functions $f_1,f_2\in \cS(\GL_n(E))$ such that $\varphi_i(x)=\int_{\GL_n(F)} f_i(hx)dh$ for $i=1,2$. Then, simple manipulations show that
$$\displaystyle (\varphi_1,\varphi_2)_{L^2(Y_n)}=\int_{\GL_n(F)} f(h)dh$$
where $f=\overline{f_2}\star f_1^\vee$ (here as usual $f_1^\vee$ denotes the function $g\mapsto f_1(g^{-1})$). Therefore, we are essentially reduced to finding a spectral expansion for the linear form $f\in \cS(\GL_n(E))\mapsto \int_{\GL_n(F)} f(h)dh$. The first step is then, by some ``local unfolding'', to prove an identity (see Proposition \ref{prop unfolding})
\begin{align}\label{eq -1 intro}
\displaystyle \int_{\GL_n(F)} f(h)dh=C_1\int_{N_n(F)\backslash P_n(F)} \int_{N_{n}(F)\backslash \GL_{n}(F)} W_f(p,h)dhdp,\;\;\; f\in \cS(\GL_n(E))
\end{align}
where $W_f$ is a certain Whittaker function associated to $f$ analog to the global Whittaker function associated to cusp forms (see Section \ref{section Plancherel Whitt} for the precise definition of $W_f$) and $C_1$ is a certain constant depending on the choice of the character $\psi$. Let $W_f=\int_{\Temp(\GL_n(E))} W_{f,\pi} d\mu_{\GL_n(E)}(\pi)$ be the spectral decomposition of the Whittaker function $W_f$ (that is the Plancherel decomposition for Whittaker functions which according to \cite[\S 6.3]{SV} can be deduced from the Plancherel formula for the group; see also Section \ref{section Plancherel Whitt} of this paper where we make this slightly more precise), where $d\mu_{\GL_n(E)}(\pi)$ denotes the Plancherel measure for $\GL_n(E)$. Then, we would like to use this expansion of $W_f$ to get a spectral decomposition for the right-hand side of \eqref{eq -1 intro}. Unfortunately, the resulting expression is not absolutely convergent and in particular we cannot switch the spectral integral with the integral over $N_{n}(F)\backslash \GL_{n}(F)$. To circumvent this difficulty, we will express the inner integral
$$\displaystyle \int_{N_{n}(F)\backslash \GL_{n}(F)} W_f(p,h)dh=\int_{Z_n(F)N_{n}(F)\backslash \GL_{n}(F)} W_{\widetilde{f}}(p,h)dh$$
where $\widetilde{f}=\int_{Z_n(F)} f(z.)dz$ ($Z_n(F)$ denoting the center of $\GL_n(F)$) as the residue at $s=0$ of some Zeta integral $Z(s,W_{\widetilde{f}}(p,.),\phi)$ associated to an auxiliary test function $\phi\in \cS(F^n)$ (see Lemma \ref{lem Zeta function}(ii) for a precise statement). The important point is that the formation of this Zeta integral ``commutes'' with the spectral decomposition of $W_{\widetilde{f}}$ (i.e. the resulting expression is absolutely convergent) when $\Re(s)>0$ so that in this range we can write
$$\displaystyle Z(s,W_{\widetilde{f}}(p,.),\phi)=\int_{\Temp(\overline{\GL_n(E)})} Z(s,W_{f,\pi}(p,.),\phi) d\mu_{\overline{\GL_n(E)}}(\pi)$$
where we have set $\overline{\GL_n(E)}=Z_n(F)\backslash \GL_n(E)$. Then, using a local functional equation for the above Zeta integrals in terms of Asai $\gamma$-factors $\gamma(s,\pi,\As,\psi')$ (see Section \ref{Section Zeta integrals}), this can be rewritten as an expression of the form
\begin{align}\label{eq 0 intro}
\displaystyle \int_{\Temp(\overline{\GL_n(E)})} \Phi_s(\pi) \gamma(s,\pi,\As,\psi')^{-1}d\mu_{\overline{\GL_n(E)}}(\pi)
\end{align}
where now the family of function $s\mapsto \Phi_s$ is analytic at $s=0$. Thus, we end up with the problem of computing the residue at $s=0$ of the above distribution. This can be achieved by a direct computation using an explicit formula for the Plancherel measure $d\mu_{\overline{\GL_n(E)}}(\pi)$ in terms of the adjoint $\gamma$-factor of $\pi$ which is essentially due to Harish-Chandra \cite{H-C} in the Archimedean case and follows from work of Shahidi \cite{Sha1} and Silberger-Zink \cite{SZ} in the $p$-adic case (see \cite{HII} for a convenient uniform reformulation of all these results). Although the computation is not very enlightening (it is done in Part \ref{Part II} of this paper), the final result is very neat: we get, up to an explicit constant, the integral of the function $\Phi=\Phi_0$ against the push-forward by base-change of the measure $\frac{\gamma^*(0,\sigma,\Ad,\psi')}{\lvert S_\sigma\rvert} d\sigma$ on $\Temp(U(n))/\stab$ (see Proposition \ref{prop1 Theo limite spectrale}). Together with the previous steps, this then very easily implies Theorem \ref{theo 2 intro}.

\subsection{Applications to the Ichino-Ikeda and formal degree conjectures for unitary groups}

In Part \ref{Part IV} of this paper we give two applications of the Plancherel formula of Theorem \ref{theo 2 intro}. Namely, we establish the formal degree conjecture of Hiraga-Ichino-Ikeda \cite{HII} for unitary groups as well as a certain comparison of local relative characters left open by W. Zhang \cite[Conjecture 4.4]{Zh3} and which has application to the so-called Ichino-Ikeda conjecture for unitary groups \cite{Ha}. We now state these two results in turn.

Let $V$ be a (finite dimensional) Hermitian space over $E$, $U(V)$ the corresponding unitary group and $\mu_{U(V)}^*(\pi)$ the ``Plancherel density'' for $U(V)(F)$ i.e. the unique function on the tempered dual $\Temp(U(V))$ of $U(V)(F)$ such that the Plancherel measure for $U(V)(F)$ is given by $d\mu_{U(V)}(\pi)=\mu_{U(V)}^*(\pi) d\pi$ where $d\pi$ is a certain natural and canonical measure on $\Temp(U(V))$ (see Section \ref{Section spectral measures}). Notice that the Plancherel measure depends on some Haar measure on $U(V)(F)$ for which, as before, there is a natural choice depending only on the nontrivial additive character $\psi'$ of $F$. Therefore, $\mu_{U(V)}^*(\pi)$ also depends on $\psi'$. With these notation we prove:

\begin{theointro}\label{theo 3 intro}
We have
$$\displaystyle \mu_{U(V)}^*(\pi)=\frac{\lvert \gamma^*(0,\pi,\Ad,\psi')\rvert}{\lvert S_\pi\rvert}$$
for almost all $\pi\in \Temp(U(V))$ where $\gamma^*(0,\pi,\Ad,\psi')$ and $S_\pi$ are defined as before. In particular, for every discrete series $\pi$ of $U(V)(F)$ we have
$$\displaystyle d(\pi)=\frac{\lvert \gamma(0,\pi,\Ad,\psi')\rvert}{\lvert S_\pi\rvert}$$
where $d(\pi)$ denotes the formal degree of $\pi$ and $S_\pi$ is the centralizer of the Langlands parameter of $\pi$.
\end{theointro}

The second part of the above theorem is precisely \cite[Conjecture 1.4]{HII} for unitary groups. Also, although we will not use it, the first part can be deduced from the second part using Langlands' normalization of standard intertwining operators (which is known for unitary groups see \cite[Proposition 3.3.1]{Mok}, \cite[Lemma 2.2.3]{KMSW}) the formal degree conjecture for $\GL_n(E)$ (which is also known by work of Silberger-Zink, see \cite[Theorem 3.1]{HII}) and the description of the Plancherel measure as in \cite{Wald1}. We also remark that Theorem \ref{theo 3 intro} is not new when $F$ is Archimedean since in that case it can be deduced from Harish-Chandra's explicit formula for the Plancherel measure \cite{H-C} (see \cite[Proposition 2.1]{HII} for the translation, at least for discrete series).

To state the second application, we need to introduce more notation. Let $(V,h)$ be a $n$-dimensional Hermitian space over $E$ and set $H=U(V)$, $G=U(V)\times U(V')$ where $V'=V\oplus Ev_0$ is equipped with the Hermitian form given by $h'(v_1+\lambda v_0,v_2+\mu v_0)=h(v_1,v_2)+\lambda \mu^c$ for all $v_1,v_2\in V$ and $\lambda,\mu\in E$ (in this paper, we denote by $c$ the non-trivial $F$-automorphism of $E$). We consider $H$ as a subgroup of $G$ through the natural diagonal embedding $H\hookrightarrow G$. We also define $G'=R_{E/F} \GL_{n,E}\times R_{E/F}\GL_{n+1,E}$ (here, as elsewhere in the paper, $R_{E/F}$ stands for Weil's restriction of scalars from $E$ to $F$) and its two subgroups $H_1=R_{E/F}\GL_{n,E}$ (diagonally embedded) and $H_2=\GL_{n,F}\times \GL_{n+1,F}$. Then, Jacquet and Rallis \cite{JR} have defined a notion of {\em transfer} between functions $f\in \cS(G(F))$ and $f'\in \cS(G'(F))$. This transfer is itself the byproduct of a natural injective correspondence between {\em regular semi-simple} orbits (or double cosets) \cite[Lemme 2.3]{Zh2}
\begin{align}\label{eq 1 intro}
\displaystyle H(F)\backslash G_{\rs}(F)/H(F) \hookrightarrow H_1(F)\backslash G'_{\rs}(F)/H_2(F).
\end{align}
where $G'_{\rs}$ (resp. $G_{\rs}$) denotes the open subset of elements $\gamma\in G'$ (resp. $\delta\in G$) whose double coset $H_1\gamma H_2$ (resp. $H\delta H$) is closed and with a trivial stabilizer in $H_1\times H_2$ (resp. in $H\times H$). For $(f,f')\in \cS(G(F))\times \cS(G'(F))$ and $(\delta,\gamma)\in G_{\rs}(F)\times G'_{\rs}(F)$ we define as usual (relative) orbital integrals $O(\delta,f)$ and $O_\eta(\gamma,f')$ (the second one being twisted by a certain quadratic character $\eta$ of $H_2(F)$; see Section \ref{Section relative orb integrals}). Then, we say that $f$ and $f'$ {\em match} (or that they are {\em transfer of each other}) if we have an identity
$$\displaystyle \Omega(\gamma)O_\eta(\gamma,f')=O(\delta,f)$$
whenever $\gamma$ and $\delta$ correspond to each other by the correspondence \eqref{eq 1 intro} and where $\Omega(\gamma)$ is a certain {\em transfer factor} (which has an explicit and elementary definition, see Section \ref{Section relative orb integrals}). It is one of the main achievement of W. Zhang \cite{Zha1} that in the $p$-adic case every function $f\in \cS(G(F))$ admits a transfer $f'\in \cS(G'(F))$ and conversely. For a partial result in that direction in the Archimedean case, which is however sufficient for applications, we refer the reader to \cite{Xue} or Section \ref{Section relative orb integrals} of this paper.

We define as in \cite{Zh3} relative characters $f\in \cS(G(F))\mapsto J_\pi(f)$ and $f'\in \cS(G'(F))\mapsto I_\Pi(f')$ for every tempered irreducible representations $\pi\in \Temp(G)$ and $\Pi\in \Temp(G')$ (we warn the reader here that our convention for $J_\pi$ and $I_\Pi$ is slightly different from the one in \cite{Zh3}, in particular we have replaced $\pi$ and $\Pi$ by their contragredient and moreover we don't use the same scalar product as in {\it loc. cit.} to normalize $I_\Pi$; these changes are nevertheless minor and don't affect the global applications: see remark below). By its very definition, the family of relative characters $\pi\mapsto J_\pi$ is supported on the set $\Temp_H(G)$ of {\em $H$-distinguished} tempered representations (i.e. the one admitting a nonzero continuous $H(F)$-invariant form). The second main theorem of Part \ref{Part IV} can now be stated as follows (see Theorem \ref{theo1 main result}).

\begin{theointro}\label{theo 4 intro}
Let $f\in \cS(G(F))$ and $f'\in \cS(G'(F))$ be matching functions. Then, for every $\pi\in \Temp_{H}(G)$, we have
$$\displaystyle \kappa_V J_\pi(f)=I_{BC(\pi)}(f')$$
where $BC(\pi)$ denotes the {\em stable} base-change of $\pi$ and $\kappa_V$ is an explicit constant which depends only on $V$ and the normalization of transfer factors.
\end{theointro}

\begin{remintro}
\begin{itemize}
\item The above result was already proved in \cite{Beu2} when $F$ is $p$-adic. The proof we give here, although using similar tools, differs at some crucial points from {\it loc. cit.} and moreover has the good feature, at least to the author's taste, of treating uniformly the Archimedean and non-Archimedean case. In particular, in this paper we make no use at all of the result of W. Zhang \cite{Zh3} on {\em truncated local expansion} for the relative characters $I_\Pi$. We need however to import from \cite{Zha1} the existence of smooth transfer as well as its compatibility with Fourier transform at the Lie algebra level. Furthermore, although it might not be so transparent, we need the Jacquet-Rallis fundamental lemma proved by Yun and Gordon \cite{Yu} in order to derive the ``weak comparison'' of relative characters given in Proposition \ref{prop 1 weak comparison} from a global comparison of (simple) Jacquet-Rallis trace formulae. 

\item The above theorem confirms \cite[Conjectuire 4.4]{Zh3}. Actually, as already pointed, our normalization for relative characters is not the same thus explaining why our constant $\kappa_V$ differs from {\it loc. cit.} However, this discrepancy is not completely afforded by the change of normalization since in \cite{Zh3} the constant up to which ``Fourier transform commutes with transfer'' was not exactly the correct one and here we used the one computed by Chaudouard in \cite{Chau}. For the global applications however (see below) this difference is inessential: all that matters is that if $V$ is a Hermitian space relative to a quadratic extension $k'/k$ of number fields then the product of the local constants $\prod_v \kappa_{V_v}$ is one.
\end{itemize}
\end{remintro}

As explained in \cite{Zh3} and \cite{Beu2}, Theorem \ref{theo 4 intro} has direct application to the Ichino-Ikeda conjecture for unitary groups. Namely, from \cite[Theorem 4.3]{Zh3} and \cite[Theorem 3.5.1]{Beu2} and the above theorem we deduce:

\begin{theointro}
Let $k'/k$ be a quadratic extension of number fields, $V$ an Hermitian space over $k'$ and define the groups $H\subset G$ as above. Let $\pi=\bigotimes'_v \pi_v$ be a cuspidal automorphic representation of $G(\bA)$ ($\bA$ being the ring of ad\`eles of $k$) satisfying:
\begin{itemize}
\item For every place $v$ of $k$, the representation $\pi_v$ is tempered;
\item There exists a non-Archimedean place $v$ of $k$ with $BC(\pi_v)$ supercuspidal (e.g. if $v$ splits in $k'$ and $\pi_v$ is itself supercuspidal).
\end{itemize}
Then, the Ichino-Ikeda conjecture as stated in \cite[Conjecture 1.2]{Ha} or \cite[Conjecture 1.1]{Zh3} is true for $\pi$.
\end{theointro}

\begin{remintro}
Of course, the above theorem also includes previous works of many other authors including W. Zhang, Z. Yun and J. Gordon as well as H. Xue. We have mainly stated it for convenience to the reader that would like to get the most up-to-date result on the Ichino-Ikeda conjecture. Moreover, the assumption that $BC(\pi_v)$ is supercuspidal for at least one place $v$ (which is the last one that needs to be removed) originates from the use of some simple version of Jacquet-Rallis trace formulae. To drop this assumption, one would need complete spectral decompositions of these trace formulae which is work in progress by Chaudouard and Zydor (see \cite{Zyd} and \cite{CZ} for partial progress in that direction).
\end{remintro}

Theorem \ref{theo 3 intro} and Theorem \ref{theo 4 intro} are proved together by comparing certain distributions on $G(F)$ and $G'(F)$. Namely, we will compare both local analogs of the aforementioned Jacquet-Rallis trace formulae as well as certain ``relatively unipotent'' orbital integrals. To be more specific, we discuss the two comparisons in turn.

First, the local analog of the unitary Jacquet-Rallis trace formula (which was already used in \cite{Beu2}) is a certain $H(F)\times H(F)$-invariant sesquilinear form $J(.,.)$ on $\cS(G(F))$ admitting the following ``geometric'' and ``spectral'' expansions (see Section \ref{unitary JR}):
\begin{align}\label{eq 2 intro}
\int_{H(F)\backslash G_{rs}(F)/H(F)} O(\delta,f_1)\overline{O(\delta,f_2)} d\delta=J(f_1,f_2)=\int_{\Temp(G)}J_\pi(f_1)\overline{J_\pi(f_2)} d\mu_{G}(\pi)
\end{align}
where $d\delta$ is a certain natural measure on the set of regular semi-simple orbits and $d\mu_{G}(\pi)$ denotes the Plancherel measure of $G(F)$. Similarly, the local analog of the linear Jacquet-Rallis trace formula is a certain $H_1(F)\times (H_2(F),\eta)$-equivariant sesquilinear form $I(.,.)$ on $\cS(G'(F))$ admitting the following ``geometric'' and ``spectral'' expansions (see Section \ref{linear JR}):
\begin{align}\label{eq 3 intro}
\displaystyle & \int_{H_1(F)\backslash G'_{rs}(F)/H_2(F)} O_\eta(\gamma,f'_1)\overline{O_\eta(\gamma,f'_2)} d\gamma=I(f_1',f_2') \\
\nonumber & =C_2\int_{\Temp(G_{\qs})/\stab} I_{BC(\pi)}(f'_1)\overline{I_{BC(\pi)}(f'_2)} \frac{\lvert \gamma^*(0,\pi,\Ad,\psi')\rvert}{\lvert S_\pi\rvert}d\pi
\end{align}
where $d\gamma$ is a certain natural measure on the set of regular semi-simple orbits, $G_{\qs}$ is a quasi-split inner form of $G$, $\Temp(G_{\qs})/\stab$ denotes the set of tempered $L$-packets for $G_{\qs}$ and $C_2$ is a certain positive constant depending on the choice of $\psi$. Here we remark that the proof of \eqref{eq 2 intro} doesn't pose any analytical difficulty and is rather easy and direct whereas to get the spectral side of \eqref{eq 3 intro} we need to use the explicit Plancherel formula of Theorem \ref{theo 2 intro}. When $f_k$ match $f'_k$ for $k=1,2$ the left-hand sides of \eqref{eq 2 intro} and \eqref{eq 3 intro} are easily seen to be equal (this is not quite correct, as we need first to sum \eqref{eq 2 intro} over some {\em pure inner forms} but we will ignore this issue in the introduction) and we deduce from this the equality of the right-hand sides. This is the first comparison that we will use.

Next, we consider the ``trivial'' orbital integral $O(1,f)=\int_{H(F)} f(h)dh$ on $\cS(G(F))$ for which we have the spectral expansion (see \eqref{unitary unipotent})
\begin{align}\label{eq 4 intro}
\displaystyle O(1,f)=\int_{\Temp(G)}J_\pi(f)d\mu_{G}(\pi)
\end{align}
as well as a certain ``regularized'' regular unipotent orbital integral $f'\mapsto O_+(f')$ on $\cS(G'(F))$ for which we have the spectral expansion (see Proposition \ref{prop 2 unipotent linear})
\begin{align}\label{eq 5 intro}
\displaystyle \gamma O_+(f')=C_3\int_{\Temp(G_{\qs})/\stab} I_{BC(\pi)}(f) \frac{\lvert \gamma^*(0,\pi,\Ad,\psi')\rvert}{\lvert S_\pi\rvert}d\pi
\end{align}
where $C_3$ is again a positive constant depending on $\psi$ and $\gamma$ a certain product of abelian local $\gamma$-factors. Incidentally (or not), the same regularized orbital integral appears in the truncated local expansions of W. Zhang \cite{Zh3} for the relative characters $I_\Pi$. Once again, the identity \eqref{eq 4 intro} comes almost for free and doesn't require much hard work whereas \eqref{eq 5 intro} is a consequence of (some form of) the explicit Plancherel formula of Theorem \ref{theo 2 intro}. Then, using analogs of \eqref{eq 4 intro} and \eqref{eq 5 intro} for Lie algebras and the fact that Fourier transform commutes, up to a constant, with transfer (\cite{Zha1},\cite{Xue}), we can show that if $f\in \cS(G(F))$ and $f'\in \cS(G'(F))$ match then $O(1,f)$ and $O_+(f')$ are equal up to some (absolute and explicit) constant (once again the correct statement involves summing \eqref{eq 4 intro} over some set of pure inner forms). This, together with \eqref{eq 4 intro} and \eqref{eq 5 intro}, implies a second spectral identity between matching functions.

Theorem \ref{theo 4 intro} and Theorem \ref{theo 3 intro} (or rather an analog of Theorem \ref{theo 3 intro} pertaining to the Plancherel measure of $G$ restricted to $\Temp_H(G)$) can then be easily deduced from these two comparisons combined with some ``weak'' identity of relative characters (Proposition \ref{prop 1 weak comparison}). We refer the reader to Section \ref{Section proof of main theorems} for all the details of the proof.

\subsection{Outline of the paper}

We now briefly describe the content of the paper. Part \ref{Part I} is devoted to fixing notation and numerous normalizations. We also collect there various results from the literature which will be needed in the sequel in particular regarding local base-change for unitary groups (Section \ref{Section LLC}), local $\gamma$-factors (Section \ref{Section gamma factors}), the Plancherel formula for the group and for Whittaker functions (Sections \ref{Section Planch} and \ref{section Plancherel Whitt}) and the local functional equation for Asai $\gamma$-factors of Rankin-Selberg type (Section \ref{Section Zeta integrals}). Part \ref{Part II} is the most technical part of the paper and is concerned with computing explicitly poles of analytic families of distributions like \eqref{eq 0 intro}. As we said, the final result is very neat but, unfortunately, the author wasn't able to find a conceptual way to derive it, therefore we have essentially reduced everything to certain computations in $\bR^n$. In Part \ref{Part III}, we establish the explicit Plancherel formula of Theorem \ref{theo 2 intro}. This is a rather easy task as the adequate preliminary results have been obtained in Part \ref{Part II}. In the last part of this paper, Part \ref{Part IV}, we prove Theorems \ref{theo 3 intro} and \ref{theo 4 intro} following the outline given above. Finally, we have added one appendix at the end of the paper containing the proof of a rather technical result related to the Harish-Chandra Plancherel formula and for which the author was unable to find a proper reference in the literature.

\subsection{Acknowledgment}

The results of this paper (in a slightly weaker form) have been announced by the author in his ``Cours Peccot'' in April-May 2017. I would like to thank the Coll\`ege de France for giving me the opportunity to give this course which has probably accelerated the present work and I apologize to anyone who was promised an earlier preprint.

I am grateful to Erez Lapid and Jean-Loup Waldspurger for helpful comments on a first version of this paper. I also thank the anonymous referee for his very accurate comments and suggestions for improvement.

The project leading to this publication has received funding from Excellence Initiative of Aix-Marseille University-A*MIDEX, a French ``Investissements d'Avenir'' programme.

\section{Preliminaries}\label{Part I}
\subsection{General notation}\label{Section general notation}
Here is a list of notation that will be used throughout in this paper:
\begin{itemize}
\item $E/F$ is a quadratic extension of local fields of characteristic zero (either Archimedean or non-Archimedean), we denote by $c$ the non-trivial $F$-automorphism of $E$, $\Tra_{E/F}:E\to F$ the trace and by $\eta_{E/F}:F^\times\to \{\pm 1 \}$ the quadratic character associated to $E/F$ by local class field theory. Also, in the non-Archimedean case we write $\cO_F$ and $\cO_E$ for the ring of integers of $F$ and $E$ respectively.
\item We fix a character $\eta'$ of $E^\times$ extending $\eta_{E/F}$ as well as nontrivial additive characters $\psi'$ and $\psi$ of $F$ and $E$ respectively with {\em $\psi$ trivial on $F$}. We denote by $\tau$ the unique element in $E$ of trace zero such that $\psi(z)=\psi'(\Tra_{E/F}(\tau z))$ for all $z\in E$. We will also denote by $\psi'_E$ the character $z\mapsto \psi'(\Tra_{E/F}(z))$ of $E$.
\item $W_F$ and $W_E$ are the Weil groups of $F$ and $E$ respectively, $q_F$, $q_E$ the cardinality of the residue fields of $F$ and $E$ in the $p$-adic case whereas in the Archimedean case we set $q_F=q_E=e^{1/2}$.
\item We fix an algebraic closure $\overline{F}$ of $F$ containing $E$ and for every finite extension $K$ of $F$ we denote by $\lvert .\rvert_K$ the normalized absolute value on $K$. Most of the time we will simply write $\lvert .\rvert$ for $\lvert .\rvert_F$.
\item If $X$ is an algebraic variety over $F$ and $K/F$ a finite extension, we denote by $X_K$ the variety obtained by base-change from $F$ to $K$. In the other direction, if $X$ is an algebraic variety over $K$, we denote by $R_{K/F}X$ the Weil's restriction of scalar of $X$ from $K$ to $F$.
\item Unless otherwise specified, all algebraic varieties will be tacitly assumed to be defined over $F$.
\item If $X$ is an algebraic variety over $F$, we will use freely the notion of {\em norms} on $X(F)$ as defined by Kottwitz \cite[Sect. 18]{Kott}. We always denote by $\lVert .\rVert_X$ such a norm and set $\sigma_X=\log(2+\lVert .\rVert_X)$ for the associated {\em log-norm}. We refer the reader to \cite[\S 1.2]{Beu1} for more details on these notions. 
\item For $K$ a compact group, we denote by $\widehat{K}$ the set of all isomorphism classes of its irreducible representations. Moreover, if $V$ is a representation of $K$ for each $\rho\in \widehat{K}$ we denote by $V[\rho]$ its $\rho$-isotypic component. Similarly, if $A$ is a locally compact abelian group, $\widehat{A}$ stands for its dual group.
\item $\Re(z)$ and $\Im(z)$ will stand for the real and imaginary part respectively of a complex number $z$.
\item $\cH=\{s\in \C\mid \Re(s)>0 \}$ and $\lim\limits_{s\to 0^+}$ means that we take the limit as $s$ goes to $0$ from $\cH$.
\item A sentence like ``$f(x)\ll g(x)$ for all $x\in X$'' means that there exists $C>0$ such that $f(x)\leqslant Cg(x)$ for all $x\in X$. Also, we write $f(x)\sim g(x)$ when $f(x)\ll g(x)$ and $g(x)\ll f(x)$.
\item Lie algebras of algebraic groups will always be denoted by the corresponding gothic letter (e.g. $\g$ for $G$ or $\h$ for $H$).
\item For each integer $n\geqslant 1$, $\fS_n$ stands for the group of permutations of $\{ 1,\ldots,n\}$.
\item A holomorphic function $\varphi: \C^n\to \C$ will be said to be of {\em moderate growth in vertical strips together with all its derivatives} if for all $a,b\in \R$ and every holomorphic differential operator with constant coefficients $D$ on $\C^n$ there exists $N\geqslant 1$ such that
$$\displaystyle \lvert (D\varphi)(z_1,\ldots,z_n)\rvert \ll (1+\lvert \Im (z_1)\rvert)^N\ldots (1+\lvert \Im (z_n)\rvert)^N$$
for ell $(z_1,\ldots,z_n)\in \C^n$ with $a<\Re(z_1),\ldots,\Re(z_n)<b$. Similarly, a meromorphic function $\varphi:\C\to \C$ will be said to be of {\em moderate growth on vertical strips away from its poles together with all its derivatives} if for all $a,b\in \R$ there exists $R\in \C(T)$ and for every $n\geqslant 0$ there exists $N\geqslant 1$  and  such that
$$\displaystyle \left\lvert \left(\frac{d}{dz}\right)^n(R(z)\varphi(z))\right\rvert\ll (1+\lvert \Im(z)\rvert)^N$$
for all $z\in \C$ with $a<\Re(z)<b$.
\item We also say that a smooth function $\varphi:\bR^n\to \C$ is of {\em moderate growth together with all its derivatives} if for every differential operator with constant coefficients $D$ on $\bR^n$ there exists $N\geqslant 1$ such that
$$\displaystyle \lvert (D\varphi)(x_1,\ldots,x_n)\rvert \ll (1+\lvert x_1\rvert)^N\ldots (1+\lvert x_n\rvert)^N$$
for all $(x_1,\ldots,x_n)\in \bR^n$. 
\item For each $n\geqslant 1$, we denote by $(e_1,\ldots,e_n)$ the standard basis of $F^n$.
\item If $f$ is a function on a group $G$, we set $f^\vee(g)=f(g^{-1})$ for all $g\in G$.
\item If a group $G$ acts on a set $X$ on the right (resp. on the left), we shall denote by $R$ (resp. $L$) the right (resp. left) regular action of $G$ on functions on $X$. This action usually extends to some space of functions on $G$. If moreover $G$ is a Lie group, $X$ is a smooth manifold and the action is differentiable, we denote by the same letter the resulting action of the Lie algebra $\g$ of $G$ and also of its enveloping algebra. If $G$ is a topological group equipped with a Haar measure and the action is on the right and continuous, for every function $f$ and $\varphi$ on $G$ and $X$ respectively we set (whenever it makes sense)
$$\displaystyle (\varphi\star f)(x)=\int_G \varphi(xg^{-1})f(g)dg,\;\;\; x\in X.$$
Notice that $\varphi\star f=R(f^\vee)\varphi$.
\end{itemize}

\subsection{Groups}\label{Section groups}

Let $G$ be a connected reductive group over $F$. We assume that $G$ is fixed until the end of Section \ref{section Plancherel Whitt}. We denote by $Z_G$ the center of $G$ and $A_G$ the maximal split torus in $Z_G$. Let $X^*(G)$ be the group of algebraic characters of $G$ defined over $F$. We set
$$\cA_G^*=X^*(G)\otimes \R=X^*(A_G)\otimes \R, \;\cA_{G,\C}^*=\cA_G^*\otimes_{\R} \C$$
and
$$\cA_G=\Hom(X^*(G),\R),\; \cA_{G,\C}=\cA_G\otimes_\R \C$$
Let $\langle.,.\rangle$ be the natural pairing between $\cA_{G,\C}^*$ and $\cA_{G,\C}$. We define a morphism $H_G:G(F)\to \cA_G$ by $\langle \chi, H_G(g)\rangle=\log \lvert \chi(g)\rvert$ for all $\chi\in X^*(G)$ and $g\in G(F)$. As usual a sentence like ``Let $P=MU$ be a parabolic subgroup of $G$'' means that $P$ is a parabolic subgroup of $G$ defined over $F$, with unipotent radical $U$ and $M$ is a Levi component of $P$. A {\em Levi subgroup} of $G$ means a Levi component of a parabolic subgroup. If $M$ is a Levi subgroup, we denote by $\mathcal{P}(M)$ the set of all parabolic subgroups with Levi component $M$ and we write $W(G,M)=\Norm_{G(F)}(M)/M(F)$ for its associated Weyl group. If $P=MU$ is a parabolic subgroup, we denote by $\delta_P$ the modular character of $P(F)$ and by $H_P:P(F)\to \cA_M$ the composition of $H_M$ with the projection from $P(F)$ to $M(F)$. If moreover a maximal compact subgroup $K$ of $G(F)$, which is special in the $p$-adic case, has been fixed (so that $G(F)=P(F)K$ by the Iwasawa decomposition) we extend $H_P$ to $G(F)$ by setting $H_P(muk)=H_M(m)$ for every $(m,u,k)\in M(F)\times U(F)\times K$.

When the group $G$ is understood from the context, we will simply write $\lVert .\rVert$ and $\sigma$ for the norm and log-norm $\lVert .\rVert_G$ and $\sigma_G$ respectively.

In the Archimedean case, we denote by $\cU(\g)$ the enveloping algebra of the complexified Lie algebra of $G(F)$ and by $\cZ(\g)$ its center.

For every integer $n\geqslant 1$, we write $G_n$ for the algebraic group $\GL_n$ (say defined over $\bZ$). We denote by $Z_n$, $B_n$, $A_n$ and $N_n$ the subgroups of scalar, resp. upper triangular, resp. diagonal, resp. unipotent upper triangular matrices in $G_n$. For $g\in G_n$, we write ${}^t g$ for its transpose and $g_{i,j}$, $1\leqslant i,j\leqslant n$, for its entries. Also, for $a\in A_n$ we simply write $a_i$ for the diagonal entry $a_{i,i}$ ($1\leqslant i\leqslant n$). We denote by $\delta_n$ and $\delta_{n,E}$ the modular characters of $B_n(F)$ and $B_n(E)$ respectively. We will always consider $G_n$ as a subgroup of $G_{n+1}$ through embedding in the upper left corner i.e. $g\in G_n\mapsto \begin{pmatrix} g & \\ & 1 \end{pmatrix}\in G_{n+1}$. We let $P_n$ be the mirabolic subgroup that is the set of elements $g\in G_n$ with last row $(0,\ldots,0,1)$, $U_n$ be the unipotent radical of $P_n$ (so that $P_n=G_{n-1}U_n$) and $N'_{n}$ be the derived subgroup of $N_n$ (i.e. the subgroup of $u\in N_n$ such that $u_{i,i+1}=0$ for $1\leqslant i\leqslant n-1$). The natural extension of $c$ to $G_n(E)$ will be written as $g\mapsto g^c$ and for every $g\in G_n(E)$ we set $g^*=({}^tg^{-1})^c$. We also let $K_n$ and $K_{n,E}$ be the standard maximal compact subgroups of $G_n(F)$ and $G_n(E)$ i.e. $K_n=G_n(\cO_F)$, $K_{n,E}=G_n(\cO_E)$ in the $p$-adic case whereas $K_n=\{g\in G_n(F)\mid g{}^t g=I_n\}$ and $K_{n,E}=\{g\in G_n(E)\mid g{}^t g^c=I_n\}$ in the Archimedean case. We will denote by $\psi'_n:N_n(F)\to \mathbb{S}^1$ and $\psi_n:N_n(E)\to \mathbb{S}^1$ the characters defined by
$$\displaystyle \psi'_n(u)=\psi'\left((-1)^n\sum_{i=1}^{n-1} u_{i,i+1}\right) \mbox{ and } \psi_n(u)=\psi\left((-1)^n\sum_{i=1}^{n-1} u_{i,i+1}\right)$$
Most of the time, we will consider $G_n(E)$ as the group of $F$-points of $R_{E/F}G_{n,E}$ (so that all constructions involving the $F$-points of a reductive group over $F$ can be applied to $G_n(E)$). We define
$$\displaystyle \overline{G_n(E)}:= G_n(E)/Z_n(F)$$
and we will always consider this group as the set of $F$-points of $(R_{E/F}G_{n,E})/Z_{n,F}$.

Recall that we have fixed a character $\eta':E^\times \to \C^\times$ such that $\eta'_{\mid F^\times}=\eta_{E/F}$. We define for every integer $k\geqslant 1$ characters $\eta_k$ and $\eta'_k$ of $G_n(F)$ and $G_n(E)$ respectively by
$$\displaystyle \eta_k(h)=\eta_{E/F}(h)^{k+1} \mbox{ for } h\in G_n(F)$$
and
$$\displaystyle \eta'_k(g)=\left\{\begin{array}{ll}
\eta'(\det g) \mbox{ if } k \mbox{ is even,} \\
 \\
1 \mbox{ if } k \mbox{ is odd,}
\end{array}\right.$$
for $g\in G_n(E)$. Note that the restriction of $\eta'_k$ to $G_n(F)$ is equal to $\eta_k$.

Whenever $\chi$ is a character of $F^\times$, for $h\in G_n(F)$ we will usually write $\chi(h)$ for $\chi(\det h)$. This in particular will be applied to the character $\eta_{E/F}$.

By a {\em Hermitian space} we will always mean a finite dimensional $E$-vector space $V$ equipped with a non-degenerate Hermitian form $h$ that is linear in the first variable. If $V$ is such a Hermitian space, we denote by $U(V)$ the corresponding unitary group thought as an algebraic group defined over $F$. For every $n\geqslant 1$, we also define $U(n)$ as the unitary group of the standard quasi-split Hermitian form on $E^n$ i.e. for every $F$-algebra $R$ we have
$$\displaystyle U(n)(R)=\left\{g\in G_n(R\otimes_F E)\mid {}^t g^c w_n g=w_n\right\} \mbox{ where } w_n=\begin{pmatrix} & & 1 \\ & \iddots & \\ 1 & & \end{pmatrix}.$$
If $G=U(V)$ is a unitary groups its Levi subgroups all have the form
$$\displaystyle M\simeq R_{E/F}G_{n_1,E}\times \ldots\times R_{E/F} G_{n_k,E}\times U(W)$$
where $n_1,\ldots,n_k$ are positive integers and $W$ is a nondegenerate subspace of $V$ such that $\dim(V)=\dim(W)+2n_1+\ldots+2n_k$. For such a Levi subgroup there is a unique identification $\cA_M\simeq \R^k$ such that $H_M(g_1,\ldots,g_k,\gamma)=(\log(\lvert \det g_1\rvert_E),\ldots, \log(\lvert \det g_k\rvert_E))$ for every $(g_1,\ldots,g_k,\gamma)\in M(F)$. There is also a natural identification of $W(G,M)$ with the subgroup of permutations $w\in \fS_{2k}$ preserving the partition $\{\{i,2k+1-i\}\mid 1\leqslant i\leqslant k \}$ of $\{ 1,\ldots,2k\}$ and such that $n_{w(i)}=n_i$ for each $1\leqslant i\leqslant k$ where we have set $n_j=n_{2k+1-j}$ for $k+1\leqslant j\leqslant 2k$. This identification is uniquely characterized by the following: for every $w\in W(G,M)$ there exists a representatives $\widetilde{w}\in Norm_{G(F)}(M)$ such that $\widetilde{w}^{-1}(g_1,\ldots,g_k,\gamma)\widetilde{w}=(g_{w(1)},\ldots,g_{w(k)},\gamma)$ for all $(g_1,\ldots,g_k,\gamma)\in M(F)$ where we have set $g_{j}=g^*_{2k+1-j}$ for $k+1\leqslant j\leqslant 2k$.

\subsection{LF spaces, completed tensor products}\label{Section Reminder on evtlcs}

In this paper, by a {\em locally convex topological vector space} (LCTVS) we always mean a Hausdorff locally convex topological vector space over $\C$. If $V$ is a LCTVS, we will denote by $V'$ its continuous dual equipped with the weak topology (i.e. the topology of pointwise convergence). Let $V$ be a LCTVS. Recall that $V$ is said to be {\em quasi-complete} if every closed bounded subset of $V$ is complete and that a {\em total} subspace $H\subseteq V'$ is a subspace such that $\bigcap_{\lambda\in H} \Ker(\lambda)=\{ 0\}$. We will use the following basic notions extensively in this paper, for details we refer the reader to \cite[\S 1.3, \S 2]{Gr2}, \cite{Bou} or \cite[Appendix A]{Beu1}.

\begin{itemize}
	\item Let $X$ be a $\sigma$-compact locally compact space and $dx$ be a Radon measure on $X$. Let $f:X\to V$ be a continuous function. We say that $f$ is {\em weakly integrable} if for every $\ell\in V'$ the function $x\in X\mapsto \ell(f(x))\in \C$ is integrable. If $f$ is weakly integrable, an integral of $f$ is an element $I\in V$ such that
	\begin{equation}\label{eq integral TVS}
	\displaystyle \ell(I)=\int_X \ell(f(x))dx
	\end{equation}
	for every $\ell \in V'$. If an integral $I$ of $f$ exists it is unique and we will use the following notation for it
	$$\displaystyle I=\int_X f(x) dx.$$
	We also say that $f$ is {\em absolutely integrable} or that the integral $\int_X f(x) dx$ {\em converges absolutely} if for every continuous seminorm $\lVert .\rVert_V$ on $V$, the function $x\in X\mapsto \lVert f(x)\rVert_V$ is integrable. An absolutely integrable function is automatically weakly integrable. Moreover if $V$ is quasi-complete, any continuous absolutely integrable function $f$ admits an integral in $V$ (cf. \cite[\S 1.3]{Gr2}) and for $I\in V$ to be the integral of $f$ it suffices to check that \eqref{eq integral TVS} is satisfied for $\ell$ in some total subspace $H\subset V'$.
	
	\item Let $M$ be a real smooth manifold and $f:M\to V$ be a function. We say that $f$ is {\em smooth} if it admits derivatives of all order in a strong sense (see \cite[\S 5.3]{Bou}). If $f$ is smooth and $D$ is a differential operator on $M$, we denote by $Df$ the function $M\to V$ obtained by applying $D$ to $f$. When $V$ is quasi-complete, $f$ is smooth if and only if it is weakly smooth i.e. for every $\ell\in V'$ the function $m\in M\mapsto \ell(f(m))\in \C$ is smooth.
	
	\item Let $U\subset \C$ be an open subset and $f:U\to V$ be a function. We say that $f$ is {\em holomorphic} (or {\em analytic}) if for every $z_0\in U$ the limit $\lim\limits_{z\to z_0} (f(z)-f(z_0))(z-z_0)^{-1}$ exists in $V$. When $V$ is quasi-complete, $f$ is holomorphic if and only if it is weakly holomorphic i.e. for every $\ell\in V'$ the function $z\in U\mapsto \ell(f(z))\in \C$ is holomorphic. We also recall the following useful criterion of holomorphicity \cite[\S 3.3.1]{Bou}:
	
	\begin{num}
	\item\label{criterion holom} Assume that $V$ is quasi-complete. Then, a function $f:U\to V$ is holomorphic if and only if it is continuous and there exists a total subspace $H\subseteq V$ such that the functions $z\in U\mapsto \ell(f(z))$ are holomorphic for every $\ell\in H$.
	\end{num}
\end{itemize}

Recall that a {\em LF space} is a locally convex topological vector space $V$ that can be written as the direct limit of an at most countable family of Fr\'echet spaces. Actually, a LF space $V$ can always be written as an increasing union $V=\bigcup_n V_n$ where the $V_n$ are Fr\'echet spaces and the inclusions $V_n\subset V_{n+1}$ are continuous. We say that a LF space $V$ is {\em strict} if it can be written as before $V=\bigcup_n V_n$ with the extra property that the inclusions $V_n\subset V_{n+1}$ are topological embeddings (i.e. $V_{n+1}$ induces on $V_n$ its original topology). All topological vector spaces considered in this paper will be LF spaces and we recall here two basic properties of LF spaces which will be used thoroughly. First, since LF spaces are barreled \cite[Corollary 3 of Proposition 33.2]{Tr} they satisfy the uniform boundedness principle (aka Banach-Steinhaus theorem). Hence, if $V$ is an LF space, $W$ any topological vector space and $T_n:V\to W$ a sequence of continuous linear maps converging pointwise to $T:V\to W$, then the $(T_n)_n$ form an equicontinuous family so that in particular $T$ is continuous and $T_n$ converges to $T$ uniformly on compact subsets of $V$. Secondly, LF spaces also satisfy the closed graph theorem \cite[Theorem 4.B]{Gr1}: to verify that a linear map $T:V\to W$ between LF spaces is continuous it suffices to check that its graph is closed or even that the functionals $\lambda\circ T$ are continuous for $\lambda$ in some total subspace $H\subset W'$.

Fr\'echet spaces and more generally strict LF spaces are quasi-complete (and even complete) but a general LF space need not be. However, if $V$ is a LF space then its dual $V'$ (equipped with the weak topology) is quasi-complete: this follows from \cite[\S 34.3, Corollary 2]{Tr} as $V$ is barreled. In particular, if $U\subset \C$ is open a function $z\in U\mapsto \ell_z\in V'$ is holomorphic if and only if for every $v\in V$ the function $z\in U\mapsto \ell_z(v)$ is holomorphic.

In the case where $F$ is $p$-adic, and in order to make uniform statements, we will endow some vector spaces with their finest locally convex topology. More precisely, if $V$ is a $\C$-vector space then the {\em finest locally convex topology} on $V$ is the topology associated to the family of {\em all} seminorms on $V$. This definition makes clear that this topology has the following universal property: it is the unique structure of LCTVS on $V$ such that for every other LCTVS $W$ any linear mapping $V\to W$ is continuous. Alternatively, this topology can also be described as the locally convex inductive limit of the natural topologies on the finite dimensional subspaces of $V$. In particular, if $V$ has countable dimension (the only case we will encounter in this paper), then this turns $V$ into a LF space. Let us also emphasize that the finest locally convex topology is different from the discrete topology (for which the multiplication map $\C\times V\to V$ isn't continuous).

Let $V$ and $W$ be LCTVS. We denote by $V\ctens W$ the completed projective tensor product of $V$ and $W$ (see \cite[Chapter 43]{Tr}). It is a complete LCTVS containing $V\otimes W$ as a dense subspace and satisfying the following universal property \cite[Proposition 43.4]{Tr}: for every complete LCTVS $U$, restriction to $V\otimes W$ induces an (algebraic) isomorphism
\begin{align}\label{iso hom bil}
\displaystyle \Hom_{\cont}(V\ctens W,U)\simeq \Bil_{\cont}(V\times W,U)
\end{align}
where the left (resp. right) hand side denotes the space of continuous linear (resp. bilinear) maps from $V\ctens W$ (resp. from $V\times W$) to $U$. In particular, if $V$ is complete we have a canonical isomorphism $V\ctens \C\simeq V$ and for every continuous linear form $\lambda:W\to \C$, we get a continuous morphism $\Id\ctens \lambda:V\ctens W\to V$. The completed projective tensor product is associative (cf. \cite{Gr1} \S I.1.4) and for any finite family $V_1,\ldots,V_n$ of LCTVS we shall denote by $V_1\ctens\ldots\ctens V_n$ or $\displaystyle \widehat{\bigotimes_{1\leqslant i\leqslant n}}V_i$ their completed projective tensor product (in any order).

In the case where $V$, $W$ are of countable dimension and endowed with their finest locally convex topology, the (completed) projective tensor product $V\ctens W$ is equal to the algebraic tensor product $V\otimes W$ equipped with its finest locally convex topology. Indeed, using the various universal properties, we see that this statement is equivalent to the fact that any bilinear map $B: V\times W\to U$, where $U$ denotes a third LCTVS, is continuous. Let $\lVert .\rVert_U$ be a continuous seminorm on $U$ and $(v_n)_{n\geqslant 0}$, $(w_m)_{m\geqslant 0}$ be basis of $V$ and $W$ respectively. Setting $z_{n,m}=\lVert B(v_n,w_m)\rVert_U$, there exist two sequences $(x_n)_n$, $(y_m)_m$ of positive real numbers such that $z_{n,m}\leqslant x_n y_m$ for every $n$, $m$. Let $\lVert .\rVert_V$ and $\lVert .\rVert_W$ be the seminorms $\sum_n \lambda_n v_n\mapsto \sum_n \lvert \lambda_n\rvert x_n$ and $\sum_m \mu_m w_m\mapsto \sum_m \lvert \mu_m\rvert y_m$ on $V$ and $W$ respectively. Then, $\lVert .\rVert_V$, $\lVert .\rVert_W$ are continuous (because $V$ and $W$ are endowed with their finest locally convex topologies) and by construction $\lVert B(v,w)\rVert_U\leqslant \lVert v\rVert_V \lVert w\rVert_W$ for $(v,w)\in V\times W$. Hence, $B$ is continuous.

The following fact will be quite useful:

\begin{num}
\item\label{eq 1 proj tensor product} For $V$, $W$ Fr\'echet spaces and $U$ a complete LCTVS, a sequence of continuous linear maps $T_n:V\ctens W\to U$ that converges pointwise on $V\otimes W$ actually converges pointwise everywhere to a continuous linear map $T:V\ctens W\to U$.
\end{num}
Indeed, let $B_n:V\times W\to U$ be the continuous bilinear map corresponding to $T_n$ via \eqref{iso hom bil}. Then, by assumption the sequence $(B_n)_n$ converges pointwise which implies, by the uniform boundedness principle, that for every $v\in V$ and $w\in W$ the families of linear mappings $(B_n(v,.))_n$ and $(B_n(.,w))_n$ are equicontinuous. From \cite[Theorem 34.1]{Tr}, we deduce that the family $(B_n)_n$ is itself equicontinuous. As a consequence, the pointwise limit of the sequence $(B_n)_n$ is a continuous bilinear mapping $B:V\times W\to U$ and, by \cite[Exercise 43.1]{Tr}, the sequence $(T_n)_n$ is also equicontinuous and converges pointwise on $V\otimes W$ to the continuous linear map $T: V\ctens W\to U$ associated to $B$ via \eqref{iso hom bil} again. Since $V\otimes W$ is dense in $V\ctens W$, this readily implies that $(T_n)_n$ converges pointwise to $T$ everywhere.

By an easy induction, we deduce a property analogous to \eqref{eq 1 proj tensor product} for the projective completed tensor product of any finite family of Fr\'echet spaces that we shall also use without further notice.

\subsection{Space of functions}\label{Section space of functions}

For $X$ a smooth variety over $F$, we will denote by $\cS(X(F))$ the {\em Schwartz space} of $X(F)$ that is:
\begin{itemize}
	\item In the $p$-adic case, the space of all locally constant and compactly supported complex-valued functions on $X(F)$;
	
	\item In the Archimedean case ($F=\R$), the Schwartz space in the sense of Aizenbud and Gourevitch \cite{AG1} roughly described as the space smooth functions on $X(F)$ which are ``rapidly decreasing together with all their derivatives''.
\end{itemize}
When $X=G$ and $F=\R$, a function $f$ belongs to $\cS(G(F))$ if and only if for every $u,v\in \cU(g)$ and $N\geqslant 1$ we have
$$\displaystyle \lvert (R(u)L(v)f)(g)\rvert\ll \lVert g\rVert^{-N},\;\;\; g\in G(F).$$
On the other hand, for $V$ a finite dimensional $F$-vector space $\cS(V)$ coincides with the usual Schwartz space. 

When $F=\R$, $\cS(X(F))$ has a natural topology making it a Fr\'echet space. When $F$ is $p$-adic, we equip $\cS(X(F))$ with its finest locally convex topology (see Section \ref{Section Reminder on evtlcs}). Since in this case $\cS(X(F))$ is of countable dimension, this makes $\cS(X(F))$ into a LF space.

If $V$ is a finite dimensional $F$-vector space with a decomposition $V=W\oplus W'$ then there is a natural topological isomorphism 
\begin{align}\label{eq 2 proj tensor product}
\displaystyle \cS(W)\ctens \cS(W')\simeq \cS(V)
\end{align}
induced from the linear map $\varphi\otimes\psi\mapsto \left(w+w'\mapsto \varphi(w)\psi(w') \right)$ (see \cite[Theorem 51.6]{Tr}). More generally, for any finite decomposition $\displaystyle V=\bigoplus_{1\leqslant i\leqslant n} W_i$ there is a natural topological isomorphism
\begin{align}\label{eq 3 proj tensor product}
\displaystyle \widehat{\bigotimes_{1\leqslant i\leqslant n}}\cS(W_i)\simeq \cS(V).
\end{align}

Let $\Xi^G$ be the spherical Xi function of Harish-Chandra for $G$ (see \cite[\S II.1]{Wald1}, \cite[\S II.8.5]{Var}). It depends on the choice of a maximal compact subgroup of $G(F)$ but two such choices yield equivalent functions and we will only use $\Xi^G$ for estimates purposes. The space of {\em tempered functions} $\cC^w(G(F))$ on $G(F)$ is defined as follows. If $F$ is $p$-adic, it is the space of functions $f:G(F)\to \C$ which are biinvariant by a compact-open subgroup and satisfying an inequality
\begin{align}\label{eq 1 Cw}
\displaystyle \lvert f(g)\rvert\ll \Xi^G(g)\sigma(g)^d,\;\;\; \forall \;g\in G(F)
\end{align}
for some $d\geqslant 1$. In the Archimedean case, it is the space of smooth functions $f:G(F)\to \C$ for which there exists $d\geqslant 1$ such that for every $u,v\in \cU(\g)$ we have an inequality
\begin{align}\label{eq 2 Cw}
\displaystyle \lvert (R(u)L(v)f)(g)\rvert\ll \Xi^G(g)\sigma(g)^d,\;\;\; \forall \; g\in G(F)
\end{align}
For every integer $d\geqslant 1$ we denote by $\cC^w_d(G(F))$ the subspace of functions $f\in \cC^w(G(F))$ satisfying the estimate \eqref{eq 1 Cw} in the $p$-adic case (resp. \eqref{eq 2 Cw} for all $u,v\in \cU(\g)$ in the Archimedean case). Then $\cC^w_d(G(F))$ has a natural topology making it a strict LF space (and even a Fr\'echet space in the Archimedean case) see \cite[\S 1.5]{Beu1}. Moreover we have $\displaystyle \cC^w(G(F))=\bigcup_{d\geqslant 1} \cC^w_d(G(F))$ so that $\cC^w(G(F))$ is naturally equipped with a structure of tame LF space (in the sense of Section \ref{Section Reminder on evtlcs}). Moreover, $\cS(G(F))$ is dense in $\cC^w(G(F))$ for this topology.

Assume now that $G$ is quasi-split, let $B=TN$ be a Borel subgroup with unipotent radical $N$ and $\xi:N(F)\to \mathbb{S}^1$ a non-degenerate (aka generic) character. We define $\cS(N(F)\backslash G(F),\xi)$ as the space of functions $f: G(F)\to \C$ satisfying $f(ug)=\xi(u)f(g)$ for every $(u,g)\in N(F)\times G(F)$ and which are:
\begin{itemize}
	\item locally constant and compactly supported modulo $N(F)$ in the $p$-adic case;
	
	\item smooth and satisfying an inequality
	$$\displaystyle \lvert (R(u)f)(g)\rvert\ll \lVert g\rVert_{N\backslash G}^{-N},\;\;\; g\in N(F)\backslash G(F)$$
	for every $N\geqslant 1$ and $u\in \cU(\g)$ in the Archimedean case.
\end{itemize}
Then, $\cS(N(F)\backslash G(F),\xi)$ has a natural structure of Fr\'echet space when $F=\R$ whereas we equip this space with its finest locally convex topology when $F$ is $p$-adic (as it is of countable dimension, it is again a LF space).

Fix a maximal compact subgroup $K$ of $G(F)$. Using the Iwasawa decomposition $G(F)=N(F)T(F)K$ we define a function $\Xi^{N\backslash G}$ on $N(F)\backslash G(F)$ by
$$\displaystyle \Xi^{N\backslash G}(tk)=\delta_B(t)^{1/2},\;\;\; t\in T(F), k\in K.$$
Let $\cC^w(N(F)\backslash G(F),\xi)$ be the space of functions $f:G(F)\to \C$ such that $f(ug)=\xi(u)f(g)$ for every $(u,g)\in N(F)\times G(F)$ and satisfying the following condition:
\begin{itemize}
	\item if $F$ is $p$-adic, $f$ is right invariant by a compact-open subgroup and there exists $d>0$ such that
	\begin{equation}\label{eq1 CwWhitt}
	\displaystyle \lvert f(g)\rvert\ll \Xi^{N\backslash G}(g)\sigma_{N\backslash G}(g)^d,\;\;\; g\in N(F)\backslash G(F);
	\end{equation}
	
	\item if $F=\R$, $f$ is smooth and there exists $d>0$ such that for every $u\in \cU(\g)$ we have
	\begin{equation}\label{eq2 CwWhitt}
	\displaystyle \lvert (R(u)f)(g)\rvert\ll \Xi^{N\backslash G}(g)\sigma_{N\backslash G}(g)^d,\;\;\; g\in N(F)\backslash G(F).
	\end{equation}
\end{itemize}
As before, for every $d>0$ we denote by $\cC_d^w(N(F)\backslash G(F),\xi)$ the subspace of functions satisfying the estimate \eqref{eq1 CwWhitt} or \eqref{eq2 CwWhitt} for the given exponent $d$. When $F=\R$, $\cC_d^w(N(F)\backslash G(F),\xi)$ has a natural structure of Fr\'echet space whereas for $F$ $p$-adic, for every compact-open subgroup $J\subset G(F)$ the subspace $\cC_d^w(N(F)\backslash G(F),\xi)^J$ of $J$-invariant functions has a natural structure of Banach space. In both cases, this makes $\cC^w(N(F)\backslash G(F),\xi)$ into a LF space by writing it as $\cC^w(N(F)\backslash G(F),\xi)=\bigcup_{d>0} \cC^w_d(N(F)\backslash G(F),\xi)$ when $F=\R$ (resp. $\cC^w(N(F)\backslash G(F),\xi)=\bigcup_{J, d>0} \cC^w_d(N(F)\backslash G(F),\xi)^J$ in the $p$-adic case).

Consider two quasi-split groups $G_1$, $G_2$ with maximal unipotent subgroups $N_1$, $N_2$ and generic characters on their $F$-points $\xi_1$, $\xi_2$ respectively. Then, we have the following lemma which will be used without further notice in this paper.

\begin{lemdef}\label{lem 0 space of functions}
Let $\lambda$ be a continuous linear form on $\cC^w(N_1(F)\backslash G_1(F),\xi_1)$. Then, the linear map
$$\displaystyle \lambda\ctens \Id:\cC^w(N_1(F)\backslash G_1(F)\times N_2(F)\backslash G_2(F),\xi_1\boxtimes \xi_2)\to \cC^w(N_2(F)\backslash G_2(F),\xi_2)$$
$$\displaystyle W\mapsto \left(g\mapsto \lambda(W(.,g)) \right)$$
is well-defined and continuous. Moreover, for $\lambda_1$ and $\lambda_2$ two continuous linear forms on $\cC^w(N_1(F)\backslash G_1(F),\xi_1)$ and $\cC^w(N_2(F)\backslash G_2(F),\xi_2)$ respectively we have $\lambda_2\circ(\lambda_1\ctens \Id)=\lambda_1\circ(\Id\ctens \lambda_2)$ and the resulting linear form will be denoted by $\lambda_1\ctens \lambda_2$. 
\end{lemdef}

\noindent\ul{Proof}: We explain the Archimedean case since the $p$-adic case is similar and easier. For the first part, it suffices to notice that for every $d>0$ and $W\in \cC_d^w(N_1(F)\backslash G_1(F)\times N_2(F)\backslash G_2(F),\xi_1\boxtimes \xi_2)$ the function
$$\displaystyle F: g\in G_2(F)\mapsto W(.,g)\in \cC_d^w(N_1(F)\backslash G_1(F),\xi_1)$$
has the following properties:
\begin{itemize}
\item it is smooth (as a function valued in $\cC_d^w(N_1(F)\backslash G_1(F),\xi_1)$);
\item $F(ug)=\xi_2(u)F(g)$ for every $g\in G_2(F)$ and $u\in N_2(F)$;
\item for every $u\in \cU(\g)$ the set
$$\left\{\Xi^{N_2\backslash G_2}(g)^{-1}\sigma_{N\backslash G}(g)^{-d}(R(u)F)(g)\mid g\in G(F) \right\}$$
is bounded in $\cC_d^w(N_1(F)\backslash G_1(F),\xi_1)$.
\end{itemize}
Therefore the same holds for the composition of $F$ with $\lambda$ showing that the linear map $\lambda\ctens \Id$ is indeed taking values in $\cC^w(N_2(F)\backslash G_2(F),\xi_2)$. The continuity follows from the closed graph theorem as the linear forms $\cC^w(N_1(F)\backslash G_1(F)\times N_2(F)\backslash G_2(F),\xi_1\boxtimes \xi_2)\to \C$, $W\mapsto \lambda(W(.,g))$ are readily seen to be continuous for every $g\in G_2(F)$. Finally $\lambda_2\circ(\lambda_1\ctens \Id)=\lambda_1\circ(\Id\ctens \lambda_2)$ since these two continuous linear forms agree on the dense subspace $\cC_d^w(N_1(F)\backslash G_1(F),\xi_1)\otimes \cC_d^w(N_2(F)\backslash G_2(F),\xi_2)$. $\blacksquare$

\begin{rem}
The notation introduced in the previous lemma might seem to enter in conflict with the one introduced in Section \ref{Section Reminder on evtlcs}. However, we can offer the following justification: the bilinear map $\cC_d^w(N_1(F)\backslash G_1(F),\xi_1)\times \cC_d^w(N_2(F)\backslash G_2(F),\xi_2)\to \cC^w(N_1(F)\backslash G_1(F)\times N_2(F)\backslash G_2(F),\xi_1\boxtimes \xi_2)$, $(\varphi_1,\varphi_2)\mapsto \left((g_1,g_2)\mapsto \varphi_1(g_1)\varphi_2(g_2) \right)$ is continuous and hence induces a continuous linear map
\begin{equation}\label{eq1 explanation tensor product}
\displaystyle \cC_d^w(N_1(F)\backslash G_1(F),\xi_1)\ctens \cC_d^w(N_2(F)\backslash G_2(F),\xi_2)\to \cC^w(N_1(F)\backslash G_1(F)\times N_2(F)\backslash G_2(F),\xi_1\boxtimes \xi_2).
\end{equation}
The content of the above lemma can be reformulated by saying that for every continuous linear forms $\lambda$ on $\cC^w(N_1(F)\backslash G_1(F),\xi_1)$, the continuous linear map $\lambda\ctens Id: \cC_d^w(N_1(F)\backslash G_1(F),\xi_1)\ctens \cC_d^w(N_2(F)\backslash G_2(F),\xi_2) \to \cC^w(N_2(F)\backslash G_2(F),\xi_2)$ factorizes through \eqref{eq1 explanation tensor product}. Note also that this factorization is necessarily unique, as the image of \eqref{eq1 explanation tensor product} is dense. Finally, we point out, although we will not need this fact in the sequel, that the map \eqref{eq1 explanation tensor product} is never an isomorphism but that nevertheless there is a canonical isomorphism
$$\displaystyle \cC_d^w(N_1(F)\backslash G_1(F),\xi_1)\ctens_\epsilon \cC_d^w(N_2(F)\backslash G_2(F),\xi_2)\simeq \cC^w(N_1(F)\backslash G_1(F)\times N_2(F)\backslash G_2(F),\xi_1\boxtimes \xi_2)$$
where $\ctens_\epsilon$ denotes the the completed $\epsilon$-tensor product \cite[Definition 43.5]{Tr}.
\end{rem}

The next lemma is standard and can be proved in much the same way as \cite[Lemma 2.4.3]{Beu3}.

\begin{lem}\label{lem 2 space of functions}
Let $W\in \cC_d^w(N_n(E)\backslash G_n(E),\psi_n)$. Then, for all $N\geqslant 1$ we have
$$\displaystyle \lvert W(ak)\rvert \ll \prod_{i=1}^{n-1}(1+\left\lvert \frac{a_i}{a_{i+1}}\right\rvert_E)^{-N} \delta_{n,E}(a)^{1/2}\sigma(a)^d,\;\;\; a\in A_n(E), k\in K_{n,E}.$$
\end{lem}

This will be many times combined with the following basic convergence result for estimates purposes and so we record it here.

\begin{lem}\label{lem basic estimates}
For every $C>0$, there exists $N\geqslant 1$ such that for all $d>0$ and $s\in \C$ with $0<\Re(s)<C$ the integral
$$\displaystyle \int_{A_n(F)} \prod_{i=1}^{n-1} (1+\left\lvert \frac{a_i}{a_{i+1}}\right\rvert)^{-N} (1+\lvert a_n\rvert)^{-N} \sigma(a)^d \lvert \det a\rvert^s da$$
converges absolutely.
\end{lem}

\subsection{Measures}\label{Section measures}

We equip $F$ with the Haar measure $dx=d_{\psi'}x$ which is autodual with respect to $\psi'$. For any $n\geqslant 1$, $F^n$ will be equipped with the $n$-fold product of this measure. We endow $F^\times$ with the Haar measure $d^\times x=d^\times_{\psi'}x=\frac{dx}{\lvert x\rvert_F}$. Similarly, we equip $E$ with the Haar measure $dz=d_{\psi'_E}z$ autodual with respect to $\psi'_E=\psi'\circ \Tra_{E/F}$ and $E^\times$ with the measure $d^\times z=\frac{dz}{\lvert z\rvert_E}$.

Let $X$ be a smooth variety over $F$. Then to any volume form $\omega$ on $X$ (i.e. a differential form of maximal degree) we can associate a measure $\lvert \omega\rvert_{\psi'}$ on $V(F)$ as follows. Let $x\in V(F)$ and choose local coordinates $x_1,\ldots,x_n$ around $x$. Then, there exists a function $f$ on some open neighborhood of $x$ such that on that neighborhood $\omega=f(x_1,\ldots,x_n)dx_1\wedge \ldots\wedge dx_n$. Then, we define $\lvert \omega\rvert_{\psi'}$ on that neighborhood to be $\lvert f(x_1,\ldots,x_n)\rvert_F d_{\psi'}x_1\ldots d_{\psi'}x_n$. These locally defined measures can be glue together to give a global measure $\lvert \omega\rvert_{\psi'}$. This standard construction actually extends to any volume form $\omega$ on $X_{\overline{F}}$. Indeed, choose a finite extension $K$ of $F$ such that $\omega$ is defined over $X_K$. Then, in local coordinates $\omega$ can be written as $\omega=f(x_1,\ldots,x_n)dx_1\wedge \ldots\wedge dx_n$ where this time $f$ takes values in $K$ and we define $\lvert \omega\rvert_{\psi'}$ locally as $\lvert f(x_1,\ldots,x_n)\rvert_K^{1/[K:F]} d_{\psi'}x_1\ldots d_{\psi'}x_n$. Once again, we can glue to get a global measure.

Let $G$ be a connected reductive group over $F$. We equip $G(F)$ with a canonical measure $dg=d_{\psi'}g$ as follows. Let $G_{\bZ}$ be the split reductive group over $\bZ$ such that $G_{\bZ,\overline{F}}=G_{\bZ}\times_{\bZ} \overline{F}\simeq G_{\overline{F}}$. We fix such an isomorphism $\alpha$. Let $\omega_{G_{\bZ}}$ be a generator of the (free of rank $1$) $\bZ$-module of $G_{\bZ}$-invariant volume form on $G_{\bZ}$ and $\omega_{G_{\bZ},\overline{F}}$ be its base change to $G_{\bZ,\overline{F}}$. We set $\omega_{G_{\overline{F}}}:=\alpha_*(\omega_{G_{\bZ},\overline{F}})$. Then, $\omega_{G_{\overline{F}}}$ is a $G_{\overline{F}}$-invariant volume form on $G_{\overline{F}}$ which depends on the various choices only up to a root of unity so that the associated Haar measure $dg=\lvert \omega_{G_{\overline{F}}}\rvert_{\psi'}$ is independent of all choices (except of course $\psi'$). We also equip $\fg(F)$ with the Haar measure $dX=d_{\psi'}X$ defined by $dX:=\lvert \omega_{G_{\overline{F}},1}\rvert_{\psi'}$ where $\omega_{G_{\overline{F}},1}$ denotes the value of $\omega_{G_{\overline{F}}}$ at $1$ (a differential form of maximal degree on $\fg$). Notice that in the case where $G=\mathbb{G}_m$, the measure so defined coincides with the one we have already fixed on $F^\times$. More generally, we obtain on $G_n(F)$ the measure $dg=\lvert \det g\rvert^{-n} \prod_{i,j} dg_{i,j}$ (where $dg_{i,j}$ corresponds to the $\psi'$-autodual Haar measure on $F$) and on $G_n(E)$ the measure $dg=\lvert \det g\rvert_E^{-n} \prod_{i,j} dg_{i,j}$ (where this time $dg_{i,j}$ corresponds to the $\psi'_E$-autodual Haar measure on $E$).

The above construction can actually be applied to any algebraic linear group $G$ over $F$ once we chose a model over $\bZ$ for $G_{\overline{F}}$ (if it exists). For example, we will apply this to $G=N_n$, $U_n$ and $N'_n$ (see Section \ref{Section groups}) with their obvious models over $\bZ$ thus equipping the unipotent groups $N_n(F)$, $U_n(F)$, $N_n'(F)$ as well as their Lie algebras with Haar measures. These Haar measures are of course very easy to describe. For example on $N_n(F)$ we get the measure $du=\prod_{1\leqslant i < j\leqslant n}du_{i,j}$. We will also apply this construction to $N_n(E)$, $U_n(E)$, $N_n'(E)$ either by considering these groups as the $F$-points of $R_{E/F}N_{n,E}$, $R_{E/F}U_{n,E}$ and $R_{E/F}N'_{n,E}$ and using the natural isomorphisms $(R_{E/F}H_E)_{\overline{F}}\simeq H_{\overline{F}}\times H_{\overline{F}}$ ($H=N_n,U_n,N'_n$) or by doing the same as before with $F$ replaced by $E$ and $\psi'$ replaced by $\psi'_E$. Finally, the same can be applied to endow the Borel subgroups $B_n(F)$, $B_n(E)$ and the mirabolic subgroups $P_n(F)$, $P_n(E)$ with {\em right} Haar measures (to be denoted simply by $db$ and $dp$ respectively) using the natural model of $B_n$ and $P_n$ over $\bZ$. Notice that the Haar measure so obtained on $P_n(F)$ is, through the decomposition $P_n(F)=U_n(F)\rtimes G_{n-1}(F)$, the product of the Haar measures on $U_n(F)$ and $G_{n-1}(F)$. We also equip $P_n(F)\backslash G_n(F)$ with the ``twisted'' measure obtained as the quotient of the Haar measure $dg$ with the (left) Haar measure $\lvert \det p\rvert^{-1}dp$ on $P_n(F)$: then the integral $\displaystyle \int_{P_n(F)\backslash G_n(F)} . dg$ is a linear form on the space $C_c(P_n(F)\backslash G_n(F),\lvert \det\rvert)$ of functions $f:G_n(F)\to \mathbb{C}$ satisfying $f(pg)=\lvert \det p\rvert f(g)$ for every $(p,g)\in P_n(F)\times G_n(F)$ and which are of compact support modulo $P_n(F)$. We have the usual integration formula
$$\displaystyle \int_{G_n(F)} f(g)dg=\int_{P_n(F)\backslash G_n(F)}\int_{P_n(F)} f(pg) \lvert \det p\rvert^{-1}dpdg,\;\;\; f\in \cS(G_n(F)).$$
Moreover the measure $\lvert \det g\rvert dg$ on $P_n(F)\backslash G_n(F)$ (which is now a true measure, although not invariant) can be identified with the Haar measure $dx_1\ldots dx_n$ on $F^n$ through the isomorphism $P_n(F)\backslash G_n(F)\simeq F^n\setminus \{ 0\}$, $g\mapsto e_ng$. Using this, we readily check the following Fourier inversion formula
\begin{align}\label{eq 1 measures}
\displaystyle \int_{P_{n-1}(F)\backslash G_{n-1}(F)} \int_{U_n(E)} \varphi(v) \psi_n(hvh^{-1})^{-1} dv \lvert \det h\rvert dh=\lvert \tau\rvert_E^{(n-1)/2}\int_{U_n(F)} \varphi(v)dv,\;\;\; \varphi\in \cS(U_n(E))
\end{align}
(of course the left-hand side is only convergent as an iterated integral). Here we recall that $\tau\in E^\times$ is of trace zero and $\psi$ is given by $\psi(z)=\psi'(\Tra_{E/F}(\tau z))$.

If $G=A$ is a split torus, then we endow $i\mathcal{A}^*:=X^*(A)\otimes i\bR$ with the unique Haar measure giving the quotient $i\mathcal{A}^*_F:=i\cA^*/(\frac{2i\pi}{\log(q_F)})X^*(A)$ volume $1$ (recall that by convention $q_F=e^{1/2}$ in the Archimedean case). Let $\widehat{A(F)}$ be the unitary dual of $A(F)$ and $d_{\psi'}\chi$ be the Haar measure on $\widehat{A(F)}$ dual to the Haar measure we have just fixed on $A(F)$. Set
\begin{align}
\displaystyle d\chi:=\gamma^*(0,\mathbf{1}_F,\psi')^{-\dim(A)} d_{\psi'}\chi
\end{align}
where $\gamma^*(0,\mathbf{1}_F,\psi')$ is the ``regularized'' value at $0$ of $\gamma(s,\mathbf{1}_F,\psi')$ as defined in Section \ref{Section gamma factors} below. Then, $d\chi$ is a measure on $\widehat{A(F)}$ which is independent of the choice of $\psi'$ and moreover for this measure the local isomorphism
$$\displaystyle i\cA^* \to \widehat{A(F)}$$
$$\displaystyle \chi\otimes \lambda\mapsto \left(a\mapsto \lvert \chi(a)\rvert_F^{\lambda} \right)$$
is locally measure preserving.

\subsection{Representations}\label{Section representations}

In this paper, all representations will be tacitly assumed to be on complex vector spaces. By a {\em representation} of $G(F)$ we will always mean a smooth representation of finite length. Here {\em smooth} has the usual meaning in the $p$-adic case (i.e. every vector has an open stabilizer) whereas in the Archimedean case it means a smooth admissible Fr\'echet representation of moderate growth in the sense of Casselman-Wallach \cite{Cas}, \cite[Sect. 11]{Wall2} or, which is the same, an admissible SF representation in the sense of \cite{BK}. We shall always abuse notation and denote by the same letter a representation and the space on which it acts. In the Archimedean case this space is always coming with a topology (it is a Fr\'echet space) whereas in the $p$-adic case it will sometimes be convenient, in order to make uniform statements, to equip this space with its finest locally convex topology (see Section \ref{Section Reminder on evtlcs}, it then becomes a LF space).

Let $\pi$ be a representation of $G(F)$. We denote by $\pi^\vee$ the contragredient representation (aka smooth dual) and by $\langle .,.\rangle$ the natural pairing between $\pi$ and $\pi^\vee$. In the Archimedean case, $\pi^\vee$ can be described as the space of linear forms on $\pi$ which are continuous with respect {\em to every continuous} $G(F)$-continuous norm on $\pi$ together with the natural $G(F)$-action on it (a norm on $\pi$ is said to be {\em $G(F)$-continuous} if the action of $G(F)$ on $\pi$ is continuous for this norm). If $G=R_{E/F} G_{n,E}$ so that $G(F)=G_n(E)$ then we write $\pi^c$ for the composition of $\pi$ with the automorphism of $G_n(E)$ induced by $c$ and we set
$$\displaystyle \pi^*=(\pi^\vee)^c$$
For $\chi$ a continuous character of $G(F)$ we write $\pi\otimes \chi$ for the twist of $\pi$ by $\chi$ and for $\lambda\in \cA_{G,\C}^*$ by $\pi_\lambda$ for the twist of $\pi$ by the character $g\mapsto e^{\langle \lambda, H_G(g)\rangle}$. When $G=G_n$, we will also write $\pi_x=\pi\otimes \lvert \det \rvert^x$ for all $x\in \C$.

If $G_1$, $G_2$ are two reductive groups over $F$ and $\pi_1$, $\pi_2$ are representations of $G_1(F)$, $G_2(F)$ respectively, we denote by $\pi_1\boxtimes \pi_2$ the tensor product representation of $G_1(F)\times G_2(F)$ where in the Archimedean case this representation is realized on the {\em completed projective tensor product} of $\pi_1$ and $\pi_2$.

For $\pi$ a representation of $G(F)$, we let $\End^\infty(\pi)$ be the space of endomorphisms $T:\pi\to \pi$ which are biinvariant by a compact-open subgroup in the $p$-adic case and which are continuous with respect to {\em every continuous} $G(F)$-continuous norm in the Archimedean case. Then, $\End^\infty(\pi)$ is a smooth representation of finite length of $G(F)\times G(F)$ in the previous sense (in particular, it is a Fr\'echet space in the Archimedean case and a LF space equipped with its finest locally convex topology in the $p$-adic case) and the natural map $\pi^\vee\otimes \pi\to \End^\infty(\pi)$ extends to an isomorphism of $G(F)\times G(F)$-representations $\pi^\vee \boxtimes \pi\simeq \End^\infty(\pi)$. The canonical pairing $\pi^\vee\otimes \pi\to \C$ extends continuously to $\End^\infty(\pi)$ and we shall denote this extension by $T\mapsto \Tr(T)$. Moreover, for every function $f\in \cS(G(F))$ the expression
$$\displaystyle \pi(f)v=\int_{G(F)} f(g)\pi(g)vdg$$
converges absolutely in (the space of) $\pi$ for all $v\in \pi$ and defines an operator $\pi(f)\in \End^\infty(\pi)$. Moreover, the map $f\in \cS(G(F))\mapsto \pi(f)\in \End^\infty(\pi)$ is continuous.

For $P=MU$ a parabolic subgroup of $G$ and $\sigma$ a representation of $M(F)$ we write $i_P^G(\sigma)$ for the corresponding normalized  parabolically induced representation: it is the right regular representation of $G(F)$ on the space of smooth functions $e:G(F)\to \sigma$ satisfying $e(mug)=\delta_P(m)^{1/2}\sigma(m)e(g)$ for all $(m,u,g)\in M(F)\times U(F)\times G(F)$. When $G=G_n$ and $M$ is of the form $M=G_{n_1}\times\ldots\times G_{n_k}$ we write
$$\displaystyle \tau_1\times\ldots \times \tau_k$$
for $i_P^G(\tau_1\boxtimes\ldots\boxtimes \tau_k)$ where $\tau_i$ is a representation of $G_{n_i}(F)$ for every $1\leqslant i \leqslant k$. Similarly when $G=U(V)$ for $V$ a Hermitian space and $M$ is of the form $M=R_{E/F}G_{n_1,E}\times\ldots\times R_{E/F}G_{n_k,E}\times U(W)$ for $W$ a non-degenerate subspace of $V$, we write
$$\displaystyle \tau_1\times\ldots \times \tau_k \rtimes \sigma_0$$
for $i_P^G(\tau_1\boxtimes\ldots\boxtimes \tau_k\boxtimes \sigma_0)$ where $\tau_i$ is a representation of $G_{n_i}(E)$ for every $1\leqslant i \leqslant k$ and $\sigma_0$ a representation of $U(W)(F)$.

We denote by $\Irr(G)$ the set of all isomorphism classes of irreducible representations of $G(F)$. Of course in the Archimedean case irreducible should be understood as {\em topologically irreducible}. For $\pi\in \Irr(G)$, we denote by $\omega_\pi:Z_G(F)\to \C^\times$ its central character and in the Archimedean case by $\chi_\pi:\cZ(\g)\to \C$ its infinitesimal character.

Let $\Temp(G)\subset \Irr(G)$ be the subset of tempered irreducible representations and $\Pi_2(G)\subset \Temp(G)$ the further subset of square-integrable representations. For every $\pi\in \Pi_2(G)$ we define its {\em formal degree} $d(\pi)$ by the relation
$$\displaystyle \int_{G(F)/A_G(F)} \langle \pi(g)v_1,v_1^\vee\rangle \langle v_2,\pi^\vee(g)v_2^\vee \rangle dg=\frac{\langle v_1,v_2^\vee\rangle \langle v_2,v_1^\vee\rangle}{d(\pi)},\;\;\; v_1,v_2\in \pi, v_1^\vee,v_2^\vee\in \pi^\vee$$
Let $\Temp_{\ind}(G)$ be the set of isomorphism classes of representations of the form $i_P^G(\sigma)$ where $P=MU$ is a parabolic subgroup of $G$ and $\sigma\in \Pi_2(M)$. The isomorphism class of $i_P^G(\sigma)$ is independent of $P\in \mathcal{P}(M)$ and we shall thus write it as $i_M^G(\sigma)$. These representations are always semi-simple and in fact even unitarizable. According to Harish-Chandra, for $M$, $M'$ two Levi subgroups of $G$ and $\sigma\in \Pi_2(M)$, $\sigma'\in \Pi_2(M')$ the two representations $i_M^G(\sigma)$ and $i_{M'}^G(\sigma')$ have a constituent in common if and only if they are isomorphic and this happens precisely when there exists $g\in G(F)$ such that $gMg^{-1}=M'$, $g\sigma g^{-1}\simeq \sigma'$. Moreover, every representation of $\Temp(G)$ embeds in a unique representation of $\Temp_{\ind}(G)$ thus yielding a map
$$\displaystyle \Temp(G)\to \Temp_{\ind}(G).$$
When $G=G_n$, the representations $i_M^G(\sigma)$ ($M$ a Levi subgroup and $\sigma\in \Pi_2(M)$) are all irreducible so that $\Temp(G)=\Temp_{\ind}(G)$.

For $\pi\in \Temp(G)$, we have the following (see \cite[(2.2.5)]{Beu1}):
\begin{num}
\item\label{eq 1 representations} The assignment $T\mapsto \left(g\in G(F)\mapsto \Tr(\pi(g)T)\right)$ defines a continuous linear map $\End^\infty(\pi)\to \cC^w_0(G(F))$. In particular for every $(v,v^\vee)\in \pi\times \pi^\vee$, the function $g\mapsto \langle \pi(g)v,v^\vee\rangle$ belongs to $\cC^w_0(G(F))$ and the bilinear map $\pi\times \pi^\vee\to \cC^w(G(F))$, $(v,v^\vee)\mapsto \left(g\mapsto  \langle \pi(g)v,v^\vee\rangle\right)$, is continuous.
\end{num}

For $M$ a Levi subgroup of $G$ and $\sigma\in \Pi_2(M)$ we set
$$\displaystyle W(G,\sigma)=\{w\in W(G,M)\mid w\sigma\simeq \sigma \}$$
Then the map $\lambda\in i\cA_M^*\mapsto i_M^G(\sigma_\lambda)\in \Temp_{\ind}(G)$ is $W(G,\sigma)$-invariant and moreover there exists a unique topology on $\Temp_{\ind}(G)$ such that the induced maps $i\cA_M^*/W(G,\sigma)\to \Temp_{\ind}(G)$ (for all $M$ and $\sigma$) are local isomorphisms near $0$. When $F=\R$ and the central character of $\sigma$ is trivial on the connected component of $A_M(F)$ the map $i\cA_M^*/W(G,\sigma)\to \Temp_{\ind}(G)$ actually induces an isomorphism between $i\cA_M^*/W(G,\sigma)$ and a connected component of $\Temp_{\ind}(G)$. In general, connected components of $\Temp_{\ind}(G)$ are always of the form
$$\displaystyle \mathcal{O}=\{i_M^G(\sigma_\lambda)\mid \lambda\in i\cA_M^* \}$$
for some Levi $M$ and $\sigma\in \Pi_2(M)$ and in the $p$-adic case these components are all compacts. For $V$ a topological vector space, we will say that a function $f:\Temp_{\ind}(G)\to V$ is {\em smooth} if for every Levi subgroup $M$ of $G$ and every $\sigma\in \Pi_2(M)$ the function $\lambda\in i\cA_M^*\mapsto f(i_M^G(\sigma_\lambda))\in V$ is smooth.

Assume one moment that $G=U(V)$ for some Hermitian space $V$ over $E$. Let $M$ be a Levi subgroup of $G$ of the form
$$\displaystyle M=R_{E/F}G_{n_1,E}\times \ldots\times R_{E/F} G_{n_k,E}\times U(W)$$
where $W\subset V$ is a nondegenerate subspace and let
$$\displaystyle \sigma=\tau_1\boxtimes\ldots\boxtimes \tau_k\boxtimes \sigma_0$$
be a square-integrable representation of $M(F)$ where $\tau_i\in \Pi_2(G_{n_i}(E))$ for all $1\leqslant i\leqslant k$ and $\sigma_0\in \Pi_2(U(W))$. Recall that the group $W(G,M)$ can be identified with the subgroup of permutations $w\in\fS_{2k}$ preserving the partition $\{\{i,2k+1-i\}\mid 1\leqslant i\leqslant k \}$ and such that $n_{w(i)}=n_i$ for all $1\leqslant i \leqslant 2k$ where we have set $n_i=n_{2k+1-i}$ for every $k+1\leqslant i\leqslant 2k$. Then, using this identification $W(G,\sigma)$ is the subgroup of element $w\in W(G,M)$ such that $\tau_{w(i)}\simeq \tau_i$ for all $1\leqslant i\leqslant 2k$ where we have set $\tau_i=\tau_{2k+1-i}^*$ for every $k+1\leqslant i\leqslant 2k$.

Assume that $F=\R$. We define a norm $\pi\mapsto N(\pi)$ on $\Temp(G)$ as in \cite[\S 2.2]{Beu1} that is: fix a maximal torus $T\subset G$ and fix a $W(G_\C,T_\C)$-invariant norm on $\mathfrak{t}(\C)^*$, then identifying the infinitesimal character $\chi_\pi$ with an element of $\mathfrak{t}(\C)^*/W(G_\C,T_\C)$ by the Harish-Chandra isomorphism we set
$$\displaystyle N(\pi)=1+\lVert \chi_\pi\rVert.$$
Since $N(\pi)$ only depends on the infinitesimal character $\chi_\pi$ it also makes sense for $\pi\in \Temp_{\ind}(G)$. The same definition also makes sense for disconnected groups and will in particular be applied to maximal compact subgroups of $G(\R)$. Note that there exists $z\in \cZ(\g)$ such that
\begin{align}\label{eq1 representations}
\displaystyle N(\pi)\ll \chi_\pi(z),\;\;\; \pi\in \Temp_{\ind}(G)
\end{align}
Indeed if $z_1,\ldots,z_n$ is a generating family of homogeneous elements of $\cZ(\g)$ then we can take $z=1+\sum_{i=1}^n z_iz_i^*$ where $u\mapsto u^*$ denotes the conjugate-linear antiautomorphism of $\cU(\g)$ sending every $X\in \g(F)$ to $-X$.
Finally, if $M\subset G$ is a Levi subgroup and $\sigma\in \Temp(M)$, we will write $N(\sigma)$ for $N(i_M^G(\sigma))$.

\subsection{Spectral measures}\label{Section spectral measures}

We define two measures $d_{\psi'}\pi$ and $d\pi$ on $\Temp_{\ind}(G)$ as follows. Let $M$ be a Levi subgroup of $G$. Then, we define $d_{\psi'}\sigma$ (resp. $d\sigma$) to be the unique Borel measure on $\Pi_2(M)$ for which the local isomorphism
\begin{align}
\displaystyle \Pi_2(M)\to \widehat{A_M(F)} \\
\nonumber \sigma\mapsto \omega_{\sigma}\mid_{A_M(F)}
\end{align}
is locally measure preserving when $\widehat{A_M(F)}$ is equipped with the Haar measure $d_{\psi'}\chi$ (resp. $d\chi$). The induction map
\begin{align}\label{eq 3 Measures}
\displaystyle i_M^G:\Pi_2(M)\to \Temp_{\ind}(G) \\
\nonumber \sigma\mapsto i_M^G(\sigma)
\end{align}
is quasi-finite and proper with image the union of certain connected components of $\Temp_{\ind}(G)$. The restrictions of $d_{\psi'}\pi$ and $d\pi$ to this image are then defined to be $\lvert W(G,M)\rvert^{-1}i^G_{M*}d_{\psi'}\sigma$ and $\lvert W(G,M)\rvert^{-1}i^G_{M*}d\sigma$ respectively where $i^G_{M*}d_{\psi'}\sigma$ and $i^G_{M*}d\sigma$ stands for the push-forward of the measures $d_{\psi'}\sigma$ and $d\sigma$ respectively. 

Near a representation $\pi_0\in \Temp_{\ind}(G)$ the measure $d\pi$ can be described more explicitly as follows. Let $M$ be a Levi subgroup of $G$ and $\sigma \in \Pi_2(M)$ such that $\pi_0\simeq i_M^G(\sigma)$. Let $\mathcal{V}$ be a sufficiently small $W(G,\sigma)$-invariant open neighborhood of $0$ in $i\cA_M^*$ such that the map $\displaystyle \lambda\in i\cA_M^*\mapsto \pi_\lambda:=i_M^G(\sigma_\lambda)$
induces a topological isomorphism between $\mathcal{V}/W(G,\sigma)$ and an open neighborhood $\mathcal{U}$ of $\pi_0$ in $\Temp_{\ind}(G)$. Then, we have the integration formula
\begin{align}\label{eq 4 Measures}
\displaystyle \int_{\mathcal{U}}\varphi(\pi)d\pi=\frac{1}{\lvert W(G,\sigma)\rvert}\int_{\mathcal{V}}\varphi(\pi_\lambda) d\lambda,\;\;\; \varphi\in C_c(\mathcal{U})
\end{align}
As said before, when $F=\R$ and the central character of $\sigma$ is trivial on the connected component of $A_M(F)$ (which we can always be arranged up to twisting $\sigma$ by an unramified character) we can take $\mathcal{V}=i\cA_M^*$ and $\mathcal{U}=\mathcal{O}$ the connected component of $\pi_0$ in $\Temp_{\ind}(G)$.

We have the following basic estimates:
\begin{num}
\item\label{basic estimates spectral measure} There exists $k\geqslant 1$ such that the integral
$$\displaystyle \int_{\Temp_{\ind}(G)}N(\pi)^{-k}d\pi$$
converges.
\end{num} 

\subsection{Whittaker models}\label{Section Whittaker models}

Assume that $G$ is quasi-split, let $B=TN$ be a Borel subgroup of $G$ and $\xi: N(F)\to \mathbb{S}^1$ be a generic character. Recall that $\pi \in \Irr(G)$ is said to be $(N,\xi)$-generic if there exists a nonzero continuous linear form $\lambda: \pi\to \C$ satisfying $\lambda\circ \pi(u)=\xi(u)\lambda$ for all $u\in N(F)$. Such a linear form, which we call a {\em $(N,\xi)$-Whittaker functional}, is always unique up to a scalar (\cite[Theorem 9.2]{CHM}, \cite[Th\'eor\`eme 2]{Rod}). If $\pi$ is $(N,\xi)$-generic, we denote by $\cW(\pi,\xi)$ its  Whittaker model i.e. the space of all functions of the form $g\in G(F)\mapsto\lambda(\pi(g)v)$ where $\lambda$ is a nonzero $(N,\xi)$-Whittaker functional on $\pi$ and $v\in \pi$. Since the uniqueness of Whittaker functionals still holds for representations parabolically induced from irreducible ones (\cite[Theorem 9.1]{CHM}, \cite[Th\'eor\`eme 4]{Rod}), all of these definitions trivially extend to any $\pi\in \Temp_{\ind}(G)$. The following lemma is \cite[Lemma 15.7.3]{Wall2} in the Archimedean case and follows from \cite[Theorem 3.1]{LM} in the $p$-adic case (the second part of the lemma is easy to obtain from the closed graph theorem).

\begin{lem}\label{lem 1 space of functions}
Let $\pi\in\Temp_{\ind}(G)$ which is $(N,\xi)$-generic. Then, we have $\cW(\pi,\xi)\subset \cC^w(N(F)\backslash G(F),\xi)$. Moreover if $v\mapsto W_v$ is an isomorphism between $\pi$ and its Whittaker model, the linear map $v\in \pi\mapsto W_v\in \cC^w(N(F)\backslash G(F),\xi)$ is continuous.
\end{lem}

\noindent Assume now that $G=R_{E/F}G_n$. We define a scalar product $(.,.)^{\Whitt}$ on $\cC^w(N_n(E)\backslash G_n(E),\psi_n)$ by
$$\displaystyle (W,W')^{\Whitt}=\int_{N_n(E)\backslash P_n(E)} W(p) \overline{W'(p)} dp,\;\;\; W,W'\in \cC^w(N_n(E)\backslash G_n(E),\psi_n)$$

This integral is readily seen to be absolutely convergent by a combination of Lemma \ref{lem 2 space of functions} and Lemma \ref{lem basic estimates}. The following result is due to Bernstein \cite{Ber1} in the $p$-adic case and to Sahi \cite{Sah} and Baruch \cite{Bar} in the Archimedean case.

\begin{theo}
For $\pi\in \Temp(G_n(E))$, the restriction of $(.,.)^{\Whitt}$ to $\cW(\pi,\psi_n)$ is $G_n(E)$-invariant.
\end{theo}

\subsection{Schwartz functions on $\Temp_{\ind}(G)$}\label{Section Schwartz functions on Temp(G)}

Let $V=\bigcup_n V_n$ be a LF space where $(V_n)_n$ is an increasing sequence of Fr\'echet spaces with continuous inclusions. We say that a function $f:\Temp_{\ind}(G)\to V$ is {\em Schwartz} if it is smooth (in the sense of Section \ref{Section representations}) and satisfies:
\begin{num}
	\item\label{eq defn Schwartz tempered} In the Archimedean case, $f$ {\em and all its derivatives are of rapid decay} that is: for every $k\geqslant 0$ there exists $n>0$ such that for every Levi subgroup $M$ of $G$, every $d>0$, every continuous seminorm $\lVert .\rVert_{V_n}$ on $V_n$ and every $D\in \Sym^k(\cA^*_{M,\C})$ that we see as a differential operator of order $k$ with constant coefficients on $i\cA_M^*$ we have
	\begin{equation*}
	\displaystyle \left\lVert D(\lambda\in i\cA_M^*\mapsto f(i_M^G(\sigma_\lambda)))_{\lambda=0}\right\rVert_{V_n} \ll N(\sigma)^{-d},\; \mbox{ for } \sigma\in \Pi_2(M).
	\end{equation*}
	
	\item\label{eq2 defn Schwartz tempered} In the $p$-adic case: $f$ has compact support and for every $k\geqslant 0$ there exists $n>0$ such that for every Levi subgroup $M$ of $G$, every $D\in \Sym^k(\cA^*_{M,\C})$ and every $\sigma\in \Pi_2(M)$, the map $\lambda\in i\cA_M^*\mapsto Df(i_M^G(\sigma_\lambda))$ takes values in $V_n$.
\end{num}
By \cite[corollaire 3, Intro IV, p.17]{Gr1}, two presentations $V=\bigcup_n V_n^1$ and $V=\bigcup_n V^2_n$ are cofinal so that the above definition is independent of the choice of the sequence of Fr\'echet spaces $(V_n)_n$. If $W=\bigcup_m W_m$ is another LF space and $T:V\to W$ a continuous linear mapping then for every $n>0$ there exists $m>0$ such that $T$ induces a continuous map $V_n\to W_m$ \cite[th\'eor\`eme A p.16]{Gr1}. Therefore, if $f:\Temp_{\ind}(G)\to V$ is a Schwartz function in the above sense then $T\circ f: \Temp_{\ind}(G)\to W$ is also a Schwartz function. Moreover, by \eqref{basic estimates spectral measure} any Schwartz function $f:\Temp_{\ind}(G)\to V$ is absolutely integrable in some $V_n$ with respect to the measure $d\pi$ and the integral $\int_{\Temp_{\ind}(G)} f(\pi) d\pi$ exists in $V$ (as $V_n$ is a Fr\'echet space hence quasi-complete).

When $V=\C$, we denote by $\cS(\Temp(G))$ the space of Schwartz functions $\Temp_{\ind}(G)\to V$. In the Archimedean case, we equip $\cS(\Temp(G))$ with the topology associated to the family of seminorms
$$\displaystyle \lVert f\rVert_{M,d,D}=\sup_{\sigma\in \Pi_2(M)} N(\sigma)^{d} \left\lvert D(\lambda\in i\cA_M^*\mapsto f(i_M^G(\sigma_\lambda)))_{\lambda=0}\right\rvert$$
where $M$ is a Levi subgroup of $G$, $d>0$ and $D\in \Sym^\bullet(\cA^*_{M,\C})$. It is then a Fr\'echet space. In the $p$-adic case, for $\cO\subset \Temp_{\ind}(G)$ a finite union of connected components we denote by $\cS_{\cO}(\Temp(G))$ the subspace of Schwartz functions supported on $\cO$ and we equip this subspace with the topology associated to the family of seminorms
$$\displaystyle \lVert f\rVert_{M,D}=\sup_{\sigma\in \Pi_2(M)} \left\lvert D(\lambda\in i\cA_M^*\mapsto f(i_M^G(\sigma_\lambda)))_{\lambda=0}\right\rvert$$
where $M$ is a Levi subgroup of $G$ and $D\in \Sym^\bullet(\cA^*_{M,\C})$. This makes $\cS_{\cO}(\Temp(G))$ into a Fr\'echet space and $\cS(\Temp(G))=\bigcup_{\cO}\cS_{\cO}(\Temp(G))$ into a strict LF space.

We write $\cS_c(\Temp(G))$ for the subspace of functions $f\in \cS(\Temp(G))$ with compact support. Thus, in the $p$-adic case we have $\cS_c(\Temp(G))=\cS(\Temp(G))$. In the Archimedean case on the other hand, $\cS_c(\Temp(G))$ is dense in $\cS(\Temp(G))$.

We will need the following lemma.

\begin{lem}\label{lem holom Schwartz functions}
Let $V$ be a LF space, $f:\Temp_{\ind}(G)\to V$ be a Schwartz function and $z\in U\mapsto \ell_z\in V'$ be a holomorphic map where $U\subset \C$ is open. Then, the map
\begin{equation}\label{eq holom map}
\displaystyle z\in U\mapsto \ell_z\circ f\in \cS(\Temp(G))
\end{equation}
is holomorphic.
\end{lem}

\noindent\ul{Proof}: We prove the lemma in the Archimedean case, the $p$-adic case being similar. Since $z\in U\mapsto \ell_z\in V'$ is holomorphic, for every $\pi\in \Temp_{\ind}(G)$, the function $z\in U\mapsto \ell_z(f(\pi))$ is holomorphic. Therefore, by the criterion \eqref{criterion holom} applied to the total subspace $H\subset \cS(\Temp(G))$ generated by the functionals of ``evaluation at a point $\pi\in \Temp_{\ind}(G)$'', it suffices to check that the map \eqref{eq holom map} is continuous. By definition of the topology on $\cS(\Temp(G))$, we need to show that for every Levi subgroup $M\subset G$, $d>0$, $D\in \Sym^\bullet(\cA^*_{M,\C})$ and $z_0\in U$, we have
\begin{equation}\label{eq2 holom map}
\displaystyle \lim\limits_{z\to z_0} \lVert \ell_z(f)-\ell_{z_0}(f)\rVert_{M,d,D}=0.
\end{equation}
The left hand side can be rewritten as
$$\displaystyle \sup_{\sigma\in \Pi_2(M)} N(\sigma)^{d} \left\lvert (\ell_z-\ell_{z_0})(Df(\sigma))\right\rvert$$
where we have set $Df(\sigma)=D(\lambda\in i\cA_M^*\mapsto f(i_M^G(\sigma_\lambda)))_{\lambda=0}\in V$. Note that by definition of a Schwartz on $\Temp_{\ind}(G)$, for every continuous seminorm $\lVert .\rVert_V$ on $V$ the map $\sigma\in \Pi_2(M)\mapsto N(\sigma)^{d} \lVert Df(\sigma)\rVert_V$ is bounded up to a constant by $N(\sigma)^{-1}$ and therefore converges to $0$ at infinity. This entails in particular that the subset
$$\displaystyle \left\{ N(\sigma)^{d} Df(\sigma)\mid \sigma\in \Pi_2(M) \right\}\cup \{ 0\}$$
of $V$ is compact. Moreover, by the uniform boundedness principle, $\ell_z$ converges to $\ell_{z_0}$ as $z\to z_0$ uniformly on compact subsets. The convergence \eqref{eq2 holom map} follows. $\blacksquare$

\subsection{Local Langlands correspondences and base-change}\label{Section LLC}

Let $W_F$ be the Weil group of $F$ and
$$\displaystyle W'_F=\left\{
    \begin{array}{ll}
        W_F\times SL_2(\C) & \mbox{ if } F \mbox{ is } p-\mbox{adic} \\
        W_F & \mbox{ if } F \mbox{ is  Archimedean}
    \end{array}
\right.
$$
the Weil-Deligne group of $F$.

For $G$ a connected reductive group over $F$, we will denote by ${}^L G=G^\vee \rtimes W_F$ the Weil form of the $L$-group of $G$ where $G^\vee$ denotes the Langlands dual of $G$ (a connected reductive complex group) and the action of $W_F$ on $G^\vee$ is by pinned automorphisms. We will actually abuse this notation slightly and denote by ${}^L G_n(E)$ and ${}^L \overline{G_n(E)}$ the $L$-groups of the natural connected reductive groups over $F$ of which $G_n(E)$ and $\overline{G_n(E)}$ are the sets of $F$-points. Thus, we have
$$\displaystyle {}^L G_n(E)=\left(\GL_n(\C)\times \GL_n(\C) \right)\rtimes W_F$$
where the action of $W_F$ on $\GL_n(\C)\times \GL_n(\C)$ is given by
$$\displaystyle \sigma\cdot (g_1,g_2)=\left\{
    \begin{array}{ll}
        (g_2,g_1) & \mbox{ if } \sigma\in W_F\setminus W_E \\
        (g_1,g_2) & \mbox{ otherwise}
    \end{array}
\right.
$$
and ${}^L \overline{G_n(E)}$ is the subgroup of elements $(g_1,g_2,\sigma)\in {}^L G_n(E)$ such that $\det(g_1)\det (g_2)=1$. If $V$ is a Hermitian space of dimension $n$ over $E$, we have ${}^L U(V)=\GL_n(\C)\rtimes W_F$ where the action of $W_F$ on $\GL_n(\C)$ is given by
$$\displaystyle \sigma\cdot g=\left\{
    \begin{array}{ll}
        J_n{}^tgJ_n^{-1} & \mbox{ if } \sigma\in W_F\setminus W_E \\
        g & \mbox{ otherwise}
    \end{array}
\right.
$$
with $J_n=\begin{pmatrix} & & 1 \\ & \iddots & \\ (-1)^{n-1} & & \end{pmatrix}$.

Recall that a {\em Langlands parameter} for $G$ (or by abuse of language for $G(F)$) is a $G^\vee$-conjugacy class of continuous homomorphisms: $\varphi: W_F'\to {}^L G$ which are algebraic when restricted to $SL_2(\C)$ (in the $p$-adic case), sending $W_F$ to semi-simple elements and commuting with the natural projections $W'_F\to W_F$, ${}^L G\to W_F$. A Langlands parameter $\varphi$ is {\em tempered} (resp. {\em discrete}) if the projection of $\varphi(W_F)$ in $G^\vee$ is bounded (resp. if it is tempered and the centralizer of $\varphi(W'_F)$ in $G^\vee$ is finite modulo $Z(G^\vee)^{W_F}$). We shall denote by $\Phi(G)$, resp. $\Phi_{\temp}(G)$, resp. $\Phi_2(G)$, the set of all, resp. all tempered, resp. all discrete, Langlands parameters for $G$.

The local Langlands correspondence postulates the existence of a finite-to-one map $\Irr(G)\to \Phi(G)$, $\pi\mapsto \varphi_\pi$ satisfying certain properties and whose fibers are usually called {\em $L$-packets} (of $G$). This correspondence has been established in some cases: for any real reductive group by Langlands \cite{La} and in the $p$-adic case for $G_n(F)$ by Harris-Taylor \cite{HT}, Henniart \cite{Hen} and Scholze \cite{Sch} and more recently for unitary groups $U(V)$ by Mok \cite{Mok} and Kaletha-Minguez-Shin-White \cite{KMSW} (following the work of Arthur \cite{Art3} for split orthogonal and symplectic groups). These correspondences naturally extend to product of such groups, including in particular all Levi subgroups of groups in the previous list, as well as groups obtained by extension of scalars or quotient by a split central torus thus including $G_n(E)$ and $\overline{G_n(E)}$ in our list. The basic properties of these correspondences that we shall use in this paper are listed below (where $M$ denotes a Levi subgroup of $G$):
\begin{itemize}
\item It sends $\Temp(G)$ (resp. $\Pi_2(G)$) to $\Phi_{\temp}(G)$ (resp. $\Phi_2(G)$). Moreover the map $\Temp(G)\to \Phi_{\temp}(G)$ factorizes through the quotient $\Temp_{\ind}(G)$;
\item If $G$ is quasi-split then $\Irr(G)\to \Phi(G)$, $\Temp_{\ind}(G)\to \Phi_{\temp}(G)$ and $\Pi_2(G)\to \Phi_2(G)$ are all surjective (in particular this applies to $U(n)$) and these are even bijections for $G_n(F)$, $G_n(E)$ and $\overline{G_n(E)}$.
\item For $\sigma\in \Temp(M(F))$ and $\pi$ a constituent of $i_M^G(\sigma)$ we have $\varphi_{\pi}=\iota\circ \varphi_\sigma$ where $\iota:{}^L M\to {}^L G$ is the natural embedding between $L$-groups. Conversely, if $\varphi_{\pi}=\iota\circ \varphi_M$ where $\varphi_M\in \Phi_2(M)$ then there exists $\sigma\in \Pi_2(M)$ with $\varphi_\sigma=\varphi_M$ such that $\pi\hookrightarrow i_M^G(\sigma)$.
\item For $\pi\in \Irr(G)$ and $\lambda\in \cA_{G,\C}^*=X^*(G)\otimes \C$, we have $\varphi_{\pi_\lambda}=\varphi \lvert .\rvert^\lambda$ where $\lvert .\rvert^\lambda$ is the character $W'_F\to (Z(G^\vee)^{W_F})^0=X^*(G)\otimes \C^\times$ obtained by composing the absolute value $\lvert .\rvert$ on $W'_F$ with the character $t\in \R_+^*\mapsto t^\lambda \in X^*(G)\otimes \C^\times$;
\item If $G=G_1\times G_2$ and $\pi=\pi_1\boxtimes \pi_2\in \Irr(G)$ then the Langlands parameter $\varphi$ of $\pi$ is the {\em product} in a suitable sense of the Langlands parameters $\varphi_1$ and $\varphi_2$ of $\pi_1$ and $\pi_2$.
\end{itemize}
We will denote by $\Temp(G)/\stab$ the quotient of $\Temp(G)$ by the relation $\pi\sim_{\stab} \pi'\Leftrightarrow \varphi_\pi=\varphi_{\pi'}$ (i.e. the set of tempered $L$-packets for $G$). This notation will actually only be used for unitary groups, as for other groups considered in this paper (general linear groups and their variants) we have $\Temp(G)/\stab=\Temp(G)$. Note that we have natural surjections $\Temp(G)\to \Temp(G)/\stab$ and $\Temp_{\ind}(G)\to \Temp(G)/\stab$ (by the first point above). If $G_{\qs}$ is a quasi-split inner form of $G$ (e.g. $G$ is a unitary groups of rank $n$ and $G_{\qs}=U(n)$), by the second point above we have a natural identification $\Temp(G_{\qs})/\stab=\Phi_{\temp}(G_{\qs})$ and since $G$ and $G_{\qs}$ share the same $L$-group $\Phi_{\temp}(G_{\qs})=\Phi_{\temp}(G)$ and the Langlands correspondence gives a map $\Temp_{\ind}(G)\to \Temp(G_{\qs})/\stab$.

\begin{lem}\label{lem 0 L parameters, LLC, basechange}
Assume that $G$ is a product of unitary groups. Then, there exists a unique topology on $\Temp(G_{\qs})/\stab$ such that $\Temp_{\ind}(G_{\qs})\to \Temp(G_{\qs})/\stab$ is a local isomorphism and for every connected component $\cO\subset \Temp_{\ind}(G)$ the map $\cO\to \Temp(G_{\qs})/\stab$ induces an isomorphism between $\cO$ and a connected component of $\Temp(G_{\qs})/\stab$.
\end{lem}

\noindent\ul{Proof}: As the Langlands correspondence is compatible with products (cf.\ last point above), we may assume that $G=U(V)$ for $V$ a Hermitian space of dimension $n$ over $E$ and $G_{\qs}=U(n)$. Let $\pi\in \Temp_{\ind}(G)$ and $\cO\subset \Temp_{\ind}(G)$ be its connected component. We can write $\pi=i_M^G(\sigma)$ where $M$ is a Levi subgroup of the form $M=R_{E/F}G_{n_1}\times \ldots\times R_{E/F} G_{n_k}\times U(W)$, for a certain non-degenerate subspace $W\subset V$ of dimension $m$, and $\sigma=\tau_1\boxtimes\ldots\boxtimes \tau_k\boxtimes \sigma_0\in \Pi_2(M)$. Let $\pi'\in \Temp_{\ind}(U(n))$ be such that $\varphi_{\pi'}=\varphi_{\pi}$ and $\cO'\subset \Temp_{\ind}(U(n))$ be its connected component. Then, by the compatibility of LLC with parabolic induction (third point above) and with products, we have $\pi'=i_L^{U(n)}(\sigma')$ where $L$ is a Levi subgroup (of $U(n)$) of the form $L=R_{E/F}G_{n_1}\times \ldots\times R_{E/F} G_{n_k}\times U(m)$ and $\sigma'=\tau_1\boxtimes\ldots\boxtimes \tau_k\boxtimes \sigma'_0$ for some $\sigma_0'\in \Pi_2(U(m))$. There exist identifications $i\cA_M^*\simeq (i\R)^k\simeq i\cA_L^*$ such that $\pi_\lambda:=i_M^G(\sigma_\lambda)=\tau_{1,\lambda_1}\times \ldots \times \tau_{k,\lambda_k}\rtimes \sigma_0$ and $\pi'_\lambda:=i_L^{U(n)}(\sigma'_\lambda)=\tau_{1,\lambda_1}\times \ldots \times \tau_{k,\lambda_k}\rtimes \sigma'_0$ for every $\lambda\in (i \R)^k$. By the compatibility of LLC with unramified twists (fourth point above) and with parabolic induction we see that $\varphi_{\pi_\lambda}=\varphi_{\pi'_\lambda}$ for every $\lambda\in (i \R)^k$. Moreover, by the precise descriptions of $W(G,\sigma)$ and $W(U(n),\sigma')$ given in Section \ref{Section representations}, there is an isomorphism $W(G,\sigma)\simeq W(U(n),\sigma')$ compatible with the previous identification. Therefore, we have a commutative diagram
$$\displaystyle \xymatrix{i\cA_M^*/W(G,\sigma) \eq[rr] \ar[d] & & i\cA_L^*/W(U(n),\sigma') \ar[d] \\ \cO \ar[rd] & & \cO' \ar[ld] \\ & \Temp(G_{\qs})/\stab & }$$
where the two vertical arrows, given by $\lambda\mapsto \pi_\lambda$ and $\lambda\mapsto \pi'_\lambda$, are surjective and local isomorphisms near $0$ (by definition of the topologies on $\Temp_{\ind}(G)$ and $\Temp_{\ind}(U(n))$). 

Using the above diagram for $G=G_{\qs}$ and $\pi$, $\pi'$ any two tempered representations lying in the same $L$-packet, we conclude that $\Temp(G_{\qs})/\stab$ has a unique topology such that $\Temp_{\ind}(G_{\qs})\to \Temp(G_{\qs})/\stab$ is a local isomorphism. Moreover, the same diagram shows that for any $G$ the map $\Temp_{\ind}(G)\to \Temp(G_{\qs})/\stab$ is a local isomorphism and that the image of any two connected components $\cO_1,\cO_2\subset \Temp_{\ind}(G)$ in $\Temp(G_{\qs})/\stab$ are either disjoint or identical.

Thus, to conclude the proof it only remains to show that $\cO\to \Temp(U(n))/stab$ is injective. If $\varphi_{\pi_\lambda}=\varphi_{\pi_\mu}$ for some $\lambda,\mu \in i\cA_M^*$, by the compatibility of LLC with parabolic induction and products we have $\pi_\mu\simeq\tau_{1,\lambda_1}\times \ldots \times \tau_{k,\lambda_k}\rtimes \sigma''_0$ for some $\sigma''_0\in \Pi_2(U(W))$. According to Harish-Chandra, this implies that the representations $\tau_{1,\mu_1}\boxtimes \ldots \boxtimes \tau_{k,\mu_k}\boxtimes \sigma_0$ and $\tau_{1,\lambda_1}\boxtimes \ldots \boxtimes \tau_{k,\lambda_k}\boxtimes \sigma''_0$ of $M(F)$ are conjugated under $W(G,M)$ hence $\sigma_0''=\sigma_0$ and finally $\pi_\mu=\pi_\lambda$. $\blacksquare$

Still assuming that $G$ is a product of unitary groups, by the previous lemma there exists a unique Borel measure on $\Temp(G_{\qs})/\stab$ for which the map $\Temp_{\ind}(G_{\qs})\to \Temp(G_{\qs})/\stab$ is locally measure preserving where we equip $\Temp_{\ind}(G_{\qs})$ with the measure $d\pi$ defined in Section \ref{Section representations} and we shall also denote this measure by $d\pi$. Arguing similarly as in the proof of previous lemma, we see that the map $\Temp_{\ind}(G)\to \Temp(G_{\qs})/\stab$ is also locally measure preserving. Therefore, we have the integration formulas
\begin{align}\label{eq 4ter Measures}
\displaystyle \int_{\Temp_{\ind}(G)}f(\pi)d\pi=\int_{\Temp(G_{\qs})/\stab}\sum_{\substack{\pi'\in \Temp_{\ind}(G) \\ \varphi_{\pi'}=\varphi_{\pi}}}f(\pi') d\pi,\;\;\; f\in C_c(\Temp_{\ind}(G))
\end{align}
\begin{align}\label{eq 4bis Measures}
\displaystyle \int_{\mathcal{U}/\stab}\varphi(\pi)d\pi=\frac{1}{\lvert W(G,\sigma)\rvert}\int_{\mathcal{V}}\varphi(\pi_\lambda) d\lambda,\;\;\; \varphi\in C_c(\mathcal{U}/\stab)
\end{align}
where $\pi$, $\sigma$, $\cU$ and $\cV$ are as in \eqref{eq 4 Measures} and $\cU/\stab$ denotes the image of $\cU$ in $\Temp(G_{qs})/\stab$ (which by Lemma \ref{lem 0 L parameters, LLC, basechange} is isomorphic to $\cU$).

For $V$ be a Hermitian space of dimension $n$ over $E$ we define a morphism $BC: {}^L U(V)\to {}^L G_n(E)$ called {\em base-change} by $(g,\sigma)\mapsto (g,J_n {}^tg^{-1}J_n^{-1},\sigma)$. This induces a map between set of Langlands parameters $\Phi_{\temp}(U(V))\to \Phi_{\temp}(G_n(E))$ and thus composing with the Langlands correspondence a map $\Temp(U(V))/\stab\to \Temp(G_n(E))$ which we shall also denote by $BC$. This map is injective and has its image contained in the subset of representations $\pi\in\Temp(G_n(E))$ satisfying $\pi\simeq \pi^*$. The base-change map satisfies the following basic properties:
\begin{num}
\item\label{eq -2 L parameters, LLC, basechange} $\omega_{BC(\sigma)}(z)=\omega_\sigma(z/z^c)$ for all $\sigma\in \Temp(U(V))$ and all $z\in E^\times$.
\item\label{eq -1 L parameters, LLC, basechange} $BC(\tau\rtimes \sigma)=\tau\rtimes BC(\sigma)$ for all $\tau\in \Temp(G_m(E))$ and $\sigma\in \Temp(U(W))$ where $W$ is a certain Hermitian space with $\dim(W)+2m=n$. 
\end{num}
At some point it will be convenient to consider certain twists of this base-change map. Recall that we have defined a character $\eta'_n$ of $G_n(E)$ in Section \ref{Section groups}. We set $BC_n(\sigma)=BC(\sigma)\otimes \eta'_n$ for all $\sigma\in \Temp(U(V))/\stab$. When $n$ is odd $BC_n=BC$ and when $n$ is even $BC_n$ is usually called {\em unstable base-change}. The following property characterizes the image of $BC_n$ from a quasi-split unitary group:
\begin{num}
\item\label{eq 0 L parameters, LLC, basechange} Let $\pi\in \Temp(G_n(E))$. Then, $\pi$ belongs to the image of $BC_n:\Temp(U(n))\to \Temp(G_n(E))$ if and only if it can be written as
$$\displaystyle \pi=\left(\bigtimes_{i=1}^k\tau_i \times \tau_i^*\right) \times \bigtimes_{j=1}^\ell \mu_j$$
where for every $1\leqslant i\leqslant k$, $\tau_i\in \Pi_2(G_{n_i}(E))$ for some $n_i\geqslant 1$ and for every $1\leqslant j\leqslant \ell$, $\mu_j\in BC_{m_j}(\Temp(U(m_j)))\cap \Pi_2(G_{m_j}(E))$ for some $m_j\geqslant 1$.
\end{num}

\subsection{Groups of centralizers}\label{Section Groups of centralizers}

Let $G$ be a general linear groups or one of its variant (e.g. such that $G(F)=G_n(E)$ or $\overline{G_n(E)}$) or a product of unitary groups. Let $\pi\in \Temp(G)$ and choose a Levi subgroup $M$ of $G$ and $\sigma\in \Pi_2(M)$ such that $\pi\hookrightarrow i_{M}^G(\sigma)$. We set $S_\pi:=S_{\varphi_{\sigma}}$ where $\varphi_{\sigma}:W'_F\to {}^L M$ is the Langlands parameter of $\sigma$ and $S_{\varphi_{\sigma}}$ denotes the centralizer of $\varphi_{\sigma}$ in $(M/A_M)^\vee$ (which is naturally a subgroup of $M^\vee$). We will actually not need the precise definition of this group. All that matter is its cardinality which can be explicitly computed using the following short list of properties that $S_\pi$ satisfies (all of them being standard or very easy):
\begin{num}
\item\label{eq 1 gp of centralizers} For $\pi\in \Temp(U(n))$ (resp. $\pi\in \Temp(G_n(E))$) which embeds in $\pi_1\times\ldots\times \pi_t\rtimes \sigma_0$ (resp. $\pi_1\times\ldots\times \pi_t$) where for every $1\leqslant i \leqslant t$, $\pi_i\in \Pi_2(G_{n_i}(E))$ for some $n_i\geqslant 1$ and $\sigma_0\in \Pi_2(U(m))$ for some $m\geqslant 0$, we have
$$\displaystyle S_\pi\simeq S_{\pi_1}\times \ldots \times S_{\pi_t}\times S_{\sigma_0} \; (\mbox{resp. }S_\pi\simeq S_{\pi_1}\times \ldots \times S_{\pi_t}).$$

\item\label{eq 2 gp of centralizers} $S_\pi\simeq \bZ/2n\bZ$ for every $\pi\in \Pi_2(G_n(E))$.

\item\label{eq 3 gp of centralizers} For $\sigma\in \Pi_2(U(n))$, we have $S_\sigma\simeq (\bZ/2\bZ)^k$ where $k$ is such that $BC(\sigma)\simeq \pi_1\times\ldots\times \pi_k$ for some $\pi_i\in \Pi_2(G_{n_i}(E))$ ($1\leqslant i\leqslant k$).
\end{num}

Finally, we note that, for $V$ a Hermitian space, the assignment $\pi\in \Temp(U(V))\mapsto S_\pi$ factorizes through $\Temp_{\ind}(U(V))$ and even through $\Temp(U(V))/\stab$ so that it can as well been considered as an assignment on these former sets.

\subsection{Local $\gamma$-factors}\label{Section gamma factors}

Let $\varphi: W'_F\to \GL(M)$ be a continuous semi-simple and algebraic when restricted to $SL_2(\C)$ finite dimensional complex representation of $W'_F$. We associate to $\varphi$ a local $L$-factor $L(s,\varphi)$ and a local $\epsilon$-factor $\epsilon(s,\varphi,\psi')$ as in \cite[\S 3]{Ta} and \cite[\S 2.2]{GR}. In the $p$-adic case, $L(s,\varphi)$ is of the form $P(q_F^{-s})$ where $P\in \C[T]$ is such that $P(0)=1$ whereas $\epsilon(s,\varphi,\psi')$ is of the form $cq^{n(s-1/2)}$ where $n\in \bZ$ and $c=\epsilon(1/2,\varphi,\psi')\in \C^\times$. In the Archimedean case, $L(s,\varphi)$ is a product of functions of the form $\pi^{-(s+s_0)/2}\Gamma((s+s_0)/2)$ for some $s_0\in \C$ and $\epsilon(s,\varphi,\psi')$ is of the form $c Q^{s-1/2}$ where $Q\in \R_+^*$ and $c=\epsilon(1/2,\varphi,\psi')\in \C^\times$. When $\varphi=\mathbf{1}_F$ is the trivial one-dimensional representation of $W'_F$, we will also write $\zeta_F(s)$ for $L(s,\mathbf{1}_F)$. We have
$$\displaystyle \zeta_F(s)=\left\{
    \begin{array}{ll}
        (1-q_F^{-s})^{-1} & \mbox{ if } F \mbox{ is } p \mbox{-adic,} \\
        \pi^{-\frac{s}{2}}\Gamma(\frac{s}{2}) & \mbox{ if } F=\R.
    \end{array}
\right.
$$
In both cases we have $\zeta_F(s)\sim_{s\to 0} (s \log(q_F))^{-1}$ (where we recall that we set $q_{\R}=e^{1/2}$). Returning to the general case, we define the local $\gamma$-factor associated to $\varphi$ as
$$\displaystyle \gamma(s,\varphi,\psi')=\epsilon(s,\varphi,\psi')\frac{L(1-s,\varphi^\vee)}{L(s,\varphi)}$$
where $\varphi^\vee$ stands for the contragredient of $\varphi$. These factors are additive: for every (semi-simple, continuous, algebraic) complex representations $\varphi$, $\varphi'$ of $W'_F$ we have
\begin{align}\label{eq -4 L parameters, LLC, basechange}
\displaystyle \gamma(s,\varphi\oplus \varphi',\psi')=\gamma(s,\varphi,\psi')\gamma(s,\varphi',\psi')
\end{align}
Of course, we can define similarly $\gamma$-factors $\gamma(s,\varphi,\psi'_E)$ for any (semi-simple, continuous and algebraic) complex representation of $W'_E$ where $\psi'_E$ is the additive character of $E$ defined by $\psi'_E(z)=\psi'(\Tra_{E/F}(z))$. Set
$$\displaystyle \lambda_{E/F}(\psi')=\gamma(\frac{1}{2},\eta_{E/F},\psi')=\epsilon(\frac{1}{2},\eta_{E/F},\psi')$$
where we consider $\eta_{E/F}$ as a character of $W'_F$ through local class field theory. Then, $\lambda_{E/F}(\psi')$ is a fourth root of unity and if $\varphi$ is a complex representation of $W'_E$ of dimension $n$, by the inductivity of $\gamma$-factors in degree $0$ (\cite[Theorem 3.4.1]{Ta}) we have
\begin{align}\label{eq -3 L parameters, LLC, basechange}
\displaystyle \gamma(s,\varphi,\psi'_E)=\lambda_{E/F}(\psi')^{-n}\gamma(s,\Ind_{W'_E}^{W'_F}\varphi,\psi')
\end{align}
where $\Ind_{W'_E}^{W'_F}\varphi$ denotes the representation induced by $\varphi$ from $W'_E$ to $W'_F$.

The following property of Archimedean $\gamma$-factors will be used repeatedly:
\begin{num}
\item \label{eq 1 L parameters, LLC, basechange} In the Archimedean case the $\gamma$-factor $\gamma(s,\varphi,\psi')$ and its inverse are of moderate growth on vertical strips away from their poles together with all their derivatives.
\end{num}
Moreover, if $\varphi$ is tempered, meaning that $\varphi(W_F)$ is bounded in $\GL(M)$, the local $\gamma$-factor $\gamma(s,\varphi,\psi')$ does not vanish in $\cH$ and has no pole in $-\overline{\cH}$. For such a $\varphi$, we will set
$$\displaystyle \gamma^*(0,\varphi,\psi')=\lim\limits_{s\to 0} \zeta_F(s)^{n_\varphi} \gamma(s,\varphi,\psi')$$
where $n_\varphi$ is the order of the zero of $\gamma(s,\varphi,\psi')$ at $s=0$.

Assume that the local Langlands correspondence is known for $G$. Let $r:{}^L G\to \GL(M)$ be a continuous complex representation which is algebraic on $G^\vee$, semisimple and {\em bounded} when restricted to $W_F$. For $\pi\in \Irr(G)$, we set $L(s,\pi,r)=L(s,r\circ \varphi_\pi)$ and $\gamma(s,\pi,r,\psi')=\gamma(s,r\circ \varphi_\pi,\psi')$. When $\pi\in\Temp(G)$, the representation $r\circ \varphi_{\pi}$ is tempered and therefore
\begin{num}
\item\label{eq 2 L parameters, LLC, basechange} $\gamma(s,\pi,r,\psi')$ does not vanish in $\cH$ and has no pole in $-\overline{\cH}$.
\end{num}
Still for $\pi\in \Temp(G)$, we set (as before)
$$\displaystyle \gamma^*(0,\pi,r,\psi')=\lim\limits_{s\to 0} \zeta_F(s)^{n_{\pi,r}} \gamma(s,\pi,r,\psi')$$
where $n_{\pi,r}$ is the order of the zero of $\gamma(s,\pi,r,\psi')$ at $s=0$. We will need the following result which is probably well-known but for which, in lack of a proper reference, we provide a quick proof.

\begin{lem}\label{lem 1 L parameters, LLC, basechange}
Assume that $F=\R$ and let, for $\pi\in \Temp(G)$, $\{s_i(\pi) \}$ be the set of zeroes (counted with multiplicity) of $\gamma(s,\pi,r,\psi')$ on the imaginary line $i\R$. Then, for every $\epsilon>0$ and $C>0$ there exists $k>0$ such that for every $n\geqslant 0$ we have
\begin{equation}\label{ineq1 gamma Arch}
\displaystyle \left\lvert \frac{d^n}{ds^n}\left(\prod_i (s-s_i(\pi))^{-1} \gamma(s,\pi,r,\psi')\right)\right\rvert \ll N(\pi)^k (1+\lvert s\rvert)^k \mbox{ for } \pi\in \Temp(G),\, -C<\Re(s)<1-\epsilon
\end{equation}
and
\begin{equation}\label{ineq1bis gamma Arch}
\left\lvert \frac{d^n}{ds^n}\left(\prod_i (s-s_i(\pi)) \gamma(s,\pi,r,\psi')^{-1}\right)\right\rvert \ll N(\pi)^k (1+\lvert s\rvert)^k \mbox{ for } \pi\in \Temp(G), \, -\frac{1}{2}+\epsilon<\Re(s)<C.
\end{equation}
In particular, for any given $s_0\in \cH$, there exists $k\geqslant 1$ such that
\begin{equation}\label{ineq2 gamma Arch}
\displaystyle \lvert \gamma^*(0,\pi,r,\psi')\rvert \ll N(\pi)^k \mbox{ and } \lvert \gamma(s_0,\pi,r,\psi')^{-1}\rvert \ll N(\pi)^k
\end{equation}
for $\pi\in \Temp(G)$.
\end{lem}

\noindent\ul{Proof}: First, we note that by Cauchy's integral formula it suffices to establish the estimates \eqref{ineq1 gamma Arch} and \eqref{ineq1bis gamma Arch} when $n=0$.
Let $M\subset G$ be a Levi subgroup. Since the functions $\sigma\in \Pi_2(M)\mapsto N(i_M^G(\sigma))$, $\sigma\in \Pi_2(M)\mapsto N(\sigma)$ are equivalent and $\gamma(s,i_M^G(\sigma),r,\psi')=\gamma(s,\sigma,r\circ \iota_M,\psi')$ where $\iota_M:{}^LM \hookrightarrow {}^L G$ is the natural embedding of $L$-groups, we are reduced to prove a similar statement for discrete series. Thus, we assume that $G(\R)$ admits discrete series. Let $T\subset G$ be a maximal elliptic torus, $B_{\C}\subset G_\C$ be a Borel subgroup containing $T_\C$ and $\rho$ be half the sum of the roots of $\mathfrak{t}_\C$ in $\mathfrak{b}_\C$. Then, there exists an embedding of $L$-groups $\iota: {}^L T\to {}^L G$ (here it is important to consider Weil forms of $L$-groups) with the property that for every $\pi\in \Pi_2(G)$ there exists $\chi\in \widehat{T(\R)}$ such that $\varphi_{\pi}=\iota\circ \varphi_\chi$. Moreover, we have the equality of infinitesimal characters $\chi_\pi=d\chi+\rho$ in $\mathfrak{t}(\C)^*/W(G_\C,T_\C)$ where $\chi_\pi$ is identified with an element of the latter set via the Harish-Chandra isomorphism. For such $\pi$ and $\chi$, we have $\gamma(s,\pi,r,\psi')=\gamma(s,\chi,r\circ \iota,\psi')$ and $N(\pi)\sim N(\chi)$ so that we are can further restrict ourself to the case where $G=T$ is a torus. By semisimplicity of $r$, we may moreover assume that it is irreducible. It is then of dimension one or two.

Assume first that $\dim(r)=2$. Then, $r$ is induced from a character $\mu$ of $T^\vee\times W_\C$ and for every $\chi\in \widehat{T(\R)}$ the representation $r\circ \varphi_{\chi}$ is induced from the representation $\mu\circ \varphi_{\chi_\C}$ of $W_\C$ where $\chi_{\C}\in \widehat{T(\C)}$ is the composite of $\chi$ with the norm map $T(\C)\to T(\R)$. Hence
$$\displaystyle \gamma(s,\chi,r,\psi')=\lambda_{\C/\R}(\psi')\gamma(s,\chi_{\C},\mu,\psi'_\C).$$
Furthermore, the restriction of $\mu$ to $T^\vee$ determines a cocharacter $\mathbb{G}_{m,\C}\to T_\C$ hence a morphism $\mu_1: W_\C=\C^\times \to T(\C)$ and denoting by $\mu_2$ the restriction of $\mu$ to $W_\C$ (a unitary character) we have
$$\displaystyle \gamma(s,\chi_{\C},\mu,\psi'_\C)=\gamma(s,(\chi_{\C}\circ\mu_1) \mu_2,\psi'_\C),\;\;\; \chi\in \widehat{T(\R)}.$$
As $N((\chi_{\C}\circ\mu_1) \mu_2)\ll N(\chi)$ for all $\chi\in \widehat{T(\R)}$, in this situation we are ultimately reduced to the case where $G=T=\mathbb{G}_{m,\C}$ and $r$ is the ``standard'' representation. For $\nu\in \widehat{\C^\times}$, $d\nu: \C\to i\R$ is of the form $z\mapsto z(\nu)z-\overline{z(\nu)z}$ for an unique $z(\nu)\in \frac{1}{2}\bZ+i\R$ and we have
$$\displaystyle \lvert \gamma(s,\nu,\psi'_\C)\rvert= C(\psi'_\C)^{\Re(s-1/2)} \left\lvert \frac{\Gamma(1-\overline{s}+z(\nu)^+)}{\Gamma(s+z(\nu)^+)}\right\rvert
$$
where $C(\psi'_\C)\in \R_+^*$ only depends on $\psi'_\C$, $\Gamma$ stands for the usual gamma function and for every complex number $u\in \C$ we have set
$$\displaystyle u^+=\left\{\begin{array}{ll}
u \mbox{ if } u\in \overline{\cH}, \\
-\overline{u} \mbox{ if } u\in -\cH.
\end{array} \right.$$
In particular, we see that the only zero of $\gamma(s,\nu,\psi'_\C)$ on the imaginary line is at $-z(\nu)$ if $z(\nu)\in i\R$ and that there are none otherwise. On the other hand, if $z(\nu)\notin i\R$ then $\lvert s+z(\nu)^+\rvert>\epsilon$ for $\Re(s)>-\frac{1}{2}+\epsilon$ and we have $N(\nu)\sim 1+\lvert z(\nu)\rvert$ for $\nu \in \widehat{\C^\times}$. From this, it follows that the estimates \eqref{ineq1 gamma Arch} and \eqref{ineq1bis gamma Arch} for $n=0$ are in this case equivalent to the following: for every $C,\epsilon>0$, there exists $k>0$ such that
$$\displaystyle \left\lvert (s+z(\nu)^+)^{-1} \frac{\Gamma(1-\overline{s}+z(\nu)^+)}{\Gamma(s+z(\nu)^+)}\right\rvert \ll (1+\lvert s\rvert)^k(1+\lvert z(\nu)\rvert)^k,\; \mbox{ for } \nu\in \widehat{\C^\times},\; -C<\Re(s)<1-\epsilon,$$
$$\displaystyle \left\lvert (s+z(\nu)^+) \frac{\Gamma(s+z(\nu)^+)}{\Gamma(1-\overline{s}+z(\nu)^+)}\right\rvert \ll (1+\lvert s\rvert)^k(1+\lvert z(\nu)\rvert)^k,\; \mbox{ for } \nu\in \widehat{\C^\times},\; -\frac{1}{2}+\epsilon<\Re(s)<C.$$
Setting $z=z(\nu)+i\Im(s)$ and $(a,b)=(1-\Re(s),\Re(s))$ in the first case, $(a,b)=(\Re(s),1-\Re(s))$ in the second case, the above estimates are consequences of the following property of the gamma function:
\begin{equation}\label{ineq gamma function}
\displaystyle \left\lvert (a+z)(b+z)^{-1}\frac{\Gamma(a+z)}{\Gamma(b+z)}\right\rvert\ll (1+\lvert z\rvert)^{\Re(a-b)}
\end{equation}
for $z\in \cH$, $a$ in a fixed compact subset of $\{\Re(s)>-1\}$ and $b$ in a fixed compact subset of $\C$ which can readily be seen using Stirling's asymptotic formula (cf.\ \cite[\S 13.6]{WW}).

Assume now that $\dim(r)=1$. Then, the restriction of $r$ to $T^\vee$ comes from a cocharacter $\mathbb{G}_{m,\R}\to T$ inducing a morphism $\mu_1: W_\R^{\ab}\simeq \R^\times \to T(\R)$ and denoting by $\mu_2$ the restriction of $r$ to $W_\R$ we have
$$\displaystyle \gamma(s,\chi,r,\psi')=\gamma(s,(\chi\circ \mu_1)\mu_2,\psi'),\;\;\; \chi\in \widehat{T(\R)}$$
As $N((\chi\circ\mu_1) \mu_2)\ll N(\chi)$ for all $\chi\in \widehat{T(\R)}$, this time we are reduced to the case where $G=T=\mathbb{G}_{m,\R}$ and $r$ is the ``standard'' representation. For $\nu\in \widehat{\R^\times}$, there exist $s_0\in i\R$ and $\delta\in \{ 0, 1\}$ such that $\nu(t)=\sgn(t)^\delta\lvert t \rvert^{s_0}$ and setting $z(\nu)=s_0+\delta\in i\R \cup i\R+1$, we have
$$\displaystyle \left\lvert \gamma(s,\nu,\psi')\right\rvert=C(\psi')^{\Re(s-1/2)} \left\lvert \frac{\Gamma(\frac{1-\overline{s}+z(\nu)}{2})}{\Gamma(\frac{s+z(\nu)}{2})}\right\rvert.$$
As $N(\nu)\sim 1+\lvert z(\nu)\rvert$, we can deduce similarly the estimates \eqref{ineq1 gamma Arch} and \eqref{ineq1bis gamma Arch} in this case from the inequality \eqref{ineq gamma function}. $\blacksquare$

In this paper, we will only consider $\gamma$-factors associated to the following representations of $L$-groups:
\begin{itemize}
\item $r=\Ad$ is the adjoint representation of ${}^L G$ on $\Lie(G^\vee)$. This notation will be used for all groups encountered in this paper with the exception of $\overline{G_n(E)}$ where we will denote this representation by $\overline{\Ad}$ (in order to distinguish it from the adjoint representation of ${}^L G_n(E)$). Note that we have
\begin{align}\label{eq 3 L parameters, LLC, basechange}
\displaystyle \gamma(s,\pi,\Ad,\psi')=\gamma(s,\mathbf{1}_F,\psi')\gamma(s,\pi,\overline{\Ad},\psi')
\end{align}
for all $\pi\in \Irr(\overline{G_n(E)})$ (the first $\gamma$-factor being defined by viewing $\pi$ as a representation of $G_n(E)$).
\item The tensor product representation $r$ of ${}^L(G_n(E)\times G_m(E))$ ($m,n\in \bN^*$). It is the induced to ${}^L(G_n(E)\times G_m(E))$ of the representation $((g_1,g_2),(g_3,g_4),\sigma)\mapsto g_1\otimes g_2\in \GL_{mn}(\C)$ of the index two subgroup $(G_n(E)^\vee\times G_m(E)^\vee)\times W_E$. In this case, for all $(\pi,\pi')\in\Irr(G_n(E))\times \Irr(G_m(E))$ we will write $\gamma(s,\pi\times \pi',\psi')$ for $\gamma(s,\pi\boxtimes \pi',r,\psi')$ (this should not be confused with the standard $\gamma$-factor of the representation $\pi\times\pi'\in \Irr(G_{m+n}(E))$, but we will never consider standard $\gamma$-factors in this paper).
\item $r=\As: {}^L G_n(E)\to \GL(\C^n\otimes \C^n)$ the Asai representation which is defined by $\As(g_1,g_2)=g_1\otimes g_2$ for $g_1,g_2\in \GL_n(\C)$ and $\As(\sigma)=\iota$ if $\sigma\in W_F\setminus W_E$ where $\iota\in \GL(\C^n\otimes \C^n)$ is given by $\iota(v\otimes w)=w\otimes v$.
\end{itemize}
The corresponding $\gamma$-factors satisfy the following properties (where we omit the additive character $\psi'$ since it is fixed throughout in this paper)
\begin{align}\label{eq 4 L parameters, LLC, basechange}
\displaystyle \gamma(s,\pi\lvert .\rvert_E^x,\Ad)=\gamma(s,\pi,\Ad), \;\; (\pi\in \Irr(G_n(E)), x\in \C).
\end{align}

\begin{align}\label{eq 5 L parameters, LLC, basechange}
\displaystyle  & \gamma(s,\pi_1\times \pi_2,\Ad)=\gamma(s,\pi_1,\Ad)\gamma(s,\pi_2,\Ad)\gamma(s,\pi_1\times \pi_2^\vee)\gamma(s,\pi^\vee_1\times \pi_2),\\
\nonumber & ((\pi_1,\pi_2)\in \Irr(G_n(E))\times \Irr(G_m(E))).
\end{align}

\begin{align}\label{eq 6 L parameters, LLC, basechange}
\displaystyle \gamma(s,\pi,\Ad) \mbox{ has a simple zero at }s=0 \mbox{ for every } \pi\in \Pi_2(G_n(E)).
\end{align}

\begin{align}\label{eq 6bis L parameters, LLC, basechange}
\displaystyle & \gamma(s,\pi_1\times (\pi_2\times \pi_3))=\gamma(s,\pi_1\times \pi_2)\gamma(s,\pi_1\times \pi_3) \\
\nonumber & ((\pi_1,\pi_2,\pi_3)\in \Temp(G_n(E))\times \Temp(G_m(E))\times \Temp(G_k(E)))
\end{align}
(here $\pi_2\times \pi_3$ denotes the parabolic induction of $\pi_2\boxtimes \pi_3$ to $G_{m+k}(E)$).
\begin{align}\label{eq 7 L parameters, LLC, basechange}
\displaystyle & \gamma(s,\pi_1\lvert .\rvert_E^x\times \pi_2 \lvert .\rvert_E^y)=\gamma(s+x+y,\pi_1\times \pi_2),\\
\nonumber & ((\pi_1,\pi_2)\in \Irr(G_n(E))\times \Irr(G_m(E)), x,y\in \C).
\end{align}

\begin{num}
\item\label{eq 8 L parameters, LLC, basechange} For every $(\pi_1,\pi_2)\in \Pi_2(G_n(E))\times \Pi_2(G_m(E))$, $\gamma(s,\pi_1\times \pi_2)$ has at most a simple zero at $s=0$ and moreover
$$\gamma(0,\pi_1\times \pi_2)=0\Leftrightarrow \pi_1\simeq \pi_2^\vee.$$
\end{num}

\begin{align}\label{eq 9 L parameters, LLC, basechange}
\displaystyle \gamma(s,\pi\lvert.\rvert_E^x,\As)=\gamma(s+2x,\pi,\As),\;\; (\pi\in \Irr(G_n(E)), x\in \C).
\end{align}

\begin{align}\label{eq 10 L parameters, LLC, basechange}
\displaystyle & \gamma(s,\pi_1\times \pi_2,\As)=\gamma(s,\pi_1,\As)\gamma(s,\pi_2,\As)\gamma(s,\pi_1\times \pi_2^c), \\
\nonumber & (\pi_1,\pi_2)\in \Irr(G_n(E))\times \Irr(G_m(E)).
\end{align}

\begin{num}
\item\label{eq 11 L parameters, LLC, basechange} For every $\pi\in \Pi_2(G_n(E))$, $\gamma(s,\pi,\As)$ has at most a simple zero at $s=0$ and moreover
$$\displaystyle \gamma(0,\pi,\As)=0\Leftrightarrow \pi\in BC_n(\Temp(U(n))).$$
\end{num}

\begin{align}\label{eq 12 L parameters, LLC, basechange}
\displaystyle \gamma(s,\sigma,\Ad)=\frac{\gamma(s,BC_n(\sigma),\Ad)}{\gamma(s,BC_n(\sigma),\As)},\;\; (\sigma\in \Temp(U(V)), n=\dim(V)).
\end{align}

\subsection{Harish-Chandra Plancherel formula}\label{Section Planch}

Let $f\in \cS(G(F))$. For all $\pi\in \Temp_{\ind}(G)$ we define a function $f_\pi$ by
$$\displaystyle f_\pi(g)=\Tr(\pi(g)\pi(f^\vee)),\;\;\; g\in G(F)$$
where $f^\vee(x)=f(x^{-1})$. Note that by \eqref{eq 1 representations}, we have $f_\pi\in \cC^w(G(F))$ and the linear map $f\in \cS(G(F))\mapsto f_\pi\in \cC^w(G(F))$ is continuous. Recall that in Section \ref{Section space of functions}, we have introduced a notion of Schwartz function $\Temp_{\ind}(G)\to \cC^w(G(F))$. The following is probably well-known but in lack of a proper reference we provide a proof in Appendix \ref{Appendice Schwartz}.
\begin{prop}\label{prop 1 Planch}
For every $f\in \cS(G(F))$, the map $\pi\in \Temp_{\ind}(G)\mapsto f_\pi\in \cC^w(G(F))$ is Schwartz.
\end{prop}

According to Harish-Chandra (see \cite{H-C}, \cite{Wald1}), there exists a (necessarily unique) Borel measure $d\mu_G$ on $\Temp(G(F))$ such that
\begin{align}\label{eq 1 Planch}
\displaystyle f(g)=\int_{\Temp(G(F))}f_\pi(g) d\mu_G(\pi)
\end{align}
for all $f\in \cS(G(F))$ and $g\in G(F)$. Moreover, $d\mu_G(\pi)=\mu_G(\pi)d_{\psi'}\pi$\footnote{Let us emphasize that, although not transparent from the notation, the measure $d\mu_G(\pi)$ depends on the choice of the additive character $\psi'$: indeed, this measure is inversely proportional to the Haar measure on $G(F)$ which, following the convention of Section \ref{Section measures}, depends itself on $\psi'$. Similarly, the density $\mu_G(\pi)$ depends on $\psi'$, a dependence which, in the case of $\GL_n$, will be made more transparent by Proposition \ref{prop Planch measure} below.} where $d_{\psi'}\pi$ is the measure on $\Temp_{\ind}(G)$ defined in Section \ref{Section measures} and for $\pi=i_M^G(\sigma)$, where $M$ is a Levi subgroup of $G$ and $\sigma\in \Pi_2(M)$, we have
$$\displaystyle \mu_G(\pi)=d(\sigma)j(\sigma)^{-1}$$
where $d(\sigma)$ is the formal degree of $\sigma$ (as defined in Section \ref{Section representations}) and $j(\sigma)$ is a certain product of standard intertwining operators that we won't need to describe precisely here. When $\pi$ is a discrete series we simply have $\mu_G(\pi)=d(\pi)$. In the Archimedean case the function $\pi\mapsto \mu_G(\pi)$ is of moderate growth i.e. there exists $k>0$ such that
\begin{align}\label{eq 2 Planch}
\displaystyle \mu_G(\pi)\ll N(\pi)^k,\;\;\; \pi\in \Temp_{\ind}(G).
\end{align}
Combining this with Proposition \ref{prop 1 Planch} and \eqref{basic estimates spectral measure}, we see that the function $\pi\in \Temp_{\ind}(G)\mapsto f_\pi\in \cC^w(G(F))$ is absolutely integrable with respect to the measure $d\mu_G(\pi)$ and by \eqref{eq 1 Planch} this implies that (where the right-hand side is a priori defined as an element of $\cC^w(G(F))$)
\begin{align}\label{eq 3 Planch}
\displaystyle f=\int_{\Temp_{\ind}(G)}f_\pi d\mu_G(\pi).
\end{align}

Recall that in Section \ref{Section Groups of centralizers} we have associated to any $\pi\in \Temp(G_n(E))$ a finite abelian group $S_\pi$. By \eqref{eq 1 gp of centralizers} and \eqref{eq 2 gp of centralizers} we immediately see that if $\pi$ is induced from a discrete series of a Levi of the form $M=R_{E/F}G_{n_1,E}\times \ldots \times R_{E/F}G_{n_k,E}$ then
\begin{align}\label{eq 3bis Planch}
\displaystyle \lvert S_\pi\rvert=2^k n_1\ldots n_k
\end{align}
Finally recall that we have set $\overline{G_n(E)}:=G_n(E)/Z_n(F)$ (i.e. the group of $F$-points of the quotient $(R_{E/F}G_{n,E})/Z_{n,F}$). The following formula for $d\mu_{\overline{G_n(E)}}$ is essentially due to Harish-Chandra \cite{H-C} in the Archimedean case and Shahidi \cite{Sha1} and Silberger-Zink \cite{SZ} in the $p$-adic case. Actually, this also incorporates a reformulation of the result of Silberger-Zink due to Hiraga-Ichino-Ikeda \cite{HII}. Because our normalization of measures does not compare obviously to theirs (in particular because we have normalized everything by considering $G_n(E)$ as an $F$-group) we provide a short proof to explain the relation.

\begin{prop}\label{prop Planch measure}
We have the following equality
$$\displaystyle d\mu_{\overline{G_n(E)}}(\pi)=\lambda_{E/F}(\psi')^{-n^2} \frac{\gamma^*(0,\pi,\overline{\Ad},\psi')}{\lvert S_\pi\rvert} d\pi$$
of measures on $\Temp(\overline{G_n(E)})$.
\end{prop}

\noindent\ul{Proof}: By Fourier inversion on $Z_n(F)$ it is easy to see that
$$\displaystyle \mu_{\overline{G_n(E)}}(\pi)=\mu_{G_n(E)}(\pi),\;\;\; \pi\in \Temp(\overline{G_n(E)}).$$
Thus in fact it suffices to show that
\begin{align}\label{eq 4 Planch}
\displaystyle d\mu_{G_n(E)}(\pi)=\lambda_{E/F}(\psi')^{-n^2} \omega_{\pi}(-1)^{n-1}\frac{\gamma^*(0,\pi,\Ad,\psi')}{\lvert S_\pi\rvert} d\pi,\;\;\; \pi\in \Temp(G_n(E)).
\end{align}
The Plancherel measure for $G_n(E)$ is completely computed by results of Shahidi \cite{Sha1} and Silberger-Zink \cite{SZ} in the $p$-adic case and by Harish-Chandra \cite{H-C} in the Archimedean case but to state their results it is actually better to consider $G_n(E)$ as the group of $E$-points of $G_n$ (and not as the group of $F$-points of $R_{E/F}G_{n,E}$ as in other parts of this paper). Using the character $\psi'_E=\psi'\circ Tr_{E/F}$ we can construct as before a measure on $G_n(E)$, which coincides with the one we already fixed, and a measure $d_{\psi'_E}\pi$ on $\Temp(G_n(E))$, which does not coincide with $d_{\psi'}\pi$. Then, we have
\begin{align}\label{eq 5 Planch}
\displaystyle d\mu_{G_n(E)}(\pi)=\mu_{G_n(E)}^E(\pi)d_{\psi'_E}\pi
\end{align}
where $\mu_{G_n(E)}^E(\pi)$ is the Plancherel density for $G_n(E)$ considered as an $E$-group. For $\pi=i_M^{G_n(E)}(\sigma)$ where $M$ is a Levi subgroup of $G_{n,E}$ and $\sigma\in \Pi_2(M(E))$ we have
\begin{align}\label{eq 6 Planch}
\displaystyle \mu_{G_n(E)}^E(\pi)=d^E(\sigma)j(\sigma)^{-1}
\end{align}
where $j(\sigma)$ is the same factor as before and $d^E(\sigma)$ stands for the formal degree of $\sigma$ when $M(E)$ is considered as an $E$-group (i.e. by integrating product of matrix coefficients on $M(E)/A_M(E)$ and not, assuming that $M$ is defined over $F$, over $M(E)/A_M(F)$). Assume that $M$ is of the form
$$\displaystyle M\simeq G_{n_1,E}\times\ldots\times G_{n_k,E}$$
and that
$$\displaystyle \sigma=\tau_1\boxtimes\ldots\boxtimes \tau_k$$
where $\tau_i\in \Pi_2(G_{n_i}(E))$ for all $1\leqslant i\leqslant k$. Then, by a result of Shahidi \cite{Sha1}, \cite{Sha2} we have
\begin{align}\label{eq 7 Planch}
\displaystyle j(\sigma)^{-1}=\left(\prod_{i=1}^k \omega_{\tau_i}(-1)^{n-n_i}\right)\gamma(0,\sigma,\Ad_{G_n/M},\psi'_E)
\end{align}
where $\Ad_{G_n/M}$ denotes the adjoint representation of $M^\vee$ on $\Lie(G_n^\vee)/\Lie(M^\vee)$. On the other hand, by a result of Silberger-Zink \cite{SZ} in the $p$-adic case and Harish-Chandra \cite[Lemma 23.1]{H-C} in the Archimedean case as reformulated by Hiraga-Ichino-Ikeda \cite{HII} and reproved and refined by Ichino-Lapid-Mao \cite[Theorem 2.1]{ILM} we have
\begin{align}\label{eq 8 Planch}
\displaystyle d^E(\sigma)=\left(\prod_{i=1}^k\omega_{\tau_i}(-1)^{n_i-1}\right)\frac{\gamma(0,\sigma,\Ad_{M/A_M},\psi'_E)}{n_1\ldots n_k}
\end{align}
where $\Ad_{M/A_M}$ stands for the adjoint representation of $M^\vee$ on $M^\vee/A_M^\vee$. On the other hand, on the connected component of $\pi$ we have
\begin{align}\label{eq 9 Planch}
\displaystyle d_{\psi'_E}\pi=\gamma^*(0,\mathbf{1}_E,\psi'_E)^kd^E\pi
\end{align}
where $d^E\pi$ is the analog of the measure $d\pi$ by viewing $G_n(E)$ as an $E$-group and
\begin{align}\label{eq 10 Planch}
\displaystyle \gamma^*(0,\mathbf{1}_E,\psi'_E)=\lim\limits_{s\to 0^+}\zeta_E(s) \gamma(s,\mathbf{1}_E,\psi'_E)=(\frac{\log(q_F)}{\log(q_E)})\lim\limits_{s\to 0^+}\zeta_F(s) \gamma(s,\mathbf{1}_E,\psi'_E)
\end{align}
Finally by the description \eqref{eq 4 Measures} of the measure $d\pi$ and its analog for $d^E\pi$ we easily check that
\begin{align}\label{eq 11 Planch}
\displaystyle d^E\pi=2^{-k}(\frac{\log(q_E)}{\log(q_F)})^kd\pi
\end{align}
(the factor $2^k$ stems from the fact that $\lvert .\rvert_E=\lvert .\rvert_F^2$ on $F^\times$). Combining \eqref{eq 5 Planch}, \eqref{eq 6 Planch}, \eqref{eq 7 Planch}, \eqref{eq 8 Planch}, \eqref{eq 9 Planch}, \eqref{eq 10 Planch} and \eqref{eq 11 Planch} we obtain
$$\displaystyle d\mu_{G_n(E)}(\pi)=\omega_{\pi}(-1)^{n-1}\gamma(0,\sigma,\Ad_{G_n/M},\psi'_E)\frac{\gamma(0,\sigma,\Ad_{M/A_M},\psi'_E)}{2^k n_1\ldots n_k} \lim\limits_{s\to 0^+} \zeta_F(s)^k \gamma(s,\mathbf{1}_E,\psi'_E)^k d\pi$$
Now, the adjoint representation of $M^\vee$ on $\Lie(G_n)^\vee$ is the direct sum of $\Ad_{G_n/M}$, $\Ad_{M/A_M}$ and $k$ copies of the trivial representation so that by additivity of $\gamma$-factors \eqref{eq -4 L parameters, LLC, basechange} and their inductivity in degree $0$ \eqref{eq -3 L parameters, LLC, basechange}
\[\begin{aligned}
\displaystyle \gamma(s,\sigma,\Ad_{G_n/M},\psi'_E)\gamma(s,\sigma,\Ad_{M/A_M},\psi'_E)\gamma(s,\mathbf{1}_E,\psi'_E)^k & =\gamma(s,\pi,\Ad,\psi'_E) \\
 & =\lambda_{E/F}(\psi')^{-n^2}\gamma(s,\pi,\Ad,\psi')
\end{aligned}\]
Together with \eqref{eq 3bis Planch} this gives precisely \eqref{eq 4 Planch} and ends the proof of the proposition. $\blacksquare$

\subsection{Plancherel formula for Whittaker functions}\label{section Plancherel Whitt}

In this section we make a little bit more precise the Plancherel formula for Whittaker functions proved by Sakellaridis and Venkatesh in \cite[\S 6.3]{SV} (in particular, our discussion includes real groups). Let us also mention that a similar Whittaker-Plancherel inversion formula has been established, via different means, by Wallach \cite[Chap. 15]{Wall2} in the Archimedean case and Delorme \cite{Del1} in the $p$-adic case.

Assume that $G$ is quasi-split, let $B=TN$ be a Borel subgroup with unipotent radical $N$ and $\xi:N(F)\to \mathbb{S}^1$ a generic character. For any $f\in \cS(G(F))$ we define a function $W_f$ on $G(F)\times G(F)$ by
$$\displaystyle W_f(g_1,g_2)=\int_{N(F)} f(g_1^{-1}ug_2) \xi(u)^{-1}du,\;\;\; g_1,g_2\in G(F)$$
Clearly, we have $W_f(u_1g_1,u_2g_2)=\xi(u_1)^{-1}\xi(u_2)W_f(g_1,g_2)$ for every $(g_1,g_2)\in G(F)^2$ and every $(u_1,u_2)\in N(F)^2$. Moreover, $W_f(g,.)\in \cS(N(F)\backslash G(F),\xi)$ for all $g\in G(F)$.

Recall from \cite[\S 7.1]{Beu1}, that the functional $\displaystyle f\in \cS(G(F))\mapsto \int_{N(F)} f(u)\xi(u)^{-1}du$ extends (necessarily uniquely) by continuity to a functional $\cC^w(G(F))\to \C$ to be denoted by
$$\displaystyle f\mapsto \int^*_{N(F)} f(u)\xi(u)^{-1}du$$
(actually in \cite{Beu1} only the case of unitary groups is considered but the proof is easily seen to work in general). This allows to extend the definition of $W_f$ to every $f\in \cC^w(G(F))$ by
$$\displaystyle W_f(g_1,g_2)=\int^*_{N(F)} f(g_1^{-1}ug_2) \xi(u)^{-1}du,\;\;\; g_1,g_2\in G(F)$$
Note that $\xi^{-1}\boxtimes \xi$ is a non-degenerate character on the $F$-points of the maximal unipotent subgroup $N\times N$ of $G\times G$. Thus, following the general construction of Section \ref{Section space of functions}, there is a corresponding space of functions $\cC^w(N(F)\times N(F)\backslash G(F)\times G(F),\xi^{-1}\boxtimes \xi)$.

\begin{lem}\label{lem 1 HC Planch}
For all $f\in \cC^w(G(F))$ we have $W_f\in \cC^w(N(F)\times N(F)\backslash G(F)\times G(F),\xi^{-1}\boxtimes \xi)$ and moreover the linear map
$$\displaystyle f\in \cC^w(G(F))\mapsto W_f\in \cC^w(N(F)\times N(F)\backslash G(F)\times G(F),\xi^{-1}\boxtimes \xi)$$
is continuous.
\end{lem}

\noindent\ul{Proof}: By the continuity of the regularized integral $\int_{N(F)}^* . \xi(u)^{-1}du$, the linear forms $f\in \cC^w(G(F))\mapsto W_f(g_1,g_2)$ are continuous for $(g_1,g_2)\in G(F)\times G(F)$. Thus, by the closed graph theorem, it suffices to show that $W_f\in \cC^w(N(F)\times N(F)\backslash G(F)\times G(F),\xi^{-1}\boxtimes \xi)$ for every $f\in \cC^w(G(F))$. Let $d>0$. By Dixmier-Malliavin \cite{DM}, any $f\in \cC_d^w(G(F))$ is a finite sum of functions of the form $\varphi_1\star f'\star \varphi_2$ where $f'\in \cC_d^w(G(F))$ and $\varphi_1,\varphi_2\in C_c^\infty(G(F))$. Thus, we may as well assume that $f$ is of this form. From the continuity of the linear form $\int^*_{N(F)} .\xi(u)^{-1}du$, we then have
$$\displaystyle W_{\varphi_1\star f'\star \varphi_2}=W_{f'}\star(\varphi_1^\vee\otimes \varphi_2)$$
where $\varphi_1^\vee\otimes \varphi_2\in C_c^\infty(G(F)\times G(F))$ is defined by $\varphi_1^\vee\otimes \varphi_2(g_1,g_2)=\varphi_1(g_1^{-1})\varphi(g_2)$. Therefore, the function $W_f$ is smooth and moreover in the Archimedean case $R(u_1\otimes u_2)W_f=W_{f'}\star(R(u_1)\varphi_1^\vee\otimes R(u_2)\varphi_2)$ for every $u_1,u_2\in \cU(\g)$. It thus suffices to check the existence of $d'>0$ such that for every $f\in \cC_d^w(G(F))$ we have
$$\displaystyle \lvert W_f(g_1,g_2)\rvert\ll \Xi^{N\backslash G}(g_1)\Xi^{N\backslash G}(g_2)\sigma_{N\backslash G}(g_1)^{d'}\sigma_{N\backslash G}(g_2)^{d'},\;\;\; g_1,g_2\in G(F).$$
Such an inequality is proved in \cite[Lemma 7.3.1(ii)]{Beu1} in the case where $G$ is a unitary group but the proof of {\it loc. cit.} actually works for any reductive group. $\blacksquare$

For $f\in \cS(G(F))$ and $\pi\in \Temp_{\ind}(G)$, we set $W_{f,\pi}:=W_{f_\pi}$ where $f_\pi\in \cC^w(G(F))$ is defined as in Section \ref{Section Planch}. By the previous lemma, we have $W_{f,\pi}\in \cC^w(N(F)\times N(F)\backslash G(F)\times G(F),\xi^{-1}\boxtimes \xi)$ and the linear map $f\in \cS(G(F))\mapsto W_{f,\pi}\in \cC^w(N(F)\times N(F)\backslash G(F)\times G(F),\xi^{-1}\boxtimes \xi)$ is continuous.

For $\pi\in \Temp_{\ind}(G)$, the linear form
$$\displaystyle T\in \End^\infty(\pi)\simeq \pi^\vee\boxtimes \pi\mapsto \int_{N(F)}^* \Tr(\pi(u)T)\xi(u)^{-1}du$$
is well-defined, continuous by \eqref{eq 1 representations} and is obviously a $(N(F)\times N(F),\xi^{-1}\boxtimes \xi)$-equivariant. Therefore, if $\pi$ is not $(N(F),\xi)$-generic we have $W_{f,\pi}=0$ whereas if $\pi$ is $(N(F),\xi)$-generic we have
\begin{align}\label{eq 1 Planch Whitt}
\displaystyle W_{f,\pi}\in \cW(\pi^\vee\boxtimes \pi,\xi^{-1}\boxtimes \xi)
\end{align}
for every $f\in \cS(G(F))$.

The following proposition is a direct consequence of the previous lemma, of Proposition \ref{prop 1 Planch} and of the Plancherel formula \eqref{eq 3 Planch}.

\begin{prop}\label{prop 1 Plancherel Whitt}
For $f\in \cS(G(F))$, the map $\pi\in \Temp_{\ind}(G)\mapsto W_{f,\pi}\in \cC^w(N(F)\times N(F)\backslash G(F)\times G(F),\xi^{-1}\boxtimes \xi)$ is Schwartz. Moreover, we have the following equality
$$\displaystyle W_f=\int_{\Temp_{\ind}(G)} W_{f,\pi}d\mu_G(\pi)$$
where the right hand side is absolutely convergent in $\cC^w(N(F)\times N(F)\backslash G(F)\times G(F),\xi^{-1}\boxtimes \xi)$.
\end{prop}

Assume now that $G=R_{E/F} G_n$ and recall that in Section \ref{Section Whittaker models} we have defined a scalar product $(.,.)^{\Whitt}$ on $\cC^w(N_n(E)\backslash G_n(E),\psi_n)$. The next proposition is \cite[Lemma 4.4]{LM2}. Notice that the extra factor $\lvert \tau\rvert_E^{-n(n-1)/2}$ is a consequence of our choice of measures which are normalized using the character $\psi'_E=\psi'\circ \Tra_{E/F}$ rather than $\psi$. Also, that the left-hand side of the proposition coincides with what is denoted $[W,W']^{\psi_n}$ in {\it loc. cit.} follows from \cite[Proposition 2.3 and Proposition 2.11]{LM2} in the $p$-adic case (this last proposition actually shows that $\displaystyle \int_{N_n(E)} . \psi_n(u)^{-1}du$ extends continuously to the bigger space $C^\infty(G_n(E))$ of all smooth functions on $G_n(E)$) and from the alternative definition of $(.,.)^{\psi_n}$ given in the middle of p.465 of {\it loc. cit.} in the Archimedean case.

\begin{prop}\label{prop 1 Planch Whitt}
For every $W,W'\in \cC^w(N_n(E)\backslash G_n(E),\psi_n)$ we have
$$\displaystyle \int^*_{N_n(E)} (R(u)W,W')^{\Whitt} \psi_n(u)^{-1}du=\lvert \tau\rvert_E^{-n(n-1)/2}W(1)\overline{W'(1)}.$$
\end{prop}

Let $f\in \cS(G_n(E))$ and $\pi\in \Temp(G_n(E))$. By the above construction, we associate to $f$ and $\pi$ a function $ W_{f,\pi}\in \cC^w(N_n(E)\backslash G_n(E)\times N_n(E)\backslash G_n(E),\psi_n^{-1}\boxtimes \psi_n)$. The $G_n(E)$-invariant scalar product $(.,.)^{\Whitt}$ on $\cW(\pi,\psi_n)$ induces an isomorphism
$$\displaystyle \End^\infty(\pi)\simeq \pi^\vee\ctens \pi\simeq \overline{\cW(\pi,\psi_n)}\ctens \cW(\pi,\psi_n)$$
where $\overline{\cW(\pi,\psi_n)}$ denotes the complex conjugate of $\cW(\pi,\psi_n)$. Let $\cB(\pi,\psi_n)$ be an orthonormal basis of (the Hilbert completion of) $\cW(\pi,\psi_n)$ for $(.,.)^{\Whitt}$ obtained by taking the union of orthonormal basis for $\cW(\pi,\psi_n)[\delta]$ for every $\delta\in \widehat{K_{n,E}}$. Then, through the above identification, we have
\begin{align}\label{eq 2 Planch Whitt}
\displaystyle \pi(f)=\sum_{W\in \cB(\pi,\psi_n)} \overline{W}\otimes R(f)W
\end{align}
where the sum converges absolutely in $\End^\infty(\pi)$ (notice that this sum is actually finite in the $p$-adic case). From this, Proposition \ref{prop 1 Planch Whitt} and \eqref{eq 1 representations}, we deduce (recall that $f_\pi(g)=\Tr(\pi(g)\pi(f^\vee))$)
\begin{align}\label{eq 3 Planch Whitt}
\displaystyle W_{f,\pi}=\lvert \tau\rvert_E^{-n(n-1)/2}\sum_{W\in \cB(\pi,\psi_n)} \overline{W}\otimes R(f^\vee)W
\end{align}
where the sum converges absolutely in $\cC^w(N_n(E)\backslash G_n(E)\times N_n(E)\backslash G_n(E),\psi_n^{-1}\boxtimes \psi_n)$.

\subsection{The linear form $\beta$}\label{Section betan}

Let $n\geqslant 1$. For every $W\in \cC^w(N_n(E)\backslash G_n(E),\psi_n)$, we set
$$\displaystyle \beta(W)=\int_{N_n(F)\backslash P_n(F)} W(p)dp \mbox{ and } \beta_n(W)=\int_{N_n(F)\backslash P_n(F)} W(p)\eta_{E/F}(\det p)^{n-1}dp.$$

\begin{lem}\label{lem continuity of beta}
The integral defining $\beta(W)$, $\beta_n(W)$ are absolutely convergent and $\beta$, $\beta_n$ are continuous linear forms on $\cC^w(N_n(E)\backslash G_n(E),\psi_n)$.
\end{lem}

\noindent\ul{Proof}: By Lemma \ref{lem 2 space of functions} and the Iwasawa decomposition $P_n(F)=N_n(F)A_{n-1}(F)K_{n-1}$, it suffices to show the existence of $N\geqslant 1$ such that for every $d>0$ the integral
$$\displaystyle \int_{A_{n-1}(F)} \prod_{i=1}^{n-2}(1+\lvert \frac{a_i}{a_{i+1}}\rvert)^{-N} (1+\lvert a_{n-1}\rvert)^{-N} \delta_{n,E}(a)^{1/2} \sigma(a)^d \delta_{n-1}(a)^{-1}da$$
converges. As $\delta_{n,E}(a)^{1/2}\delta_{n-1}(a)^{-1}=\delta_{n}(a)\delta_{n-1}(a)^{-1}=\lvert \det a\rvert$ this follows from Lemma \ref{lem basic estimates}. $\blacksquare$

Let $\sigma\in \Temp(U(n))$ and set $\pi=BC_n(\sigma)$. Then, by slightly reformulating the main results of \cite{Har} and \cite{Mat} we have that $\pi$ is {\em $G_n(F)$-distinguished} (i.e. it admits a nonzero continuous $G_n(F)$-invariant linear form). On the other hand, the form $\beta$ is obviously $P_n(F)$-invariant. By a result of Ok and Offen (\cite{Ok}, \cite{Off}) in the $p$-adic case and of Kemarsky (\cite{Kem}) in the Archimedean case, this implies
\begin{num}
\item\label{prop0 form betan} The restriction of $\beta$ to $\cW(\pi,\psi_n)$ is $G_n(F)$-invariant.
\end{num}
Following the reasoning of \cite[(2) p.263]{LM3}, this also implies the existence of a sign $c_1(\pi)\in \{\pm 1\}$ such that, setting $\widetilde{W}(g)=W(w_n{}^tg^{-1})$ for every $g\in G_n(E)$, we have
\begin{align}\label{prop form betan}
\displaystyle \beta(\widetilde{W})=c_1(\pi)\beta(W)
\end{align}
for every $W\in \cW(\pi,\psi_n)$. By \cite[Corollary 7.2]{Off} we have $c_1(\pi)=1$ in the $p$-adic case. Actually, it should be possible to prove that $c_1(\pi)=1$ in general using Corollary \ref{cor lim spectrale} but we decided not do so since it would increase the length of this paper and we don't need this result.

\subsection{Rankin-Selberg local functional equation for Asai $\gamma$-factors}\label{Section Zeta integrals}

The following variants of the local Rankin-Selberg zeta integrals have been introduced and studied by Flicker \cite{Fli3} and Kable \cite{Kab}. For every $W\in \cC^w(N_n(E)\backslash G_n(E),\psi_n)$, $\phi\in \cS(F^n)$ and $s\in \cH$ set
$$\displaystyle Z(s,W,\phi):=\int_{N_n(F)\backslash G_n(F)} W(h) \phi(e_nh) \lvert \det h\rvert^s dh$$
where $e_n=(0,\ldots,0,1)\in F^n$. For all $\phi\in \cS(F^n)$ we shall also denote by
$$\displaystyle \widehat{\phi}(y)=\int_{F^n} \phi(x)\psi'(y{}^t x)dx,\;\;\; y\in F^n$$
the Fourier transform of $\phi$ (recall that $dx$ denotes the autodual measure with respect to $\psi'$).

\begin{lem}\label{lem Zeta function}
\begin{enumerate}[(i)]
\item For $s\in \cH$, the integral defining $Z(s,W,\phi)$ is absolutely convergent and defines a continuous linear form on $\cC^w(N_n(E)\backslash G_n(E),\psi_n)$. Moreover, for all $\phi\in \cS(F^n)$ the map $s\in \cH\mapsto Z(s,.,\phi)\in \cC^w(N_n(E)\backslash G_n(E),\psi_n)'$ is analytic.

\item For all $W\in \cS(Z_n(F)N_n(E)\backslash G_n(E),\psi_n)$ and $\phi\in \cS(F^n)$ we have
$$\displaystyle \lim\limits_{s\to 0^+} n\gamma(s,\mathbf{1}_F,\psi')Z(s,W,\phi)=\phi(0)\int_{Z_n(F)N_n(F)\backslash G_n(F)} W(h)dh.$$
\end{enumerate}
\end{lem}

\noindent\ul{Proof}:
\begin{enumerate}[(i)]
\item Let $\phi\in \cS(F^n)$. By Lemma \ref{lem 2 space of functions} and the Iwasawa decomposition $G_n(F)=N_n(F)A_n(F)K_n$, it suffices to show that for all $C>0$ and $d>0$ there exists $N\geqslant 1$ such that the following integral converges uniformly on compact subsets of $\left\{s\in \C\mid 0<\Re(s)<C \right\}$
$$\displaystyle \int_{A_n(F)} \prod_{i=1}^{n-1} (1+\lvert \frac{a_i}{a_{i+1}}\rvert)^{-N} \delta_{n,E}(a)^{1/2}\sigma(a)^d \lvert \phi(e_na)\rvert \delta_n(a)^{-1} \lvert \det a\rvert^{\Re(s)}da.$$
Since $\delta_{n,E}(a)^{1/2}=\delta_n(a)$ and $\lvert \phi(e_na)\rvert\ll (1+\lvert a_n\rvert)^{-N}$ for all $a\in A_n(F)$, this follows from Lemma \ref{lem basic estimates}.

\item Let $W\in \cS(Z_n(F)N_n(E)\backslash G_n(E),\psi_n)$ and $\phi\in \cS(F^n)$. Then, again by the Iwasawa decomposition, for all $s\in \cH$ we have
\[\begin{aligned}
\displaystyle & \gamma(ns,\mathbf{1}_F,\psi') Z(s,W,\phi)= \\
& \int_{A_{n-1}(F)}\int_{K_n} W(ak) \gamma(ns,\mathbf{1}_F,\psi')\left(\int_{Z_n(F)} \phi(e_nzk) \lvert \det z\rvert^s dz\right) \, dk\lvert \det a\rvert^{s}\delta_n(a)^{-1}da
\end{aligned}\]
Set $\phi_k(x)=\phi(xe_nk)$ for all $x\in F$ and $k\in K$. Then $\phi_k\in \cS(F)$ and by Tate's thesis for every $k\in K_n$ we have
$$\displaystyle \gamma(ns,\mathbf{1}_F,\psi')\int_{Z_n(F)} \phi(e_nzk) \lvert \det z\rvert^s dz=\int_F \widehat{\phi_k}(x) \lvert x\rvert^{-ns} dx$$
for all $s\in \cH$ sufficiently close to $0$. Since $\{\phi_k\mid k\in K_n\}$ is a bounded subset of $\cS(F^n)$ so does $\{\widehat{\phi_k}\mid k\in K_n\}$ and therefore $\displaystyle (k,s)\mapsto \int_F \widehat{\phi_k}(x) \lvert x\rvert^{-ns} dx$ is uniformly bounded for $s$ in a neighborhood of $0$. By dominated convergence we deduce that
\[\begin{aligned}
\displaystyle & \lim\limits_{s\to 0^+}\gamma(ns,\mathbf{1}_F,\psi')Z(s,W,\phi)=\int_{A_{n-1}(F)} \int_{K_n} W(ak)\left(\int_F \widehat{\phi_k}(x)dx\right)\, dk\lvert \det a\rvert^{s}\delta_n(a)^{-1}da \\
 & =\phi(0)\int_{A_{n-1}(F)\times K_n} W(ak) \delta_n(a)^{-1}dadk=\phi(0)\int_{Z_n(F)N_n(F)\backslash G_n(F)} W(h)dh
\end{aligned}\]
As $\gamma(ns,\mathbf{1}_F,\psi')\sim_{s\to 0} n\gamma(s,\mathbf{1}_F,\psi')$ this ends the proof of the lemma. $\blacksquare$
\end{enumerate}

Let $\pi\in \Temp(G_n(E))$. Then we have $\cW(\pi,\psi_n)\subset \cC^w(N_n(E)\backslash G_n(E),\psi_n)$ and therefore for $W\in \cW(\pi,\psi_n)$ and $\phi\in \cS(F^n)$ the expression $Z(s,W,\phi)$ is well-defined for all $s\in \cH$. Recall that for $W\in \cW(\pi,\psi_n)$, $\widetilde{W}$ denotes the function defined by $\widetilde{W}(g)=W(w_n{}^t g^{-1})$ for all $g\in G_n(E)$. Note that $\widetilde{W}\in \cW(\pi^\vee,\psi_n^{-1})$ and that up to replacing $\psi$ by $\psi^{-1}$ the previous lemma ensures that $Z(s,\widetilde{W},\phi)$ is well-defined for all $s\in \cH$ and $\phi\in \cS(F^n)$. Finally, recall that $\phi\in \cS(F^n)\mapsto \widehat{\phi}\in \cS(F^n)$ denotes the Fourier transform with respect to $\psi'$ and the corresponding autodual measure. The following result is \cite[Theorem 3.4.1]{Beu3}.

\begin{theo}\label{theo1 Zeta function}
For every $W\in \cW(\pi,\psi_n)$ and $\phi\in \cS(F^n)$ and $s\in \cH$ with $\Re(s)<1$, we have 
$$\displaystyle Z(1-s,\widetilde{W},\widehat{\phi})=\omega_\pi(\tau)^{n-1}\lvert \tau\rvert_E^{\frac{n(n-1)}{2}(s-1/2)}\lambda_{E/F}(\psi')^{-\frac{n(n-1)}{2}}\gamma(s,\pi,\As,\psi')Z(s,W,\phi)$$
\end{theo}

Recall that to any representation $\pi\in \Temp(G_n(E))$ of the form $\pi=BC_n(\sigma)$, $\sigma\in \Temp(U(n))$, we have associated in the previous section a sign $c_1(\pi)$ satisfying \eqref{prop form betan}. We will need the following. 

\begin{lem}\label{lem 2 Zeta function}
Let $\sigma\in \Temp(U(n))$ and set $\pi=BC_n(\sigma)$. Then, we have
$$\displaystyle Z(1,\widetilde{W},\widehat{\phi})=\phi(0)c_1(\pi)\beta(W)$$
for all $W\in \cW(\pi,\psi_n)$ and $\phi\in \cS(F^n)$.
\end{lem}

\noindent\ul{Proof}: Let $W\in \cW(\pi,\psi_n)$ and $\phi\in \cS(F^n)$. We have
\[\begin{aligned}
\displaystyle Z(1,\widetilde{W},\widehat{\phi})=\int_{P_n(F)\backslash G_n(F)}\int_{N_n(F)\backslash P_n(F)} \widetilde{W}(ph)dp \, \widehat{\phi}(e_nh) \lvert \det h\rvert \, dh =\int_{P_n(F)\backslash G_n(F)} \beta(R(h)\widetilde{W})\, \widehat{\phi}(e_nh) \lvert \det h\rvert \, dh.
\end{aligned}\]
By \eqref{prop form betan}, the linear form $W\in \cW(\pi,\psi_n)\mapsto \beta(\widetilde{W})$ is invariant by $P_n(F)$ and its transpose hence by $G_n(F)$. Therefore,
\[\begin{aligned}
\displaystyle Z(1,\widetilde{W},\widehat{\phi})=\beta(\widetilde{W})\int_{P_n(F)\backslash G_n(F)} \widehat{\phi}(e_nh) \lvert \det h\rvert \, dh=\beta(\widetilde{W})\int_{F^n} \widehat{\phi}(x)\, dx=\phi(0)\beta(\widetilde{W})
\end{aligned}\]
and a new appeal to \eqref{prop form betan} ends the proof of the lemma. $\blacksquare$

\section{Computation of certain spectral distributions}\label{Part II}

The final goal of this chapter is Corollary \ref{cor lim spectrale} which in essence furnishes a spectral expansion for the period integral $\int_{N_n(F)\backslash G_n(F)} W(h) dh$ of a rapidly decaying Whittaker function $W$ on $G_n(E)$. The proof roughly goes as follows: using Lemma \ref{lem Zeta function} we are able to express this period as the residue of certain zeta integrals. For an argument of positive real part, these zeta integrals can in turn be spectrally decomposed using the Whittaker-Plancherel formula of Section \ref{section Plancherel Whitt} and applying the local functional equation recalled in Section \ref{Section Zeta integrals}, we roughly end up with computing the residue of an analytic family of Schwartz distributions on $\Temp(G_n(E))$ given in terms of gamma factors. By localization, this computation essentially reduces to the study of certain families of distributions on real vector spaces. The object of Sections \ref{Section PV} to \ref{Section distributions} is to calculate explicitely these residual distributions that are kind of multidimensional generalizations of a principal value integral. More precisely, Section \ref{Section PV} records some basic properties of such distributions, in Section \ref{Section polynomials} we write explicit decompositions of rational functions that show up in the computation and in Section \ref{Section distributions} we combine these identities with the simple lemmas proved in the first section to calculate the relevant residues explicitely. The link with the aforementioned distribution on $\Temp(G_n(E))$ given in terms of gamma factors is made in Section \ref{Section spectral limit}. Finally, we prove the main result in Section \ref{Section Corollary} as a corollary of the previous (rather long) computations.

\subsection{Principal values and poles of certain distributions}\label{Section PV}

Set
$$\displaystyle PV(\int_{i\R} \frac{\varphi(x)}{x}dx):=\lim\limits_{\epsilon\to 0^+} \int_{\lvert x\rvert>\epsilon} \frac{\varphi(x)}{x}dx$$
for every $\varphi\in \cS(i\R)$ (the limit is well-known, and actually easily seen, to exist). We record the following standard formula (see \cite[Eq. 4.4(6)]{GS})
\begin{align}\label{formula PV}
\displaystyle \lim\limits_{\substack{s\to 0^+}} \int_{i\R} \frac{\varphi(x)}{x+s}dx=PV(\int_{i\R} \frac{\varphi(x)}{x}dx)+\pi \varphi(0),\;\;\; \varphi\in \cS(i\R).
\end{align}
The next two results will be needed in the proof of Proposition \ref{prop 1 Distributions}.

\begin{lem}\label{lem cont family in Schwartz}
Let $W\subset V$ be real vector spaces and $\lambda_1,\ldots,\lambda_k$ be (real) linear forms on $V$ whose restrictions to $W$ are linearly independent. Fix Haar measures on $iW$ and $iV$. Let $\varphi\in \cS(iV)$ and set
$$\displaystyle \varphi_s(v)=\int_{iW} \frac{\varphi(v+w)}{\prod\limits_{i=1}^k (\lambda_i(v+w)+s)}dw$$
for all $v\in iV/iW$ and $s\in \cH$. Then, $\varphi_s\in \cS(iV/iW)$ for all $s\in \cH$ and the map $s\mapsto \varphi_s$ extends (uniquely) to a continuous map $s\in\cH\cup\{ 0\}\mapsto \varphi_s\in\cS(iV/iW)$. Moreover, if the function $v\mapsto \frac{\varphi(v)}{\prod_{i=1}^k \lambda_i(v)}$ belongs to $\cS(iV)$ we have
$$\displaystyle \varphi_0(0)=\int_{iW} \frac{\varphi(w)}{\prod\limits_{i=1}^k (\lambda_i(w))}dw.$$
\end{lem}

\noindent\ul{Proof}: Clearly, $\varphi_s\in \cS(iV/iW)$ and the linear map $\varphi\in \cS(iV)\mapsto \varphi_s\in \cS(iV/iW)$ is continuous for every $s\in \cH$. Let $X$ be a complement subspace of $W$ in $V$ on which $\lambda_1,\ldots,\lambda_k$ are trivial. As $\cS(iV)=\cS(iX)\widehat{\otimes} \cS(iW)$ (see \eqref{eq 2 proj tensor product}), by \eqref{eq 1 proj tensor product} for the first part of the lemma it suffices to show that for every $\varphi\in \cS(iW)$ the limit
$$\displaystyle \lim\limits_{s\to 0^+}\int_{iW} \frac{\varphi(w)}{\prod\limits_{i=1}^k (\lambda_i(w)+s)}dw$$
exists. For this, without loss of generality, we may assume that $W=\R^n$ and $\lambda_i(x_1,\ldots,x_n)=x_i$ for every $1\leqslant i\leqslant k$. Then, as $\cS(i\R^n)\simeq \widehat{\bigotimes\limits_{1\leqslant j\leqslant n}}\cS(i\R)$, by \eqref{eq 1 proj tensor product} again we may further assume that $\varphi(x_1,\ldots,x_n)=\varphi_1(x_1)\ldots\varphi_n(x_n)$ for some $\varphi_1,\ldots,\varphi_n\in \cS(i\R)$. We are then reduced to the case $n=k=1$ where the limit exists by \eqref{formula PV}. The second part of the lemma is an easy consequence of dominated convergence. $\blacksquare$

\begin{prop}\label{prop 1 PV}
Let $n\geqslant 1$ and $h_1,\ldots,h_n>0$. Let $(i\R^n)_0$ be the subspace of vectors $(x_1,\ldots,x_n)\in i\R^n$ with $x_1+\ldots+x_n=0$ and equip it with the unique Haar measure such that if $i\R^n$ is endowed with the Lebesgue measure then the quotient measure on $i\R^n/(i\R^n)_0\simeq i\R$ (where the isomorphism is given by $(x_1,\ldots,x_n)\mapsto x_1+\ldots+x_n$) is the Lebesgue measure on $i\R$. Then, we have
$$\displaystyle \lim\limits_{\substack{s\to 0^+ \\ s\in \cH}} s\int_{(i\R^n)_0}\frac{\varphi(\ux)}{\prod\limits_{i=1}^n (\frac{x_i}{h_i}+s)}d\ux=\frac{\prod\limits_{i=1}^n h_i}{\sum\limits_{i=1}^n h_i} (2\pi)^{n-1} \varphi(0)$$
for every $\varphi\in \cS(i\R^n)$.
\end{prop}

\noindent\ul{Proof}: First, since
$$\displaystyle \int_{(i\R^n)_0}\frac{\varphi(\ux)}{\prod\limits_{i=1}^n (\frac{x_i}{h_i}+s)}d\ux=\prod_{i=1}^nh_i  \int_{(i\R^n)_0}\frac{\varphi(\ux)}{\prod\limits_{i=1}^n (x_i+h_is)}d\ux$$
for every $s\in\cH$, it suffices to prove
\begin{align}\label{eq transformee}
\displaystyle \lim\limits_{\substack{s\to 0^+ \\ s\in \cH}} s\int_{(i\R^n)_0}\frac{\varphi(\ux)}{\prod\limits_{i=1}^n (x_i+h_is)}d\ux=\frac{(2\pi)^{n-1}}{\sum\limits_{i=1}^n h_i} \varphi(0)
\end{align}
for every $\varphi\in \cS(i\R^n)$. We do this by induction on $n$. The case $n=1$ is trivial but we treat also the case $n=2$ which will be needed for the general case. When $n=2$, we want
\begin{align}\label{eq cas n=2}
\displaystyle \lim\limits_{s\to 0^+} s\int_{i\R} \frac{\varphi(x,-x)}{(x+h_1s)(-x+h_2s)}dx=\frac{2\pi}{h_1+h_2}\varphi(0,0),\;\;\; \varphi\in\cS(i\R^2).
\end{align}
Fix $\varphi\in\cS(i\R^2)$. Then,
$$\displaystyle (h_1+h_2)s\int_{i\R} \frac{\varphi(x,-x)}{(x+h_1s)(-x+h_2s)}dx=\int_{i\R} \frac{\varphi(x,-x)}{x+h_1s}dx+\int_{i\R} \frac{\varphi(x,-x)}{-x+h_2s}dx$$
for every $s\in\cH$ and by \eqref{formula PV},
$$\displaystyle\lim\limits_{s\to 0^+} \int_{i\R} \frac{\varphi(x,-x)}{x+h_1s}dx=PV\left( \int_{i\R} \frac{\varphi(x,-x)}{x} dx\right) +\pi \varphi(0,0),$$
\[\begin{aligned}
\displaystyle\lim\limits_{s\to 0^+} \int_{i\R} \frac{\varphi(x,-x)}{-x+h_2s}dx & =PV\left( \int_{i\R} -\frac{\varphi(x,-x)}{x} dx\right) +\pi \varphi(0,0) \\
& =-PV\left( \int_{i\R} \frac{\varphi(x,-x)}{x} dx\right) +\pi \varphi(0,0).
\end{aligned}\]
The equality \eqref{eq cas n=2} follows.

We now treat the general case, assuming the result is known for $n-1$. Since $\cS(i\R^n)=\cS(i\R^{n-1})\widehat{\otimes}\cS(i\R)$, by \eqref{eq 1 proj tensor product} we may assume that there exist $\varphi'\in \cS(i\R^{n-1})$ and $\varphi''\in \cS(i\R)$ such that
$$\displaystyle \varphi(x_1,\ldots,x_n)=\varphi'(x_1,\ldots,x_{n-1})\varphi''(x_n),\;\;\; (x_1,\ldots,x_n)\in i\R^n$$
For all $a\in i\R$, set
$$\displaystyle (i\R^{n-1})_a:=\{(x_1,\ldots,x_{n-1})\in i\R^{n-1}\mid x_1+\ldots+x_{n-1}=a \}.$$
We equip this affine subspace with the invariant measure transferred from the one on $(i\R^{n-1})_0$. Then, we have
\begin{align}\label{eq PV1}
\displaystyle \int_{(i\R^n)_0}\frac{\varphi(\ux)}{\prod\limits_{i=1}^n (x_i+h_is)}d\ux=\int_{i\R} \int_{(i\R^{n-1})_a} \frac{\varphi'(\ux)}{\prod\limits_{i=1}^{n-1} (x_i+h_is)}d\ux \frac{\varphi''(-a)}{-a+h_ns}da
\end{align}
for every $s\in\cH$. Set
$$\displaystyle \varphi'_s(a):=(a+(\sum_{i=1}^{n-1}h_i)s)\int_{(i\R^{n-1})_a} \frac{\varphi'(\ux)}{\prod\limits_{i=1}^{n-1} (x_i+h_is)}d\ux$$
for every $a\in i\R$ and $s\in \cH$. Noticing that
$$\displaystyle \varphi'_s(a)=\sum_{i=1}^{n-1}\int_{(i\R^{n-1})_a} \frac{\varphi'(\ux)}{\prod\limits_{\substack{1\leqslant j\leqslant n-1 \\ j\neq i}} (x_i+h_is)}d\ux$$
and since for all $1\leqslant i\leqslant n-1$ the linear forms $\ux\mapsto x_j$ for $1\leqslant j\leqslant n-1$ and $j\neq i$ have linearly independent restrictions to $i\R^{n-1}$, by Lemma \ref{lem cont family in Schwartz} we see that the family $s\mapsto \varphi'_s$ extends to a continuous map $s\in \cH\cup\{ 0\}\mapsto \varphi'_s\in \cS(i\R)$. Moreover, by the induction hypothesis we have $\varphi'_0(0)=(2\pi)^{n-2} \varphi'(0)$. By \eqref{eq PV1}, we have
$$\displaystyle \int_{(i\R^n)_0}\frac{\varphi(\ux)}{\prod\limits_{i=1}^n (x_i+h_is)}d\ux=\int_{i\R} \frac{\varphi'_s(a)\varphi''(-a)}{(a+(\sum_{i=1}^{n-1}h_i)s)(-a+h_ns)}da.$$
Finally, by the $n=2$ case that we already treated and the uniform boundedness principle, we get
$$\displaystyle \lim\limits_{s\to 0^+} s\int_{i\R} \frac{\varphi'_s(a)\varphi''(-a)}{(a+(\sum_{i=1}^{n-1}h_i)s)(-a+h_ns)}da=\frac{2\pi}{\sum_{i=1}^{n}h_i}\varphi'_0(0)\varphi''(0)=\frac{(2\pi)^{n-1}}{\sum_{i=1}^{n}h_i}\varphi'(0)\varphi''(0)$$
As $\varphi'(0)\varphi''(0)=\varphi(0)$, this finishes the proof of \eqref{eq transformee} and thus of the proposition by induction. $\blacksquare$

\subsection{Some polynomial identities}\label{Section polynomials}

\begin{prop}\label{prop 1 Polynomial identities}
\begin{enumerate}[(i)]
\item Let $m\geqslant n$ be two nonnegative integers and define $P_{m,n}, S_{m,n}\in \mathbb{Q}(X_1,\ldots,X_m,X_1^*,\ldots,X_n^*)$ by
$$\displaystyle P_{m,n}:=\frac{\prod\limits_{1\leqslant j\neq k\leqslant m}(X_j-X_k)\prod\limits_{1\leqslant j\neq k\leqslant n} (X_j^*-X_k^*)}{\prod\limits_{\substack{1\leqslant j\leqslant m \\ 1\leqslant k\leqslant n}}(X_j+X_k^*)}$$
and
$$\displaystyle S_{m,n}:=\frac{\prod\limits_{1\leqslant j<k\leqslant n} (X_j^*-X_k^*)\prod\limits_{1\leqslant j<k\leqslant m}(X_j-X_k)\prod\limits_{n+1\leqslant k<j\leqslant m} (X_j-X_k)}{\prod\limits_{1\leqslant k\leqslant n} (X_k+X_k^*)}$$
Let $(w,P)\mapsto w\cdot P$ denote the natural action of $W=\fS_m\times \fS_n$ on $\mathbb{Q}(X_1,\ldots,X_m,X_1^*,\ldots,X_n^*)$ where $\fS_m$ acts by permutation of the variables $X_1,\ldots,X_m$ whereas $\fS_n$ acts by permutation of the variables $X_1^*,\ldots,X^*_n$. Let $W':=\fS_{m-n}\times \fS_n^{\diag}$ be the subgroup of elements $(\sigma_1,\sigma_2)\in W$ such that $\sigma_2^{-1}\sigma_1$ fixes $\{1,\ldots,n \}$ point-wise. Then, we have the identity
\begin{align}\label{eq Pol Id 1}
\displaystyle P_{m,n}=\frac{1}{\lvert W'\rvert}\sum_{w\in W} w\cdot S_{m,n}.
\end{align}
Moreover, set
$$\displaystyle P_{m,n}^*:=\left(\prod\limits_{1\leqslant k\leqslant n} (X_k+X_k^*) \right) P_{m,n},\; S_{m,n}^*:=\left(\prod\limits_{1\leqslant k\leqslant n} (X_k+X_k^*)\right) S_{m,n}$$
and
$$\displaystyle P_{m,n,s}(\underline{x}):=P_{m,n}(x_1+\frac{s}{2},\ldots,x_m+\frac{s}{2},x_1^*+\frac{s}{2},\ldots,x_n^*+\frac{s}{2})$$
$$\displaystyle S_{m,n,s}(\underline{x}):=S_{m,n}(x_1+\frac{s}{2},\ldots,x_m+\frac{s}{2},x_1^*+\frac{s}{2},\ldots,x_n^*+\frac{s}{2})$$
for every $s\in \cH$ and $\underline{x}=(x_1,\ldots,x_m,x_1^*,\ldots,x_n^*)\in i\R^m\times i\R^n$. Finally, let $V$ be the subspace $\{\underline{x}\in i\R^m\times i\R^n\mid x_k=-x_k^*, \;\forall \; 1\leqslant k\leqslant n \}$ of $i\R^m\times i\R^n$. Then, the function $\ux\in V\mapsto P_{m,n}^*(\ux)$ (which is a priori only well-defined on a Zariski open subset) extends to a polynomial function on $V$ and we have
\begin{align}\label{eq Pol Id 2}
\displaystyle \lim\limits_{s\to 0^+}s^nP_{m,n,s}(\underline{x})=\lim\limits_{s\to 0^+} s^nS_{m,n,s}(\underline{x})=P_{m,n}^*(\underline{x})=S_{m,n}^*(\underline{x})
\end{align}
for almost all $\underline{x}\in V$.

\item Let $p$ be a positive integer and define $Q_p,T_p\in \mathbb{Q}(Y_1,\ldots,Y_p)$ by
$$\displaystyle Q_p:=\frac{\prod\limits_{1\leqslant j\neq k\leqslant p} (Y_j-Y_k)}{\prod\limits_{1\leqslant j<k\leqslant p}(Y_j+Y_k)}$$
and
$$\displaystyle T_p:=(-1)^\epsilon \frac{\prod\limits_{1\leqslant j<k\leqslant p}(Y_j-Y_k)\prod\limits_{1\leqslant k\leqslant \floor*{\frac{p}{2}}}(Y_{p+1-k}-Y_k)}{\prod\limits_{1\leqslant k\leqslant \floor*{\frac{p}{2}}}(Y_k+Y_{p+1-k})}$$
where $\epsilon:=\frac{1}{2}(\frac{p(p-1)}{2}-\floor*{\frac{p}{2}})$. Let $(w,P)\mapsto w\cdot P$ denote the natural action of $W:=\mathfrak{S}_p$ on $\mathbb{Q}(Y_1,\ldots,Y_p)$ and  $W':=\fS_{\floor*{\frac{p}{2}}}\ltimes (\bZ/2\bZ)^{\floor*{\frac{p}{2}}}$ be the subgroup of $W$ preserving the partition \\
$\left\{\{\ell, p+1-\ell\}\mid 1\leqslant \ell\leqslant \ceil*{\frac{p}{2}} \right\}$ of $\{1,\ldots,p \}$. Then, we have the identity
\begin{align}\label{eq Pol Id 3}
\displaystyle Q_p=\frac{1}{\lvert W'\rvert} \sum_{w\in W}w\cdot T_p.
\end{align}
Moreover, set
$$\displaystyle Q_p^*:=\left(\prod\limits_{1\leqslant k\leqslant \floor*{\frac{p}{2}}} (Y_k+Y_{p+1-k}) \right) Q_p,\; T_p^*:=\left(\prod\limits_{1\leqslant k\leqslant \floor*{\frac{p}{2}}} (Y_k+Y_{p+1-k})\right) T_p$$
and
$$\displaystyle Q_{p,s}(\underline{y}):=Q_{p}(y_1+\frac{s}{2},\ldots,y_p+\frac{s}{2}),\; T_{p,s}(\underline{y}):=T_{p}(y_1+\frac{s}{2},\ldots,y_p+\frac{s}{2})$$
for every $s\in \cH$ and $\underline{y}=(y_1,\ldots,y_p)\in i\R^p$. Finally, let $V$ be the subspace $\{\underline{y}\in i\R^p\mid y_k=-y_{p+1-k}, \;\forall \; 1\leqslant k\leqslant \floor*{\frac{p}{2}} \}$ of $i\R^p$. Then, the function $\uy\in V\mapsto Q_p^*(\uy)$ (which is a priori only well-defined on a Zariski open subset) extends to a polynomial function on $V$ and we have
\begin{align}\label{eq Pol Id 4}
\displaystyle \lim\limits_{s\to 0^+}s^{\floor*{\frac{p}{2}}}Q_{p,s}(\underline{y})=\lim\limits_{s\to 0^+} s^{\floor*{\frac{p}{2}}}T_{p,s}(\underline{y})=Q_{p}^*(\underline{y})=T_{p}^*(\underline{y})
\end{align}
for almost all $\underline{y}\in V$.

\item Let $q$ be a positive integer and define $R_q,U_q\in \mathbb{Q}(Z_1,\ldots,Z_q)$ by
$$\displaystyle R_q:=\frac{\prod\limits_{1\leqslant j\neq k\leqslant q} (Z_j-Z_k)}{\prod\limits_{1\leqslant j\leqslant k\leqslant q}(Z_j+Z_k)}$$
and
$$\displaystyle U_q:=(-1)^\epsilon\frac{\prod\limits_{1\leqslant j<k\leqslant q}(Z_j-Z_k)}{\prod\limits_{1\leqslant k\leqslant \floor*{\frac{q}{2}}}(Z_k+Z_{q+1-k})\prod\limits_{1\leqslant k\leqslant \ceil*{\frac{q}{2}}}2Z_k}$$
where $\epsilon:=\frac{1}{2}(\frac{q(q-1)}{2}-\floor*{\frac{q}{2}})$. Let $(w,P)\mapsto w\cdot P$ denote the natural action of $W:=\mathfrak{S}_q$ on $\mathbb{Q}(Z_1,\ldots,Z_q)$ and  $W':=\fS_{\floor*{\frac{q}{2}}}\ltimes (\bZ/2\bZ)^{\floor*{\frac{q}{2}}}$ be the subgroup of $W$ preserving the partition \\
$\left\{\{\ell, q+1-\ell\}\mid 1\leqslant \ell\leqslant \ceil*{\frac{q}{2}} \right\}$ of $\{1,\ldots,q \}$. Then, we have the identity
\begin{align}\label{eq Pol Id 5}
\displaystyle R_q=\frac{1}{\lvert W'\rvert} \sum_{w\in W}w\cdot U_q.
\end{align}
Moreover, set
$$\displaystyle R_q^*:=\left(\prod\limits_{1\leqslant k\leqslant \ceil*{\frac{q}{2}}} (Z_k+Z_{q+1-k}) \right) R_q,\; U_q^*:=\left(\prod\limits_{1\leqslant k\leqslant \ceil*{\frac{q}{2}}} (Z_k+Z_{q+1-k})\right) U_q$$
and
$$\displaystyle R_{q,s}(\underline{z}):=R_{q}(z_1+\frac{s}{2},\ldots,z_q+\frac{s}{2}),\; U_{q,s}(\underline{z}):=U_{q}(z_1+\frac{s}{2},\ldots,z_q+\frac{s}{2})$$
for every $s\in \cH$ and $\underline{z}=(z_1,\ldots,z_q)\in i\R^q$. Finally, let $V$ be the subspace $\{\underline{z}\in i\R^q\mid z_k=-z_{q+1-k}, \;\forall \; 1\leqslant k\leqslant \ceil*{\frac{q}{2}} \}$ of $i\R^q$. Then, the function $\uz\in V\mapsto R_q^*(\uz)$ and $\uz\in V\mapsto U_q^*(\uz)$ (which are a priori only well-defined on a Zariski open subset) extends to polynomial functions on $V$ and we have
\begin{align}\label{eq Pol Id 6}
\displaystyle \lim\limits_{s\to 0^+}s^{\ceil*{\frac{q}{2}}}R_{q,s}(\underline{z})=\lim\limits_{s\to 0^+} s^{\ceil*{\frac{q}{2}}}U_{q,s}(\underline{z})=R_{q}^*(\underline{z})=U_{q}^*(\underline{z})
\end{align}
for almost all $\underline{z}\in V$.
\end{enumerate}
\end{prop}

\vspace{2mm}

\noindent\ul{Proof}:
\begin{enumerate}[(i)]
\item Set
$$\displaystyle \Delta:=\prod_{1\leqslant j<k\leqslant m} (X_j-X_k),\;\;\; \Delta^*:=\prod_{1\leqslant j<k\leqslant n} (X^*_j-X^*_k)$$
and $\sgn: W\to \{\pm 1 \}$ be the product of the sign characters on $\fS_m$ and $\fS_n$. Then,
$$\displaystyle w\cdot(\Delta\Delta^*)=\sgn(w)\Delta\Delta^*$$
for every $w\in W$. Hence,
$$\displaystyle \sum_{w\in W} w\cdot S_{m,n}=\Delta\Delta^*\sum_{w\in W}\sgn(w)w\cdot\left( \frac{\prod\limits_{n+1\leqslant k<j\leqslant m} (X_j-X_k)}{\prod\limits_{1\leqslant k\leqslant n} (X_k+X_k^*)}\right)$$
can be written as $\displaystyle \Delta\Delta^*\frac{Q}{\prod\limits_{\substack{1\leqslant j\leqslant m \\ 1\leqslant k\leqslant n}}(X_j+X_k^*)}$ where $Q\in \Q[X_1,\ldots,X_m,X_1^*,\ldots,X_n^*]$ is homogeneous of total degree
$$\displaystyle \frac{(m-n)(m-n-1)}{2}-n+mn=\frac{m(m-1)}{2}+\frac{n(n-1)}{2}=\deg(\Delta)+\deg(\Delta^*)$$
and satisfies $w\cdot Q=\sgn(w)Q$ for every $w\in W$. It follows that $Q$ is a scalar multiple of $\Delta\Delta^*$ and therefore there exists $c\in \Q$ such that
$$\displaystyle cP_{m,n}=\frac{1}{\lvert W'\rvert} \sum_{w\in W} w\cdot S_{m,n}$$
Now, noticing that $W'$ is the stabilizer of $S_{m,n}$ in $W$ and that for every $w\in W\setminus W'$ the restriction of
$$\displaystyle \left(\prod\limits_{1\leqslant k\leqslant n} (X_k+X_k^*) \right) w\cdot S_{m,n}$$
to $V$ vanishes almost everywhere, we also have
$$\displaystyle cP^*_{m,n}(\ux)=S_{m,n}^*(\ux)$$
for almost all $\ux\in V$. On the other hand, it is clear that
$$\displaystyle \lim\limits_{s\to 0^+}s^nP_{m,n,s}(\ux)=P_{m,n}^*(\ux) \mbox{ and } \lim\limits_{s\to 0^+}s^nS_{m,n,s}(\ux)=S_{m,n}^*(\ux)$$
for almost all $\ux\in V$. Hence, to get \eqref{eq Pol Id 1} and \eqref{eq Pol Id 2} it only remains to prove that $P^*_{m,n}(\ux)=S_{m,n}^*(\ux)$ almost everywhere on $V$. For a generic point $\ux=(x_1,\ldots,x_m,x_1^*,\ldots,x_n^*)\in V$, we have
\[\begin{aligned}
\displaystyle P_{m,n}^*(\ux) & =S_{m,n}^*(\ux)\frac{\prod\limits_{1\leqslant j<k\leqslant n}(x_k^*-x_j^*)\prod\limits_{1\leqslant k<j\leqslant n}(x_j-x_k)\prod\limits_{\substack{1\leqslant k\leqslant n \\ n+1\leqslant j\leqslant m}}(x_j-x_k)}{\prod\limits_{\substack{1\leqslant j\leqslant m \\ 1\leqslant k\leqslant n \\ j\neq k}}(x_j+x_k^*)} \\
 & =S_{m,n}^*(\ux)\frac{\prod\limits_{1\leqslant j<k\leqslant n}(x_k^*+x_j)\prod\limits_{1\leqslant k<j\leqslant n}(x_j+x^*_k)\prod\limits_{\substack{1\leqslant k\leqslant n \\ n+1\leqslant j\leqslant m}}(x_j+x^*_k)}{\prod\limits_{\substack{1\leqslant j\leqslant m \\ 1\leqslant k\leqslant n \\ j\neq k}}(x_j+x_k^*)} \\
 & =S_{m,n}^*(\ux)
\end{aligned}\]
and the claim follows.
\item This time we set $\displaystyle \Delta:=\prod_{1\leqslant j<k\leqslant p}(Y_j-Y_k)$ and let $\sgn:W\to \{\pm 1 \}$ be the sign character. Then, $w\cdot \Delta=\sgn(w)\Delta$ for every $w\in W$. Therefore,
$$\displaystyle \sum_{w\in W}w\cdot T_p=(-1)^\epsilon\Delta\sum_{w\in W}\sgn(w) w\cdot\left(\frac{\prod\limits_{1\leqslant k\leqslant \floor*{\frac{p}{2}}}(Y_{p+1-k}-Y_k)}{\prod\limits_{1\leqslant k\leqslant \floor*{\frac{p}{2}}}(Y_k+Y_{p+1-k})} \right)$$
can be written as $\displaystyle \Delta \frac{P}{\prod\limits_{1\leqslant j< k\leqslant p}(Y_j+Y_k)}$ where $P\in \Q[Y_1,\ldots,Y_p]$ is homogeneous of degree $\floor*{\frac{p}{2}}-\floor*{\frac{p}{2}}+\frac{p(p-1)}{2}=\deg(\Delta)$ and satisfies $w\cdot P=\sgn(w)P$ for every $w\in W$. It follows that $P$ is a scalar multiple of $\Delta$ and thus there exists $c\in \Q$ such that
$$\displaystyle cQ_{p}=\frac{1}{\lvert W'\rvert} \sum_{w\in W} w\cdot T_p$$
Noticing that $W'$ is the stabilizer of $T_{p}$ in $W$ and that for all $w\in W\setminus W'$ the restriction of
$$\displaystyle \left(\prod\limits_{1\leqslant k\leqslant \floor*{\frac{p}{2}}} (Y_k+Y_{p+1-k}) \right) w\cdot T_{p}$$
to $V$ vanishes almost everywhere, we also have
$$\displaystyle cQ^*_{p}(\uy)=T_{p}^*(\uy)$$
for almost all $\uy\in V$. On the other hand, it is clear that
$$\displaystyle \lim\limits_{s\to 0^+}s^{\floor*{\frac{p}{2}}}Q_{p,s}(\uy)=Q_{p}^*(\uy) \mbox{ and } \lim\limits_{s\to 0^+}s^{\floor*{\frac{p}{2}}}T_{p,s}(\uy)=T_{p}^*(\uy)$$
for almost all $\uy\in V$. Hence, to get \eqref{eq Pol Id 3} and \eqref{eq Pol Id 4} it only remains to prove that $Q^*_{p}(\uy)=T_{p}^*(\uy)$ almost everywhere on $V$. For a generic point $\uy=(y_1,\ldots,y_p)\in V$, we have
\[\begin{aligned}
\displaystyle Q_p^*(\uy) & =(-1)^\epsilon T_p^*(\uy)\frac{\prod\limits_{\substack{1\leqslant k<j\leqslant p \\ j\neq p+1-k}}(y_j-y_k)}{\prod\limits_{\substack{1\leqslant j< k\leqslant p \\ k\neq p+1-j}} (y_j+y_k)} \\
 & =(-1)^\epsilon T_p^*(\uy)\frac{\prod\limits_{\substack{1\leqslant k<j\leqslant p \\ j< p+1-k}}(y_j-y_k)\times \prod\limits_{\substack{1\leqslant k<j\leqslant p \\ k> p+1-j}}(y_j-y_k)}{\prod\limits_{\substack{1\leqslant j< k\leqslant p \\ k\neq p+1-j}} (y_j+y_k)} \\
 & =(-1)^\epsilon T_p^*(\uy)\frac{\prod\limits_{\substack{1\leqslant j<k\leqslant p \\ j> p+1-k}}(y_j+y_k)\times \prod\limits_{\substack{1\leqslant j<k\leqslant p \\ k< p+1-j}}(-y_j-y_k)}{\prod\limits_{\substack{1\leqslant j< k\leqslant p \\ k\neq p+1-j}} (y_j+y_k)} \\
 & =(-1)^\epsilon T_p^*(\uy)\times (-1)^\epsilon=T_p^*(\uy)
\end{aligned}\]
and the claim follows.

\item Notice that
\begin{align}\label{eq Pol Id 7bis}
\displaystyle R_q(Z_1,\ldots,Z_q)=\frac{Q_q(Z_1,\ldots,Z_q)}{\prod\limits_{1\leqslant j\leqslant q}2Z_j}
\end{align}
Therefore, as $\prod\limits_{1\leqslant j\leqslant q} 2Z_j$ is $W$-invariant, by (ii) we have
\begin{align}\label{eq Pol Id 7}
\displaystyle R_q(Z_1,\ldots,Z_q)=\frac{1}{\lvert W'\rvert}\sum_{w\in W}w\cdot\left(\frac{T_q(Z_1,\ldots,Z_q)}{\prod\limits_{1\leqslant j\leqslant q}2Z_j} \right)
\end{align}
Moreover,
\begin{align}\label{eq Pol Id 8}
\displaystyle \frac{T_q(Z_1,\ldots,Z_q)}{\prod\limits_{1\leqslant j\leqslant q}2Z_j}=(-1)^{\epsilon}\left(\prod\limits_{1\leqslant j\leqslant \floor*{\frac{q}{2}}}\frac{Z_{q+1-j}-Z_j}{4Z_jZ_{q+1-j}}\right)\times \frac{\Delta}{\prod\limits_{1\leqslant j\leqslant \ceil*{\frac{q}{2}}}(Z_j+Z_{q+1-j})}
\end{align}
where we have set $\displaystyle \Delta:=\prod\limits_{1\leqslant j< k\leqslant q}(Z_j-Z_k)$ and
\begin{align}\label{eq Pol Id 9}
\displaystyle \prod\limits_{1\leqslant j\leqslant \floor*{\frac{q}{2}}}\frac{Z_{q+1-j}-Z_j}{4Z_jZ_{q+1-j}} & =2^{-\floor*{\frac{q}{2}}}\prod\limits_{1\leqslant j\leqslant \floor*{\frac{q}{2}}}(\frac{1}{2Z_j}-s_j\cdot \frac{1}{2Z_j}) \\
\nonumber & =2^{-\floor*{\frac{q}{2}}}\sum_{s\in (\bZ/2\bZ)^{\floor*{\frac{q}{2}}}} \sgn(s) s\cdot\left(\frac{1}{\prod\limits_{1\leqslant j\leqslant \floor*{\frac{q}{2}}}2Z_j}\right)
\end{align}
where $s_j$ denotes the transposition $(j,q+1-j)$, we have identified $(\bZ/2\bZ)^{\floor*{\frac{q}{2}}}$ with the subgroup of $W$ generated by the $s_j$ for $1\leqslant j\leqslant \floor*{\frac{q}{2}}$ and $\sgn:W\to \{\pm 1 \}$ denotes the sign character. Since
$$\displaystyle s\cdot\left( \frac{\Delta}{\prod\limits_{1\leqslant j\leqslant \ceil*{\frac{q}{2}}}(Z_j+Z_{q+1-j})}\right)=\sgn(s) \frac{\Delta}{\prod\limits_{1\leqslant j\leqslant \ceil*{\frac{q}{2}}}(Z_j+Z_{q+1-j})}$$
for every $s \in (\bZ/2\bZ)^{\floor*{\frac{q}{2}}}$, we obtain from \eqref{eq Pol Id 8} and \eqref{eq Pol Id 9} that
\[\begin{aligned}
\displaystyle \frac{T_q(Z_1,\ldots,Z_q)}{\prod\limits_{1\leqslant j\leqslant q}2Z_j} & =2^{-\floor*{\frac{q}{2}}}\sum_{s \in (\bZ/2\bZ)^{\floor*{\frac{q}{2}}}}(-1)^\epsilon s\cdot\left( \frac{\Delta}{\prod\limits_{1\leqslant j\leqslant \ceil*{\frac{q}{2}}}(Z_j+Z_{q+1-j})\prod\limits_{1\leqslant j\leqslant \floor*{\frac{q}{2}}}2Z_j}\right) \\
 & =2^{-\floor*{\frac{q}{2}}}\sum_{s \in (\bZ/2\bZ)^{\floor*{\frac{q}{2}}}} s\cdot U_q
\end{aligned}\]
Hence, by \eqref{eq Pol Id 7},
\[\begin{aligned}
\displaystyle R_q(Z_1,\ldots,Z_q)=\frac{2^{-\floor*{\frac{q}{2}}}}{\lvert W'\rvert}\sum_{w\in W}\sum_{s\in (\bZ/2\bZ)^{\floor*{\frac{q}{2}}}} ws\cdot U_q=\frac{1}{\lvert W'\rvert}\sum_{w\in W} w\cdot U_q
\end{aligned}\]
and this proves \eqref{eq Pol Id 5}.

Once again, it is obvious that
$$\displaystyle \lim\limits_{s\to 0^+}s^{\ceil*{\frac{q}{2}}}R_{q,s}(\underline{z})=R_{q}^*(\underline{z}) \mbox{ and } \lim\limits_{s\to 0^+} s^{\ceil*{\frac{q}{2}}}U_{q,s}(\underline{z})=U_{q}^*(\underline{z})$$
for almost all $\uz\in V$. By \eqref{eq Pol Id 7bis}, we get
$$\displaystyle R^*_q(Z_1,\ldots,Z_q)=\frac{Q^*_q(Z_1,\ldots,Z_q)}{\prod\limits_{1\leqslant j\leqslant \floor*{\frac{q}{2}}}4Z_jZ_{q+1-j}}$$
Similarly, we have
$$\displaystyle T_q^*(Z_1,\ldots,Z_q)=\left(\prod\limits_{1\leqslant j\leqslant \floor*{\frac{q}{2}}} 2Z_j(Z_{q+1-j}-Z_j)\right) U_q^*(Z_1,\ldots,Z_q)$$
Therefore, by (ii) we obtain
\[\begin{aligned}
\displaystyle R_q^*(\uz) & =\frac{Q^*_q(\uz)}{\prod\limits_{1\leqslant j\leqslant \floor*{\frac{q}{2}}}4z_jz_{q+1-j}}=\frac{T^*_q(\uz)}{\prod\limits_{1\leqslant j\leqslant \floor*{\frac{q}{2}}}4z_jz_{q+1-j}} \\
 & =\left(\prod\limits_{1\leqslant j\leqslant \floor*{\frac{q}{2}}}\frac{2z_j(z_{q+1-j}-z_j)}{4z_jz_{q+1-j}}\right)U_q^*(\uz) \\
 & =\left(\prod\limits_{1\leqslant j\leqslant \floor*{\frac{q}{2}}}\frac{2z_j\times 2z_{q+1-j}}{4z_jz_{q+1-j}}\right)U_q^*(\uz)=U_q^*(\uz)
\end{aligned}\]
for almost all $\uz \in V$ and \eqref{eq Pol Id 6} follows. We show similarly that
$$\displaystyle U_q^*(\uz)=(-1)^\epsilon \prod_{\substack{1\leqslant j< k\leqslant q \\ k\neq q+1-j}}(z_j-z_k)$$
for almost all $\uz \in V$ and consequently $U_q^*(\uz)$ extends to a polynomial function on $V$. $\blacksquare$
\end{enumerate}

\subsection{Explicit computation of certain residual distributions}\label{Section distributions}

In this section, we fix the following data:
\begin{itemize}
\item $r,s,t\in \bN$ three nonnegative integers;
\item $(m_i,n_i)\in \bN^*\times \bN$ such that $m_i\geqslant n_i$ and $d_i\in \bN^*$ for all $1\leqslant i\leqslant r$;
\item $p_j\in \bN^*$ and $e_j\in \bN^*$ for all $1\leqslant j\leqslant s$;
\item $q_k\in \bN^*$ and $f_k\in \bN^*$ for all $1\leqslant k\leqslant t$;
\end{itemize}
and we set
$$\displaystyle \cA:=\prod\limits_{1\leqslant i\leqslant r} (i\R)^{m_i}\times (i\R)^{n_i}\times \prod\limits_{1\leqslant j\leqslant s} (i\R)^{p_j}\times \prod\limits_{1\leqslant k\leqslant t} (i\R)^{q_k}.$$
For any $\lambda\in \cA$, we will write $\ux_i(\lambda)$ ($1\leqslant i\leqslant r$), $\uy_j(\lambda)$ ($1\leqslant j\leqslant s$) and $\uz_k(\lambda)$ ($1\leqslant k\leqslant t$) for the projections of $\lambda$ onto $(i\R)^{m_i}\times (i\R)^{n_i}$, $(i\R)^{p_j}$ and $(i\R)^{q_k}$ respectively. We will further write the coordinate of these vectors as follows:
$$\displaystyle \ux_i(\lambda)=(x_{i,1}(\lambda),\ldots,x_{i,m_i}(\lambda),x_{i,1}^*(\lambda),\ldots,x_{i,n_i}^*(\lambda))\in (i\R)^{m_i}\times (i\R)^{n_i}$$
$$\displaystyle \uy_j(\lambda)=(y_{j,1}(\lambda),\ldots,y_{j,p_j}(\lambda))\in (i\R)^{p_j}$$
$$\displaystyle \uz_k(\lambda)=(z_{k,1}(\lambda),\ldots,z_{k,q_k}(\lambda))\in (i\R)^{q_k}$$
We equip $\cA$ with the product of the Lebesgue measure on $i\R$. Let $\Sigma:\cA\to i\R$ be the linear form given by
$$\displaystyle \Sigma(\lambda)=\sum_{i=1}^{r}\left(\sum_{\ell=1}^{m_i}x_{i,\ell}(\lambda)+\sum_{\ell=1}^{n_i}x_{i,\ell}^*(\lambda)\right)+\sum_{j=1}^s \sum_{\ell=1}^{p_j}y_{j,\ell}(\lambda)+\sum_{k=1}^t \sum_{\ell=1}^{q_k}z_{k,\ell}(\lambda)$$
and $\cA_0:=\Ker(\Sigma)$. We equip $\cA_0$ with the unique Haar measure such that the quotient measure on $\cA/\cA_0\simeq i\R$ (the isomorphism being induced by $\Sigma$) is the Lebesgue measure.

For all $1\leqslant i\leqslant r$ (resp. $1\leqslant i\leqslant r$, $1\leqslant j\leqslant s$ and $1\leqslant k\leqslant t$), we let $\fS_{m_i}$ (resp. $\fS_{n_i}$, $\fS_{p_j}$ and $\fS_{q_k}$) act on $(i\R)^{m_i}$ (resp. $(i\R)^{n_i}$, $(i\R)^{p_j}$ and $(i\R)^{q_k}$) by permutation of the coordinates and we set
$$\displaystyle W:=\prod_{i=1}^r \left(\fS_{m_i}\times \fS_{n_i}\right)\times \prod_{j=1}^s \fS_{p_j}\times \prod_{k=1}^t \fS_{q_k}.$$
Then $W$ acts on $\cA$ by the product of the previous actions and therefore also on the Schwartz space $\cS(\cA)$. Denote by $\cS(\cA)^W$ the subspace of $W$-invariant functions. Using the notation of Proposition \ref{prop 1 Polynomial identities}, for every $s\in \cH$ we define a distribution $D_s\in \cS(\cA)'$ by
$$\displaystyle D_s(\varphi):=\int_{\cA_0}\varphi(\lambda) \prod_{i=1}^r P_{m_i,n_i,s}(\frac{\ux_i(\lambda)}{d_i})\times \prod_{j=1}^s Q_{p_j,s}(\frac{\uy_j(\lambda)}{e_j})\times \prod_{k=1}^t R_{q_k,s}(\frac{\uz_k(\lambda)}{f_k})d\lambda,\;\;\; \varphi\in \cS(\cA).$$

Let $\mathcal{A}'$ be the subspace of $\cA$ defined by the relations
\begin{itemize}
\item $x_{i,\ell}(\lambda)+x_{i,\ell}^*(\lambda)=0$ for every $1\leqslant i\leqslant r$ and $1\leqslant \ell\leqslant n_i$;
\item $y_{j,\ell}(\lambda)+y_{j,p_j+1-\ell}(\lambda)=0$ for every $1\leqslant j\leqslant s$ and $1\leqslant \ell\leqslant \floor*{\frac{p_j}{2}}$;
\item $z_{k,\ell}(\lambda)+z_{k,q_k+1-\ell}(\lambda)=0$ for every $1\leqslant k\leqslant t$ and $1\leqslant \ell\leqslant \ceil*{\frac{q_k}{2}}$.
\end{itemize}
The map $\lambda\mapsto \left( (x_{i,\ell}(\lambda))_{\substack{1\leqslant i\leqslant r\\ 1\leqslant \ell\leqslant n_i}}, (y_{j,\ell}(\lambda))_{\substack{1\leqslant j\leqslant s\\ 1\leqslant \ell\leqslant \ceil*{\frac{p_j}{2}}}}, (z_{k,\ell}(\lambda))_{\substack{ 1\leqslant k\leqslant t \\ 1\leqslant \ell \leqslant \floor*{\frac{q_k}{2}}}}\right)$ induces an isomorphism
$$\displaystyle \cA'\simeq \prod_{i=1}^r (i\R)^{n_i}\times \prod_{j=1}^s (i\R)^{\ceil*{\frac{p_j}{2}}}\times \prod_{k=1}^t (i\R)^{\floor*{\frac{q_k}{2}}}$$
and we equip $\cA'$ with the measure which transfer to the Lebesgue measure via this isomorphism. Define the following subgroup of $W$:
$$\displaystyle W':=\prod_{i=1}^r \left(\fS_{n_i}^{\diag}\times \fS_{m_i-n_i}\right) \times \prod_{j=1}^s \left(\fS_{\floor*{\frac{p_j}{2}}}\ltimes (\bZ/2\bZ)^{\floor*{\frac{p_j}{2}}}\right) \times \prod_{k=1}^t \left(\fS_{\floor*{\frac{q_k}{2}}}\ltimes (\bZ/2\bZ)^{\floor*{\frac{q_k}{2}}}\right)$$
where for all integers $m\geqslant n$ and $p$, $\fS_n^{\diag}\times \fS_{m-n}$ denotes the subgroup of elements $(\sigma,\tau)\in\fS_m\times\fS_n$ such that $\tau^{-1}\sigma$ fixes $\{1,\ldots,n \}$ point-wise and we identify $\fS_{\floor*{\frac{p}{2}}}\ltimes (\bZ/2\bZ)^{\floor*{\frac{p}{2}}}$ with the subgroup of $\fS_p$ preserving the partition $\left\{\{\ell, p+1-\ell\}\mid 1\leqslant \ell\leqslant \ceil*{\frac{p}{2}} \right\}$ of $\{1,\ldots,p \}$. It is easy to check that $\cA'$ is $W'$-invariant. Finally, we define a distribution $D'\in \cS(\cA)'$ by
$$\displaystyle D'(\varphi):=\frac{D}{n}(2\pi)^{N-1}2^{1-c}\int_{\cA'} \varphi(\mu) \lim\limits_{s\to 0^+} s^N\prod_{i=1}^r P_{m_i,n_i,s}(\frac{\ux_i(\mu)}{d_i})\times \prod_{j=1}^s Q_{p_j,s}(\frac{\uy_j(\mu)}{e_j})\times \prod_{k=1}^t R_{q_k,s}(\frac{\uz_k(\mu)}{f_k})d\mu$$
for all $\varphi\in \cS(\cA)$, where
\begin{itemize}
\item $\displaystyle n:=\sum_{i=1}^r (n_i+m_i)d_i+\sum_{j=1}^s p_j e_j+\sum_{k=1}^t q_kf_k$;
\item $\displaystyle D:=\prod_{i=1}^r d_i^{n_i} \times \prod_{j=1}^s e_j^{\frac{p_j}{2}}\times \prod_{k=1}^t f_k^{\ceil*{\frac{q_k}{2}}}$;
\item $\displaystyle c:=\left\lvert \{1\leqslant k\leqslant t\mid q_k\equiv 1\; [2] \}\right\rvert$;
\item $\displaystyle N:=\sum_{i=1}^r n_i+\sum_{j=1}^s \floor*{\frac{p_j}{2}}+\sum_{k=1}^t \ceil*{\frac{q_k}{2}}$.
\end{itemize}
Note that by Proposition \ref{prop 1 Polynomial identities}, we have
\begin{align}\label{eq Distr0}
\displaystyle \lim\limits_{s\to 0^+} s^N & \prod_{i=1}^r P_{m_i,n_i,s}(\frac{\ux_i(\mu)}{d_i})\times \prod_{j=1}^s Q_{p_j,s}(\frac{\uy_j(\mu)}{e_j})\times \prod_{k=1}^t R_{q_k,s}(\frac{\uz_k(\mu)}{f_k})= \\
\nonumber & \prod_{i=1}^r S^*_{m_i,n_i}(\frac{\ux_i(\mu)}{d_i})\times \prod_{j=1}^s T_{p_j}^*(\frac{\uy_j(\mu)}{e_j})\times \prod_{k=1}^t U_{q_k}^*(\frac{\uz_k(\mu)}{f_k})
\end{align}
for almost all $\mu\in \cA'$ and that this extends to a polynomial function on $\cA'$ so that the distribution $D'$ is well-defined.

\begin{prop}\label{prop 1 Distributions}
For every $\varphi\in \cS(\cA)^W$, we have
$$\displaystyle \lim\limits_{s\to 0^+} sD_s(\varphi)=\left\{
    \begin{array}{ll}
        \frac{\lvert W\rvert}{\lvert W'\rvert}D'(\varphi) & \mbox{if } m_i=n_i \mbox{ for all } 1\leqslant i\leqslant r \mbox{ and } p_j \mbox{ is even for all } 1\leqslant j\leqslant s, \\
         \\
        0 & \mbox{otherwise.}
    \end{array}
\right.
$$

\end{prop}

\noindent\ul{Proof}: Let $\varphi\in \cS(\cA)^W$. By Proposition \ref{prop 1 Polynomial identities}, we have
\[\begin{aligned}
\displaystyle & \prod_{i=1}^r P_{m_i,n_i,s}(\frac{\ux_i(\lambda)}{d_i})\times \prod_{j=1}^s Q_{p_j,s}(\frac{\uy_j(\lambda)}{e_j})\times \prod_{k=1}^t R_{q_k,s}(\frac{\uz_k(\lambda)}{f_k})= \\
 & \frac{1}{\lvert W'\rvert}\sum_{w\in W} \prod_{i=1}^r S_{m_i,n_i,s}(\frac{\ux_i(w\lambda)}{d_i})\times \prod_{j=1}^s T_{p_j,s}(\frac{\uy_j(w\lambda)}{e_j})\times \prod_{k=1}^t U_{q_k,s}(\frac{\uz_k(w\lambda)}{f_k})
\end{aligned}\]
for every $\lambda\in \cA$ and $s\in \cH$. Therefore, as $\varphi$ is $W$-invariant, setting
$$\displaystyle \widetilde{D}_s(\varphi)=\int_{\cA_0}\varphi(\lambda) \prod_{i=1}^r S_{m_i,n_i,s}(\frac{\ux_i(\lambda)}{d_i})\times \prod_{j=1}^s T_{p_j,s}(\frac{\uy_j(\lambda)}{e_j})\times \prod_{k=1}^t U_{q_k,s}(\frac{\uz_k(\lambda)}{f_k})d\lambda,\;\;\; s\in \cH$$
we need to show that
\begin{align}\label{eq Distr1}
\displaystyle \lim\limits_{s\to 0^+} s\widetilde{D}_s(\varphi)=\left\{
    \begin{array}{ll}
        D'(\varphi) & \mbox{if } m_i=n_i \mbox{ for all } 1\leqslant i\leqslant r \mbox{ and } p_j \mbox{ is even for all } 1\leqslant j\leqslant s \\
        0 & \mbox{otherwise.}
    \end{array}
\right.
\end{align}

Set
$$\displaystyle \psi(\lambda)=\varphi(\lambda)\prod_{i=1}^r S^*_{m_i,n_i}(\frac{\ux_i(\lambda)}{d_i})\prod_{j=1}^s T_{p_j}^*(\frac{\uy_j(\lambda)}{e_j})\prod_{k=1}^t \left(U_{q_k}^*(\frac{\uz_k(\lambda)}{f_k})\prod\limits_{1\leqslant \ell\leqslant \floor*{\frac{q_k}{2}}} \frac{2z_{k,\ell}(\lambda)}{f_k}  \right),\;\;\; \lambda\in \cA$$
Notice that $\psi$ being the product of $\varphi$ with a polynomial function on $\cA$ we have $\psi\in \cS(\cA)$. Moreover,
\[\begin{aligned}
\displaystyle \widetilde{D}_s(\varphi)=\int_{\cA_0}\psi(\lambda)  & \prod_{i=1}^r \prod\limits_{1\leqslant \ell\leqslant n_i} \left(\frac{x_{i,\ell}(\lambda)+x_{i,\ell}^*(\lambda)}{d_i}+s \right)^{-1} \prod_{j=1}^s \prod\limits_{1\leqslant \ell\leqslant \floor*{\frac{p_j}{2}}} \left( \frac{y_{j,\ell}(\lambda)+y_{j,p_j+1-\ell}(\lambda)}{e_j}+s\right)^{-1} \\
 & \prod_{k=1}^t \prod\limits_{1\leqslant \ell\leqslant \ceil*{\frac{q_k}{2}}} \left(\frac{z_{k,\ell}(\lambda)+z_{k,q_k+1-\ell}(\lambda)}{f_k}+s \right)^{-1} \prod\limits_{1\leqslant \ell\leqslant \floor*{\frac{q_k}{2}}} \left(\frac{2z_{k,\ell}(\lambda)}{f_k}+s \right)^{-1} d\lambda.
\end{aligned}\]

Let $\cA_0':=\cA'\cap \cA_0$. Fixing a Haar measure on $\cA_0'$ which coincides with the one we fixed on $\cA'$ if $\cA_0'=\cA'$, we define
$$\displaystyle \psi_s(\lambda)=\int_{\cA_0'} \psi(\lambda+\mu)\prod_{k=1}^t\prod\limits_{1\leqslant \ell\leqslant \floor*{\frac{q_k}{2}}} \left(\frac{2z_{k,\ell}(\lambda+\mu)}{f_k}+s \right)^{-1} d\mu$$
for every $s\in \cH$ and $\lambda\in \cA_0/\cA'_0$. Then, equipping $\cA_0/\cA_0'$ with the quotient measure, we have
\begin{align}\label{eq Distr2}
\displaystyle \widetilde{D}_s(\varphi)=\int_{\cA_0/\cA_0'}\psi_s(\lambda)  & \prod_{i=1}^r \prod\limits_{1\leqslant \ell\leqslant n_i} \left(\frac{x_{i,\ell}(\lambda)+x_{i,\ell}^*(\lambda)}{d_i}+s \right)^{-1} \prod_{j=1}^s \prod\limits_{1\leqslant \ell\leqslant \floor*{\frac{p_j}{2}}} \left( \frac{y_{j,\ell}(\lambda)+y_{j,p_j+1-\ell}(\lambda)}{e_j}+s\right)^{-1} \\
\nonumber & \prod_{k=1}^t \prod\limits_{1\leqslant \ell\leqslant \ceil*{\frac{q_k}{2}}} \left(\frac{z_{k,\ell}(\lambda)+z_{k,q_k+1-\ell}(\lambda)}{f_k}+s \right)^{-1} d\lambda
\end{align}
for every $s\in \cH$.

We readily check that the linear forms $\lambda\mapsto z_{k,\ell}(\lambda)$ ($1\leqslant k\leqslant t$, $1\leqslant \ell\leqslant \floor*{\frac{q_k}{2}}$) are linearly independent on $\cA_0'$. Hence, by Lemma \ref{lem cont family in Schwartz}, the family $s\mapsto \psi_s$ extends to a continuous map $\cH\cup\{ 0\}\to \cS(\cA_0/\cA'_0)$. If there exists $1\leqslant i\leqslant r$ such that $m_i\neq n_i$ or $1\leqslant j\leqslant s$ such that $p_j\notin 2\bN$, then the linear forms $x_{i,\ell}(.)+x_{i,\ell}^*(.)$ ($1\leqslant i\leqslant r$, $1\leqslant \ell\leqslant n_i$), $y_{j,\ell}(.)+y_{j,p_j+1-\ell}(.)$ ($1\leqslant j\leqslant s$, $1\leqslant \ell\leqslant \floor*{\frac{p_j}{2}}$) and $z_{k,\ell}(.)+z_{k,q_k+1-\ell}(.)$ ($1\leqslant k\leqslant t$, $1\leqslant \ell\leqslant \ceil*{\frac{q_k}{2}}$) restricted to $\cA_0$ are linearly independent. Hence, by \eqref{eq Distr2}, Lemma \ref{lem cont family in Schwartz} and the uniform boundedness principle, $\widetilde{D}_s(\varphi)$ admits a limit as $s\to 0^+$ and thus $\lim\limits_{s\to 0^+} s\widetilde{D}_s(\varphi)=0$. This shows \eqref{eq Distr1} in this case. 

Assume from now on that $m_i=n_i$ for all $1\leqslant i\leqslant r$ and $p_j$ is even for all $1\leqslant j\leqslant s$. Then, we have $\cA_0'=\cA'$ and moreover the map sending $\lambda\in \cA_0/\cA'$ to
$$\displaystyle \left((x_{i,\ell}(\lambda)+x_{i,\ell}^*(\lambda))_{\substack{ 1\leqslant i\leqslant r \\ 1\leqslant \ell\leqslant n_i}}, (y_{j,\ell}(\lambda)+y_{j,p_j+1-\ell}(\lambda))_{\substack{1\leqslant j\leqslant s\\ 1\leqslant \ell\leqslant \frac{p_j}{2}}}, (z_{k,\ell}(\lambda)+z_{k,q_k+1-\ell}(\lambda))_{\substack{1\leqslant k\leqslant t\\ 1\leqslant \ell\leqslant \floor*{\frac{q_k}{2}}}}, (z_{k,\ceil*{\frac{q_k}{2}}})_{\substack{1\leqslant k\leqslant t\\ q_k\equiv 1[2]}} \right)$$ 
induces an isomorphism of vector spaces
$$\displaystyle \cA_0/\cA'\simeq (i\R^{N})_0$$
which sends the measure on $\cA_0/\cA'$ to the measure on $(i\R^{N})_0$ appearing in Proposition \ref{prop 1 PV}. Thus, by \eqref{eq Distr2} viewing $\psi_s$ as a function on $(i\R^N)_0$ via this isomorphism, we have
$$\displaystyle \widetilde{D}_s(\varphi)=\int_{(i\R^N)_0} \frac{\psi_s(t_1,\ldots,t_N)}{\prod\limits_{\ell=1}^N (\frac{t_\ell}{h_\ell}+s)}d\underline{t}$$
where the sequence $(h_\ell)_{1\leqslant \ell\leqslant N}$ is the concatenation of the sequences $(d_i)_{\substack{ 1\leqslant i\leqslant r \\ 1\leqslant \ell\leqslant n_i}}$, $(e_j)_{\substack{1\leqslant j\leqslant s\\ 1\leqslant \ell\leqslant \frac{p_j}{2}}}$, $(f_k)_{\substack{1\leqslant k\leqslant t\\ 1\leqslant \ell\leqslant \floor*{\frac{q_k}{2}}}}$ and $(\frac{f_k}{2})_{\substack{1\leqslant k\leqslant t\\ q_k\equiv 1[2]}}$. Hence, by Proposition \ref{prop 1 PV} and the uniform boundedness principle, we get
\begin{align}\label{eq Distr3}
\displaystyle \lim\limits_{s\to 0^+} s\widetilde{D}_s(\varphi)=\frac{\prod\limits_{i=1}^N h_\ell}{\sum_{\ell=1}^Nh_\ell} (2\pi)^{N-1}\psi_0(0)=\frac{D}{n} (2\pi)^{N-1} 2^{1-c} \psi_0(0)
\end{align}
where the last equality follows from a painless computation. Finally, by the last part of Lemma \ref{lem cont family in Schwartz}, and since the function
$$\displaystyle \mu\in \cA'\mapsto \psi(\mu)\prod_{k=1}^t \prod_{1\leqslant \ell\leqslant \floor*{\frac{q_k}{2}}} \left( \frac{2z_{k,\ell}(\mu)}{f_k} \right)^{-1}=\varphi(\mu) \prod_{i=1}^r S^*_{m_i,n_i}(\frac{\ux_i(\mu)}{d_i})\times \prod_{j=1}^s T_{p_j}^*(\frac{\uy_j(\mu)}{e_j})\times \prod_{k=1}^t U_{q_k}^*(\frac{\uz_k(\mu)}{f_k}),$$
being the product of a Schwartz function by a polynomial (by Proposition \ref{prop 1 Polynomial identities}), belongs to $\cS(\cA')$, we have
\[\begin{aligned}
\displaystyle \psi_0(0)=\int_{\cA'} \varphi(\mu) \prod_{i=1}^r S^*_{m_i,n_i}(\frac{\ux_i(\mu)}{d_i})\times \prod_{j=1}^s T_{p_j}^*(\frac{\uy_j(\mu)}{e_j})\times \prod_{k=1}^t U_{q_k}^*(\frac{\uz_k(\mu)}{f_k})d\mu
\end{aligned}\]
Combining this with \eqref{eq Distr0} and \eqref{eq Distr3}, we obtain \eqref{eq Distr1} and this ends the proof of the proposition. $\blacksquare$

\subsection{A spectral limit}\label{Section spectral limit}

Recall that we have set $BC_n(\sigma)=BC(\sigma)\otimes \eta_n'$ for every $\sigma\in \Temp(U(n))$ (see Section \ref{Section LLC}) and that $\cS(\Temp(\overline{G_n(E)}))$ denotes the space of Schwartz functions $\Temp(\overline{G_n(E)})\to \C$ (see Section \ref{Section Schwartz functions on Temp(G)}).

\begin{prop}\label{prop1 Theo limite spectrale}
For every $\Phi\in \cS(\Temp(\overline{G_n(E)}))$, we have
\begin{align}\label{eq 1 Theo limite spectrale}
\displaystyle & \lim\limits_{s\to 0^+} n\gamma(s,\mathbf{1}_F,\psi')\int_{\Temp(\overline{G_n(E)})} \Phi(\pi) \gamma(s,\pi,\As,\psi')^{-1} d\mu_{\overline{G_n(E)}}(\pi)= \\
\nonumber & \lambda_{E/F}(\psi')^{-n^2} \int_{\Temp(U(n))/\stab} \Phi(BC_n(\sigma)) \frac{\gamma^*(0,\sigma,\Ad,\psi')}{\lvert S_\sigma\rvert} d\sigma
\end{align}
where the right-hand side is absolutely convergent and so does the left hand side for any $s\in \cH$.
\end{prop}

\vspace{2mm}

\noindent\ul{Proof}: The convergence of the right hand side and the left hand side of \eqref{eq 1 Theo limite spectrale} for $s\in \cH$ follows directly from Lemma \ref{lem 1 L parameters, LLC, basechange}, \eqref{eq 2 Planch} and \eqref{basic estimates spectral measure}. By Proposition \ref{prop Planch measure}, we have
\begin{align}\label{eq 2 Theo limite spectrale}
\displaystyle \int_{\Temp(\overline{G_n(E)})} \Phi(\pi) \gamma(s,\pi,\As,\psi')^{-1} d\mu_{\overline{G_n(E)}}(\pi)= \lambda_{E/F}(\psi')^{-n^2} \int_{\Temp(\overline{G_n(E)})} \Phi(\pi) \frac{\gamma^*(0,\pi,\overline{\Ad},\psi')}{\lvert S_\pi\rvert \gamma(s,\pi,\As,\psi')}d\pi
\end{align}
for every $s\in \cH$. Fix $\pi\in \Temp(\overline{G_n(E)})$. We can write it as
$$\displaystyle \pi=\left[\bigtimes\limits_{i=1}^r \tau_i^{\times m_i}\times (\tau_i^*)^{\times n_i} \right] \times \left[ \bigtimes\limits_{j=1}^s \mu_j^{\times p_j}\right] \times \left[\bigtimes\limits_{k=1}^t \nu_k^{\times q_k} \right]$$
(Recall that $\tau^*$ stands for $(\tau^c)^\vee$) where
\begin{itemize}
\item For all $1\leqslant i\leqslant r$, $\tau_i\in \Pi_2(G_{d_i}(E))$ for some $d_i\in \bN^*$ is such that $\tau_i\not\simeq \tau_i^*$ and $m_i,n_i\in \bN^*$ are such that $m_i\geqslant n_i$. Moreover, $\tau_i\not{\simeq} \tau_j$ and $\tau_i\not{\simeq} \tau_j^*$ for all $1\leqslant i<j\leqslant r$.

\item For all $1\leqslant j\leqslant s$, $\mu_j\in \Pi_2(G_{e_j}(E))$ for some $e_j\in \bN^*$ is such that $\mu_j\simeq \mu_j^*$ but $\gamma(0,\mu_j,\As,\psi')\neq 0$ and $p_j\in \bN^*$. Moreover, $\mu_i\not{\simeq}\mu_j$ for all $1\leqslant i<j\leqslant s$.

\item For all $1\leqslant k\leqslant t$, $\nu_k\in \Pi_2(G_{f_k}(E))$ for some $f_k\in \bN^*$ is such that $\nu_k\simeq \nu_k^*$ and $\gamma(0,\nu_k,\As,\psi')=0$ and $q_k\in \bN^*$. Moreover, $\nu_i\not{\simeq} \nu_j$ for all $1\leqslant i<j\leqslant t$.
\end{itemize}
In other words, we have written $\pi=i_{M(F)}^{\overline{G_n(E)}}(\Pi)$ where
$$\displaystyle M(F)=\left(\prod_{i=1}^r G_{d_i}(E)^{m_i+n_i}\times \prod_{j=1}^s G_{e_j}(E)^{p_j} \times \prod_{k=1}^t G_{f_k}(E)^{q_k} \right)/Z_n(F)$$
is a Levi subgroup of $\overline{G_n(E)}$ and
\begin{equation}\label{eq 2bis Theo limite spectrale}
\displaystyle \Pi=\bigboxtimes\limits_{i=1}^r \tau_i^{\boxtimes m_i}\boxtimes (\tau_i^*)^{\boxtimes n_i} \boxtimes \bigboxtimes\limits_{j=1}^s \mu_j^{\boxtimes p_j} \boxtimes \bigboxtimes\limits_{k=1}^t \nu_k^{\boxtimes q_k}
\end{equation}
is a certain square-integrable representation of $M(F)$.

In the $p$-adic case, using a partition of unity, we may assume that $\Phi$ is supported in a small enough neighborhood $\cU \subset \Temp(\overline{G_n(E)})$ of $\pi$.

In the Archimedean case (i.e. when $F=\R$), we will assume for the moment that $\Phi$ is supported in the connected component $\cO\subset \Temp(\overline{G_n(E)})$ of $\pi$ and we will explain at the end of the proof how to extend the result to all Schwartz functions on $\Temp(\overline{G_n(E)})$. Up to twisting $\Pi$, we will also assume that the restriction of the central character of $\Pi$ to $A_M(F)$ has finite order.

Let $\cA_0$ be as in Section \ref{Section distributions}. Then, there exists a unique isomorphism of vector spaces
\begin{align}\label{eq 3 Theo limite spectrale}
\displaystyle \cA_0\simeq i\cA_M^*
\end{align}
which when composed with the map $\lambda\in i\cA_M^*\mapsto \pi_\lambda:=i_{M(F)}^{\overline{G_n(E)}}(\Pi_\lambda)$ becomes (with the notation of Section \ref{Section distributions})
\begin{align}\label{eq 4 Theo limite spectrale}
\displaystyle \lambda\in \cA_0\mapsto \pi_\lambda:= & \left(\bigtimes\limits_{i=1}^r \bigtimes\limits_{\ell=1}^{m_i} \tau_i\otimes \lvert \det\rvert_E^{\frac{x_{i,\ell}(\lambda)}{d_i}} \times \bigtimes\limits_{\ell=1}^{n_i} \tau^*_i\otimes \lvert \det\rvert_E^{\frac{x^*_{i,\ell}(\lambda)}{d_i}}\right) \times \left(\bigtimes\limits_{j=1}^s \bigtimes\limits_{\ell=1}^{p_j} \mu_j\otimes \lvert \det \rvert_E^{\frac{y_{j,\ell}(\lambda)}{e_j}} \right) \\
\nonumber & \times \left( \bigtimes\limits_{k=1}^t \bigtimes\limits_{\ell=1}^{q_k} \nu_k\otimes \lvert \det \rvert_E^{\frac{z_{k,\ell}(\lambda)}{f_k}}\right)
\end{align}

Define $W$ as in Section \ref{Section distributions}. Then, there exists an isomorphism $W\simeq W(\overline{G_n(E)},\Pi)$ such that \eqref{eq 3 Theo limite spectrale} transports the action of $W$ on $\cA_0$ to the action of $W(\overline{G_n(E)}, \Pi)$ on $i\cA_M^*$. In the $p$-adic case, up to shrinking $\cU$ we may assume, which we do in what follows, that there exists a small open neighborhood $\cV\subset \cA_0$ of $0$ such that $\lambda\mapsto \pi_\lambda$ induces a topological isomorphism $\cV/W\simeq \cU$. The inverse image by \eqref{eq 3 Theo limite spectrale} of $iX^*(A_M)$ is $\frac{1}{2}\Lambda_0$ where $\Lambda_0$ denotes the intersection of $\cA_0$ with the lattice
$$\displaystyle \Lambda:=\prod\limits_{1\leqslant i\leqslant r} (i\bZ)^{m_i}\times (i\bZ)^{n_i}\times \prod\limits_{1\leqslant j\leqslant s} (i\bZ)^{p_j}\times \prod\limits_{1\leqslant k\leqslant t} (i\bZ)^{q_k}$$
of $\cA$ (the factor $\frac{1}{2}$ stemming from the fact that $\lvert x\rvert_E=\lvert x\rvert_F^2$ for every $x\in F^\times$). Recall that we have fixed a Haar measure on $\cA_0$ in Section \ref{Section distributions}. By definition of this Haar measure, we have $\vol(\cA_0/\Lambda_0)=1$. Hence, by the definition of the Haar measure on $i\cA_M^*$ (see Section \ref{Section measures}) the isomorphism \eqref{eq 3 Theo limite spectrale} sends this Haar measure to $\left(\frac{\pi}{\log(q_F)}\right)^{1-S}$ times the Haar measure on $\cA_0$ where $S=\dim(\cA_0)+1=\sum_{i=1}^r m_i+n_i+\sum_{j=1}^s p_j+\sum_{k=1}^t q_k$.

Therefore, by \eqref{eq 2 Theo limite spectrale} and \eqref{eq 4 Measures} we have
\begin{align}\label{eq 5 Theo limite spectrale}
\displaystyle & \int_{\Temp(\overline{G_n(E)})} \Phi(\pi) \gamma(s,\pi,\As,\psi')^{-1} d\mu_{\overline{G_n(E)}}(\pi)= \\
\nonumber & \frac{\lambda_{E/F}(\psi')^{-n^2}}{\lvert W\rvert}\left(\frac{\pi}{\log(q_F)} \right)^{1-S} \int_{\cA_0} \varphi(\lambda) \frac{\gamma^*(0,\pi_\lambda,\overline{\Ad},\psi')}{\lvert S_{\pi_\lambda}\rvert \gamma(s,\pi_\lambda,\As,\psi')}d\lambda
\end{align}
where we have set $\varphi(\lambda)=\Phi(\pi_\lambda)$ in the Archimedean case and
$$\displaystyle \varphi(\lambda)=\left\{
    \begin{array}{ll}
        \Phi(\pi_\lambda) & \mbox{ if } \lambda\in \cV \\
        0 & \mbox{ otherwise}
    \end{array}
\right.
$$
in the $p$-adic case. Notice that in both cases we have $\varphi\in \cS(\cA_0)^W$.

By \eqref{eq 1 gp of centralizers}, \eqref{eq 2 gp of centralizers} and \eqref{eq 3 gp of centralizers}, we readily check that
\begin{align}\label{eq 6 Theo limite spectrale}
\displaystyle \lvert S_{\pi_\lambda}\rvert=2^S P
\end{align}
for every $\lambda\in \cA_0$ where we have set $\displaystyle P=\prod_{i=1}^r d_i^{m_i+n_i}\prod_{j=1}^s e_j^{p_j} \prod_{k=1}^t f_k^{q_k}$.

From \eqref{eq 1 L parameters, LLC, basechange}, \eqref{eq 2 L parameters, LLC, basechange}, \eqref{eq 4 L parameters, LLC, basechange}, \eqref{eq 5 L parameters, LLC, basechange}, \eqref{eq 6bis L parameters, LLC, basechange}, \eqref{eq 7 L parameters, LLC, basechange} and \eqref{eq 8 L parameters, LLC, basechange}, we infer that there exists a function $F\in C^\infty(\cA_0)^W$ which is of moderate growth together with all its derivatives in the Archimedean case such that
\begin{align}\label{eq 7 Theo limite spectrale}
\displaystyle & \gamma^*(0,\pi_\lambda,\overline{\Ad},\psi')=\left(\prod_{i=1}^r \prod_{1\leqslant \ell \neq \ell'\leqslant m_i} (\frac{x_{i,\ell}(\lambda)-x_{i,\ell'}(\lambda)}{d_i}) \prod_{1\leqslant \ell \neq \ell'\leqslant n_i} (\frac{x^*_{i,\ell}(\lambda)-x^*_{i,\ell'}(\lambda)}{d_i})\right) \times \\
\nonumber & \left(\prod_{j=1}^s \prod_{1\leqslant \ell \neq \ell'\leqslant p_j} (\frac{y_{j,\ell}(\lambda)-y_{j,\ell'}(\lambda)}{e_j}) \right) \times \left(\prod_{k=1}^t \prod_{1\leqslant \ell \neq \ell'\leqslant q_k} (\frac{z_{k,\ell}(\lambda)-z_{k,\ell'}(\lambda)}{f_k}) \right) F(\lambda)
\end{align}
for almost all $\lambda\in \cA_0$.

Set $\displaystyle \cA_\C:=\cA\otimes_{\R} \C=\prod_{i=1}^r \C^{m_i+n_i} \prod_{j=1}^s \C^{p_j}\prod_{k=1}^t \C^{q_k}\simeq \C^S$ and embed $\C$ in $\cA_\C$ by $s\mapsto s\lambda_0$ where $\lambda_0\in \cA_\C$ has coordinates
$$\displaystyle x_{i,\ell}(\lambda_0)=x^*_{i,\ell'}(\lambda_0)=d_i \;(1\leqslant i\leqslant r, 1\leqslant \ell\leqslant m_i, 1\leqslant \ell'\leqslant n_i),\; y_{j,\ell}(\lambda_0)=e_j \;(1\leqslant j\leqslant s, 1\leqslant \ell\leqslant p_j),$$
$$\displaystyle z_{k,\ell}(\lambda_0)=f_k \;(1\leqslant k\leqslant t, 1\leqslant \ell\leqslant q_k).$$
From \eqref{eq 1 L parameters, LLC, basechange}, \eqref{eq 2 L parameters, LLC, basechange}, \eqref{eq 6bis L parameters, LLC, basechange}, \eqref{eq 7 L parameters, LLC, basechange}, \eqref{eq 8 L parameters, LLC, basechange}, \eqref{eq 9 L parameters, LLC, basechange}, \eqref{eq 10 L parameters, LLC, basechange} and \eqref{eq 11 L parameters, LLC, basechange}, we infer similarly that, up to shrinking $\cU$ in the $p$-adic case, there exists a $W$-invariant meromorphic function $G$ on $\cA_\C$ whose polar divisors are affine subspaces not meeting $(-\epsilon+\cH)^S$ for some $\epsilon>0$ and which is of moderate growth on vertical strips with all its derivatives there in the Archimedean case, resp. whose polar divisors are affine subspaces disjoint from $(\cH\cup\{ 0\})^S+\cV$ in the $p$-adic case, such that
\begin{align}\label{eq 8 Theo limite spectrale}
\displaystyle & \gamma(s,\pi_\lambda,\As,\psi')^{-1}=\left(\prod_{i=1}^r \prod_{\substack{1\leqslant \ell \leqslant m_i \\ 1\leqslant \ell'\leqslant n_i}} (s+\frac{x_{i,\ell}(\lambda)+x^*_{i,\ell'}(\lambda)}{d_i})^{-1}\right) \times \\
\nonumber & \left(\prod_{j=1}^s \prod_{1\leqslant \ell < \ell'\leqslant p_j} (s+\frac{y_{j,\ell}(\lambda)+y_{j,\ell'}(\lambda)}{e_j})^{-1} \right) \times \left(\prod_{k=1}^t \prod_{1\leqslant \ell \leqslant \ell'\leqslant q_k} (s+\frac{z_{k,\ell}(\lambda)+z_{k,\ell'}(\lambda)}{f_k})^{-1} \right) G(\frac{s}{2}+\lambda)
\end{align}
for all $\lambda\in \cA_0$ and $s\in \cH$.

From \eqref{eq 7 Theo limite spectrale} and \eqref{eq 8 Theo limite spectrale} it follows that, with the notation of Section \ref{Section distributions}, there exists a continuous family $s\in \cH\cup\{0\}\mapsto \varphi_s\in \cS(\cA_0)^W$ such that
\begin{align}\label{eq 8bis Theo limite spectrale}
\displaystyle \varphi(\lambda)\frac{\gamma^*(0,\pi_\lambda,\overline{\Ad},\psi')}{\gamma(s,\pi_\lambda,\As,\psi')}=\varphi_s(\lambda) \prod_{i=1}^r P_{m_i,n_i,s}(\frac{\ux_i(\lambda)}{d_i})\times \prod_{j=1}^s Q_{p_j,s}(\frac{\uy_j(\lambda)}{e_j})\times \prod_{k=1}^t R_{q_k,s}(\frac{\uz_k(\lambda)}{f_k})
\end{align}
for all $s\in \cH$ and almost all $\lambda\in \cA_0$. Therefore, using again notation from Section \ref{Section distributions}, by \eqref{eq 5 Theo limite spectrale}, \eqref{eq 6 Theo limite spectrale} and since $\gamma(s,\mathbf{1}_F,\psi')\sim_{s\to 0} s\log(q_F) \gamma^*(0,\mathbf{1}_F,\psi')$, we have
\begin{align}\label{eq 9 Theo limite spectrale}
\displaystyle & n\gamma(s,\mathbf{1}_F,\psi')\int_{\Temp(\overline{G_n(E)})} \Phi(\pi) \gamma(s,\pi,\As,\psi')^{-1} d\mu_{\overline{G_n(E)}}(\pi)\sim_{s\to 0^+} \\
\nonumber & n\frac{\lambda_{E/F}(\psi')^{-n^2}}{\lvert W\rvert}\left(\frac{2\pi}{\log(q_F)} \right)^{-S}\pi P^{-1}\gamma^*(0,\mathbf{1}_F,\psi') sD_s(\varphi_s).
\end{align}

Assume first that $\pi$ cannot be written as $BC_n(\sigma)$ for some $\sigma\in \Temp(U(n))$. This means that (by \eqref{eq 0 L parameters, LLC, basechange} and \eqref{eq 11 L parameters, LLC, basechange}) either there exists $1\leqslant i\leqslant r$ such that $m_i\neq n_i$ or there exists $1\leqslant j\leqslant s$ such that $p_j$ is odd. In the Archimedean case, by the assumption that the central character of the representation $\Pi$ \eqref{eq 2bis Theo limite spectrale} restricted to $A_M(F)$ has finite order, it follows that $BC_n(\Temp(U(n)))$ does not meet $\cO$ at all. In the $p$-adic case on the other hand, up to shrinking $\cU$ if necessary, we may assume that  $BC_n(\Temp(U(n)))\cap\cU=\emptyset$. Then, in both cases, the right-hand side of \eqref{eq 1 Theo limite spectrale} turns out to be just zero whereas by \eqref{eq 9 Theo limite spectrale} and Proposition \ref{prop 1 Distributions} so does the left-hand side. This proves the proposition in this case.

Assume now that there exists $\sigma\in \Temp(U(n))/\stab$ such that $\pi=BC_n(\sigma)$. Then $m_i=n_i$ for all $1\leqslant i\leqslant r$, $p_j$ is even for all $1\leqslant j\leqslant s$ and we can write
$$\displaystyle \sigma=\left[\bigtimes_{i=1}^r \tau_i^{\times m_i}\bigtimes_{j=1}^s \mu_j^{\times \frac{p_j}{2}} \bigtimes_{k=1}^t \nu_k^{\times \floor*{\frac{q_k}{2}}}\right]\rtimes \sigma_0$$
where $\sigma_0$ is a discrete series of some $U(m)$ with
$$\displaystyle BC_m(\sigma_0)=\bigtimes_{\substack{1\leqslant k\leqslant t \\ q_k\equiv 1 \mod{2}}} \nu_k.$$
In other words, we have $\sigma=i_{L}^{U(n)}(\Sigma)$ where
$$\displaystyle L=\prod_{i=1}^r R_{E/F}G_{d_i,E}^{m_i}\times \prod_{j=1}^s R_{E/F}G_{e_j,E}^{\frac{p_j}{2}}\times \prod_{k=1}^t R_{E/F}G_{f_k,E}^{\floor*{\frac{q_k}{2}}}\times U(m),$$
for a suitable integer $m$, is a Levi subgroup of $U(n)$ and $\Sigma$ is a certain square-integrable representation of $L(F)$. Let $\cA'$ be as in Section \ref{Section distributions}. Then, there exists a unique isomorphism
\begin{align}\label{eq 10 Theo limite spectrale}
\displaystyle \cA'\simeq i\cA_L^*
\end{align}
which when composed with the map $\mu\in i\cA_L^*\mapsto \sigma_\mu:=i_{L(F)}^{U(n)}(\Sigma_\mu)$ becomes
\[\begin{aligned}
\displaystyle \mu\in \cA'\mapsto \sigma_\mu=\left[\bigtimes_{i=1}^r \bigtimes_{\ell=1}^{m_i}\tau_i\otimes \lvert \det\rvert_E^{\frac{x_{i,\ell}(\mu)}{d_i}} \times \bigtimes_{j=1}^s \bigtimes_{\ell=1}^{\frac{p_j}{2}} \mu_j\otimes \lvert \det\rvert_E^{\frac{y_{j,\ell}(\mu)}{e_j}}\times  \bigtimes_{k=1}^t \bigtimes_{\ell=1}^{\floor*{\frac{q_k}{2}}}\nu_k\otimes \lvert \det\rvert_E^{\frac{z_{k,\ell}(\mu)}{f_k}}\right]\rtimes \sigma_0
\end{aligned}\]
By \eqref{eq -1 L parameters, LLC, basechange} and \eqref{eq 0 L parameters, LLC, basechange}, in the Archimedean case, for every $\lambda\in \cA_0$ we have $\pi_\lambda\in BC_n(\Temp(U(n)))$ if and only if $\lambda\in \cA'$ in which case $\pi_\lambda=BC_n(\sigma_\lambda)$. Similarly, in the $p$-adic case and up to shrinking $\cU$ if necessary, for every $\lambda\in \cV$ we have $\pi_\lambda\in BC_n(\Temp(U(n)))$ if and only if $\lambda\in \cA'$ in which case we again have $\pi_\lambda=BC_n(\sigma_\lambda)$. It follows that the function $\Phi\circ BC_n$ on $\Temp(U(n))/\stab$ is supported in the connected component $\{\sigma_\mu\mid \mu\in \cA' \}$ of $\sigma$ in the Archimedean case and even in the open neighborhood $\{\sigma_\mu\mid \mu\in \cA'\cap \cV \}$ of $\sigma$ in the $p$-adic case. We check easily, as before, that the isomorphism \eqref{eq 10 Theo limite spectrale} sends the Haar measure on $i\cA_L^*$ to $\left(\frac{\pi}{\log(q_F)}\right)^{N-S}$ times the Haar measure on $\cA'$ that was defined in Section \ref{Section distributions} where
$$\displaystyle N=S-\dim(\cA')=\sum_{i=1}^r n_i+\sum_{j=1}^s \frac{p_j}{2}+\sum_{k=1}^t \ceil*{\frac{q_k}{2}}$$
(using again notation from Section \ref{Section distributions}). Moreover, using the precise description of $W(U(n),\Sigma)$ given in Section \ref{Section representations}, we see that there exists an isomorphism $W'\simeq W(U(n),\Sigma)$ (where $W'$ is defined as in Section \ref{Section distributions}) such that \eqref{eq 10 Theo limite spectrale} transports the action of $W'$ on $\cA'$ to the action of $W(U(n),\Sigma)$ on $i\cA_L^*$. Therefore by \eqref{eq 4bis Measures}, assuming $\cU$ sufficiently small, we have
\begin{align}\label{eq 11 Theo limite spectrale}
\displaystyle & \int_{\Temp(U(n))/\stab} \Phi(BC_n(\sigma)) \frac{\gamma^*(0,\sigma,\Ad,\psi')}{\lvert S_\sigma\rvert} d\sigma= \\
\nonumber & \lvert W'\rvert^{-1}\left( \frac{\pi}{\log(q_F)}\right)^{N-S}\int_{\cA'} \varphi(\mu) \frac{\gamma^*(0,\sigma_\mu,\Ad,\psi')}{\lvert S_{\sigma_\mu}\rvert}d\mu
\end{align}
Moreover, using \eqref{eq 1 gp of centralizers}, \eqref{eq 2 gp of centralizers} and \eqref{eq 3 gp of centralizers} we readily check that
$$\displaystyle \lvert S_{\sigma_\mu}\rvert=2^{c+S-N}\frac{P}{D}$$
for all $\mu\in \cA'$ where $c=\lvert \{1\leqslant k\leqslant t\mid q_k\equiv 1 \mod{2} \}\rvert$ and $\displaystyle D=\prod_{i=1}^r d_i^{n_i} \times \prod_{j=1}^s e_j^{\frac{p_j}{2}}\times \prod_{k=1}^t f_k^{\ceil*{\frac{q_k}{2}}}$ (same notation as in Section \ref{Section distributions}) so that $\displaystyle c+S-N=\sum_{i=1}^r m_i+\sum_{j=1}^s \frac{p_j}{2}+\sum_{k=1}^t \ceil*{\frac{q_k}{2}}$ and $\displaystyle \frac{P}{D}=\prod_{i=1}^r d_i^{m_i} \times \prod_{j=1}^s e_j^{\frac{p_j}{2}}\times \prod_{k=1}^t f_k^{\floor*{\frac{q_k}{2}}}$. Consequently, \eqref{eq 11 Theo limite spectrale} can be rewritten
\begin{align}\label{eq 12 Theo limite spectrale}
\displaystyle & \int_{\Temp(U(n))/\stab} \Phi(BC_n(\sigma)) \frac{\gamma^*(0,\sigma,\Ad,\psi')}{\lvert S_\sigma\rvert} d\sigma= \\
\nonumber & \lvert W'\rvert^{-1}\left( \frac{2\pi}{\log(q_F)}\right)^{N-S} 2^{-c}\frac{D}{P}\int_{\cA'} \varphi(\mu) \gamma^*(0,\sigma_\mu,\Ad,\psi')d\mu
\end{align}
On the other hand, by \eqref{eq 9 Theo limite spectrale}, Proposition \ref{prop 1 Distributions} and the uniform boundedness principle\footnote{\label{footnote}More precisely, the uniform boundedness principle implies that, for any sequence $s_n\in \cH$ converging to $0$, the sequence of tempered distributions $(s_nD_{s_n})_n$ is equicontinuous hence that, by continuity of $s\mapsto \varphi_s\in \cS(\cA_0)$, $s_nD_{s_n}(\varphi_{s_n})-s_nD_{s_n}(\varphi_0)$ converges to $0$ as $n\to \infty$ and we can then apply Proposition \ref{prop 1 Distributions} to the single Schwartz function $\varphi_0$.} we have
\[\begin{aligned}
\displaystyle & \lim\limits_{s\to 0^+} n\gamma(s,\mathbf{1}_F,\psi')\int_{\Temp(\overline{G_n(E)})} \Phi(\pi) \gamma(s,\pi,\As,\psi')^{-1} d\mu_{\overline{G_n(E)}}(\pi)= \\
& \frac{\lambda_{E/F}(\psi')^{-n^2}}{\lvert W'\rvert}\left(\frac{2\pi}{\log(q_F)} \right)^{-S}(2\pi)^N 2^{-c} \frac{D}{P}\gamma^*(0,\mathbf{1}_F,\psi')\times \\
& \int_{\cA'} \lim\limits_{s\to 0^+} s^N \varphi_s(\mu) \prod_{i=1}^r P_{m_i,n_i,s}(\frac{\ux_i(\mu)}{d_i})\times \prod_{j=1}^s Q_{p_j,s}(\frac{\uy_j(\mu)}{e_j})\times \prod_{k=1}^t R_{q_k,s}(\frac{\uz_k(\mu)}{f_k})d\mu
\end{aligned}\]
From \eqref{eq 8bis Theo limite spectrale} and the fact that $\zeta_F(s)\sim_{s\to 0} (s\log(q_F))^{-1}$ this can be rewritten as
\[\begin{aligned}
\displaystyle & \lim\limits_{s\to 0^+} n\gamma(s,\mathbf{1}_F,\psi')\int_{\Temp(\overline{G_n(E)})} \Phi(\pi) \gamma(s,\pi,\As,\psi')^{-1} d\mu_{\overline{G_n(E)}}(\pi)= \\
& \frac{\lambda_{E/F}(\psi')^{-n^2}}{\lvert W'\rvert}\left(\frac{2\pi}{\log(q_F)} \right)^{N-S} 2^{-c} \frac{D}{P}\gamma^*(0,\mathbf{1}_F,\psi')\int_{\cA'} \varphi(\mu) \lim\limits_{s\to 0^+} \zeta_F(s)^{-N}\frac{\gamma^*(0,\pi_\mu,\overline{\Ad},\psi')}{\gamma(s,\pi_\mu,\As,\psi')} d\mu
\end{aligned}\]
From \eqref{eq 8 Theo limite spectrale} it is easy to infer that $s\mapsto \gamma(s,\pi_\mu,\As,\psi')$ has a zero of order $N$ at $s=0$ for almost all $\mu\in \cA'$. Therefore by \eqref{eq 3 L parameters, LLC, basechange} and \eqref{eq 12 L parameters, LLC, basechange} the above equality is equivalent to
\[\begin{aligned}
\displaystyle & \lim\limits_{s\to 0^+} n\gamma(s,\mathbf{1}_F,\psi')\int_{\Temp(\overline{G_n(E)})} \Phi(\pi) \gamma(s,\pi,\As,\psi')^{-1} d\mu_{\overline{G_n(E)}}(\pi)= \\
& \frac{\lambda_{E/F}(\psi')^{-n^2}}{\lvert W'\rvert}\left(\frac{2\pi}{\log(q_F)} \right)^{N-S} 2^{-c} \frac{D}{P}\int_{\cA'} \varphi(\mu) \gamma^*(0,\sigma_\mu,\Ad,\psi') d\mu
\end{aligned}\]
Comparing this with \eqref{eq 12 Theo limite spectrale} we get the identity of the proposition.

It remains to explain how to remove the assumption that $\Phi$ is supported on one component of $\Temp(\overline{G_n(E)})$ in the Archimedean case. The result obviously extend to the case where $\Phi$ is supported on a finite number of components hence to the subspace $\cS_c(\Temp(\overline{G_n(E)}))$ of compactly supported Schwartz functions. Let $\mathbb{D}_s(\Phi)$ be the left hand side of \eqref{eq 1 Theo limite spectrale} for a given $s\in \cH$. As we already argue, $\mathbb{D}_s$ is a continuous linear form on $\cS(\Temp(\overline{G_n(E)}))$. Since the right hand side of \eqref{eq 1 Theo limite spectrale} is also continuous and $\cS_c(\Temp(\overline{G_n(E)}))$ is dense in $\cS(\Temp(\overline{G_n(E)}))$, by the uniform boundedness principle it suffices to check that for any sequence $(s_n)_n\in \cH$ converging to $0$, the sequence $(\mathbb{D}_{s_n})_n$ is pointwise bounded (as it will then be an equicontinuous family of linear forms). Let $(s_n)_n\in \cH$ be such a sequence and $\Phi\in \cS(\Temp(\overline{G_n(E)})$. For every connected component $\cO\subset \Temp(\overline{G_n(E)})$, let us write $\Phi_{\cO}$ for the restriction of $\Phi$ to $\cO$. Then, the series $\sum_{\cO} \Phi_{\cO}$ is readily seen to converge absolutely towards $\Phi$ in the Schwartz space $\cS(\Temp(\overline{G_n(E)})$. Therefore, we just need to show the existence of a continuous seminorm $\lVert .\rVert$ on $\cS(\Temp(\overline{G_n(E)}))$ such that $\lvert \mathbb{D}_{s_n}(\Phi_{\cO})\rvert\leqslant \lVert \Phi_{\cO}\rVert$ for every $n$ and $\cO$. We can as well forgot the index $\cO$ and just assume that $\Phi$ is supported in one connected component $\cO\subset \Temp(\overline{G_n(E)})$. Using the above notation, we have then shown that, up to some immaterial constant, $\mathbb{D}_{s_n}(\Phi)$ is equal to $s_n D_{s_n}(\varphi_{s_n})$ where $D_{s_n}$ is the Schwartz distribution on $\cA_0\simeq i\cA_M^*$ introduced in Section \ref{Section distributions} and the function $\varphi_{s_n}$ is defined as above. Of course, the distribution $D_{s_n}$ depends on auxilliary parameters associated to the component $\cO$ (more precisely: the Levi $M$ and integers $r$, $s$, $t$, $(m_i)$, $(n_i)$, $(p_j)$ and $(q_k)$) but these parameters can only take finitely many values and for the purpose of uniformly bounding $\mathbb{D}_{s_n}(\Phi)$ we can as well assume these parameters as fixed. By Proposition \ref{prop 1 Distributions} and the uniform boundedness principle, the family $(s_nD_{s_n})_n$ is itself equicontinuous. Thus, it only remains to show that for every continuous seminorm $\lVert .\rVert_{\cS}$ on $\cS(\cA_0)=\cS(i\cA_M^*)$ there exists a continuous seminorm $\lVert .\rVert$ on $\cS(\Temp(\overline{G_n(E)}))$ such that $\lVert \varphi_{s_n}\rVert_{\cS}\leqslant \lVert \Phi\rVert$ for every $n$ and $\Phi$ (supported in one component $\cO\subset \Temp(\overline{G_n(E)})$). Returning to the definition of $\varphi_{s_n}$, we see that it is a product of the function $\lambda\mapsto \Phi(\pi_\lambda)$ by the function
$$\displaystyle \lambda\mapsto \left(\prod_i (t-t_i(\pi_\lambda))^{-1} \gamma(t,\pi_\lambda,\overline{\Ad},\psi')\right)\mid_{t=0} \prod_j (s_n-s_j(\pi_\lambda))\gamma(s_n,\pi_\lambda,\As,\psi')^{-1}$$
where $\{ t_i(\pi_\lambda)\}$ (resp. $\{s_j(\pi_\lambda) \}$) stand for the zeroes (counted with multiplicities) of $\gamma(t,\pi_\lambda,\overline{\Ad},\psi')$ (resp. of $\gamma(s,\pi_\lambda,\As,\psi')$) on the imaginary line $i\R$. Therefore, by Lemma \ref{lem 1 L parameters, LLC, basechange} and the definition of $\cS(\Temp(\overline{G_n(E)}))$, we see that for every $\partial\in \Sym^\bullet(\cA^*_{M,\C})$ and $k>0$ there exists a continuous seminorm $\lVert .\rVert_{\partial, k}$ on $\cS(\Temp(\overline{G_n(E)}))$ such that
$$\displaystyle \sup_{\lambda\in i\cA_M^*}N(\pi_\lambda)^k\left\lvert \partial \varphi_{s_n}(\lambda)\right\rvert\leqslant \lVert \Phi\rVert_{\partial,k}$$
for every $n$, $\lambda\in i\cA_M^*$ and $\Phi$. Finally, as $N(\pi_\lambda)\gg 1+\lVert \lambda\rVert$, any continuous seminorm on $\cS(i\cA_M^*)$ is bounded up to a constant by the left hand side of the above inequality for some $\partial$ and $k$. This shows the desired pointwise boundedness and therefore ends the proof of the proposition. $\blacksquare$

\subsection{A corollary to Proposition \ref{prop1 Theo limite spectrale}}\label{Section Corollary}

For every $\sigma\in \Temp(U(n))/\stab$, set
$$\displaystyle c(\sigma):=\lambda_{E/F}(\psi')^{-\frac{n(n+1)}{2}}c_1(\pi)\omega_\sigma(-1)^{1-n}\eta_{E/F}(-1)^{\frac{n(n-1)^2}{2}}$$
where $\pi=BC_n(\sigma)$ and $c_1(\pi)$ is the constant defined by \eqref{prop form betan}. Notice that $c(\sigma)$ is just a certain root of unity.

Recall from Section \ref{Section betan} the continuous linear form $\beta:\cC^w(N_n(E)\backslash G_n(E),\psi_n)\to \C$ and also that in Section \ref{section Plancherel Whitt} we have associated to any function $f\in \cS(G_n(E))$ and $\pi\in \Temp(G_n(E))$ a function $W_{f,\pi}\in \cC^w(N_n(E)\times N_n(E)\backslash G_n(E)\times G_n(E),\psi_n^{-1}\boxtimes \psi_n)$. In particular, we have  $W_{f,\pi}(g,.)\in \cC^w(N_n(E)\backslash G_n(E),\psi_n)$ for all $g\in G_n(E)$.

\begin{cor}\label{cor lim spectrale}
For every $f\in \cS(G_n(E))$ and $g\in G_n(E)$, we have
\begin{align}\label{eq 1 Cor to lim spectrale}
\displaystyle & \int_{N_n(F)\backslash G_n(F)} W_f(g,h) dh= \\
\nonumber & \lvert \tau\rvert_E^{\frac{n(n-1)}{4}} \int_{\Temp(U(n))/\stab} \beta(W_{f,BC_n(\sigma)}(g,.)) \frac{\gamma^*(0,\sigma,\Ad,\psi')}{\lvert S_\sigma\rvert} c(\sigma) d\sigma
\end{align}
where the right-hand side is an absolutely convergent integral.
\end{cor}

\noindent\ul{Proof}: First note that up to replacing $f$ by $L(g)f$ we may assume that $g=1$. Let $f\in \cS(G_n(E))$ and define $\widetilde{f}\in \cS(\overline{G_n(E)})$ by $\displaystyle \widetilde{f}(g)=\int_{Z_n(F)}f(zg)dz$. Then, we have
$$\displaystyle \int_{N_n(F)\backslash G_n(F)} W_f(1,h) dh=\int_{Z_n(F)N_n(F)\backslash G_n(F)} W_{\widetilde{f}}(1,h) dh$$
Choose $\phi\in \cS(F^n)$ such that $\phi(0)=1$. Since $\widetilde{f}_\pi=f_\pi$ for every $\pi\in \Temp(\overline{G_n(E)})$, by Lemma \ref{lem Zeta function}(i)-(ii) and Proposition \ref{prop 1 Plancherel Whitt}, we have
\[\begin{aligned}
\displaystyle \int_{Z_n(F)N_n(F)\backslash G_n(F)} W_{\widetilde{f}}(1,h) dh & =\lim\limits_{s\to 0^+} n\gamma(s,\mathbf{1}_F,\psi') Z(s,W_{\widetilde{f}}(1,.),\phi) \\
 & =\lim\limits_{s\to 0^+} n\gamma(s,\mathbf{1}_F,\psi') \int_{\Temp(\overline{G_n(E)})}Z(s,W_{f,\pi}(1,.),\phi) d\mu_{\overline{G_n(E)}}(\pi).
\end{aligned}\]
By \eqref{eq 1 Planch Whitt}, applying the functional equation of Theorem \ref{theo1 Zeta function}, this becomes
\begin{align}\label{eq 2 Cor to lim spectrale}
\displaystyle & \int_{N_n(F)\backslash G_n(F)} W_f(1,h) dh = \lim\limits_{s\to 0^+} n\gamma(s,\mathbf{1}_F,\psi') \lvert\tau\rvert_E^{\frac{n(n-1)}{2}(\frac{1}{2}-s)} \lambda_{E/F}(\psi')^{\frac{n(n-1)}{2}} \times \\
\nonumber &  \int_{\Temp(\overline{G_n(E)})}Z(1-s,\widetilde{W_{f,\pi}}(1,.),\widehat{\phi})\omega_\pi(\tau)^{1-n}\gamma(s,\pi,\As,\psi')^{-1}d\mu_{\overline{G_n(E)}}(\pi).
\end{align}
Set $\Phi_s(\pi)=Z(1-s,\widetilde{W_{f,\pi}}(1,.),\widehat{\phi})\omega_\pi(\tau)^{1-n}$ for every $\pi\in \Temp(\overline{G_n(E)})$. Then, by Proposition \ref{prop 1 Plancherel Whitt}, Lemma \ref{lem Zeta function}(i) and Lemma \ref{lem holom Schwartz functions}, the map $s\mapsto \Phi_s\in \cS(\Temp(\overline{G_n(E)}))$ is holomorphic for $\Re(s)<1$. Therefore, by Proposition \ref{prop1 Theo limite spectrale} and the uniform boundedness principle\footnote{See footnote \ref{footnote}.}, we have
\begin{align}\label{eq 3 Cor to lim spectrale}
\displaystyle & \lim\limits_{s\to 0^+} n\gamma(s,\mathbf{1}_F,\psi') \lvert\tau\rvert_E^{\frac{n(n-1)}{2}(\frac{1}{2}-s)} \lambda_{E/F}(\psi')^{\frac{n(n-1)}{2}} \times \\
\nonumber & \int_{\Temp(\overline{G_n(E)})}Z(1-s,\widetilde{W_{f,\pi}}(1,.),\widehat{\phi})\omega_\pi(\tau)^{1-n}\gamma(s,\pi,\As,\psi')^{-1}d\mu_{\overline{G_n(E)}}(\pi)= \\
\nonumber & \lambda_{E/F}(\psi')^{-\frac{n(n+1)}{2}}\lvert\tau\rvert_E^{\frac{n(n-1)}{4}}\int_{\Temp(U(n))/\stab}Z(1,\widetilde{W}_{f,BC_n(\sigma)}(1,.),\widehat{\phi}) \omega_{BC_n(\sigma)}(\tau)^{1-n} \frac{\gamma^*(0,\sigma,\Ad,\psi')}{\lvert S_\sigma\rvert}d\sigma.
\end{align}
Let $\sigma\in \Temp(U(n))$ and $\pi=BC_n(\sigma)$. By \eqref{eq -2 L parameters, LLC, basechange}, we have
$$\displaystyle \omega_{\pi}(\tau)^{n-1}=\omega_{\pi}(\tau)^{n-1}\eta'_n(\tau)^{(n-1)}=\omega_{\sigma}(-1)^{n-1} \eta_{E/F}(-1)^{\frac{n(n-1)^2}{2}}$$
whereas by Lemma \ref{lem 2 Zeta function} and the fact that $\phi(0)=1$, we have
$$\displaystyle Z(1,\widetilde{W}_{f,BC_n(\sigma)}(1,.),\widehat{\phi})=c_1(\pi)\beta(W_{f,BC_n(\sigma)}(1,.)).$$
Combined with \eqref{eq 2 Cor to lim spectrale} and \eqref{eq 3 Cor to lim spectrale} this gives \eqref{eq 1 Cor to lim spectrale} and ends the proof of the corollary. $\blacksquare$

\section{A Plancherel formula for $G_n(F)\backslash G_n(E)$}\label{Part III}

\subsection{Pr\'ecis on Plancherel decompositions}

Let $G$ be a reductive algebraic group over $F$ which acts on the right of a smooth $F$-algebraic variety $X$. For every $\pi\in \Irr(G(F))$, we denote by $\cS(X(F))_\pi$ the {\em $\pi$-isotypic quotient of $\cS(X(F))$} i.e. the quotient of $\cS(X(F))$ by the intersection of kernels of all continuous $G(F)$-equivariant linear maps $\cS(X(F))\to \pi$. Let $\Irr_{\unit}(G)\subset \Irr(G)$ be the subset of {\em unitarizable} irreducible representations that we equip with the Fell topology \cite[Sect. 18.1]{Dix}. Assume moreover that $X(F)$ is equipped with a $G(F)$-invariant measure. We denote by $L^2(X(F))$ the Hilbert space of square-integrable functions for this measure and by $(.,.)_X$ the corresponding inner product. The action by right translation induces a unitary representation of $G(F)$ on $L^2(X(F))$. Since $G(F)$ is a postliminal locally compact group \cite[13.11.12]{Dix}, \cite{Ber0}, by \cite[Theorem 8.6.6]{Dix} there exists a unique (in a suitable sense) direct integral decomposition of $L^2(X(F))$
$$\displaystyle L^2(X(F))=\int_{\Irr_{\unit}(G)}^{\oplus}\cH_\pi d\mu_X(\pi).$$
Here $\mu_X$ is a Borel measure on $\Irr_{\unit}(G)$, and $\pi\mapsto \cH_\pi$ is a measurable field of continuous unitary representations of $G(F)$ (in the sense of \cite[Sect. 18.7]{Dix}) with $\cH_\pi$ isomorphic for $\mu_X$-almost all $\pi$ to an Hilbert direct sum of copies of (an Hilbert completion of) $\pi$. Such a decomposition is usually called a {\em Plancherel decomposition} for $L^2(X(F))$.

According to Gelfand-Kostyuchenko and Bernstein \cite{Ber3} (see also \cite[Sect. 3.3]{Li}), a Plancherel decomposition for $L^2(X(F))$ is equivalent to the following set of data:
\begin{itemize}
\item A Borel measure $\mu_X$ on $\Irr_{\unit}(G)$;
\item For $\mu_X$-almost all $\pi\in \Irr_{\unit}(G)$, a continuous Hermitian form $(.,.)_{X,\pi}$ on $\cS(X(F))$ which is $G(F)$-invariant positive semi-definite (i.e. $(\varphi,\varphi)_{X,\pi}\geqslant 0$ for every $\varphi\in \cS(X(F))$; from now on we will call such an Hermitian form a {\em semi-positive scalar product}) and factorizes through the quotient $\cS(X(F))\to \cS(X(F))_\pi$;
\end{itemize}
such that the following condition is satisfied: for every $\varphi_1,\varphi_2\in \cS(X(F))$, the function $\pi\mapsto (\varphi_1,\varphi_2)_{X,\pi}$ is $\mu_X$-integrable and we have
\begin{align}\label{eq 1 Planch formula}
\displaystyle (\varphi_1,\varphi_2)_{X}=\int_{\Irr_{\unit}(G)} (\varphi_1,\varphi_2)_{X,\pi} d\mu_X(\pi)
\end{align}
It is under this form that we will describe the Plancherel decomposition for $L^2(G_n(F)\backslash G_n(E))$ in the next section.

Actually, we don't even need to assume $(.,.)_{X,\pi}$ to be a semi-positive scalar product: once \eqref{eq 1 Planch formula} is satisfied for some continuous $G(F)$-invariant Hermitian forms $(.,.)_{X,\pi}$ factorizing through $\cS(X(F))\to \cS(X(F))_\pi$ then these forms are automatically semi-positive scalar products (for $\mu_X$-almost all $\pi$). This is the content of \cite[Proposition 6.1.1]{SV} in the $p$-adic case. The proof of {\em loc. cit.} actually also applies in the Archimedean case provided we replace each occurrence of ``the Bernstein center'' in it by ``Arthur's multiplier algebra'' (\cite{Art2}, \cite{Del}). We will actually need this result in a slightly different form but for which the proof of {\em loc. cit.} is still valid and we content ourself to state it referring however the interested reader to Lemma \ref{lem separation} for a similar kind of argument.

\begin{prop}\label{prop Planch meas}
Assume given
\begin{itemize}
\item A Borel measure $\mu_X$ on $\Irr_{\unit}(G)$;
\item For $\mu_X$-almost all $\pi\in \Irr_{\unit}(G)$, a continuous sesquilinear form $(.,.)_{X,\pi}$ on $\cS(X(F))$ which is $G(F)$-invariant and factorizes through the quotient $\cS(X(F))\to \cS(X(F))_\pi$
\end{itemize}
such that for every $\varphi_1,\varphi_2\in \cS(X(F))$, the function $\pi\mapsto (\varphi_1,\varphi_2)_{X,\pi}$ is $\mu_X$-integrable and \eqref{eq 1 Planch formula} is satisfied. Then, for $\mu_X$-almost all $\pi$ the form $(.,.)_{X,\pi}$ is a semi-positive scalar product so that the data $(\mu_X, (.,.)_{X,\pi})$ induce a Plancherel decomposition for $L^2(X(F))$.
\end{prop}

Finally we remark that the restriction of the Fell topology to $\Temp(G_n(E))\subset \Irr_{\unit}(G_n(E))$ coincides with the natural topology on $\Temp(G_n(E))$ that was described in Section \ref{Section representations}. This can be seen as follows: let $(\pi_n)_n$ be a sequence in $\Temp(G_n(E))$ converging to $\pi\in \Temp(G_n(E))$ for the topology defined in Section \ref{Section representations}. Then, there exists a Levi subgroup $M\subseteq G_n$, $\sigma\in \Pi_2(M(E))$ and a sequence $(\chi_n)_n$ of unramified characters of $M(E)$ converging to the trivial character such that $\pi_n=i_{M(E)}^{G_n(E)}(\sigma\otimes \chi_n)$ for $n$ large enough and $\pi=i_{M(E)}^{G_n(E)}(\sigma)$. Identifying each $\pi_n$ and $\pi$ with a compact model $i_{K_P}^{K}(\sigma)$ (where $K=K_{n,E}$, $P\in \mathcal{P}(M)$ and $K_P=K\cap P(E)$) and equipping this last space with the usual $K$-invariant scalar product we have $(\pi_n(g)e,e)\to (\pi(g)e,e)$ uniformly on compact subsets for every $e\in i_{K_P}^{K}(\sigma)$. Therefore, by definition of the Fell topology, $\pi_n\to \pi$ in $\Irr_{\unit}(G_n(E))$. Conversely, if $(\pi_n)_n$ is a sequence of tempered representations converging to $\pi\in \Temp(G_n(E))$ for the Fell topology then by \cite[Theorem 2.2(i)]{Tad} and \cite{BDix} the infinitesimal character $\chi_{\pi_n}$ of $\pi_n$ converges point-wise to $\chi_\pi$ where in the $p$-adic case by ``infinitesimal character'' we mean the corresponding character of the Bernstein center. However, the map sending $\pi\in \Temp(G_n(E))$ to $\chi_\pi$ is proper (for the topology of Section \ref{Section representations} on $\Temp(G_n(E))$ and the topology of point-wise convergence on the set of all infinitesimal characters): this follows from \cite[Th\'eor\`eme VIII.1.2]{Wald1} in the $p$-adic case and \cite[Corollary 35 p.81]{H-C0} in the Archimedean case. Since the sequence $(\pi_n)_n$ can only have $\pi$ as a limit point in $\Temp(G_n(E))$ (because, as we just saw, a limit point in $\Temp(G_n(E))$ would also be a limit point for the Fell topology) this shows that $\pi_n\to \pi$ in $\Temp(G_n(E))$ and thus the claim follows.

\subsection{The theorem}\label{Section the Theorem}

Let $n\geqslant 1$ and set $Y_n=G_n(F)\backslash G_n(E)$ that we equip with the quotient measure. There is a surjection $\cS(G_n(E))\to \cS(Y_n)$ given by
$$\displaystyle f\mapsto \varphi_f(x)=\int_{G_n(F)} f(hx)dh$$
which identifies $\cS(Y_n)$ with the space of $G_n(F)$-coinvariants of $\cS(G_n(E))$ for the left regular action (see \cite[Corollary D.5]{AL}). For $\pi\in \Temp(G_n(E))$, we define a semi-positive scalar product on $\cS(G_n(E))$ by
$$\displaystyle (f_1,f_2)_{Y_n,\pi}=\sum_{W\in \cB(\pi,\psi_n)} \beta(R(f_1)W)\overline{\beta(R(f_2)W)},\; f_1,f_2\in \cS(G_n(E))$$
where $\cB(\pi,\psi_n)$ is an orthonormal basis of $\cW(\pi,\psi_n)$ as in Section \ref{section Plancherel Whitt} and $\beta$ is the linear form defined in Section \ref{Section betan}.
\begin{lem}\label{lem scalar products}
The expression defining $(.,.)_{Y_n,\pi}$ is absolutely convergent and does not depend on the choice of the basis $\cB(\pi,\psi_n)$. Moreover, $(.,.)_{Y_n,\pi}$ is continuous $G_n(E)$-invariant (for the action by right translation) and if $\pi\in BC_n(\Temp(U(n)))$ it factorizes through the quotient $\cS(Y_n)_{\pi^\vee}$. Finally, the function $\pi\in \Temp(G_n(E))\mapsto (f_1,f_2)_{Y_n,\pi}$ belongs to $\cS(\Temp(G_n(E)))$ for every $f_1,f_2\in \cS(G_n(E))$ and the sesquilinear map
$$\displaystyle \cS(G_n(E))\times \cS(G_n(E))\to \cS(\Temp(G_n(E)))$$
$$\displaystyle (f_1,f_2)\mapsto \left(\pi\mapsto (f_1,f_2)_{Y_n,\pi} \right)$$
is continuous. 
\end{lem}

\noindent\ul{Proof}: Let $f_1,f_2\in \cS(G_n(E))$. By \eqref{eq 3 Planch Whitt} we have
$$\displaystyle W_{\overline{f_2}\star f_1^\vee,\pi}=R_1(\overline{f_2})W_{f^\vee_1,\pi}=\lvert \tau\rvert_E^{-n(n-1)/2}\sum_{W\in \cB(\pi,\psi_n)} \overline{R(f_2)W}\otimes R(f_1)W$$
the sum being absolutely convergent in $\cC^w(N_n(E)\backslash G_n(E)\times N_n(E)\backslash G_n(E),\psi_n^{-1}\boxtimes \psi_n)$ (where we have denoted by $R_1(\overline{f_2})W_{f^\vee_1,\pi}$ the right regular action of $\overline{f_2}$ on the first variable of $W_{f^\vee_1,\pi}$). Applying the continuous functional $\beta\ctens \beta$ to this decomposition we get that $(f_1,f_2)_{Y_n,\pi}$ is indeed defined by an absolutely convergent expression and that
\begin{align}\label{eq 0 enonce Planch}
\displaystyle (f_1,f_2)_{Y_n,\pi}=\lvert \tau\rvert_E^{n(n-1)/2}(\beta\ctens \beta)(W_{\overline{f_2}\star f_1^\vee,\pi}).
\end{align}
This also shows that $(f_1,f_2)_{Y_n,\pi}$ does not depend on the choice of the basis $\cB(\pi,\psi_n)$ and that $(.,.)_{Y_n,\pi}$ is right $G_n(E)$-invariant and continuous (by continuity of convolution and of the map $f\in \cS(G(F))\mapsto W_{f,\pi}\in \cC^w(N_n(E)\backslash G_n(E)\times N_n(E)\backslash G_n(E),\psi_n^{-1}\boxtimes \psi_n)$). Moreover, it follows from Proposition \ref{prop 1 Plancherel Whitt} that the function $\pi\in \Temp(G_n(E))\mapsto (f_1,f_2)_{Y_n,\pi}$ belongs to $\cS(\Temp(G_n(E)))$. The continuity of the sesquilinear map
$$\displaystyle \cS(G_n(E))\times \cS(G_n(E))\to \cS(\Temp(G_n(E)))$$
$$\displaystyle (f_1,f_2)\mapsto \left(\pi\in \Temp(G_n(E))\mapsto (f_1,f_2)_{Y_n,\pi} \right)$$
then follows from the closed graph theorem. Since $(f_1,f_2)_{Y_n,\pi}$ only depends on $\pi(f_1)$ and $\pi(f_2)$, we have that $(.,.)_{Y_n,\pi}$ factorizes through the maximal $\pi^\vee$-isotypic quotient of $\cS(G_n(E))$ (for the right $G_n(E)$-action). Finally, by \eqref{prop0 form betan}, if $\pi\in BC_n(\Temp(U(n)))$, the form $(.,.)_{Y_n,\pi}$ is $G_n(F)$-invariant for the left regular action so that it factorizes through $\cS(Y_n)$, hence through $\cS(Y_n)_{\pi^\vee}$, in both variables. $\blacksquare$

By the above lemma, for every $\pi\in BC_n(\Temp(U(n)))$, $(.,.)_{Y_n,\pi}$ induces a semi-positive scalar product on $\cS(Y_n)$ that we denote the same way. Recall that $(.,.)_{Y_n}$ stands for the $L^2$-inner product on $L^2(Y_n)$. We can now state the main theorem of this chapter.

\begin{theo}\label{theo Plancherel}
For every $\varphi_1,\varphi_2\in \cS(Y_n)$, we have
$$\displaystyle (\varphi_1,\varphi_2)_{Y_n}=\int_{\Temp(U(n))/stab} (\varphi_1,\varphi_2)_{Y_n,BC_n(\sigma)} \frac{\lvert \gamma^*(0,\sigma,\Ad,\psi')\rvert}{\lvert S_\sigma\rvert} d\sigma$$
where the right-hand side is absolutely convergent. Moreover, we have
$$\displaystyle c(\sigma)\gamma^*(0,\sigma,\Ad,\psi')=\lvert \gamma^*(0,\sigma,\Ad,\psi')\rvert$$
for every $\sigma\in \Temp(U(n))$, where $c(\sigma)$ is the constant defined in Section \ref{Section Corollary}.
\end{theo}

\noindent First, the convergence of the right-hand side of the proposition follows directly from Lemma \ref{lem 1 L parameters, LLC, basechange}, Lemma \ref{lem scalar products} and \eqref{basic estimates spectral measure}.

\noindent Let $\varphi_1,\varphi_2\in \cS(Y_n)$ and choose $f_1,f_2\in \cS(G_n(E))$ such that $\varphi_i=\varphi_{f_i}$ for $i=1,2$. Then, we have
\begin{align}\label{eq 1 enonce Planch}
\displaystyle (\varphi_1,\varphi_2)_{Y_n}=\int_{Y_n}\int_{G_n(F)\times G_n(F)} f_1(h_1x)\overline{f_2(h_2x)}dh_1dh_2dx=\int_{G_n(F)} f(h)dh
\end{align}
where we have set $f=\overline{f_2}\star f_1^\vee$. Moreover, by \eqref{eq 0 enonce Planch} we also have
\begin{align}\label{eq 2 enonce Planch}
\displaystyle (\varphi_1,\varphi_2)_{Y_n,\pi}=\lvert \tau\rvert_E^{n(n-1)/2}(\beta\ctens \beta)(W_{f,\pi})
\end{align}
for every $\pi\in BC_n(\Temp(U(n)))$.

\subsection{A local unfolding identity}

Recall that in Section \ref{section Plancherel Whitt} we have associated to any $f\in \cS(G_n(E))$ a function $W_f\in \cC^w(N_n(E)\backslash G_n(E)\times N_n(E)\backslash G_n(E),\psi_n^{-1}\boxtimes \psi_n)$.

\begin{prop}\label{prop unfolding}
For every $f\in \cS(G_n(E))$, we have
$$\displaystyle \int_{G_n(F)} f(h)dh=\lvert \tau\rvert_E^{n(n-1)/4}\int_{N_n(F)\backslash P_n(F)} \int_{N_{n}(F)\backslash G_{n}(F)} W_f(p,h)dhdp$$
where the right-hand side is given by an absolutely convergent expression.
\end{prop}

\vspace{2mm}

\noindent\ul{Proof}: First, we check the absolute convergence of the right-hand side. Let $\nu$ and $r$ be the functions on $G_n(E)$ defined by $\nu(g)=\lvert \det g\rvert_E^{-1/4}$ and $r(g)=(1+\lVert e_ng\rVert^{1/2})$ where $e_n=(0,\ldots,0,1)$ and $\lVert .\rVert$ denotes the following norm on $E^n$: $(x_1,\ldots,x_n)\mapsto \max(\lvert x_1\rvert_E,\ldots,\lvert x_n\rvert_E)$ in the $p$-adic case, $(x_1,\ldots,x_n)\mapsto\left(\lvert x_1\rvert_E^2+\ldots + \lvert x_n\rvert_E^2\right)^{1/2}$ in the Archimedean case. Then, for every integer $N\geqslant 1$, $r^N\nu f\in \cS(G_n(E))$ and moreover we easily check that
$$W_{r^N \nu f}(a_{n-1}k_{n-1},a_nk_n)=(1+\lvert a_{n,n}\rvert)^{N}\lvert \det a^{-1}_{n-1}a_n\rvert^{-1/2}W_f(a_{n-1}k_{n-1},a_nk_n)$$
for every $(a_{n-1},a_n,k_{n-1},k_n)\in A_{n-1}(F)\times A_n(F)\times K_{n-1}\times K_n$. Therefore by Lemma \ref{lem 2 space of functions} (or rather its obvious analog for $G_n(E)\times G_n(E)$) applied to $W_{r^N\nu f}$ and the Iwasawa decomposition $G_n(F)=N_n(F)A_n(F)K_n$, $P_n(F)=N_n(F)A_{n-1}(F)K_{n-1}$, it suffices to show the existence of $N\geqslant 1$ such that for every $d>0$ the following expression converges:
\[\begin{aligned}
\displaystyle & \int_{A_n(F)\times A_{n-1}(F)} (1+\lvert a_{n,n}\rvert)^{-N}\lvert \det a^{-1}_{n-1}a_n\rvert^{1/2}\prod_{i=1}^{n-1} \left(1+\lvert \frac{a_{n-1,i}}{a_{n-1,i+1}}\rvert \right)^{-N}\prod_{i=1}^{n-1} \left(1+\lvert \frac{a_{n,i}}{a_{n,i+1}}\rvert \right)^{-N}   \\
 &  \delta_{n,E}(a_n)^{1/2}\delta_{n,E}(a_{n-1})^{1/2} \sigma(a_n)^{d}\sigma(a_{n-1})^d\delta_n(a_n)^{-1}\delta_{n-1}(a_{n-1})^{-1}da_nda_{n-1} \\
 & =\int_{A_n(F)} (1+\lvert a_{n,n}\rvert)^{-N} \prod_{i=1}^{n-1} \left(1+\lvert \frac{a_{n,i}}{a_{n,i+1}}\rvert \right)^{-N} \sigma(a_n)^d \lvert \det a_n\rvert^{1/2} da_n \\
 & \times \int_{A_{n-1}(F)} (1+\lvert a_{n-1,n-1}\rvert)^{-N} \prod_{i=1}^{n-2} \left(1+\lvert \frac{a_{n-1,i}}{a_{n-1,i+1}}\rvert \right)^{-N} \sigma(a_{n-1})^d \lvert \det a_{n-1}\rvert^{1/2} da_{n-1}.
\end{aligned}\]
But, by Lemma \ref{lem basic estimates} both integrals above are convergent for $N\gg 1$ (and any $d$).

We now show the identity of the proposition by induction on $n$. For $n=1$, this equality is a tautology. Assume from now on that $n\geqslant 2$ and the result is known for smaller values of $n$. Then, we can write
\begin{align}\label{eq 1 Prop unfolding}
\displaystyle & \int_{N_n(F)\backslash P_n(F)}\int_{N_{n}(F)\backslash G_{n}(F)} W_f(p,h)dhdp= \\
\nonumber & \int_{P_{n}(F)\backslash G_{n}(F)}\int_{P_{n-1}(F)U_n(F)\backslash P_n(F)}\int_{N_{n-1}(F)\backslash P_{n-1}(F)} \int_{N_{n-1}(F)\backslash G_{n-1}(F)} W_f(p'p,h'h)\lvert \det(p'h')\rvert^{-1}  dh'dp'dpdh
\end{align}
Let $(p,h)\in P_n(F)\times G_{n}(F)$ and set $f':=L(p)R(h)f$, $\displaystyle \widetilde{f'}(g):=\lvert \det g\rvert_E^{-1/2} \int_{U_n(E)} f'(vg)\psi_n(v)^{-1}dv$ for every $g\in G_{n-1}(E)$. Then, $\widetilde{f'}\in \cS(G_{n-1}(E))$ and we have
\[\begin{aligned}
\displaystyle W_f(p'p,h'h) & =W_{f'}(p',h')=\int_{N_{n-1}(E)} \int_{U_n(E)} f'({p'}^{-1}vuh')\psi_n(v)^{-1}dv\psi_n(u)^{-1}du \\
 & =\lvert \det p'\rvert^2 \int_{N_{n-1}(E)} \int_{U_n(E)} f'(v{p'}^{-1}uh')\psi_n(v)^{-1}dv\psi_n(u)^{-1}du \\
 & =\lvert \det(p'h')\rvert \int_{N_{n-1}(E)} \widetilde{f'}({p'}^{-1}uh')\psi_{n-1}(\epsilon_{n-1}u\epsilon_{n-1}^{-1})^{-1}du \\
 & =\lvert \det(p'h')\rvert W_{\widetilde{f'}}(\epsilon_{n-1}p',\epsilon_{n-1}h')
\end{aligned}\]
for every $(p',h')\in P_{n-1}(F)\times G_{n-1}(F)$ where we have set
$$\displaystyle \epsilon_{n-1}=\begin{pmatrix} (-1)^{n-2} \\ & \ddots \\ & & 1\end{pmatrix}\in A_{n-1}(F)\cap P_{n-1}(F).$$
By the induction hypothesis, we thus get
\[\begin{aligned}
\displaystyle & \int_{N_{n-1}(F)\backslash P_{n-1}(F)}\int_{N_{n-1}(F)\backslash G_{n-1}(F)} W_f(p'p,h'h)\lvert \det(p'h')\rvert^{-1}  dh'dp' \\
 & =\int_{N_{n-1}(F)\backslash P_{n-1}(F)}\int_{N_{n-1}(F)\backslash G_{n-1}(F)} W_{\widetilde{f'}}(p',h') dh'dp' \\
 & =\lvert \tau\rvert_E^{(n-1)(n-2)/4}\int_{G_{n-1}(F)} \widetilde{f'}(h') dh'=\lvert \tau\rvert_E^{(n-1)(n-2)/4}\int_{G_{n-1}(F)} \lvert \det h'\rvert^{-1}\int_{U_n(E)} f(p^{-1}vh'h)\psi_n(v)^{-1}dvdh'.
\end{aligned}\]
Combining this with \eqref{eq 1 Prop unfolding}, we obtain
\begin{align}\label{eq 3 Prop unfolding}
\displaystyle & \lvert \tau\rvert_E^{-(n-1)(n-2)/4}\int_{N_n(F)\backslash P_n(F)}\int_{N_{n}(F)\backslash G_{n}(F)} W_f(p,h)dhdp \\
\nonumber & =\int_{P_{n}(F)\backslash G_{n}(F)}\int_{P_{n-1}(F)U_n(F)\backslash P_n(F)}\int_{G_{n-1}(F)} \lvert \det h'\rvert^{-1} \int_{U_n(E)}f(p^{-1} vh'h)\psi_n(v)^{-1}dvdh' dpdh \\
\nonumber & =\int_{P_{n}(F)\backslash G_{n}(F)}\int_{P_{n-1}(F)U_n(F)\backslash P_n(F)}\int_{G_{n-1}(F)} \lvert \det p^{-1}h'\rvert^{-1} \int_{U_n(E)}f(vp^{-1}h' h )\psi_n(pvp^{-1})^{-1}dvdh' \lvert \det p\rvert dpdh \\
\nonumber & =\int_{P_{n}(F)\backslash G_{n}(F)}\int_{P_{n-1}(F)U_n(F)\backslash P_n(F)}\int_{G_{n-1}(F)} \lvert \det h' \rvert^{-1} \int_{U_n(E)}f(vh'h)\psi_n(pvp^{-1})^{-1}dvdh' \lvert \det p\rvert dpdh
\end{align}
Let $h\in G_n(F)$ and define $\varphi\in \cS(U_n(E))$ by $\displaystyle \varphi(v):=\int_{G_{n-1}(F)} (R(h)f)(vh') \lvert \det h'\rvert^{-1} dh'$. Then, by \eqref{eq 1 measures} we have
\[\begin{aligned}
\displaystyle & \int_{P_{n-1}(F)U_n(F)\backslash P_n(F)}\int_{G_{n-1}(F)} \lvert \det h' \rvert^{-1} \int_{U_n(E)}f(vh'h)\psi_n(pvp^{-1})^{-1}dvdh' \lvert \det p\rvert dp \\
 & =\int_{P_{n-1}(F)\backslash G_{n-1}(F)} \int_{U_n(E)} \varphi(v) \psi_n(h'' v{h''}^{-1})^{-1}dv \lvert \det h'' \rvert dh''\\
 & =\lvert \tau\rvert_E^{(n-1)/2}\int_{U_n(F)} \varphi(v)dv=\lvert \tau\rvert_E^{(n-1)/2}\int_{U_n(F)}\int_{G_{n-1}(F)} f(vh'h)\lvert \det h'\rvert^{-1} dh'dv
\end{aligned}\]
Plugging this into \eqref{eq 3 Prop unfolding}, we obtain
\[\begin{aligned}
\displaystyle & \int_{N_n(F)\backslash P_n(F)}\int_{N_{n}(F)\backslash G_{n}(F)} W_f(p,h)dhdp \\
 & =\lvert \tau\rvert_E^{n(n-1)/4}\int_{P_n(F)\backslash G_n(F)} \int_{U_n(F)}\int_{G_{n-1}(F)} f(vh'h)\lvert \det h'\rvert^{-1} dh'dvdh \\
 & =\lvert \tau\rvert_E^{n(n-1)/4} \int_{G_n(F)} f(h) dh
\end{aligned}\]
ending the proof of the proposition. $\blacksquare$

\subsection{End of the proof of Theorem \ref{theo Plancherel}}

Let $\varphi_1,\varphi_2\in \cS(Y_n)$ and $f_1,f_2,f\in \cS(G_n(E))$ be as in the end of Section \ref{Section the Theorem}. By Proposition \ref{prop unfolding} and Corollary \ref{cor lim spectrale}, we have
\begin{align}\label{eq 1 proof of theo Planch}
\displaystyle \int_{G_n(F)} f(h)dh=\lvert \tau\rvert_E^{n(n-1)/2}\int_{N_n(F)\backslash P_n(F)} \int_{\Temp(U(n))/\stab} \beta(W_{f,BC_n(\sigma)}(p,.)) \frac{\gamma^*(0,\sigma,\Ad,\psi')}{\lvert S_\sigma\rvert} c(\sigma) d\sigma dp
\end{align}
From Lemma \ref{lem continuity of beta} and Lemma \ref{lem 0 space of functions}, we deduce that the linear map
\begin{align}\label{eq 2 proof of theo Planch}
\displaystyle \cC^w(N_n(E)\backslash G_n(E)\times N_n(E)\backslash G_n(E),\psi_n^{-1}\boxtimes \psi_n)\to \cC^w(N_n(E)\backslash G_n(E),\psi^{-1}_n)
\end{align}
$$W\mapsto \left( g\mapsto \beta(W(g,.))\right)$$
is continuous. Therefore, by Lemma \ref{lem continuity of beta} again, Proposition \ref{prop 1 Plancherel Whitt}, Lemma \ref{lem 1 L parameters, LLC, basechange} and \eqref{basic estimates spectral measure}, we see that the right-hand side of \eqref{eq 1 proof of theo Planch} is an absolutely convergent expression and therefore
\[\begin{aligned}
\displaystyle \int_{G_n(F)} f(h)dh & = \lvert \tau\rvert_E^{n(n-1)/2}\int_{\Temp(U(n))/\stab} \int_{N_n(F)\backslash P_n(F)}\beta(W_{f,BC_n(\sigma)}(p,.))dp \frac{\gamma^*(0,\sigma,\Ad,\psi')}{\lvert S_\sigma\rvert} c(\sigma) d\sigma \\
 & =\lvert \tau\rvert_E^{n(n-1)/2}\int_{\Temp(U(n))/\stab} (\beta\ctens \beta)(W_{f,BC_n(\sigma)})\frac{\gamma^*(0,\sigma,\Ad,\psi')}{\lvert S_\sigma\rvert} c(\sigma) d\sigma \\
 & =\int_{\Temp(U(n))/stab} (\varphi_1,\varphi_2)_{Y_n,BC_n(\sigma)} \frac{\gamma^*(0,\sigma,\Ad,\psi')}{\lvert S_\sigma\rvert} c(\sigma)d\sigma
\end{aligned}\]
where in the last equality we have used \eqref{eq 2 enonce Planch}. Together with \eqref{eq 1 enonce Planch} this shows
$$\displaystyle (\varphi_1,\varphi_2)_{Y_n}=\int_{\Temp(U(n))/stab} (\varphi_1,\varphi_2)_{Y_n,BC_n(\sigma)} \frac{\gamma^*(0,\sigma,\Ad,\psi')}{\lvert S_\sigma\rvert} c(\sigma)d\sigma$$
From this, Lemma \ref{lem scalar products} and Proposition \ref{prop Planch meas}, we deduce that $c(\sigma)\gamma^*(0,\sigma,\Ad,\psi')=\lvert \gamma^*(0,\sigma,\Ad,\psi')\rvert$ for almost all $\sigma\in \Temp(U(n))/stab$, hence for all. Given this, the above identity is precisely the content of Theorem \ref{theo Plancherel}. $\blacksquare$

\section{Applications to the Ichino-Ikeda and formal degree conjectures for unitary groups}\label{Part IV}

In this chapter, we prove the main results of this paper stated as Theorem \ref{theo 3 intro} and Theorem \ref{theo 4 intro} in the introduction. As explained there, these will be obtain through comparisons of local relative trace formulas on the one hand and of relative unipotent orbital integrals on ther other hand via the Jacquet-Rallis transfer. The main result of the previous chapter (Theorem \ref{theo Plancherel}) is needed to get spectral expansions of these distributions on the linear side of the Jacquet-Rallis transfer. Besides this, we also need the existence of the Jacquet-Rallis transfer as well as its compatibility with Fourier transform on the level of Lie algebra. Both of these results are due to Zhang \cite{Zha1} (in the $p$-adic case) and Xue \cite{Xue} (in the Archimedean case). Moreover, to start the local comparison, we require a weak comparison of relative characters (see Proposition \ref{prop 1 weak comparison}). This was established in \cite{Beu2} using a comparison of global trace formulas and therefore involves the Jacquet-Rallis fundamental lemma of Yun and Gordon \cite{Yu}. Let us mention however that the author has recently given an alternative proof of this fundamental lemma \cite{Beu4} only based on techniques of harmonic analysis as well as the previous work of Zhang \cite{Zha1}.

\subsection{Notation, matching of orbits}\label{Section Notation, matching of orbits}

Fix an integer $n\geqslant 1$. In this chapter we will consider the following groups and subgroups:
\begin{itemize}
\item $G'=R_{E/F} G_{n,E}\times R_{E/F}G_{n+1,E}$ with its two subgroups $H_1=R_{E/F}G_{n,E}$ ({\em diagonally} embedded) and $H_2=G_{n,F}\times G_{n+1,F}$. We also equip $H_2(F)=G_n(F)\times G_{n+1}(F)$ with the character $\eta=\eta_n\boxtimes \eta_{n+1}$ where we recall that $\eta_k$ stands for the character of $G_k(F)$ given by $\eta_k(g)=\eta_{E/F}(\det g)^k$.
\item For $V$ a Hermitian space of dimension $n$ (with underlying Hermitian form $h$), we set $H^V=U(V)$ and $G^V=U(V)\times U(V')$ where $V'=V\oplus Ev_0$ is equipped with the Hermitian form $h'$ given by $h'(v_1+\lambda v_0,v_2+\mu v_0)=h(v_1,v_2)+\lambda \mu^c$ for all $v_1,v_2\in V$ and $\lambda,\mu\in E$. We consider $H^V$ as a subgroup of $G^V$ through the natural {\em diagonal} embedding $H^V \hookrightarrow G^V$.
\end{itemize}

We let $H_1\times H_2$ (resp. $H^V\times H^V$) act on $G'$ (resp. $G^V$) by $(h_1,h_2)\cdot g=h_1gh_2^{-1}$. A geometric point $g\in G'$ (resp. $g\in G^V$) is said to be {\em regular semi-simple} if its stabilizer in $H_1\times H_2$ (resp. $H^V\times H^V$) for this action is trivial and its orbit $H_1gH_2$ (resp. $H^VgH^V$) is Zariski closed. We denote by $G'_{\rs}\subset G'$ and $G^V_{\rs}\subset G^V$ the open subsets of regular semi-simple elements. These are not empty (\cite[\S 2.1]{Zh2}) and the actions of $H_1\times H_2$ and $H^V\times H^V$ on $G'_{\rs}$ and $G^V_{\rs}$ respectively are free. Let $\cB$ and $\cB^V$ be the geometric quotients $H_1\backslash G'/H_2$ and $H^V\backslash G^V/H^V$ respectively. Geometric points in $\cB$ and $\cB^V$ correspond bijectively to closed geometric orbits in $G'$ and $G^V$ respectively. Let $\cB_{\rs}\subset \cB$ and $\cB^V_{\rs}\subset \cB^V$ be the open subsets corresponding to regular semi-simple elements (or orbits) i.e. $\cB_{\rs}=H_1\backslash G'_{\rs}/H_2$ and $\cB^V_{\rs}=H^V\backslash G^V_{\rs}/H^V$. Then, there is a natural isomorphism $\cB\simeq \cB^V$ which restricts to $\cB_{\rs}\simeq \cB^V_{\rs}$ (\cite[Lemme 15.1.4.1]{CZ}). Moreover, when taking $F$-points this isomorphism induces a bijection (\cite[Lemme 2.3]{Zh2})
\begin{align}\label{eq 1 matching of orbits}
\displaystyle H_1(F)\backslash G'_{\rs}(F)/H_2(F)\simeq \bigsqcup_V H^V(F)\backslash G^V_{\rs}(F)/H^V(F)
\end{align}
where the disjoint union of the right-hand side runs over a set of representatives of isomorphism classes of Hermitian spaces of dimension $n$. Two regular semi-simple elements $\gamma\in G'_{\rs}(F)$ and $\delta\in G^V_{\rs}(F)$ whose orbits correspond to each other by the above bijection will be said to {\em match} and we will abbreviate this by the notation $\gamma \leftrightarrow \delta$.

Since the action is free, the quotient of the (restriction to $G'_{\rs}(F)$ of the) Haar measure on $G'(F)$ by the Haar measure on $H_1(F)\times H_2(F)$ defines a measure on $H_1(F)\backslash G'_{\rs}(F)/H_2(F)$ that we shall denote by $d\gamma$. Similarly for every Hermitian space $V$ of dimension $n$, the quotient of the Haar measure on $G^V(F)$ by the Haar measure on $H^V(F)\times H^V(F)$ gives rise to a measure $d\delta=d^V \delta$ on $H^V(F)\backslash G^V_{\rs}(F)/H^V(F)$.

\begin{lem}\label{lem 1 matching of orbits}
The bijection \eqref{eq 1 matching of orbits} is measure preserving, i.e. it sends the measure $d\gamma$ to the sum of the measures $d^V\delta$.
\end{lem}

\noindent\ul{Proof}: Let $V$ be a Hermitian space of dimension $n$. Then, by our normalization of measures (see Section \ref{Section measures}), $d\gamma=\lvert \omega'\rvert_{\psi'}$ and $d^V\delta=\lvert \omega^V\rvert_{\psi'}$ where $\omega'$ and $\omega^V$ are certain volume forms on $\cB_{\rs,\overline{F}}=\cB^V_{\rs,\overline{F}}$. More precisely, $\omega^V$ is obtained as the quotient of the volume form $\omega_{G^V_{\oF}}$ by $\omega_{H^V_{\oF}\times H^V_{\oF}}$ fixed in Section \ref{Section measures} whereas $\omega'$ is the quotient of $\omega_{G'_{\oF}}$ by $\omega_{H_{1,\oF}\times H_{2,\oF}}$. Obviously, it suffices to show that $\omega'$ and $\omega^V$ are equal up to a sign. We have $\cB_{\rs,\overline{F}}=H_{1,\overline{F}}\backslash G'_{\rs,\oF}/H_{2,\oF}$ and using the $\overline{F}$-algebra isomorphism $E\otimes_F \oF\simeq \oF\times \oF$, $x\otimes y\mapsto (xy,x^cy)$ we get identifications
$$\displaystyle H_{1,\oF}\simeq G_{n,\oF}\times G_{n,\oF},\; G'_{\oF}\simeq (G_{n,\oF}\times G_{n,\oF})\times (G_{n+1,\oF}\times G_{n+1,\oF}),\; H_{2,\oF}\simeq G_{n,\oF}\times G_{n+1,\oF}$$
such that the inclusion $H_{2,\oF}\subset G'_{\oF}$ is the product of the diagonal embeddings $G_{k,\oF}\hookrightarrow G_{k,\oF}\times G_{k,\oF}$ for $k=n,n+1$ whereas the inclusion $H_{1,\oF}\subset G'_{\oF}$ is the identity in the first component and the natural embedding $G_{n,\oF}\times G_{n,\oF}\hookrightarrow G_{n+1,\oF}\times G_{n+1,\oF}$ in the second. Similarly, using the same $\oF$-algebra isomorphism and a linear bijection $V\simeq E^n$ that we extend to $V'\simeq E^{n+1}$ by sending $v_0$ to $(0,\ldots,0,1)$ we get identifications
$$\displaystyle H^V_{\oF}\simeq G_{n,\oF},\; G^V_{\oF}\simeq G_{n,\oF}\times G_{n+1,\oF}$$
such that the inclusion $H^V_{\oF}\subset G^V_{\oF}$ is the identity in the first component and the natural embedding $G_{n,\oF}\hookrightarrow G_{n+1,\oF}$ in the second. Using these identifications, we get an isomorphism
\begin{align}\label{eq 2 matching of orbits}
\displaystyle H_{1,\oF}\simeq H^V_{\oF}\times H^V_{\oF}
\end{align}
and the projection onto the first components yields another isomorphism
\begin{align}\label{eq 3 matching of orbits}
\displaystyle G'_{\oF}/H_{2,\oF}\simeq G^V_{\oF} 
\end{align}
which is $G'_{\oF}\simeq G^V_{\oF}\times G^V_{\oF}$-equivariant thus inducing an identification $\cB_{\oF}=H_{1,\oF}\backslash G'_{\oF}/H_{2,\oF}\simeq H^V_{\oF}\backslash G^V_{\oF}/H^V_{\oF}=\cB^V_{\oF}$. This last isomorphism is by its very definition the same as before (base-changed to $\oF$). Moreover, both \eqref{eq 2 matching of orbits} and \eqref{eq 3 matching of orbits} obviously extends to the natural split forms over $\bZ$ of all the groups under consideration. Therefore (and since the group of outer automorphisms of $G_{n}$ has order $2$), \eqref{eq 2 matching of orbits} sends $\omega_{H_{1,\oF}}$ to $\pm \omega_{H^V_{\oF}\times H^V_{\oF}}$ and \eqref{eq 3 matching of orbits} sends the quotient of the volume form $\omega_{G'_{\oF}}$ by $\omega_{H_{2,\oF}}$ to $\pm \omega_{G^V_{\oF}}$. From this, it immediately follows that the isomorphism $\cB_{\rs, \oF}\simeq \cB^V_{\rs, \oF}$ sends $\omega'$ to $\pm \omega^V$. $\blacksquare$

\subsection{Linearization and Fourier transform}

At some point we will need to ``linearize'' certain expressions and thus we introduce the following extra notation:
\begin{itemize}
\item Let $S$ be the symmetric space $S=\{g\in R_{E/F}G_{n+1} \mid gg^c=1 \}$ and $\fs$ be its tangent space at the identity i.e. $\fs=\{Y\in R_{E/F}\gl_{n+1}\mid Y+Y^c=0 \}$.
\item For every Hermitian space $V$ of dimension $n$, let $\fu^V=\mathfrak{u}(V')$ be the Lie algebra of $U(V')$.
\end{itemize}

We let $G_n$ (resp. $U(V)$) act on $\fs$ (resp. $\fu^V$) by conjugation. As before, a geometric point $X\in \fs$ (resp. $X\in \fu^V$) is said to be {\em regular semi-simple} if its stabilizer in $G_n$ (resp. $U(V)$) for this action is trivial and the corresponding orbit is closed. We denote by $\fs_{\rs}\subset \fs$ and $\fu^V_{\rs}\subset \fu^V$ the open subsets of regular semi-simple elements. These are not empty (\cite[\S 3 \& 4]{JR}) and the actions of $G_n$ and $U(V)$ on $\fs_{\rs}$ and $\fu^V_{\rs}$ respectively are free. Moreover, there is also a natural isomorphism between geometric quotients $\fs/G_n\simeq \fu^V/U(V)$ which restricts to $\fs_{\rs}/G_n\simeq \fu^V_{\rs}/U(V)$ (see \cite[Proposition 2.2.2.1]{Chau}) and gives rise, when taking $F$-points, to a bijection (\cite[Lemmes 2.1.5.1, 2.1.5.3 \& Proposition 2.2.4.1]{Chau})
\begin{align}\label{eq 4 matching of orbits}
\displaystyle \fs_{\rs}(F)/G_n(F)\simeq \bigsqcup_V \fu^V_{\rs}(F)/U(V)(F)
\end{align}
where the disjoint union of the right-hand side runs, once again, over a set of representatives of isomorphism classes of Hermitian spaces of dimension $n$. Two regular semi-simple elements $Y\in \fs_{\rs}(F)$ and $X\in \fu^V_{\rs}(F)$ whose orbits correspond by the above bijection will be said to {\em match} and we will abbreviate this by $Y \leftrightarrow X$.

There is the following $G_n$-equivariant isomorphism
$$\displaystyle \nu: R_{E/F}G_{n+1}/G_{n+1}\to S$$
$$g\mapsto g(g^{-1})^c.$$
We use $\nu$ to transfer the measure on $G_{n+1}(E)/G_{n+1}(F)$ to a measure on $S(F)$ which in turn induces a Haar measure on $\fs(F)$ by taking its fiber at $1$ (this makes sense since the measure on $S(F)$ thus obtained is in the natural class of measures on this $F$-analytic manifold). Similarly, the Haar measure on $U(V')(F)$ induces one on $\fu^V(F)$. Now, since the actions are free, the quotient of the measure on $\fs(F)$ (resp. $\fu^V(F)$) by the Haar measure on $G_n(F)$ (resp. $U(V)(F)$) defines a measure on $\fs_{\rs}(F)/G_n(F)$ (resp. $\fu_{\rs}^V(F)/U(V)(F)$) that we shall denote by $dY$ (resp. $dX=d^V X$). By essentially the same arguments as for Lemma \ref{lem 1 matching of orbits} we have:

\begin{lem}\label{lem 2 matching of orbits}
The bijection \eqref{eq 4 matching of orbits} is measure preserving, i.e. it sends the measure $dY$ to the sum of the measures $d^VX$.
\end{lem}

Let $\fc:\fs\dashrightarrow S$ (resp. $\fc^V:\fu^V\dashrightarrow U(V')$) be the birational isomorphism (henceforth called {\em Cayley maps}) sending $X$ to $\frac{1+X/2}{1-X/2}$. Note that $\fc$ (resp. $\fc^V$) is $G_n$-equivariant (resp. $U(V)$-equivariant). By abuse of notation we will also write $\fc$ for $\fc^V$, hoping that it will not create any confusion for the attentive reader. We fix once and for all a small open neighborhood $\cU$ of $0$ in $(\fs/G_n)(F)=(\fu^V/U(V))(F)$ on the inverse image of which (both in $\fs(F)$ and $\fu^V(F)$ for every Hermitian space $V$ of dimension $n$) $\fc$ is well-defined as well as a cut-off function $\alpha\in C_c^\infty(\cU)$ (here we remark that the ``base'' $\fs/G_n$ is smooth, and even an affine space, see \cite[Proposition 2.1.5.2]{Chau}) such that $\alpha=1$ on some neighborhood of $0$. This allows to define two applications $\Phi'\in \cS(S(F))\mapsto \Phi'_\natural\in \cS(\fs(F))$ and $\Phi^V\in \cS(U(V')(F))\mapsto \Phi^V_\natural\in \cS(\fu^V(F))$ by
$$\displaystyle \Phi'_\natural(Y)=\left\{
    \begin{array}{ll}
        \eta'(\det(1-Y))^{-n}\alpha(Y)\Phi'(\fc(Y)) & \mbox{ if } Y\in \cU \\
        0 & \mbox{ otherwise}
    \end{array}
\right.$$
and
$$\displaystyle \displaystyle \Phi^V_\natural(X)=\left\{
    \begin{array}{ll}
        \alpha(X)\Phi^V(\fc(X)) & \mbox{ if } X\in \cU \\
        0 & \mbox{ otherwise}
    \end{array}
\right.$$
for every $Y\in \fs(F)$ and $X\in \fu^V(F)$. The presence of the extra factor $\eta'(\det(1-Y))^{-n}$ is justified a posteriori by \eqref{eq 1 matching of functions}. Note that in the $p$-adic case, up to shrinking $\cU$ we can make this factor identically $1$ on $\cU$.

We also define two applications $f'\in \cS(G'(F))\mapsto \widetilde{f'}\in \cS(S(F))$ and $f^V\in \cS(G^V(F))\mapsto \widetilde{f^V}\in \cS(U(V')(F))$ by
$$\displaystyle \widetilde{f^V}(g)=\int_{H^V(F)}f(h(1,g))dh,\;\;\; g\in U(V')(F)$$
and
$$\displaystyle \widetilde{f'}(s)=\int_{H_1(F)\times G_{n+1}(F)} f'(h_1(1,\nu^{-1}(s)h_{n+1})) \eta'_{n+1}(\nu^{-1}(s)h_{n+1}) dh_{n+1}dh_1,\;\;\; s\in S(F)$$
where $\nu$ is the isomorphism $G_{n+1}(E)/G_{n+1}(F)\simeq S(F)$ defined above and we recall that $\eta'_{n+1}$ is a character of $G_{n+1}(E)$ extending the character $\eta_{n+1}$ of $G_{n+1}(F)$ (see \ref{Section groups}).

We shall denote the composition of the two maps $f'\mapsto \widetilde{f'}$ and $\Phi'\mapsto \Phi'_\natural$ (resp. $f^V\mapsto \widetilde{f^V}$ and $\Phi^V\mapsto \Phi^V_\natural$) by $f'\in \cS(G'(F))\mapsto \widetilde{f'_\natural}\in \cS(\fs(F))$ (resp. $f^V\in \cS(G^V(F))\mapsto \widetilde{f^V_\natural}\in \cS(\fu^V(F))$). Thus, we have
\begin{align}\label{def linearization}
\displaystyle \widetilde{f'_\natural}(Y)=\left\{
\begin{array}{ll}
\eta'(\det(1-Y))^{-n}\alpha(Y)\widetilde{f'}(\fc(Y)) & \mbox{ if } Y\in \cU \\
0 & \mbox{ otherwise}
\end{array}
\right. , \\
\nonumber \widetilde{f^V_\natural}(X)=\left\{
\begin{array}{ll}
\alpha(X)\widetilde{f^V}(\fc(X)) & \mbox{ if } X\in \cU \\
0 & \mbox{ otherwise}
\end{array}
\right.
\end{align}
for every $f'\in \cS(G'(F))$, $f^V\in \cS(G^V(F))$, $Y\in \fs(F)$ and $X\in \fu^V(F)$.

We define two non-degenerate $G_n(F)$- and $U(V)(F)$-invariant symmetric bilinear forms $\langle .,.\rangle: \fs(F)\times \fs(F)\to F$ and $\langle .,.\rangle: \fu^V(F)\times \fu^V(F)\to F$ by
$$\displaystyle \langle X,Y\rangle=\Tr(XY)$$
for all $X,Y\in \fs(F)$ or $X,Y\in \fu^V(F)$. This allows to define Fourier transforms $\varphi'\in \cS(\fs(F))\mapsto \cF\varphi'\in \cS(\fs(F))$ and $\varphi^V\in \cS(\fu^V(F))\mapsto \cF\varphi^V\in \cS(\fu^V(F))$ by
$$\displaystyle \cF\varphi'(Y)=\int_{\fs(F)}\varphi'(Y')\psi'(\langle Y',Y\rangle) dY' \left(\mbox{resp. } \cF\varphi^V(X)=\int_{\fu^V(F)}\varphi^V(X')\psi'(\langle X',X\rangle) dX'\right)$$
for every $Y\in \fs(F)$ and $X\in \fu^V(F)$. By our choice of Haar measures, we have $\cF(\cF\varphi')(Y)=\varphi'(-Y)$ and $\cF(\cF\varphi^V)(X)=\varphi^V(-X)$ for every $\varphi'\in \cS(\fs'(F))$ and $\varphi^V\in \cS(\fu^V(F))$.

Finally, we fix once an for all a set $\cV$ of representatives of the isomorphism classes of Hermitian spaces of dimension $n$ and set $G=\bigsqcup_{V\in \cV} G^V$, $\fu=\bigsqcup_{V\in \cV} \fu^V$ so that
$$\displaystyle \cS(G(F))=\bigoplus_{V\in \cV} \cS(G^V(F)) \mbox{ and } \cS(\fu(F))=\bigoplus_{V\in \cV} \cS(\fu^V(F)).$$
We extend the maps $f^V\mapsto \widetilde{f^V_\natural}$ by linearity to
$$\displaystyle \cS(G(F))\to \cS(\fu(F))$$
$$\displaystyle f=(f^V)_V\mapsto \widetilde{f_\natural}=(\widetilde{f^V_\natural})_V.$$
We could also extend the Fourier transforms $\varphi^V\mapsto \cF\varphi^V$ by linearity to $\cS(\fu(F))$ but it turns out to be a better choice to use the following convention (see in particular Theorem \ref{theo transfer FT})
$$\displaystyle \cF\varphi=\left(\eta_{E/F}(\disc(V))^n\cF\varphi^V\right)_V$$
for every $\varphi=(\varphi^V)_V\in \cS(\fu(F))$ where $\disc(V)$ stands for the {\em discriminant} of the Hermitian space $V$ (that is the determinant of the matrix representing the Hermitian form, in any basis, seen as an element of $F^\times/N(E^\times)$).

\subsection{Relative orbital integrals and matching of functions}\label{Section relative orb integrals}

We keep the notation of the previous section. Let $V\in \cV$. For $f'\in \cS(G'(F))$ (resp. $f^V\in \cS(G^V(F))$) and $\gamma\in G'_{\rs}(F)$ (resp. $\delta\in G^V_{\rs}(F)$), we define a {\em relative orbital integral} $O_\eta(\gamma,f')$ (resp. $O(\delta,f^V)$) by
$$\displaystyle O_\eta(\gamma,f')=\int_{H_1(F)\times H_2(F)} f'(h_1\gamma h_2)\eta(h_2)dh_2dh_1 \left(\mbox{resp. } O(\delta,f^V)=\int_{H^V(F)\times H^V(F)} f^V(h_1\delta h_2)dh_1dh_2 \right)$$
Similarly, for $\varphi'\in \cS(\fs(F))$ (resp. $\varphi^V\in \cS(\fu^V(F))$) and $Y\in \fs_{\rs}(F)$ (resp. $X\in \fu^V_{\rs}(F)$) we define a relative orbital integral
$$\displaystyle O_\eta(Y,\varphi')=\int_{G_n(F)} \varphi'(hYh^{-1}) \eta_{E/F}(h)dh \left(\mbox{resp. } O(X,\varphi^V)=\int_{U(V)(F)} \varphi^V(hXh^{-1})dh \right)$$

We define {\em transfer factors}
$$\displaystyle \Omega:G'_{\rs}(F)\to \C^\times \mbox{ and } \omega:\fs_{\rs}(F)\to \C^\times$$
by
\begin{equation}\label{def omega}
\displaystyle \omega(Y)=\eta'\left( \det(e_{n+1},e_{n+1}Y,\ldots,e_{n+1}Y^n)\right),\;\;\; Y\in \fs_{\rs}(F)
\end{equation}
and
\begin{equation}\label{def Omega}
\displaystyle \Omega((g_n,g_{n+1}))=\eta'(g_n^{-1}g_{n+1})^{-n}\eta'\left( \det(e_{n+1},e_{n+1}s,\ldots,e_{n+1}s^n)\right),\;\;\; (g_n,g_{n+1})\in G'_{\rs}(F)
\end{equation}
where $e_{n+1}=(0,\ldots,0,1)$ and in the last equality we have used the notation $s=\nu(g_n^{-1}g_{n+1})$ for simplicity.

We say that two functions $f=(f^V)_V\in \cS(G(F))$ and $f'\in \cS(G'(F))$ {\em match} or that they are {\em transfer of each other} if
$$\displaystyle \Omega(\gamma)O_\eta(\gamma,f')=O(\delta,f^V)$$
for every $V\in \cV$ and every pair $(\gamma,\delta)\in G'_{\rs}(F)\times G^V_{\rs}(F)$ with matching orbits. Similarly, we say that two functions $\varphi=(\varphi^V)_V\in \cS(\fu(F))$ and $\varphi'\in \cS(\fs(F))$ {\em match} or that they are {\em transfer of each other} if
$$\displaystyle \omega(Y)O_\eta(Y,\varphi')=O(X,\varphi^V)$$
for every $V\in \cV$ and every pair $(Y,X)\in \fs_{\rs}(F)\times \fu^V_{\rs}(F)$ with matching orbits.

The following follows easily from the definitions and a painless computation:
\begin{num}
\item\label{eq 1 matching of functions} If $f'\in \cS(G'(F))$ and $f\in \cS(G(F))$ match then $\widetilde{f'_\natural}$ and $\widetilde{f_\natural}$ match.
\end{num}

Finally, we recall the following deep results from \cite{Zha1} and \cite{Xue}.
\begin{theo}[Zhang, Xue]\label{theo existence transfer}
\begin{enumerate}[(i)]
\item Assume that $F$ is a $p$-adic field. Then, for every $f\in \cS(G(F))$ (resp. $\varphi\in \cS(\fu(F))$) there exists $f'\in \cS(G'(F))$ (resp. $\varphi'\in \cS(\fs(F))$) such that $f$ and $f'$ match (resp. $\varphi$ and $\varphi'$ match). Conversely, for every $f'\in \cS(G'(F))$ (resp. $\varphi'\in \cS(\fs(F))$) there exists $f\in \cS(G(F))$ (resp. $\varphi\in \cS(\fu(F))$) such that $f$ and $f'$ match (resp. $\varphi$ and $\varphi'$ match).
\item Assume that $F$ is Archimedean. Then, there exists dense subspaces $\cS(G(F))_{\trans}\subset \cS(G(F))$, $\cS(G'(F))_{\trans}\subset \cS(G'(F))$, $\cS(\fu(F))_{\trans}\subset \cS(\fu(F))$ and $\cS(\fs(F))_{\trans}\subset \cS(\fs(F))$ satisfying the following: for every $f\in \cS(G(F))_{\trans}$ (resp. $\varphi\in \cS(\fu(F))_{\trans}$) there exists $f'\in \cS(G'(F))$ (resp. $\varphi'\in \cS(\fs(F))$) such that $f$ and $f'$ match (resp. $\varphi$ and $\varphi'$ match) and conversely, for every $f'\in \cS(G'(F))_{\trans}$ (resp. $\varphi'\in \cS(\fs(F))_{\trans}$) there exists $f\in \cS(G(F))$ (resp. $\varphi\in \cS(\fu(F))$) such that $f$ and $f'$ match (resp. $\varphi$ and $\varphi'$ match).
\end{enumerate}
\end{theo}

\begin{theo}[Zhang, Xue]\label{theo transfer FT}
Let $\varphi\in \cS(\fu(F))$ and $\varphi'\in \cS(\fs(F))$. Then, if $\varphi$ and $\varphi'$ match so do $\eta_{E/F}(-1)^{n(n+1)/2}\lambda_{E/F}(\psi')^{n(n+1)/2}\cF\varphi$ and $\cF\varphi'$ (where the Fourier transform $\cF\varphi$ of $\varphi$ is as defined in the end of Section \ref{Section Notation, matching of orbits}).
\end{theo}

\begin{rem}\label{rem transfer FT}
We remark that the constant $\eta_{E/F}(-1)^{n(n+1)/2}\lambda_{E/F}(\psi')^{n(n+1)/2}$ appearing in the above theorem is not exactly the same as the one we can extract from the computations of \cite[Sect. 4]{Zha1}. We refer the reader to \cite[Theorem 3.4.2.1]{Chau}  for a precise determination of this constant (which fits precisely the one given above) in the $p$-adic case. In the Archimedean case, the above theorem corresponds precisely what is stated in \cite[\S 9]{Xue}.
\end{rem}

\subsection{Relative characters}

Let $V\in \cV$ and $\pi\in \Temp(G^V)$. Then, for every $f\in \cS(G^V(F))$ we set
$$\displaystyle J_\pi(f)=\int_{H^V(F)} \Tr(\pi(h)\pi(f^\vee))dh=\int_{H^V(F)} f_{\pi}(h)dh$$
the integral being absolutely convergent by the following lemma.

\begin{lem}\label{lem 1 relative characters}
For every $f\in \cC^w(G^V(F))$ the integral
$$\displaystyle \int_{H^V(F)} f(h)dh$$
converges absolutely and defines a continuous linear form on $\cC^w(G^V(F))$.
\end{lem}

\noindent\ul{Proof}: This follows easily from \cite[Lemma 6.5.1(i)]{Beu1}. $\blacksquare$

Set $N'=R_{E/F}N_n\times R_{E/F} N_{n+1}$ (a maximal unipotent subgroup of $G'$) and $\psi_{N'}=\psi_n\boxtimes \psi_{n+1}$ (a generic character of $N'(F)$). We define continuous linear forms
$$\displaystyle \beta': \cC^w(N'(F)\backslash G'(F),\psi_{N'})\to \C$$
by 
$$\displaystyle \beta'(W)=\int_{N_2(F)\backslash P_2(F)} W(p)\eta(p)dp$$
and
$$\displaystyle \lambda(W)=\int_{N_1(F)\backslash H_1(F)} W(h_1)dh_1,\;\;\; W\in \cW(\Pi,\psi)$$
where $P_2=H_2\cap P'$, $N_2=H_2\cap N'$ and $N_1=H_1\cap N'$. Note that $\beta'=\beta_n\widehat{\otimes}\beta_{n+1}$ where $\beta_n$ and $\beta_{n+1}$ are defined in Section \ref{Section betan}. That $\lambda$ is absolutely convergent and continuous is a consequence of the following lemma.

\begin{lem}\label{lem 2 relative characters}
\begin{enumerate}[(i)]
\item For every $W\in \cC^w(N'(F)\backslash G'(F),\psi_{N'})$ the integral
$$\displaystyle \lambda(W)=\int_{N_1(F)\backslash H_1(F)} W(h_1)dh_1$$
is absolutely convergent and defines a continuous linear form on $\cC^w(N'(F)\backslash G'(F),\psi_{N'})$.

\item For every $\phi\in \cC^w(G'(F))$ the integral
$$\displaystyle \int_{H_1(F)} \phi(h_1)dh_1$$
is absolutely convergent and defines a continuous linear form on $\cC^w(G'(F))$.
\end{enumerate}
\end{lem}

\noindent\ul{Proof}:
\begin{enumerate}[(i)]
\item The proof is completely similar to the proof of Lemma \ref{lem continuity of beta} and left to the reader.

\item By definition of $\phi\in \cC^w(G'(F))$ it suffices to show that for every $d>0$ the integral
$$\displaystyle \int_{H_1(F)} \Xi^{G'}(h_1)\sigma(h_1)^ddh_1$$
converges. This is the content of \cite[4.1(3)]{Wald3} in the $p$-adic case but the proof works verbatim in the Archimedean case. $\blacksquare$
\end{enumerate}

Let $\Pi\in \Temp(G')$. We can write $\Pi=\Pi_n\boxtimes \Pi_{n+1}$ where $\Pi_k\in \Temp(G_k(E))$, $k=n,n+1$. Let $\cW(\Pi,\psi)=\cW(\Pi_n,\psi_n)\ctens \cW(\Pi_{n+1},\psi_{n+1})$ be the Whittaker model of $\Pi$ with respect to the character $\psi_{N'}$. We equip $\cW(\Pi,\psi)$ with the $G'(F)$-invariant scalar product (see Section \ref{Section Whittaker models})
$$\displaystyle (W,W')^{\Whitt}_{G'}=\int_{N'(F)\backslash P'(F)} W(p)\overline{W'(p)} dp,\;\;\; W,W'\in \cW(\Pi,\psi)$$
where $P'=R_{E/F} P_n\times R_{E/F} P_{n+1}$. Let $\cB(\Pi,\psi)$ be an orthonormal basis of (the Hilbert completion of) $\cW(\Pi,\psi)$ for this scalar product obtained by taking the union of orthonormal basis for $\cW(\Pi,\psi)[\delta]$ for every $\delta\in \widehat{K'}$ where $K'=K_{n,E}\times K_{n+1,E}$. Then, for every $f\in \cS(G'(F))$ we set
$$\displaystyle I_\Pi(f)=\sum_{W\in \cB(\Pi,\psi)} \overline{\lambda(W)}\beta'(\Pi(f^\vee)W).$$
We also define on $\cS(G'(F))$ the following semi-positive scalar product
$$\displaystyle (f_1,f_2)_{X_2,\Pi}=\sum_{W\in \cB(\Pi,\psi)} \beta'(\Pi(f_1^\vee)W) \overline{\beta'(\Pi(f_2^\vee)W)},\;\;\; f_1,f_2\in \cS(G'(F))$$
where as usual $f_k^\vee(g)=f_k(g^{-1})$ for $k=1,2$.

\begin{prop}\label{prop 1 relative characters}
The expressions defining $I_\Pi(f)$ and $(f_1,f_2)_{X_2,\Pi}$ are absolutely convergent and do not depend on the choice of the basis $\cB(\Pi,\psi)$. The functions
$$\displaystyle \Pi\in \Temp(G'(F))\mapsto I_\Pi(f) \mbox{ and } \Pi\in \Temp(G'(F))\mapsto (f_1,f_2)_{X_2,\Pi}$$
belong to $\cS(\Temp(G'(F)))$ and the linear (resp. sesquilinear) map
$$\displaystyle f\in \cS(G'(F))\mapsto (\Pi\mapsto I_\Pi(f))\in \cS(\Temp(G'))$$
$$\left(\mbox{resp. } (f_1,f_2)\in \cS(G'(F))^2\mapsto (\Pi\mapsto (f_1,f_2)_{X_2,\Pi})\in \cS(\Temp(G'))\right)$$
is continuous. Moreover, we have
\begin{align}\label{eq 0 relative characters}
\displaystyle \int_{H_1(F)} (L(h_1)f_1,f_2)_{X_2,\Pi}dh_1=\lvert \tau\rvert_E^{-n(n-1)/2}I_\Pi(f_1) \overline{I_\Pi(f_2)}
\end{align}
for every $f_1,f_2\in \cS(G'(F))$.
\end{prop}

\noindent\ul{Proof}: The proof of the first part of the proposition is completely similar to the proof of the first part of Lemma \ref{lem scalar products}. We prove the last part (i.e. identity \eqref{eq 0 relative characters}). Let $f_1,f_2\in \cS(G'(F))$. Then, we have (at least formally)
\[\begin{aligned}
\displaystyle & \int_{H_1(F)} (L(h_1)f_1,f_2)_{X_2,\Pi}dh_1 =\int_{H_1(F)} \sum_{W\in \cB(\Pi,\psi)} \beta'(\Pi(f_1^\vee)\Pi(h_1)W)\overline{\beta'(\Pi(f_2^\vee)W)}dh_1 \\
 & =\int_{H_1(F)} \sum_{W\in \cB(\Pi,\psi)}\sum_{W'\in \cB(\Pi,\psi)}(\Pi(h_1)W,W')^{\Whitt}_{G'}\beta'(\Pi(f_1^\vee)W') \overline{\beta'(\Pi(f_2^\vee)W)}dh_1
\end{aligned}\]
We claim
\begin{num}
\item\label{eq 1 relative characters} The above expression is absolutely convergent.
\end{num}
In the $p$-adic case, the sums in $W$ and $W'$ are actually finite and the result follows from \eqref{eq 1 representations} and Lemma \ref{lem 2 relative characters} (ii). In the Archimedean case, by \eqref{eq 2 Planch Whitt} the series
$$\displaystyle \sum_{W'\in \cB(\Pi,\psi)} \beta'(\Pi(f_1^\vee)W') W' \mbox{ and } \sum_{W\in \cB(\Pi,\psi)} \overline{\beta'(\Pi(f_2^\vee)W)} W$$
converge absolutely in $\cW(\Pi,\psi)$ whereas by \eqref{eq 1 representations}, the sesquilinear map
$$\displaystyle \cW(\Pi,\psi)\times \cW(\Pi,\psi)\to \cC^w(G'(F))$$
$$\displaystyle (W,W')\mapsto \left(g\mapsto (R(g)W,W')_{G'}^{\Whitt} \right)$$
is continuous. It follows that the series of functions
$$\displaystyle g\in G'(F)\mapsto \sum_{W\in \cB(\Pi,\psi)}\sum_{W'\in \cB(\Pi,\psi)}(\Pi(g)W,W')^{\Whitt}_{G'}\beta'(\Pi(f_1^\vee)W') \overline{\beta'(\Pi(f_2^\vee)W)}$$
converges absolutely in $\cC^w(G'(F))$ and the claim now follows from Lemma \ref{lem 2 relative characters} (ii).

By \eqref{eq 1 relative characters}, we can write
\begin{equation}\label{eq 2bis relative characters}
\displaystyle \int_{H_1(F)} (L(h_1)f_1,f_2)_{X_2,\Pi}dh_1=\sum_{W\in \cB(\Pi,\psi)}\sum_{W'\in \cB(\Pi,\psi)}\int_{H_1(F)}(\Pi(h_1)W,W')^{\Whitt}_{G'}dh_1\beta'(\Pi(f_1^\vee)W') \overline{\beta'(\Pi(f_2^*)W)}.
\end{equation}
We will prove the following:
\begin{num}
\item\label{eq 2 relative characters} For every $W,W'\in \cW(\Pi,\psi)$ we have 
$$\displaystyle \int_{H_1(F)}(\Pi(h_1)W,W')^{\Whitt}_{G'}dh_1=\lvert \tau\rvert_E^{-n(n-1)/2}\lambda(W)\overline{\lambda(W')}$$
\end{num}
First, we explain how this claim implies the proposition. Indeed, by \eqref{eq 2bis relative characters}, \eqref{eq 2 relative characters} gives
\[\begin{aligned}
\displaystyle \int_{H_1(F)} (L(h_1)f_1,f_2)_{X_2,\Pi}dh_1 & =\lvert \tau\rvert_E^{-n(n-1)/2}\sum_{W\in \cB(\Pi,\psi)}\sum_{W'\in \cB(\Pi,\psi)}\overline{\lambda(W')}\lambda(W)\beta'(\Pi(f_1^\vee)W') \overline{\beta'(\Pi(f_2^\vee)W)} \\
 & =\lvert \tau\rvert_E^{-n(n-1)/2} I_\Pi(f_1)\overline{I_\Pi(f_2)}
\end{aligned}\]
hence the result.

It only remains to prove \eqref{eq 2 relative characters}. Since, by \eqref{eq 1 representations} and Lemma \ref{lem 2 relative characters} (ii) again, the sesquilinear form
$$\displaystyle (W,W')\in \cW(\Pi,\psi)^2\mapsto \int_{H_1(F)}(\Pi(h_1)W,W')^{\Whitt}_{G'}dh_1$$
is continuous, we just need to show \eqref{eq 2 relative characters} when $W=W_n\otimes W_{n+1}$ and $W'=W'_n\otimes W'_{n+1}$ where $W_k,W'_k\in \cW(\Pi_k,\psi_k)$, $k=n,n+1$. Then, returning to the definitions, we have
$$\displaystyle \int_{H_1(F)} (\Pi(h_1)W,W')^{\Whitt}_{G'}dh_1=\int_{G_n(E)}\int_{N_n(E)\backslash G_n(E)} (R(h)W_n,W'_n)^{\Whitt} W_{n+1}(gh)\overline{W'_{n+1}(g)}dgdh$$
where $(.,.)^{\Whitt}$ is the scalar products on $\cC^w(N_n(E)\backslash G_n(E),\psi_n)$ defined in Section \ref{Section Whittaker models}. Then, by formal manipulations we get
\[\begin{aligned}
\displaystyle & \int_{H_1(F)} (\Pi(h_1)W,W')^{\Whitt}_{G'}dh_1=\int_{N_n(E)\backslash G_n(E)}\int_{G_n(E)} (R(h)W_n,W'_n)^{\Whitt} W_{n+1}(gh)\overline{W'_{n+1}(g)}dhdg \\
 & =\int_{N_n(E)\backslash G_n(E)}\int_{G_n(E)} (R(h)W_n,R(g)W'_n)^{\Whitt} W_{n+1}(h)\overline{W'_{n+1}(g)}dhdg \\
 & =\int_{N_n(E)\backslash G_n(E)}\int_{N_n(E)\backslash G_n(E)}\int_{N_n(E)}^* (R(uh)W_n,R(g)W'_n)^{\Whitt}\psi_n(u)^{-1}du W_{n+1}(h)\overline{W'_{n+1}(g)}dhdg \\
 & =\lvert \tau\rvert_E^{-n(n-1)/2}\int_{N_n(E)\backslash G_n(E)}\int_{N_n(E)\backslash G_n(E)} W_n(h) \overline{W'_n(g)} W_{n+1}(h)\overline{W'_{n+1}(g)}dhdg=\lvert \tau\rvert_E^{-n(n-1)/2}\lambda(W')\overline{\lambda(W)}
\end{aligned}\]
where in the fourth equality we have used Proposition \ref{prop 1 Planch Whitt}. If these formal manipulations were justified this would prove \eqref{eq 2 relative characters}. Unfortunately, the above expression is certainly not absolutely convergent in general. However, by \eqref{eq 1 representations} the function
$$\phi: g\in G_n(E)\mapsto (R(g)W_n,W'_n)^{\Whitt}$$
belongs to $\cC^w(G_n(E))$ and we recall that to every function $\phi\in \cC^w(G_n(E))$ we have associated in Section \ref{section Plancherel Whitt} a function
$$W_\phi\in \cC^w(N_n(E)\backslash G_n(E)\times N_n(E)\backslash G_n(E),\psi_n^{-1}\boxtimes \psi_n).$$
Then that the result of the above formal manipulations is indeed correct follows from:
\begin{num}
\item\label{eq 3 relative characters} For every $\phi\in \cC^w(G_n(E))$ we have
$$\displaystyle \int_{G_n(E)} \phi(h) (R(h)W_{n+1},W'_{n+1})^{\Whitt}dh=\int_{\left(N_n(E)\backslash G_n(E)\right)^2} W_\phi(g,h) W_{n+1}(h)\overline{W'_{n+1}(g)}dhdg$$
\end{num}
By Lemma \ref{eq 1 representations}, Lemma \ref{lem 1 HC Planch} and Lemma \ref{lem 2 relative characters} both sides of \eqref{eq 3 relative characters} are absolutely convergent and define continuous linear forms on $\cC^w(G_n(E))$. Therefore, it suffices to check the equality for $\phi\in \cS(G_n(E))$ where the same formal manipulations as before are now justified due to the absolute convergence of the relevant expressions. This shows \eqref{eq 3 relative characters} and ends the proof of the proposition. $\blacksquare$

\subsection{Statement of the main theorems}

Let $V\in \cV$. Recall that a representation $\pi\in \Temp(G^V)$ is said to be {\em $H^V$-distinguished} if there exists a non-zero (continuous) $H^V(F)$-invariant linear form on (the space of) $\pi$. By \cite[Theorem 7.2.1]{Beu1}, $\pi$ is $H^V$-distinguished if and only if the relative character $J_\pi$ is not identically zero. We denote by $\Temp_{H^V}(G^V)$ the subset of irreducible $H^V$-distinguished tempered representations of $G^V(F)$. By \cite[Corollary 7.6.1]{Beu1}, $\Temp_{H^V}(G^V)$ is a union of connected components of $\Temp(G^V)$.

The first main theorem of this chapter is the following one which, as explained in the introduction, has direct applications to the global Ichino-Ikeda conjecture for unitary groups.

\begin{theo}\label{theo1 main result}
Let $f=(f^V)_V\in \cS(G(F))$ and $f'\in \cS(G'(F))$ be matching functions. Then, for every $V\in \cV$ and every $\pi\in \Temp_{H^V}(G^V)$, we have
$$\displaystyle \kappa_V J_\pi(f^V)=I_{BC(\pi)}(f')$$
where
$$\displaystyle \kappa_V=\left(\eta'((-1)^{n+1}\tau)\lambda_{E/F}(\psi') \right)^{n(n+1)/2}\lvert \tau\rvert_E^{n(n-1)/4} \eta_{E/F}(\disc(V))^n.$$
\end{theo}

\begin{rem}
The constant that appears in the theorem above differs slightly from \cite[Conjecture 4.4]{Zh3}. For this, we offer the following explanation. First our normalization of the relative characters $J_\pi$ and $I_\Pi$ is not the same as in {\it loc. cit.} since we have replaced the representations $\pi$ and $\Pi$ by their contragredient and we have used the Whittaker model of $\Pi$ with respect to $\psi_n$ rather than $\psi'_E$. This last point explain the discrepancy for the powers of $\lvert \tau\rvert_E$. The difference for the exponents of $\eta_{E/F}(-1)$ and $\eta_{E/F}(\disc(V))$ seems for its part to originate from the precise computation of the constant up to which ``Fourier transform and transfer commute'' (see Remark \ref{rem transfer FT}).
\end{rem}

To state the second main result of this chapter, we introduce the following notation: for $G$ a connected reductive group over $F$, we let $\pi\in \Temp(G)\mapsto \mu^*_{G}(\pi)$ be the function (or {\em density}) such that $d\mu_{G}(\pi)=\mu^*_{G}(\pi)d\pi$. Note that by definition of $d\pi$, $\mu^*_{G}(\pi)$ differs from $\mu_{G}(\pi)$ (defined in Section \ref{Section Planch}) by an integral power of $\gamma^*(0,\mathbf{1}_F,\psi')$.

\begin{theo}\label{theo2 main result}
We have
$$\displaystyle \mu^*_{G^V}(\pi)=\frac{\lvert \gamma^*(0,\pi,\Ad,\psi')\rvert}{\lvert S_\pi\rvert}$$
for almost all $\pi\in \Temp_{H^V}(G^V)$.
\end{theo}

We note the following interesting corollary which is a particular case of a general conjecture of Hiraga-Ichino-Ikeda (\cite{HII}).

\begin{cor}\label{cor1 main result}
For every Hermitian space $W$ over $E$, we have
$$\displaystyle \mu^*_{U(W)}(\sigma)=\frac{\lvert \gamma^*(0,\sigma,\Ad,\psi')\rvert}{\lvert S_\sigma\rvert}$$
for almost all $\sigma\in \Temp(U(W))$. In particular, for every $\sigma\in \Pi_2(U(W))$ we have the following formula for its formal degree:
$$\displaystyle d(\sigma)=\frac{\lvert \gamma(0,\sigma,\Ad,\psi')\rvert}{\lvert S_\sigma\rvert}.$$
\end{cor}

\noindent\ul{Proof}: We proceed by induction on $\dim(W)$, the case where $\dim(W)=0$ being evident. Assume now that $\dim(W)\geqslant 1$ and set $n=\dim(W)-1$. Up to scaling the Hermitian form on $W$, we may assume that $W=V'$ for some $V\in\cV$. Let $\cO\subset \Temp(U(W))=\Temp(U(V'))$ be a connected component. Then, we claim:
\begin{num}
\item\label{eq 1 main result} there exists a connected component $\cO_0\subset \Temp(U(V))$ such that $\cO_0\boxtimes \cO \subset \Temp_{H^V}(G^V)$.
\end{num}
Indeed, since $\Temp_{H^V}(G^V)$ is an union of connected components of $\Temp(G^V)$ it suffices to find $\cO_0$ such that $\cO_0\boxtimes \cO \cap \Temp_{H^V}(G^V)\neq \emptyset$. Let $\sigma\in \cO$ and $f_\sigma$ be a matrix coefficient of $\sigma$ whose restriction to $U(V)(F)$ is nonzero. By \cite[(6.5.1)]{Beu1}, this restriction is in the Harish-Chandra Schwartz space of $U(V)(F)$ and by the Harish-Chandra Plancherel formula (which still holds for this kind of functions), we see that there exists $\sigma_0\in \Temp(U(V))$ such that $\sigma_0(f_\sigma)\neq 0$. This readily implies that $J_{\sigma_0\boxtimes \sigma}\neq 0$ so that $\sigma_0\boxtimes \sigma\in \Temp_{H^V}(G^V)$. Therefore, the connected component of $\sigma_0$ has the desired property.

Let $\cO_0$ be as in \eqref{eq 1 main result}. Then, by Theorem \ref{theo2 main result} we have
$$\displaystyle \mu^*_{U(V)}(\sigma_0)\mu^*_{U(V')}(\sigma)=\mu^*_{G^V}(\sigma_0\boxtimes \sigma)=\frac{\lvert \gamma^*(0,\sigma_0,\Ad,\psi')\rvert}{\lvert S_{\sigma_0}\rvert}\frac{\lvert \gamma^*(0,\sigma,\Ad,\psi')\rvert}{\lvert S_\sigma\rvert}$$
for almost all $(\sigma_0,\sigma)\in \cO_0\times \cO$. By the induction hypothesis, we also have $\mu^*_{U(V)}(\sigma_0)=\frac{\lvert \gamma^*(0,\sigma_0,\Ad,\psi')\rvert}{\lvert S_{\sigma_0}\rvert}$ for almost all $\sigma_0\in \cO_0$ and moreover this term is almost everywhere nonzero. Therefore
$$\displaystyle \mu^*_{U(V')}(\sigma)=\frac{\lvert \gamma^*(0,\sigma,\Ad,\psi')\rvert}{\lvert S_\sigma\rvert}$$
for almost all $\sigma\in \cO$ and since $\cO$ was arbitrary the same holds for almost all $\sigma\in \Temp(U(W))$. $\blacksquare$

\subsection{Local Jacquet-Rallis trace formulas}

\subsubsection{The unitary case}\label{unitary JR}

Let $V\in \cV$ and let $f_1,f_2\in \cS(G^V(F))$. Consider the following expression
$$\displaystyle J(f_1,f_2)=\int_{H^V(F)}\int_{H^V(F)}\int_{G^V(F)} f_1(h_1gh_2)\overline{f_2(g)}dgdh_2dh_1$$
which is absolutely convergent by \cite[Lemma A.4]{Zha1}. By definition of the measure on $H^V(F)\backslash G^V_{rs}(F)/H^V(F)$ and since the complement of $G^V_{rs}(F)$ in $G^V(F)$ is of measure $0$ we have
\begin{align}\label{eq 1 local JR unitary}
\displaystyle J(f_1,f_2)=\int_{H^V(F)\backslash G^V_{rs}(F)/H^V(F)} O(\delta,f_1)\overline{O(\delta,f_2)} d\delta.
\end{align}
On the other hand by \cite[Lemma 7.2.2 (v)]{Beu1}, we also have
\begin{align}\label{eq 2 local JR unitary}
\displaystyle J(f_1,f_2)=\int_{\Temp_{\ind}(G^V)}J_\pi(f_1)\overline{J_\pi(f_2)} d\mu_{G^V}(\pi).
\end{align}
where we have extended the definition of $J_\pi$ to $\pi\in \Temp_{\ind}(G^V)$ by linearity. Notice that the above integral is absolutely convergent by \eqref{eq 3 Planch} and Lemma \ref{lem 1 relative characters}. Combining \eqref{eq 1 local JR unitary} and \eqref{eq 2 local JR unitary} we get
\begin{align}\label{eq 3 local JR unitary}
\int_{H^V(F)\backslash G^V_{rs}(F)/H^V(F)} O(\delta,f_1)\overline{O(\delta,f_2)} d\delta=\int_{\Temp_{\ind}(G^V)}J_\pi(f_1)\overline{J_\pi(f_2)} d\mu_{G^V}(\pi).
\end{align}

\subsubsection{The linear case}\label{linear JR}

Let $f_1,f_2\in \cS(G'(F))$. Consider the following expression
$$\displaystyle I(f_1,f_2)=\int_{H_1(F)}\int_{H_2(F)}\int_{G'(F)} f_1(h_1gh_2)\overline{f_2(g)}dg\eta(h_2)dh_2dh_1.$$
We claim
\begin{num}
\item\label{eq 1 local JR linear} The above expression is absolutely convergent and defines a continuous sesquilinear form on $\cS(G'(F))$.
\end{num}

\noindent\ul{Proof}: For every $d>0$ the above integral is, up to continuous norms in $f_1$ and $f_2$, bounded by
$$\displaystyle \int_{H_1(F)}\int_{H_2(F)}\int_{G'(F)} \Xi^{G'}(h_1gh_2)\Xi^{G'}(g)\sigma(h_1gh_2)^{-d}\sigma(g)^{-2d}dgdh_2dh_1$$
hence by
$$\displaystyle \int_{H_1(F)}\int_{H_2(F)}\int_{G'(F)} \Xi^{G'}(h_1gh_2)\Xi^{G'}(g)\sigma(g)^{-d}dg\sigma(h_2)^{-d}dh_2\sigma(h_1)^ddh_1.$$
Assuming (as we may) that the log-norm $\sigma$ is $K'$-bi-invariant this last expression equals
$$\displaystyle \int_{H_1(F)}\int_{H_2(F)}\int_{G'(F)} \int_{K'\times K'}\Xi^{G'}(h_1k_1gk_2h_2)dk_2dk_1\Xi^{G'}(g)\sigma(g)^{-d}dg\sigma(h_2)^{-d}dh_2\sigma(h_1)^ddh_1$$
and by the ``doubling principle'' (\cite[Lemme II.1.3]{Wald1}, \cite[Proposition 16(iii) p.329]{Var}) we have
$$\displaystyle \int_{K'\times K'}\Xi^{G'}(h_1k_1gk_2h_2)dk_2dk_1=\Xi^{G'}(h_1)\Xi^{G'}(g)\Xi^{G'}(h_2).$$
To conclude, it suffices to remark that for every $d>0$ the integral
$$\displaystyle \int_{H_1(F)}\Xi^{G'}(h_1)\sigma(h_1)^ddh_1$$
is convergent by \cite[4.1(3)]{Wald3} whereas for $d>0$ sufficiently large the two integrals
$$\displaystyle \int_{G'(F)} \Xi^{G'}(g)^2 \sigma(g)^{-d}dg \mbox{ and } \int_{H_2(F)} \Xi^{G'}(h_2)\sigma(h_2)^{-d}dh_2$$
are convergent by \cite[Lemme II.1.5]{Wald1} and \cite[Proposition 31 p.340]{Var} noting that $\Xi^{G'}_{\mid H_2}\ll (\Xi^{H_2})^2\sigma^{d'}$ for some $d'>0$. $\blacksquare$

By \eqref{eq 1 local JR linear}, the definition of the measure on $H_1(F)\backslash G'_{rs}(F)/H_2(F)$ and the fact that the complement of $G'_{rs}(F)$ in $G'(F)$ is of measure $0$ we have
\begin{align}\label{eq 2 local JR linear}
\displaystyle I(f_1,f_2)=\int_{H_1(F)\backslash G'_{rs}(F)/H_2(F)} O_\eta(\gamma,f_1)\overline{O_\eta(\gamma,f_2)} d\gamma.
\end{align}
In order to get a ``spectral'' expression for $I(f_1,f_2)$ we will have to use Theorem \ref{theo Plancherel} in a slightly disguised form. There exists a unique $V_0\in \cV$ such that $G^{V_0}$ is quasi-split and we set $G_{\qs}=G^{V_0}$. Let
$$\displaystyle BC: \Temp(G_{\qs})/\stab\to \Temp(G')$$
be the ``tensor product'' of the two base-change maps $\Temp(U(V_0))/\stab\to \Temp(G_n(E))$ and $\Temp(U(V_0'))/\stab\to \Temp(G_{n+1}(E))$.

\begin{prop}\label{prop 1 local JR linear}
For every $f_1,f_2\in \cS(G'(F))$ we have
$$\displaystyle \int_{H_2(F)}\int_{G'(F)} f_1(gh_2)\overline{f_2(g)}dg\eta(h_2)dh_2=\int_{\Temp(G_{\qs})/\stab} (f_1,f_2)_{X_2,BC(\pi)} \frac{\lvert \gamma^*(0,\pi,\Ad,\psi')\rvert}{\lvert S_\pi\rvert}d\pi$$
where the right-hand side is absolutely convergent.
\end{prop}

\noindent\ul{Proof}: By Proposition \ref{prop 1 relative characters} together with Lemma \ref{lem 1 L parameters, LLC, basechange} and \eqref{basic estimates spectral measure} we see that the right-hand side is absolutely convergent and defines a continuous sesquilinear form on $\cS(G'(F))$. Obviously, so does the left-hand side and therefore it suffices to establish the proposition when $f_1,f_2$ belong to the dense subspace $\cS(G_n(E))\otimes \cS(G_{n+1}(E))$. Thus, we assume that $f_k=f_{k,n}\otimes f_{k,n+1}$ where $f_{k,n}\in \cS(G_n(E))$ and $f_{k,n+1}\in \cS(G_{n+1}(E))$ for $k=1,2$. Set $f'_{k,l}(g)=f_{k,l}(g^{-1})\eta'_l(g^{-1})$ for $k=1,2$ and $l=n,n+1$ where $\eta'_l$ is the character of $G_l(E)$ defined in Section \ref{Section groups} which extends $\eta_l$. Then, the left-hand side of the proposition can be rewritten as
\begin{align}\label{eq 3 local JR linear}
\displaystyle \prod_{l=n,n+1} \int_{Y_l}\varphi_{1,l}'(x)\overline{\varphi_{2,l}'(x)}dx
\end{align}
where we have set $\displaystyle \varphi_{k,l}'(x)=\int_{G_l(F)} f_{k,l}'(hx)dh$  and $Y_l=G_l(F)\backslash G_l(E)$ for all $k\in \{1,2\}$, $l\in \{n,n+1 \}$.

On the other hand, for $\Pi=\Pi_n\boxtimes \Pi_{n+1}\in \Temp(G')$, if we denote by $\cB(\Pi_n,\psi_n)$ and $\cB(\Pi_{n+1},\psi_{n+1})$ orthonormal basis of $\cW(\Pi_n,\psi_n)$ and $\cW(\Pi_{n+1},\psi_{n+1})$ as in Section \ref{section Plancherel Whitt}, we have (by definition and Proposition \ref{prop 1 relative characters})
\[\begin{aligned}
\displaystyle (f_1,f_2)_{X_2,\Pi}=\prod_{l=n,n+1}\sum_{W\in \cB(\Pi_l,\psi_l)} \beta_l(R(f_{1,l}^\vee)W) \overline{\beta_l(R(f_{2,l}^\vee)W)}
\end{aligned}\]
Since $W\in \cW(\Pi_l,\psi_l)\mapsto W':=\eta'_l W\in \cW(\Pi_l\otimes \eta'_l,\psi_l)$ is an isomorphism preserving the scalar products and $\beta_l(R(f_{k,l}^\vee)W)=\beta(R(f_{k,l}')W')$ for $k=1,2$ and $l=n,n+1$, we see that the above expression equals
$$\displaystyle \prod_{l=n,n+1}(f_{1,l}',f'_{2,l})_{Y_l,\Pi_l\otimes \eta_l'}$$
where $(.,.)_{Y_l,\Pi_l\otimes \eta_l'}$, $l=n,n+1$, are the semi-positive scalar products defined in Section \ref{Section the Theorem}. Moreover, for $\pi=\pi_n\boxtimes \pi_{n+1}\in \Temp(G_{\qs})$, we have
$$\displaystyle S_\pi\simeq S_{\pi_n}\times S_{\pi_{n+1}}$$
$$\displaystyle \gamma^*(0,\pi,\Ad,\psi')=\gamma^*(0,\pi_n,\Ad,\psi')\gamma^*(0,\pi_{n+1},\Ad,\psi')$$
All in all, we conclude that the right-hand side of the proposition equals
\begin{align}\label{eq 4 local JR linear}
\displaystyle \prod_{l=n,n+1}\int_{\Temp(U(l))/\stab} (\varphi_{1,l}',\varphi_{2,l}')_{Y_l,BC_l(\pi)} \frac{\lvert \gamma^*(0,\pi,\Ad,\psi')\rvert}{\lvert S_\pi\rvert}d\pi
\end{align}
where we recall that $BC_l(\pi)=BC(\pi)\otimes \eta'_l$ for $l=n,n+1$. The equality of \eqref{eq 3 local JR linear} and \eqref{eq 4 local JR linear} now follows directly from Theorem \ref{theo Plancherel}. $\blacksquare$

By Proposition \ref{prop 1 local JR linear}, we have
\begin{align}\label{eq 5 local JR linear}
\displaystyle I(f_1,f_2)=\int_{H_1(F)}\int_{\Temp(G_{\qs})/\stab} (L(h_1)f_1,f_2)_{X_2,BC(\pi)} \frac{\lvert \gamma^*(0,\pi,\Ad,\psi')\rvert}{\lvert S_\pi\rvert}d\pi dh_1.
\end{align}
We claim
\begin{num}
\item\label{eq 6 local JR linear} The above expression is absolutely convergent.
\end{num}
Indeed, by \cite[Theorem 2]{CHH} in the $p$-adic case, \cite[Theorem 1.2]{Sun} in the Archimedean case, and since for every $\Pi\in \Temp(G')$ the function $g\in G'(F)\mapsto (L(g)f_1,f_2)_{X_2,\Pi}$ is a matrix coefficient of $\Pi^\vee$, there exist $f_1',f_2'\in \cS(G'(F))$ such that
$$\displaystyle \lvert (L(g)f_1,f_2)_{X_2,\Pi}\rvert\leqslant \lVert f_1'\rVert_{X_2,\Pi} \lVert f_2'\rVert_{X_2,\Pi} \Xi^{G'}(g)$$
for every $\Pi\in \Temp(G')$ and $g\in G'(F)$ where we have set $\lVert f_k'\rVert_{X_2,\Pi}=(f_k',f_k')_{X_2,\Pi}^{1/2}$ for $k=1,2$. By Proposition \ref{prop 1 local JR linear} and Cauchy-Schwartz inequality we see that the integral
$$\displaystyle \int_{\Temp(G_{\qs})/\stab} \lVert f_1'\rVert_{X_2,BC(\pi)} \lVert f_2'\rVert_{X_2,BC(\pi)} \frac{\lvert \gamma^*(0,\pi,\Ad,\psi')\rvert}{\lvert S_\pi\rvert}d\pi$$
is convergent. Therefore, the expression \eqref{eq 5 local JR linear} where we replace the integrand by its absolute value is bounded up to a constant by
$$\displaystyle \int_{H_1(F)}\Xi^{G'}(h_1)dh_1.$$
This last expression is convergent by \cite[4.1(3)]{Wald3} therefore proving \eqref{eq 6 local JR linear}.

By \eqref{eq 6 local JR linear} and Proposition \ref{prop 1 relative characters}, we obtain
$$\displaystyle I(f_1,f_2)=\lvert \tau\rvert_E^{-n(n-1)/2}\int_{\Temp(G_{\qs})/\stab} I_{BC(\pi)}(f_1)\overline{I_{BC(\pi)}(f_2)} \frac{\lvert \gamma^*(0,\pi,\Ad,\psi')\rvert}{\lvert S_\pi\rvert}d\pi$$
Combined with \eqref{eq 2 local JR linear}, this gives the identity
\begin{align}\label{eq 7 local JR linear}
\displaystyle & \int_{H_1(F)\backslash G'_{rs}(F)/H_2(F)} O_\eta(\gamma,f_1)\overline{O_\eta(\gamma,f_2)} d\gamma= \\
\nonumber & \lvert \tau\rvert_E^{-n(n-1)/2}\int_{\Temp(G_{\qs})/\stab} I_{BC(\pi)}(f_1)\overline{I_{BC(\pi)}(f_2)} \frac{\lvert \gamma^*(0,\pi,\Ad,\psi')\rvert}{\lvert S_\pi\rvert}d\pi
\end{align}

\subsection{Spectral expansions of certain unipotent relative orbital integrals}

\subsubsection{The unitary case}\label{unitary unipotent}
Let $V\in \cV$. For each $f\in \cS(G^V(F))$ we set
$$\displaystyle O(1,f)=\int_{H^V(F)}f(h)dh$$
By \eqref{eq 3 Planch} and Lemma \ref{lem 1 relative characters}, we have
\begin{align}\label{eq 1 unipotent unitary}
\displaystyle O(1,f)=\int_{\Temp_{\ind}(G^V)}J_\pi(f)d\mu_{G^V}(\pi)
\end{align}
the integral being absolutely convergent. On the other hand, by definition of $\widetilde{f_\natural}$ (see \eqref{def linearization}), we have
$$\displaystyle O(1,f)=\widetilde{f_\natural}(0)$$
By Fourier inversion and the choice of the measure on $\fu_{\rs}^V(F)/U(V)(F)$ this can be rewritten as
\begin{align}\label{eq 2 unipotent unitary}
\displaystyle O(1,f)=\int_{\fu^V(F)} \cF\widetilde{f_\natural}(X)dX=\int_{\fu_{\rs}^V(F)/U(V)(F)} O(X,\cF\widetilde{f_\natural}) dX
\end{align}

\subsubsection{The linear case Lie algebra version}

In this section and the next, we introduce and study certain (relative) nilpotent/unipotent orbital integrals on the linear side. In contrast to the unitary case, these distributions are not supported on the smallest orbit/coset (that of $0\in \mathfrak{s}(F)$ or $1\in G'(F)$)\footnote{Note that in the linear case we need to consider distributions that are $H_1(F)\times (H_2(F),\eta)$-equivariant, and there is no such distribution supported on the trivial coset (because $\eta$ is nontrivial on $H_1(F)\cap H_2(F)$).} but rather on certain regular \footnote{where here ``regular'' means ``with trivial stabilizer''.} nilpotent/unipotent orbits. We will see that, in some sense, these relative orbital integrals correspond to the trivial one on the unitary side via the Jacquet-Rallis transfer (see equation \eqref{eq 1ter two comparisons}). This however, can only be established a posteriori after the careful analysis of this section. On the other hand, an a priori reason for considering these nilpotent/unipotent relative orbital integrals is their homogeneity property: after linearization, these are the only orbital integrals which, like the trivial orbital integral on the unitary side, are homogeneous of degree $0$.

Set
$$\displaystyle \xi_{+}=(-1)^{n}\begin{pmatrix} 0 & \tau^{-1} & 0 & \ldots & 0 \\ \vdots & \ddots & \ddots & \ddots & \vdots \\ \vdots & & \ddots & \ddots & 0 \\ \vdots & & & \ddots & \tau^{-1} \\ 0 & \ldots & \ldots & \ldots & 0
\end{pmatrix}\in \fs(F),\;\; \xi_{-}=\tau^2{}^t\xi_{+}=(-1)^{n}\begin{pmatrix} 0 & \ldots & \ldots & \ldots & 0 \\ \tau & \ddots & & & \vdots \\ 0 & \ddots & \ddots & & \vdots \\ \vdots & \ddots & \ddots & \ddots & \vdots \\ 0 & \ldots & 0 & \tau & 0
\end{pmatrix}\in \fs(F)$$
Let $\varphi\in \cS(\fs(F))$. For $s\in \C$ we define
$$\displaystyle O_s(\xi_{+},\varphi)=\int_{G_n(F)} \varphi(h\xi_{+}h^{-1}) \lvert \det h\rvert^s \eta_{E/F}(h)dh$$
whenever the integral is convergent. Note that the formula defining $\omega(Y)$ (see \eqref{def omega}) for $Y\in \fs_{\rs}(F)$ still makes sense for $\xi_{-}$ and that we have
\begin{align}\label{eq 0 unipotent linear}
\displaystyle \omega(\xi_{-})=\eta'((-1)^{n-1}\tau)^{n(n+1)/2}
\end{align}
The goal of this section is to show the following:

\begin{prop}\label{prop 1 unipotent linear}
For $\Re(s)>1-1/n$, $O_s(\xi_{+},\varphi)$ is defined by an absolutely convergent expression. Moreover, the function $s\mapsto O_s(\xi_{+},\varphi)$ admits a meromorphic continuation to $\C$ with no pole at $s=0$ and setting $O(\xi_{+},\varphi)=O_0(\xi_{+},\varphi)$, we have
$$\displaystyle \gamma\omega(\xi_{-})O(\xi_{+},\varphi)=\int_{\fs(F)} \omega(Y)(\cF\varphi)(Y) dY=\int_{\fs_{\rs}(F)/G_n(F)}\omega(Y)O_\eta(Y,\cF \varphi) dY$$
where
$$\displaystyle \gamma=\prod_{k=1}^n \gamma(1-k,\eta_{E/F}^k,\psi')$$
\end{prop}

\begin{rem}\label{rem 1 unipotent linear}
The factor $\gamma$ is non-zero: this boils down to the fact that $L(s,\eta_{E/F}^k)$ has no pole at $s=1-k$ for every $1\leqslant k\leqslant n$. This last fact is easy to check in the $p$-adic case whereas in the Archimedean case it follows from the fact that $L(s,\eta_{E/F}^k)$ has the same poles as $\Gamma(\frac{s+1}{2})$ if $k$ is odd, $\Gamma(\frac{s}{2})$ if $k$ is even.
\end{rem}

For the proof of Proposition \ref{prop 1 unipotent linear} we need some preparation. Let $N_S$ be the image of $R_{E/F}N_{n+1,E}$ by $\nu$ and $\n_S=\fs\cap R_{E/F}\n_{n+1,E}$ be its tangent space at the origin. Then, $\nu$ induces isomorphisms $N_S(F)\simeq N_{n+1}(E)/N_{n+1}(F)$ and $\n_S(F)\simeq \n_{n+1}(E)/\n_{n+1}(F)$ and we equip $N_S(F)$, $\n_S(F)$ with the quotient measures. Set
$$\displaystyle I^1(\varphi,s)=\lvert \tau\rvert_E^{\dim(N_S)/2}\int_{G_n(F)/N_n(F)} \int_{\n_S(F)}\varphi(hYh^{-1}) \psi'(\langle \xi_{-},Y\rangle)dY \lvert \det h\rvert^s \eta_{E/F}(h)dh$$
for every $s\in \C$ for which this expression makes sense. First we show:

\begin{lem}\label{lem 1 unipotent linear}
The expression defining $O_s(\xi_{+},\varphi)$ converges absolutely for $\Re(s)>1-1/n$ and the expression defining $I^1(\varphi,s)$ converges as an iterated integral for $\Re(s)<1$. Moreover, $O_s(\xi_{+},\varphi)$ and $I^1(\varphi,s)$ are holomorphic functions in their region of convergence and in the range $1-1/n<\Re(s)<1$ we have
$$\displaystyle I^1(\varphi,s)=\prod_{k=1}^n \gamma(ks-(k-1),\eta_{E/F}^k,\psi') O_s(\xi_{+},\varphi)$$
\end{lem}

Recall that $N'_{n+1}$ denotes the derived subgroup of $N_{n+1}$. Let $N'_S$ be the image of $R_{E/F}N'_{n+1,E}$ by $\nu$ and $\n_S'$ be its tangent space at the origin. As before, $\nu$ induces an isomorphism $\n_S'(F)\simeq \n_{n+1}'(E)/\n_{n+1}'(F)$ and we equip this space with the quotient of the Haar measures we fixed on $\n_{n+1}'(E)$ and $\n_{n+1}'(F)$ (see Section \ref{Section measures}). As an intermediate step for the proof of Lemma \ref{lem 1 unipotent linear} we need the following lemma:

\begin{lem}\label{lem 2 unipotent linear}
For every $f\in \cS(\n_S(F))$, we have
$$\displaystyle \int_{N_n(F)} f(h\xi_{+}h^{-1})dh=\lvert \tau\rvert_E^{\dim(N'_S)/2}\int_{\n'_S(F)}f(\xi_{+}+Y)dY$$
\end{lem}

\noindent\ul{Proof}: The isomorphism $\n_{n+1}(F)\simeq \n_S(F)$, $Y\mapsto (-1)^{n}\tau^{-1}Y$, sends $\xi'_+=(-1)^{n}\tau \xi_{+}$ to $\xi_{+}$ and the Haar measure on $\n_{n+1}'(F)$ to $\lvert \tau\rvert_E^{\dim(N'_S)/2}$ times the Haar measure on $\n_S'(F)$. Therefore, it suffices to show that
$$\displaystyle N_n(F)\to \xi'_++\n_{n+1}'(F)$$
$$\displaystyle h\mapsto h\xi'_+h^{-1}$$
is an isomorphism preserving measures. Given the definition of the Haar measures on $N_n(F)$ and $\n_{n+1}'(F)$, it even suffices to show that
$$\displaystyle \iota_n: N_n\to \xi'_++\n_{n+1}'$$
$$\displaystyle h\mapsto h\xi'_+h^{-1}$$
is an isomorphism over $\bZ$. This last statement is easy to show by induction on $n$, noting that
$$\displaystyle \iota_n(vu)=\begin{pmatrix} v\iota_{n-1}(u)v^{-1} & \begin{matrix} v_{1,n} \\ \vdots \\ v_{n-1,n} \\ 1 \end{matrix} \\ \begin{matrix} 0 & \ldots & 0 \end{matrix} & 0\end{pmatrix}$$
for every $(u,v)\in N_{n-1}\times U_n$. $\blacksquare$

\noindent\ul{Proof}(of Lemma \ref{lem 1 unipotent linear}): By Lemma \ref{lem 2 unipotent linear} and the Iwasawa decomposition $G_n(F)=K_nA_n(F)N_n(F)$, provided everything is absolutely convergent, we have
\begin{align}\label{eq 1 unipotent linear}
\displaystyle & \lvert \tau\rvert_E^{-\dim(N'_S)/2}O_s(\xi_{+},\varphi) =\int_{G_n(F)/N_n(F)} \int_{\n'_{S}(F)} \varphi(h(\xi_++Y)h^{-1})dY \lvert \det h\rvert^s \eta_{E/F}(h) dh \\
\nonumber & =\int_{K_n}\int_{A_n(F)} \int_{\n_{S}'(F)} \varphi(ka(\xi_++Y)a^{-1}k^{-1})dY \lvert \det a\rvert^s \eta_{E/F}(a) \delta_n(a) da \eta_{E/F}(k) dk \\
\nonumber & =\int_{K_n}\int_{A_n(F)} \int_{\n_{S}'(F)} \varphi(k(a\xi_+a^{-1}+Y)k^{-1})dY \lvert \det a\rvert^s \eta_{E/F}(a) \delta_n(a) \delta_{\n'_S}(a)^{-1} da \eta_{E/F}(k) dk
\end{align}
for every $s\in \C$ where we have set $\delta_{\n'_S}(a)=\lvert \det(\Ad(a)_{\mid \n'_{S}}) \rvert$. Set
$$\displaystyle f_\varphi(x_1,\ldots,x_n):=\int_{K_n} \int_{\n_{S}'(F)} \varphi(k(\xi_+(x_1,\ldots,x_n)+Y)k^{-1})dY\eta_{E/F}(k)dk$$
for every $(x_1,\ldots,x_n)\in F^n$ where by definition
$$\xi_+(x_1,\ldots,x_n):=(-1)^{n}\begin{pmatrix} 0 & x_1\tau^{-1} & 0 & \ldots & 0 \\ \vdots & \ddots & \ddots & \ddots & \vdots \\ \vdots & & \ddots & \ddots & 0 \\ \vdots & & & \ddots & x_n\tau^{-1} \\ 0 & \ldots & \ldots & \ldots & 0
\end{pmatrix}.$$
Then $f_\varphi$ is a Schwartz function on $F^n$ and by \eqref{eq 1 unipotent linear} we have (provided everything converges)
$$\displaystyle \lvert \tau\rvert_E^{-\dim(N'_S)/2}O_s(\xi_{+},\varphi) =\int_{A_n(F)} f_\varphi(a_1/a_2,\ldots,a_{n-1}/a_n,a_n) \lvert \det a \rvert^s \eta_{E/F}(a) \delta_n(a) \delta_{\n'_S}(a)^{-1} da.$$
A direct and painless computation shows that $\delta_n(a) \delta_{\n'_S}(a)^{-1}=\lvert a_2\ldots a_n\rvert^{-1}$ and therefore the above expression is equal to
$$\displaystyle \int_{(F^\times)^n} f_\varphi(a_1/a_2,\ldots,a_{n-1}/a_n,a_n) \eta_{E/F}(a_1\ldots a_n) \lvert a_1\rvert^s  \lvert a_2\ldots a_n\rvert^{s-1}d^\times a_1\ldots d^\times a_n.$$
By the change of variables $a_1\mapsto a_1a_2,\ldots, a_{n-1}\mapsto a_{n-1}a_n$ this can further be rewritten
\begin{align}\label{eq 2 unipotent linear}
\displaystyle \lvert \tau\rvert_E^{-\dim(N'_S)/2}O_s(\xi_{+},\varphi) =\int_{(F^\times)^n} f_\varphi(a_1,\ldots,a_n) \prod_{k=1}^n \eta_{E/F}(a_k)^k \lvert a_k\rvert^{ks-(k-1)}d^\times a_1\ldots d^\times a_n
\end{align}
which is clearly a convergent integral when $\Re(s)>1-1/n$ showing the convergence of $O_s(\xi_+,\varphi)$ in this range (indeed, as we easily see by the same computations where we replace the integrand by its absolute value).

By similar manipulations, and noticing that the isomorphism $F^n\simeq \n_S(F)/\n'_S(F)$, $(x_1,\ldots,x_n)\mapsto \xi_+(x_1,\ldots,x_n)$, sends the Haar measure on $F^n$ to $\lvert \tau\rvert_E^{(\dim(N_S)-\dim(N'_S))/2}$ times the Haar measure on $\n_S(F)/\n'_S(F)$, we have
\[\begin{aligned}
\displaystyle \lvert \tau\rvert_E^{-\dim(N'_S)/2}I^1(\varphi,s) & =\int_{A_n(F)}\int_{F^n} f_\varphi(x_1,\ldots,x_n) \psi'(\frac{a_2}{a_1}x_1+\ldots+\frac{a_n}{a_{n-1}}x_{n-1}+a^{-1}_nx_n)dx_1\ldots dx_n \\
 & \delta_{n+1}(a)^{-1}\delta_n(a) \lvert \det a\rvert^s \eta_{E/F}(a) da \\
 & =\int_{(F^\times)^n} \widehat{f_\varphi}(\frac{a_2}{a_1},\ldots,\frac{a_n}{a_{n-1}},a^{-1}_n) \lvert a_1\ldots a_n\rvert^{s-1} \eta_{E/F}(a_1\ldots a_n)d^\times a_1\ldots d^\times a_n
\end{aligned}\]
where we recall that $f\in \cS(F^n)\mapsto \widehat{f_\varphi}\in \cS(F^n)$ denotes the Fourier transform for $\psi'$ and the corresponding autodual measure. By the same change of variables as before, this becomes
\begin{align}\label{eq 3 unipotent linear}
\displaystyle \lvert \tau\rvert_E^{-\dim(N'_S)/2}I^1(\varphi,s) & =\int_{(F^\times)^n} \widehat{f_\varphi}(a^{-1}_1,\ldots,a^{-1}_n) \prod_{k=1}^n \eta_{E/F}(a_k)^k \lvert a_k\rvert^{ks-k}d^\times a_1\ldots d^\times a_n \\
\nonumber & =\int_{(F^\times)^n} \widehat{f_\varphi}(a_1,\ldots,a_n) \prod_{k=1}^n \eta_{E/F}(a_k)^k \lvert a_k\rvert^{k-ks}d^\times a_1\ldots d^\times a_n.
\end{align}
This shows that when $\Re(s)<1$, $I^1(\varphi,s)$ is defined by a convergent expression. Moreover the identity of the lemma follows from \eqref{eq 2 unipotent linear}, \eqref{eq 3 unipotent linear} and Tate's thesis. $\blacksquare$

\noindent\ul{Proof of Proposition} \ref{prop 1 unipotent linear}: By Lemma \ref{lem 1 unipotent linear} and Remark \ref{rem 1 unipotent linear}, we already know that $s\mapsto O_s(\xi_+,\varphi)$ has meromorphic continuation to $\C$ with no pole at $s=0$ and that
\begin{align}\label{eq 4 unipotent linear}
\displaystyle \lvert \tau\rvert_E^{-\dim(N_S)/2}\gamma O(\xi_+,\varphi)=\int_{G_n(F)/N_n(F)}\int_{\n_S(F)}\varphi(hYh^{-1}) \psi'(\langle \xi_{-},Y\rangle)dY \eta_{E/F}(h)dh.
\end{align}
Let $B_S$ be the image of $R_{E/F}B_{n+1}$ by $\nu$, $\fb_S$ its tangent space at the origin. We equip as before $\fb_S(F)$ with the transfer of the quotient measure on $\fb_{n+1}(E)/\fb_{n+1}(F)$ through $\nu$. Then, by Fourier inversion we have
$$\displaystyle \int_{\n_S(F)}f(Y)\psi'(\langle \xi_{-},Y\rangle)dY=\int_{\fb_S(F)}\cF f(\xi_{-}+Y) dY$$
for every $f\in \cS(\fs(F))$. Together with \eqref{eq 4 unipotent linear}, this gives
$$\displaystyle \lvert \tau\rvert_E^{-\dim(N_S)/2}\gamma O(\xi_+,\varphi)=\int_{G_n(F)/N_n(F)}\int_{\fb_S(F)}\cF \varphi(h(\xi_{-}+Y)h^{-1})dY \eta_{E/F}(h)dh.$$
Noticing that $\omega(h(\xi_{-}+Y)h^{-1})=\eta_{E/F}(h)\omega(\xi_{-})$ for every $h\in G_n(F)$ and $Y\in \fb_S(F)$, the above can be rewritten as
$$\displaystyle \lvert \tau\rvert_E^{-\dim(N_S)/2}\gamma \omega(\xi_{-})O(\xi_+,\varphi)=\int_{G_n(F)/N_n(F)}\int_{\fb_S(F)}\cF \varphi(h(\xi_{-}+Y)h^{-1})\omega(h(\xi_{-}+Y)h^{-1})dY dh.$$
Therefore, the proposition would follow if we can show:
\begin{num}
\item\label{eq 5 unipotent linear} For every $f\in L^1(\fs(F))$, we have the integration formula
$$\displaystyle \int_{\fs(F)}f(Y)dY=\lvert \tau\rvert_E^{\dim(N_S)/2}\int_{G_n(F)/N_n(F)}\int_{\fb_S(F)} f(h(\xi_{-}+Y)h^{-1}) dY dh.$$
\end{num}
As in the proof of Lemma \ref{lem 2 unipotent linear}, we are easily reduced to the same statement with $\fs$ and $\fb_S$ replaced by $\gl_{n+1}$ and $\fb_{n+1}$, $\xi_{-}$ replaced by $\xi_{-}'=(-1)^{n}\tau^{-1}\xi_{-}$ and without the factor $\lvert \tau\rvert_E^{\dim(N_S)/2}$. Then, given the definition of our Haar measures, it suffices to show that the morphism
$$\displaystyle G_n\times^{N_n} (\xi'_{-}+\fb_{n+1})\to \gl_{n+1}$$
$$\displaystyle (h,Y)\mapsto hYh^{-1}$$
is an open immersion over $\bZ$ where we have denoted by $G_n\times^{N_n} (\xi'_{-}+\fb_{n+1})$ the quotient of $G_n\times (\xi'_{-}+\fb_{n+1})$ by the free $N_n$-action given by $(h,Y)\cdot u=(hu,u^{-1}Yu)$. Set
$$\displaystyle \gl_{n+1}^{r}=\{Y\in \gl_{n+1}\mid \det(e_{n+1},e_{n+1}Y,\ldots,e_{n+1}Y^n)\neq 0 \}$$
Then $\gl_{n+1}^{r}$ is an open subscheme of $\gl_{n+1}$ and clearly  the above morphism factorizes through it. We claim that the induced map
$$\displaystyle G_n\times^{N_n} (\xi'_{-}+\fb_{n+1})\to \gl^{r}_{n+1}$$
is an isomorphism. Indeed, it suffices to show that for every commutative ring $R$ the map $G_n(R)\times^{N_n(R)} (\xi'_{-}+\fb_{n+1}(R))\to \gl^{r}_{n+1}(R)$ is a bijection. This in turn amounts to establishing the two following facts:
\begin{itemize}
\item For every $Y\in \gl^{r}_{n+1}(R)$, there exists $h\in G_n(R)$ such that $h^{-1}Yh\in \xi'_{-}+\fb_{n+1}(R)$ i.e. such that $e_{n+1}Y^kh\in e_{n+1-k}+\langle e_{n+2-k},\ldots,e_{n+1}\rangle_R$ for every $1\leqslant k\leqslant n$ (here we have denoted by $(e_1,\ldots,e_{n+1})$ the standard basis of $R^{n+1}$ and $\langle S\rangle_R$ stands for the $R$-submodule generated by $S$);
\item For every $Y\in \xi'_{-}+\fb_{n+1}(R)$ and $h\in G_n(R)$, if $hYh^{-1}\in \xi'_{-}+\fb_{n+1}(R)$ then $h\in N_n(R)$.
\end{itemize}
For the first point, denoting by $\overline{x}$ the image of $x\in R^{n+1}$ in $R^n=R^{n+1}/Re_{n+1}$, it suffices to choose for $h$ the unique element of $G_n(R)$ sending the basis $(\overline{e_{n+1}Y},\ldots,\overline{e_{n+1}Y^n})$ of $R^n$ to the basis $(\overline{e_n},\ldots,\overline{e_1})$. The second point follows by noticing that if $hYh^{-1}\in \xi'_{-}+\fb_{n+1}(R)$ then $h$ preserves the submodule $V_k=\langle e_{n+1-k},\ldots,e_{n+1}\rangle_R$ and acts trivially on the quotient $V_k/V_{k-1}$ for every $1\leqslant k\leqslant n$ and this readily implies that $h\in N_n(R)$. $\blacksquare$

\subsubsection{The linear case group version}\label{linear unipotent}

Let $f\in \cS(G'(F))$. We set
$$\displaystyle O_+(f)=O(\xi_+,\widetilde{f_\natural}).$$
We also recall that by the general construction of Section \ref{section Plancherel Whitt}, we associate to $f$ a function
$$\displaystyle W_f\in \cC^w(N'(F)\backslash G'(F)\times N'(F)\backslash G'(F),\psi_{N'}^{-1}\boxtimes \psi_{N'}).$$

\begin{lem}\label{lem 3 unipotent linear}
We have the equality
$$\displaystyle \gamma O_+(f)=\lvert \tau\rvert_E^{\dim(N_S)/2}\int_{N_1(F)\backslash H_1(F)}\int_{N_2(F)\backslash H_2(F)} W_f(h_1,h_2)\eta(h_2)dh_2dh_1$$
where $\gamma$ is the same constant as in Proposition \ref{prop 1 unipotent linear} and the right-hand side is absolutely convergent.
\end{lem}

\noindent\ul{Proof}: That the right hand side is absolutely convergent follows from Lemma \ref{lem continuity of beta} and Lemma \ref{lem 2 relative characters}. Unfolding the definitions, we have
\[\begin{aligned}
\displaystyle & \int_{N_1(F)\backslash H_1(F)}\int_{N_2(F)\backslash H_2(F)} W_f(h_1,h_2)\eta(h_2)dh_2dh_1= \\
 & \int_{N_n(F)\backslash G_n(F)} \int_{N_{n+1}(F)\backslash G_{n+1}(F)}\int_{N_n(E)\backslash G_n(E)} \int_{N_n(E)\times N_{n+1}(E)} f(g_n^{-1}u_nh_n,g_n^{-1}u_{n+1}h_{n+1}) \\
 & \psi_n(u_n)^{-1}\psi_{n+1}(u_{n+1})^{-1}du_ndu_{n+1}dg_n\eta_{n+1}(h_{n+1})dh_{n+1}\eta_n(h_n) dh_n
\end{aligned}\]
By the change of variable $u_{n+1}\mapsto u_nu_{n+1}$ and merging the integrals over $N_n(E)\backslash G_n(E)$ and $N_n(E)$ this becomes
\[\begin{aligned}
\displaystyle & \int_{N_n(F)\backslash G_n(F)} \int_{N_{n+1}(F)\backslash G_{n+1}(F)}\int_{G_n(E)} \int_{N_{n+1}(E)}f(g_nh_n,g_nuh_{n+1})\psi_{n+1}(u)^{-1}du \\
 & dg_n\eta_{n+1}(h_{n+1})dh_{n+1}\eta_n(h_n) dh_n
\end{aligned}\]
Note that the triple inner integral is absolutely convergent (this follows from the fact that $H_1\times R_{E/F}N_{n+1}\times^{N_{n+1}}G_{n+1}\to G'$, $(h_1,u,h_{n+1})\mapsto h_1(1,uh_{n+1})$, is a closed embedding). Therefore, by breaking the integral over $N_{n+1}(E)$ into one over $N_{n+1}(F)$ followed by one over $N_{n+1}(E)/N_{n+1}(F)$ and the change of variable $g_n\mapsto g_nh_n^{-1}$ we see that the above expression equals
\[\begin{aligned}
\displaystyle & \int_{N_n(F)\backslash G_n(F)} \int_{N_{n+1}(E)/N_{n+1}(F)}\int_{G_n(E)\times G_{n+1}(F)} f(g_n,g_nh_n^{-1}uh_{n+1})\eta'_{n+1}(h_n^{-1}uh_{n+1})dh_{n+1} \\
 & dg_n\psi_{n+1}(u)^{-1}du \eta_{E/F}(h_n)dh_n
\end{aligned}\]
By definition of $\widetilde{f}$ and of the Haar measure on $N_S(F)$, this last expression can be rewritten as
$$\displaystyle \int_{G_n(F)/N_n(F)} \int_{N_S(F)} \widetilde{f}(hvh^{-1})\psi_S(v)^{-1}dv\eta_{E/F}(h)dh$$
where we have denoted by $\psi_S$ the (unique) factorization of $\psi_{n+1}$ through $\nu: N_{n+1}(E)\to N_S(F)$. Finally, noting that the Cayley map $\fc$ induces a $G_n(F)$-equivariant isomorphism between $\n_S(F)$ and $N_S(F)$ preserving measures and sending the character $Y\mapsto \psi'(\langle \xi_{-},Y\rangle)$ to $\psi^{-1}_S$, we obtain that
\[\begin{aligned}
\displaystyle \int_{N_1(F)\backslash H_1(F)}\int_{N_2(F)\backslash H_2(F)} & W_f(h_1,h_2)\eta(h_2)dh_2dh_1= \\
 & \int_{G_n(F)/N_n(F)} \int_{\n_S(F)} \widetilde{f_\natural}(hYh^{-1})\psi'(\langle \xi_{-},Y\rangle)dY \eta_{E/F}(h)dh
\end{aligned}\]
The lemma now follows readily from Lemma \ref{lem 1 unipotent linear}. $\blacksquare$

Recall that $G_{\qs}$ stands for the unique quasi-split group of the form $G^{V_0}$ where $V_0\in \cV$.

\begin{prop}\label{prop 2 unipotent linear}
We have
$$\displaystyle \gamma O_+(f)=\lvert \tau\rvert_E^{-n(n-1)/4}\int_{\Temp(G_{\qs})/\stab} I_{BC(\pi)}(f) \frac{\lvert \gamma^*(0,\pi,\Ad,\psi')\rvert}{\lvert S_\pi\rvert}d\pi$$
where the right-hand side is absolutely convergent.
\end{prop}

\noindent\ul{Proof}: That the right-hand side is absolutely convergent and defines a continuous linear form on $\cS(G'(F))$ follows from Proposition \ref{prop 1 relative characters} together with Lemma \ref{lem 1 L parameters, LLC, basechange} and \eqref{basic estimates spectral measure}. On the other hand, by Lemma \ref{lem 3 unipotent linear}, Lemma \ref{lem 2 relative characters} and Lemma \ref{lem continuity of beta} the left hand side also defines a continuous linear form on $\cS(G'(F))$. Therefore, it suffices to establish the proposition when $f=f_n\otimes f_{n+1}$ where $f_k\in \cS(G_k(E))$ for $k=n,n+1$. Then, by Corollary \ref{cor lim spectrale} applied to the functions $f'_k=f_k \eta'_k$ and Theorem \ref{theo Plancherel}, we readily see that (see the proof of Proposition \ref{prop 1 local JR linear} for a similar argument)
\[\begin{aligned}
\displaystyle & \int_{N_2(F)\backslash H_2(F)} W_f(g,h_2)\eta(h_2)dh_2 \\
 & =\lvert \tau\rvert_E^{\frac{n(n-1)}{4}+\frac{n(n+1)}{4}} \int_{\Temp(G_{\qs})/\stab}\beta'(W_{f,BC(\pi)}(g,.)) \frac{\lvert \gamma^*(0,\pi,\Ad,\psi')\rvert}{\lvert S_\pi\rvert}d\pi
\end{aligned}\]
for every $g\in G'(F)$. From this, Lemma \ref{lem 3 unipotent linear}, Proposition \ref{prop 1 Plancherel Whitt} and Lemma \ref{lem 2 relative characters}, it follows that
\[\begin{aligned}
\displaystyle & \gamma O_+(f)=\lvert \tau\rvert_E^{\frac{n(n-1)}{4}+\frac{n(n+1)}{4}+\frac{\dim(N_S)}{2}}\int_{N_1(F)\backslash H_1(F)}\int_{\Temp(G_{\qs})/\stab}\beta'(W_{f,BC(\pi)}(h_2,.)) \frac{\lvert \gamma^*(0,\pi,\Ad,\psi')\rvert}{\lvert S_\pi\rvert}d\pi dh_1 \\
 & =\lvert \tau\rvert_E^{\frac{n(n-1)}{4}+\frac{n(n+1)}{4}+\frac{\dim(N_S)}{2}} \int_{\Temp(G_{\qs})/\stab}(\lambda\ctens\beta')(W_{f,BC(\pi)}) \frac{\lvert \gamma^*(0,\pi,\Ad,\psi')\rvert}{\lvert S_\pi\rvert}d\pi.
\end{aligned}\]
Moreover, by definition of $I_\Pi(f)$ and \eqref{eq 3 Planch Whitt}, we have
$$\displaystyle (\lambda \ctens \beta')(W_{f,\Pi})=\lvert \tau\rvert_E^{-\frac{n(n-1)}{2}-\frac{n(n+1)}{2}}I_\Pi(f)$$
for every $\Pi\in \Temp(G')$ and the identity of the proposition follows since
$$\displaystyle \frac{n(n-1)}{4}+\frac{n(n+1)}{4}+\frac{\dim(N_S)}{2}-\frac{n(n-1)}{2}-\frac{n(n+1)}{2}=-\frac{n(n-1)}{4}.$$ $\blacksquare$

\subsection{Proof of Theorems \ref{theo1 main result} and \ref{theo2 main result}}\label{Section proof of main theorems}

\subsubsection{Weak comparison of relative characters}

We recall the following result from \cite[Proposition 4.2.1]{Beu2}.

\begin{prop}\label{prop 1 weak comparison}
For every $V\in \cV$ and every $\pi\in \Temp_{H^V}(G^V)$, there exists $\kappa(\pi)\in \C^\times$ such that
$$\displaystyle J_\pi(f^V)=\kappa(\pi)I_{BC(\pi)}(f')$$
for all matching functions $f=(f^V)_V\in \cS(G(F))$ and $f'\in \cS(G'(F))$. Moreover, the function $\pi\in \Temp_{H^V}(G^V)\mapsto \kappa(\pi)$ is continuous.
\end{prop}

Here we remark that in {\it loc. cit.} only the $p$-adic case was considered. However, the proof extends readily to the Archimedean case the main points being that the globalization result \cite[Proposition 3.6.1]{Beu2} (which was borrowed from \cite{ILM}) still holds for Archimedean places. Looking closer into the proof we see that everything works equally well in the Archimedean situation except that the proof of \cite[Lemma 3.6.2]{Beu2} has to be slightly modified. Indeed, rather than appealing to M\oe{}glin-Tadic's classification of discrete series for classical $p$-adic groups, we should use the description by Harish-Chandra of discrete series for real reductive groups (and in particular of their infinitesimal character). Except from this modification,, the proof of \cite[Proposition 3.6.1]{Beu2} also works for Archimedean places.

\subsubsection{Comparison of local trace formulas}\label{Section two comparisons}

Let $f_1,f_2\in \cS(G(F))$ and $f_1',f_2'\in \cS(G'(F))$ such that $f_k$ and $f'_k$ match for $k=1,2$. Since the transfer factors have absolute value $1$, by Lemma \ref{lem 1 matching of orbits} we have
\[\begin{aligned}
\displaystyle \sum_{V\in \cV} \int_{H^V(F)\backslash G^V_{\rs}(F)/H^V(F)} O(\delta,f_1^V)\overline{O(\delta,f^V_2)}d\delta=\int_{H_1(F)\backslash G'_{rs}(F)/H_2(F)} O_\eta(\gamma,f_1') \overline{O_\eta(\gamma,f_2')}d\gamma.
\end{aligned}\]
By \eqref{eq 3 local JR unitary} and \eqref{eq 7 local JR linear}, this gives
\[\begin{aligned}
\displaystyle & \sum_{V\in \cV}\int_{\Temp_{\ind}(G^V)}J_\pi(f^V_1)\overline{J_\pi(f^V_2)} d\mu_{G^V}(\pi)= \\
 & \lvert \tau\rvert_E^{-n(n-1)/2}\int_{\Temp(G_{\qs})/\stab} I_{BC(\pi)}(f'_1)\overline{I_{BC(\pi)}(f'_2)} \frac{\lvert \gamma^*(0,\pi,\Ad,\psi')\rvert}{\lvert S_\pi\rvert}d\pi.
\end{aligned}\]
Using Proposition \ref{prop 1 weak comparison} and \eqref{eq 4ter Measures}, this can be rewritten
\begin{align}\label{eq 1 two comparisons}
\displaystyle & \int_{\Temp(G_{\qs})/\stab} \left(\sum_{V\in \cV}\sum_{\substack{\pi'\sim_{\stab}\pi \\ \pi'\in \Temp_{H^V}(G^V)}} \lvert \kappa(\pi')\rvert^2 \mu^*_{G^V}(\pi')\right) I_{BC(\pi)}(f'_1)\overline{I_{BC(\pi)}(f'_2)}d\pi= \\
\nonumber & \lvert \tau\rvert_E^{-n(n-1)/2}\int_{\Temp(G_{\qs})/\stab} I_{BC(\pi)}(f'_1)\overline{I_{BC(\pi)}(f'_2)} \frac{\lvert \gamma^*(0,\pi,\Ad,\psi')\rvert}{\lvert S_\pi\rvert}d\pi
\end{align}
where $\pi'\sim_{\stab}\pi$ means that $\pi'$ and $\pi$ share the same Langlands parameter. We remark that both sides of the above identity are absolutely convergent by \eqref{eq 6 local JR linear} and the convergence of \eqref{eq 2 local JR unitary}.

\subsubsection{Comparison of unipotent orbital integrals}

Let $f\in \cS(G(F))$ and $f'\in \cS(G'(F))$ be matching functions. By \eqref{eq 1 matching of functions} and Theorem \ref{theo transfer FT}, the functions $(\epsilon_V \cF \widetilde{f^V_\natural})_V$ and $\cF \widetilde{f'_\natural}$ match where we have set
\begin{align}\label{eq 1bis two comparisons}
\displaystyle \epsilon_V=\eta_{E/F}(-1)^{n(n+1)/2}\lambda_{E/F}(\psi')^{n(n+1)/2} \eta_{E/F}(\disc V)^n ,\;\;\; V\in \cV.
\end{align}
Therefore, by Lemma \ref{lem 2 matching of orbits}, we have
\[\begin{aligned}
\displaystyle \sum_{V\in \cV}\epsilon_V \int_{\fu_{\rs}^V(F)/U(V)(F)} O(X,\cF \widetilde{f_\natural^V})dX=\int_{\fs_{\rs}(F)/G_n(F)}\omega(Y) O_\eta(Y,\cF \widetilde{f'_\natural})dY.
\end{aligned}\]
By \eqref{eq 2 unipotent unitary} and Proposition \ref{prop 1 unipotent linear}, this implies
\begin{equation}\label{eq 1ter two comparisons}
\displaystyle \sum_{V\in \cV}\epsilon_VO(1,f^V)=\gamma\omega(\xi_{-})O(\xi_+,\widetilde{f'_\natural})=\gamma\omega(\xi_{-})O_+(f').
\end{equation}
Then, by \eqref{eq 1 unipotent unitary} and Proposition \ref{prop 2 unipotent linear}, this gives
\[\begin{aligned}
\displaystyle \sum_{V\in \cV}\epsilon_V \int_{\Temp_{\ind}(G^V)}J_\pi(f^V)d\mu_{G^V}(\pi)=\omega(\xi_{-})\lvert \tau\rvert_E^{-n(n-1)/4}\int_{\Temp(G_{\qs})/\stab} I_{BC(\pi)}(f') \frac{\lvert \gamma^*(0,\pi,\Ad,\psi')\rvert}{\lvert S_\pi\rvert}d\pi.
\end{aligned}\]
Finally, using Proposition \ref{prop 1 weak comparison} and \eqref{eq 4ter Measures}, this can be rewritten
\begin{align}\label{eq 2 two comparisons}
\displaystyle & \int_{\Temp(G_{\qs})/\stab} \left(\sum_{V\in \cV}\epsilon_V \sum_{\substack{\pi'\sim_{\stab}\pi \\ \pi'\in \Temp_{H^V}(G^V)}} \kappa(\pi') \mu^*_{G^V}(\pi')\right) I_{BC(\pi)}(f')d\pi= \\
\nonumber & \omega(\xi_{-})\lvert \tau\rvert_E^{-n(n-1)/4}\int_{\Temp(G_{\qs})/\stab} I_{BC(\pi)}(f') \frac{\lvert \gamma^*(0,\pi,\Ad,\psi')\rvert}{\lvert S_\pi\rvert}d\pi.
\end{align}
Note that both sides of the above identity are absolutely convergent by Proposition \ref{prop 2 unipotent linear} and the convergence of \eqref{eq 1 unipotent unitary}.

\subsubsection{How to separate each spectral contribution}

In order to finish the proofs of Theorems \ref{theo1 main result} and \ref{theo2 main result}, we need to separate each spectral contribution in \eqref{eq 1 two comparisons} and \eqref{eq 2 two comparisons} (i.e. to deduce from these identities similar equalities for each $\pi\in \Temp(G_{\qs})/\stab$). As usual, the argument ultimately rests upon the Stone-Weierstrass theorem by using some multipliers algebra. It is actually quite standard in the $p$-adic case (see e.g. \cite[Proof of Proposition 6.1.1]{SV}) but is slightly more subtle in the Archimedean case since using the center of the enveloping algebra as a substitute for the Bernstein center is not enough. Therefore, we explain carefully the proof here.

Consider the following general situation: $G$ is a connected reductive group over $F$, $\mu$ is a Borel measure on the set $\Irr_{\unit}(G)$ of unitary (or rather {\em unitarizable}) irreducible representations of $G(F)$ (equipped with the Fell topology) and we are given for $\mu$-almost every $\Pi\in \Irr_{\unit}(G)$ a continuous linear form $L_\Pi:\cS(G(F))\to \C$ which factorizes through the map $f\mapsto \Pi^\vee(f)$. In the applications we have in mind, we will take $G=G'$, for $\mu$ the pushforward of the measure $\frac{\lvert \gamma^*(0,\pi,\Ad,\psi')\rvert}{\lvert S_\pi\rvert}d\pi$ to $\Temp(G_{\qs})/\stab$ by $BC$ and for $L_\Pi$ a certain multiple of $I_\Pi$. We assume moreover that the following holds:
\begin{itemize}
\item For every $f\in \cS(G(F))$, the function $\Pi\in \Irr_{\unit}(G)\mapsto L_\Pi(f)$ is $\mu$-integrable and we have
$$\displaystyle \int_{\Irr_{\unit}(G)}L_\Pi(f) \mu(\Pi)=0.$$
\end{itemize}

\begin{lem}\label{lem separation}
Under the above assumptions, we have $L_\Pi=0$ for $\mu$-almost all $\Pi\in \Irr_{\unit}(G)$.
\end{lem}

\noindent\ul{Proof}:(the proof is inspired from \cite[Proof of Proposition 6.1.1]{SV}) Since $\cS(G(F))$ is separable (it is even of countable dimension in the $p$-adic case), up to multiplying $\mu$ by
$$\displaystyle \Pi\mapsto \sum_n \frac{\lvert L_\Pi(f_n)\rvert}{n^2\lVert L(f_n)\rVert_{L^1(\mu)}}$$
where $(f_n)_{n\geqslant 1}$ is a dense sequence in $\cS(G(F))$, $\lVert L(f)\rVert_{L^1(\mu)}:=\int_{\Irr_{\unit}(G)}\lvert L_\Pi(f)\rvert \mu(\Pi)$ and dividing the family $\Pi\mapsto L_\Pi$ by the same function, we may assume that $\mu$ is finite. Let $\cZ(G)$ denote the Bernstein center in the $p$-adic case or the center of the enveloping algebra $\cU(\g)$ in the Archimedean case. Then, we have a continuous map with finite fibers
$$\displaystyle p:\Irr_{\unit}(G)\to \widehat{\cZ(G)}$$
which associates to $\Pi$ its ``infinitesimal character'' $\chi_\Pi$. Let $X$ be the image of this map and $\overline{\mu}$ be the push-forward of $\mu$ to $X$. Then, by the disintegration of measures, there exists a measurable family $\chi\in X\mapsto \mu_\chi$ of finite measures on $\Irr_{\unit}(G)$ such that for $\overline{\mu}$-almost every $\chi\in X$, $\mu_\chi$ is supported on the finite set $p^{-1}(\chi)$ and moreover
$$\displaystyle \int_{\Irr_{\unit}(G)} f(\Pi) \mu(\Pi)=\int_X \sum_{\Pi\in p^{-1}(\chi)} f(\Pi) \mu_\chi(\Pi) \overline{\mu}(\chi)$$
for every $\mu$-integrable function $f$. By the hypothesis, we therefore have
\begin{align}\label{eq 1 separation}
\displaystyle \int_X \sum_{\Pi\in p^{-1}(\chi)} L_\Pi(f) \mu_\chi(\Pi) \overline{\mu}(\chi)=0
\end{align}
for every $f\in \cS(G(F))$. Let $K\subset G(F)$ be a maximal compact subgroup and $\cH(G(F))\subset \cS(G(F))$ be the corresponding ``Hecke algebra'' of $G(F)$, that is the space of smooth compactly supported and bi-$K$-finite functions on $G(F)$ (of course in the $p$-adic case we simply have $\cH(G(F))=\cS(G(F))$). We will now use the existence of a suitable ``algebra of multipliers'' $\cM(G)$ on $\cH(G(F))$ i.e. an algebra of endomorphisms $z$ of $\cH(G(F))$ for which there exists a function on $\widehat{\cZ(G)}$ (to be denoted by the same letter) such that $\Pi^\vee(zf)=z(\chi_\Pi)\Pi^\vee(f)$ for every $\Pi\in \Irr(G)$ and $f\in \cH(G(F))$. Assuming the existence of such a multiplier algebra, by applying \eqref{eq 1 separation} to $zf$ (and by the hypothesis made on $L_\Pi$) we get
\begin{align}\label{eq 2 separation}
\displaystyle \int_X z(\chi)\sum_{\Pi\in p^{-1}(\chi)} L_\Pi(f) \mu_\chi(\Pi) \overline{\mu}(\chi)=0
\end{align}
for every $f\in \cH(G(F))$ and $z\in \cM(G)$. For each finite set $S\subset \widehat{K}$, let $X_S$ be the set of infinitesimal characters of representations $\Pi\in \Irr_{\unit}(G)$ which are generated by their $\rho$-isotypic component for some $\rho\in S$. Assume that the following holds: 
\begin{num}
\item\label{eq 2bis separation} for every $z\in \cM(G)$ and $S\subset\widehat{K}$, the function $\chi\in X_S\mapsto z(\chi)$ converges to zero at infinity and the algebra of functions $\{z_{\mid X_S}\mid z\in \cM(G) \}$ is stable by complex conjugation, separates points and does not vanish identically anywhere.
\end{num}
Then, by the Stone-Weierstrass theorem applied to the one-point compactification of $X_S$, this algebra is dense in the Banach space $C_0(X_S)$ of continuous functions on $X_S$ tending to zero at infinity. Since for each $f\in \cH(G(F))$ the function $\chi\mapsto \sum_{\Pi\in p^{-1}(\chi)} L_\Pi(f) \mu_\chi(\Pi)$ is supported on $X_S$ for some $S$ (again by the hypothesis made on $L_\Pi$), this together with \eqref{eq 2 separation} implies that for every $f\in \cH(G(F))$
$$\displaystyle \sum_{\Pi\in p^{-1}(\chi)} L_\Pi(f) \mu_\chi(\Pi)=0$$
for $\overline{\mu}$-almost every $\chi\in X$. As $\cH(G(F))$ contains a sequence which is dense in $\cS(G(F))$, by continuity of the linear forms $L_\Pi$, the above equality actually holds for $\overline{\mu}$-almost every $\chi\in X$ and every $f\in \cS(G(F))$. Since the map
$$\displaystyle \cS(G(F))\to \bigoplus_{\Pi\in p^{-1}(\chi)} \End^\infty(\Pi^\vee),\; f\mapsto (\Pi^\vee(f))_\Pi$$
is surjective (indeed, the image is a closed and $G(F)\times G(F)$-invariant subspaces, the $G(F)\times G(F)$-representations $\End^\infty(\Pi^\vee)$ are irreducible and pairwise non-equivalent and for each $\Pi$ the map $f\in \cS(G(F))\mapsto \Pi^\vee(f)\in \End^\infty(\Pi^\vee)$ is surjective), by the hypothesis made on linear forms $L_\Pi$ this implies $\mu_\chi(\Pi)L_\Pi=0$ for $\overline{\mu}$-almost every $\chi\in X$ and every $\Pi\in p^{-1}(\chi)$ i.e. $L_\Pi=0$ for $\mu$-almost every $\Pi\in \Irr_{\unit}(G)$.

Thus, it only remains to show the existence of an algebra of multipliers satisfying \eqref{eq 2bis separation}. Notice that if $f\mapsto zf$ is a multiplier then so is $f\mapsto z^*f:=(zf^*)^*$ (where we recall that $f^*(g)=\overline{f(g^{-1})}$) and that $z^*(\chi_\Pi)=\overline{z(\chi_\Pi)}$ whenever $\Pi$ is a unitary representation. Therefore, we only need to find an algebra of multipliers tending to $0$ at infinity on $X_S$ for every $S\subset \widehat{K}$ and separating points (including infinity). In the $p$-adic case, we just take $\cM(G)=\cZ(G)$ the Bernstein center. Clearly, it separates points and does not vanish identically anywhere on $\widehat{\cZ(G)}$. Since for each $S\subset \widehat{K}$, $X_S$ is compact (\cite[Theorem 2.5]{Tad}) this algebra has all the desired properties. In the Archimedean case, we will use a certain subalgebra of Arthur's algebra of multipliers (\cite{Art2}, \cite{Del}). To be more precise, we need to introduce more notation. Let $T\subset G$ be a maximal torus. Harish-Chandra's isomorphism gives an identification
$$\displaystyle \widehat{\cZ(G)}=\mathfrak{t}(\C)^*/W$$
where $W=W(G_\C,T_\C)$. Let $\mathfrak{t}_{\bR}\subset \mathfrak{t}(\C)$ be the $\bR$-points of the split form of $\mathfrak{t}$ and for $\varphi\in C_c^\infty(\mathfrak{t}_{\bR})$ write
$$\displaystyle \widehat{\varphi}(\lambda)=\int_{\mathfrak{t}_{\bR}}\varphi(X)e^{\lambda(X)}dX,\;\; \lambda\in \mathfrak{t}(\C)^*,$$
for its Fourier transform (where $dX$ is a fixed Haar measure on $\mathfrak{t}_{\bR}$). Then, $W$ preserves $\mathfrak{t}_{\bR}$ and by \cite[Theorem 4.2]{Art2}, \cite[Theorem 3]{Del} there exists an algebra of multipliers $\cM(G)$ whose associated set of functions on $\widehat{\cZ(G)}$ is precisely $\widehat{C_c^\infty(\mathfrak{t}_{\bR})^W}$. This algebra has all the desired properties the only non-trivial point being that $z$ tends to $0$ at infinity on $X_S$ for any $z\in \cM(G)$ and $S\subset \widehat{K}$ but this follows from the fact that $X_S$ has compact image in $i\mathfrak{t}_{\bR}^*\backslash \mathfrak{t}(\C)^*/W$ (see \cite[Theorem 5.2]{BW} and  \cite[Corollary 7.7.3]{Wall}) and usual properties of the Fourier transform. $\blacksquare$

\subsubsection{End of the proof}

First we show that \eqref{eq 1 two comparisons} holds for every $f_1',f_2'\in \cS(G'(F))$, both sides being absolutely convergent. In the $p$-adic case, this already follows from the computations of \S \ref{Section two comparisons} since every function in $\cS(G'(F))$ a admits a transfer to $\cS(G(F))$ (Theorem \ref{theo existence transfer} (i)). In the Archimedean case, by Theorem \ref{theo existence transfer} (ii), we know at least that it holds for $f_1'$, $f_2'$ in a certain dense subspace $\cS(G'(F))_{\trans}$ of $\cS(G'(F))$. Therefore, it suffices to show that both side of \eqref{eq 1 two comparisons} are absolutely convergent for every $f_1',f_2'\in \cS(G'(F))$ and that they define continuous sesquilinear forms on $\cS(G'(F))$. For the right-hand side, this follows readily from Proposition \ref{prop 1 relative characters} together with Lemma \ref{lem 1 L parameters, LLC, basechange} and \eqref{basic estimates spectral measure}. By Cauchy-Schwartz inequality, it then suffices to prove that the left-hand side is always less or equal to the right-hand side whenever $f_1'=f_2'$. That it is indeed the case is a consequence of Fatou's lemma together with the fact that the linear forms $I_\Pi$, $\Pi\in \Temp(G')$, are continuous.

Thus, we can now apply Lemma \ref{lem separation} to \eqref{eq 1 two comparisons}, giving us the identity
\[\begin{aligned}
\displaystyle & \left(\sum_{V\in \cV}\sum_{\substack{\pi'\sim_{\stab}\pi \\ \pi'\in \Temp_{H^V}(G^V)}} \lvert \kappa(\pi')\rvert^2 \mu^*_{G^V}(\pi')\right) I_{BC(\pi)}(f'_1)\overline{I_{BC(\pi)}(f'_2)}= \\
 & \lvert \tau\rvert_E^{-n(n-1)/2}I_{BC(\pi)}(f'_1)\overline{I_{BC(\pi)}(f'_2)} \frac{\lvert \gamma^*(0,\pi,\Ad,\psi')\rvert}{\lvert S_\pi\rvert}
\end{aligned}\]
for almost every $\pi \in \Temp(G_{\qs})/\stab$ and every $f_1',f_2'\in \cS(G'(F))$. By the local Gan-Gross-Prasad conjecture (\cite[Theorem 12.4.1]{Beu1}), for every $\pi \in \Temp(G_{\qs})/\stab$ there exists exactly one $V\in \cV$ and one representation $\pi'\in \Temp_{H^V}(G^V)$ such that $\pi'\sim_{\stab}\pi$. Therefore, as the linear form $I_{BC(\pi)}$ is non-zero (this can for example be deduced from Proposition \ref{prop 1 weak comparison} since $J_\pi\neq 0 \Leftrightarrow \pi\in \Temp_{H^V}(G^V)$), the above identity can be rewritten
\begin{align}\label{eq 1 end of proof}
\displaystyle \lvert \kappa(\pi)\rvert^2 \mu^*_{G^V}(\pi)=\lvert \tau\rvert_E^{-n(n-1)/2}\frac{\lvert \gamma^*(0,\pi,\Ad,\psi')\rvert}{\lvert S_\pi\rvert}
\end{align}
for almost every $\displaystyle \pi\in \bigsqcup_{V\in \cV} \Temp_{H^V}(G^V)$. As a first consequence, multiplying both sides by $\mu^*_{G^V}(\pi)$ and using \eqref{eq 2 Planch} together with Lemma \ref{lem 1 L parameters, LLC, basechange}, we see that in the Archimedean case there exists $k>0$ such that
\begin{align}\label{eq 2 end of proof}
\displaystyle \lvert \kappa(\pi)\rvert \mu^*_{G^V}(\pi)\ll N(\pi)^k
\end{align}
for almost every $\displaystyle \pi\in \bigsqcup_{V\in \cV} \Temp_{H^V}(G^V)$.

We now show that \eqref{eq 2 two comparisons} holds for every $f'\in \cS(G'(F))$. Once again, in the $p$-adic case there is nothing to say and in the Archimedean case it suffices to show that both sides are absolutely convergent and define continuous linear forms on $\cS(G'(F))$. But this follows readily from Proposition \ref{prop 1 relative characters} together with \eqref{eq 2 end of proof} (for the left-hand side), Lemma \ref{lem 1 L parameters, LLC, basechange} (for the right-hand side) and \eqref{basic estimates spectral measure}.

Thus, we can apply Lemma \ref{lem separation} to \eqref{eq 2 two comparisons} and using the non-vanishing of $I_{BC(\pi)}$ this gives
$$\displaystyle \sum_{V\in \cV}\epsilon_V \sum_{\substack{\pi'\sim_{\stab}\pi \\ \pi'\in \Temp_{H^V}(G^V)}} \kappa(\pi') \mu^*_{G^V}(\pi')=\omega(\xi_{-})\lvert \tau\rvert_E^{-n(n-1)/4}\frac{\lvert \gamma^*(0,\pi,\Ad,\psi')\rvert}{\lvert S_\pi\rvert}$$
for almost every $\pi \in \Temp(G_{\qs})/\stab$. Applying again the local Gan-Gross-Prasad conjecture, this can be rewritten as
\begin{align}\label{eq 3 end of proof}
\displaystyle \epsilon_V \kappa(\pi) \mu^*_{G^V}(\pi)=\omega(\xi_{-})\lvert \tau\rvert_E^{-n(n-1)/4}\frac{\lvert \gamma^*(0,\pi,\Ad,\psi')\rvert}{\lvert S_\pi\rvert}
\end{align}
for every $V\in \cV$ and almost every $\pi\in \Temp_{H^V}(G^V)$.

Squaring the module of \eqref{eq 3 end of proof} and dividing it by \eqref{eq 1 end of proof}, we obtain (as $\mu^*_{G^V}(\pi)\geqslant 0$ for every $\pi$)
$$\displaystyle \mu^*_{G^V}(\pi)=\frac{\lvert \gamma^*(0,\pi,\Ad,\psi')\rvert}{\lvert S_\pi\rvert}$$
for every $V\in \cV$ and almost every $\pi\in \Temp_{H^V}(G^V)$ thus proving Theorem \ref{theo2 main result}. Moreover, dividing \eqref{eq 3 end of proof} by the above identity gives
$$\displaystyle \epsilon_V \kappa(\pi)=\omega(\xi_{-})\lvert \tau\rvert_E^{-n(n-1)/4}$$
By \eqref{eq 0 unipotent linear} and \eqref{eq 1bis two comparisons} we easily check that $\kappa_V=\epsilon_V \omega(\xi_{-})^{-1}\lvert \tau\rvert_E^{n(n-1)/4}$ (where $\kappa_V$ is the constant defined in the statement of Theorem \ref{theo1 main result}) and therefore we have
$$\displaystyle \kappa(\pi)=\kappa_V^{-1}$$
for every $V\in \cV$ and almost every $\pi\in \Temp_{H^V}(G^V)$. By continuity of $\pi\mapsto \kappa(\pi)$ (Proposition \ref{prop 1 weak comparison}) this last equality is true for every $\pi\in \Temp_{H^V}(G^V)$ and this proves Theorem \ref{theo1 main result}.

\appendix

\section{Proof of Proposition \ref{prop 1 Planch}}\label{Appendice Schwartz}

First, we recall the following elementary lemma (see \cite[Proposition]{CHH}).
\begin{lem}[Sobolev lemma for compact groups]\label{Sobolev}
Let $H$ be a Hilbert space, $K$ be a compact group and $C(K,H)$ be the space of continuous functions from $K$ to $H$. Consider $C(K,H)$ as a representation of $K$ through the right regular action and for $\rho\in \widehat{K}$ denote by $C(K,H)[\rho]$ its $\rho$-isotypic component. Then, for every $\rho\in \widehat{K}$ and $\varphi\in C(K,H)[\rho]$ we have
$$\displaystyle \sup_{k\in K}\lVert \varphi(k)\rVert\leqslant \dim(\rho)\left(\int_K \lVert \varphi(k)\rVert^2 dk\right)^{1/2}.$$
\end{lem}

Now we proceed to the proof of Proposition \ref{prop 1 Planch}.

\vspace{2mm}

Let $f\in \cS(G(F))$. The function $\pi\in \Temp_{\ind}(G)\mapsto f_\pi\in \cC^w(G(F))$ is smooth by \cite[Lemma 2.3.1(ii)]{Beu1} together with \cite[\S 3]{Art1} in the Archimedean case and \cite[Proposition VII.1.3]{Wald1} in the $p$-adic case. Moreover, in the $p$-adic case the function $\pi\mapsto f_\pi$ is compactly supported by \cite[Th\'eor\`eme VIII.1.2]{Wald1}. It remains to check that $f$ satisfies condition \eqref{eq defn Schwartz tempered} in the Archimedean case and condition \eqref{eq2 defn Schwartz tempered} in the $p$-adic case. We will concentrate on the Archimedean case, the $p$-adic case being similar and actually easier to handle. By definition of the topology on $\cC^w(G(F))$, it suffices to establish the following: for every Levi subgroup $M\subset G$, $D\in \Sym^\bullet(\cA_{M,\C}^*)$, $u,v\in \cU(\g)$ and $k\geqslant 0$ we have
\begin{align}\label{eq 0 Planch}
\displaystyle N(\pi)^k\left\lvert D(\lambda\mapsto (R(u)L(v)f_{\pi_\lambda})(g))_{\lambda=0}\right\rvert\ll \Xi^G(g)\sigma(g)^{\deg(D)}
\end{align}
for $\tau\in \Pi_2(M)$ and $g\in G(F)$ where $\deg(D)$ stands for the degree of $D$ and we have set $\pi_\lambda=i_M^G(\tau_\lambda)$ and $\pi=\pi_0$ for $\lambda\in i\cA_M^*$.

As $R(u)L(v)f_\pi=(R(u)L(v)f)_\pi$, up to replacing $f$ by $R(u)L(v)f$ we only need to establish \eqref{eq 0 Planch} when $u=v=1$. Assume proved the slightly weaker inequality (for any $D\in \Sym^\bullet(\cA_{M,\C}^*)$)
\begin{align}\label{eq 0ter Planch}
\displaystyle \left\lvert D(\lambda\mapsto f_{\pi_\lambda}(g))_{\lambda=0}\right\rvert\ll \Xi^G(g)\sigma(g)^{\deg(D)},\;\;\; \tau\in \Pi_2(M), g\in G(F).
\end{align}
Then, we will show that \eqref{eq 0 Planch} holds for any $k\geqslant 0$ (and $u=v=1$). We do this by induction on $\deg(D)$. Let $z\in \cZ(\g)$ be such that \eqref{eq1 representations} is satisfied. Then we may as well assume that $N(\pi)=\chi_{\pi}(z)$. Since $(z^kf)_\pi=\chi_\pi(z)^kf_\pi$ for every $\pi\in \Temp_{\ind}(G)$ the result in degree $0$ just follows from replacing $f$ by $z^k f$ in \eqref{eq 0ter Planch}. In the general case the difference between
\begin{align}\label{eq 0quad Planch}
\displaystyle D(\lambda\mapsto (z^kf)_{\pi_\lambda}(g))_{\lambda=0}=D(\lambda\mapsto \chi_{\pi_\lambda}(z)^kf_{\pi_\lambda}(g))_{\lambda=0}
\end{align}
and
$$\displaystyle \chi_\pi(z)^k D(\lambda\mapsto f_{\pi_\lambda}(g))_{\lambda=0}$$
can be written as a finite sum
$$\displaystyle \sum_{i=1}^nD_i(\lambda\mapsto \chi_{\pi_\lambda}(z)^k)_{\lambda=0} D_i'(\lambda\mapsto f_{\pi_\lambda}(g))_{\lambda=0}$$
where $D_i,D'_i\in \Sym^\bullet(\cA_{M,\C})$ for all $1\leqslant i\leqslant n$ are of degree strictly less than $D$. Since the terms $D_i(\lambda\mapsto \chi_{\pi_\lambda}(z)^k)_{\lambda=0}$ are all essentially bounded by a power of $N(\pi)$, the above sum can be controlled by the induction hypothesis whereas \eqref{eq 0quad Planch} is controlled by \eqref{eq 0ter Planch} applied to $z^k f$. This shows \eqref{eq 0 Planch} assuming \eqref{eq 0ter Planch}.

Fix $D\in \Sym^\bullet(\cA_{M,\C}^*)$, a parabolic subgroup $P$ with Levi component $M$, a maximal compact subgroup $K$ of $G(F)$ in good position relative to $M$ and set $K_M=K\cap M(F)$. Then, by restriction to $K$ we can identify $\pi_\lambda=i_P^G(\tau_\lambda)$ with $\pi_K=i_{K_M}^K(\tau_{\mid K_M})$ as a $K$-representation for every $\tau\in \Pi_2(M)$ and $\lambda\in i\cA_M^*$. Choosing an invariant scalar product $(.,.)$ on $\tau$, we endow $\pi_K$ with the $K$-invariant scalar product
$$\displaystyle (e,e')=\int_K (e(k),e'(k))dk.$$
Finally, choose for any $\rho\in \widehat{K}$ an orthonormal basis $\cB_{\tau}(\rho)$ of the $\rho$-isotypic component $\pi_K[\rho]$ of $\pi_K$. Then, we have
\begin{align}\label{eq 0cinq Planch}
\displaystyle f_{\pi_\lambda}(g)=\sum_{\rho\in \widehat{K}}\sum_{e\in \cB_{\tau}(\rho)} (\pi_\lambda(g)\pi_\lambda(f^\vee)e,e)
\end{align}
for all $\tau\in \Pi_2(M)$, $\lambda\in i\cA_M^*$ and $g\in G(F)$. Assume now that we can show the existence of $r>0$ such that
\begin{align}\label{eq 0six Planch}
\displaystyle \left\lvert D\left(\lambda\mapsto (\pi_\lambda(g)\pi_\lambda(f^\vee)e,e)\right)_{\lambda=0}\right\rvert \ll N(\rho)^r \Xi^G(g)\sigma(g)^{\deg(D)}
\end{align}
for all $\tau\in \Pi_2(M)$, $g\in G(F)$, $\rho\in \widehat{K}$ and $e\in \cB_{\tau}(\rho)$. Let $z_K\in \cZ(\mathfrak{k})$ be such that \eqref{eq1 representations} is satisfied for $G(F)=K$. Then we may as well assume that $N(\rho)=\rho(z_K)$ for all $\rho\in \widehat{K}$ and hence up to replacing $f$ by $L(z_K)^{k+r}f$ in \eqref{eq 0six Planch} we obtain the same inequality with $N(\rho)^r$ replaced by $N(\rho)^{-k}$. Therefore, by \eqref{eq 0cinq Planch}, we would get for any $k>0$ an inequality
$$\displaystyle \left\lvert D(\lambda\mapsto f_{\pi_\lambda}(g))_{\lambda=0}\right\rvert\ll \Xi^G(g)\sigma(g)^{\deg(D)}\sum_{\rho\in \widehat{K}} \frac{\dim(\pi_K[\rho])}{N(\rho)^k},\;\;\; \tau\in \Pi_2(M),\; g\in G(F).$$
As for $k$ large enough the sum $\sum_{\rho\in \widehat{K}} \frac{\dim(\pi_K[\rho])}{N(\rho)^k}$ converges and is bounded independently of $\pi$ (see e.g. \cite[(2.2.2)]{Beu1}). For such a $k$ the above inequality implies \eqref{eq 0ter Planch}.

Thus, it only remains to establish \eqref{eq 0six Planch}. Actually, it suffices to show the existence of $r>0$ such that
\begin{align}\label{eq 0sept Planch}
\displaystyle \left\lvert D\left(\lambda\mapsto (\pi_\lambda(g)e,e)\right)_{\lambda=0}\right\rvert \ll N(\rho)^r \Xi^G(g)\sigma(g)^{\deg(D)}
\end{align}
for all $\tau\in \Pi_2(M)$, $g\in G(F)$, $\rho\in \widehat{K}$ and $e\in \cB_{\tau}(\rho)$. Indeed, if this is the case we would get
\[\begin{aligned}
\displaystyle & \left\lvert D\left(\lambda\mapsto (\pi_\lambda(g)\pi_\lambda(f^\vee)e,e)\right)_{\lambda=0}\right\rvert =\left\lvert D\left(\lambda\mapsto\int_G f^\vee(\gamma)  (\pi_\lambda(g\gamma)e,e)d\gamma\right)_{\lambda=0}\right\rvert \\
 & =\left\lvert \int_G f^\vee(\gamma) D\left(\lambda\mapsto (\pi_\lambda(g\gamma)e,e)\right)_{\lambda=0}d\gamma\right\rvert \ll N(\rho)^r\int_G \lvert f^\vee(\gamma)\rvert \Xi^G(g\gamma)\sigma(g\gamma)^{\deg(D)}d\gamma \\
 & \ll N(\rho)^r \sigma(g)^{\deg(D)} \int_G \sup_{k\in K}\lvert f^\vee(k\gamma)\rvert \int_K\Xi^G(gk\gamma)dk\sigma(\gamma)^{\deg(D)}d\gamma \\
 & =N(\rho)^r \Xi^G(g)\sigma(g)^{\deg(D)}\int_G \sup_{k\in K}\lvert f^\vee(k\gamma)\rvert \Xi^G(\gamma)\sigma(\gamma)^{\deg(D)}d\gamma
\end{aligned}\]
for all $\tau\in \Pi_2(M)$, $g\in G(F)$, $\rho\in \widehat{K}$ and $e\in \cB_{\tau}(\rho)$, where the differentiation under the integral sign is justified by the absolute convergence of the resulting expression and in the last line we have used the well-known `doubling formula' $\int_K \Xi^G(gk\gamma)dk=\Xi^G(g)\Xi^G(\gamma)$ (see \cite[Proposition 16.(iii) p.329]{Var}).

Fix a decomposition $\gamma=m_P(\gamma)u_P(\gamma)k_P(\gamma)$ for every $\gamma\in G(F)$ where $m_P(\gamma)\in M(F)$, $u_P(\gamma)\in U_P(F)$ and $k_P(\gamma)\in K$ and set $H_P(\gamma):=H_M(m_P(\gamma))$. Then to prove \eqref{eq 0sept Planch} we first note that
$$\displaystyle (\pi_\lambda(g)e,e)=\int_K \delta_P(m_P(kg))^{1/2} e^{\langle \lambda, H_P(kg)\rangle} (\tau(m_P(kg))e(k_P(kg)),e(k)) dk$$
so that (once again differentiation under the integral is easily justified)
$$\displaystyle D\left(\lambda\mapsto (\pi_\lambda(g)e,e)\right)_{\lambda=0}=\int_K \delta_P(m_P(kg))^{1/2} D(H_P(kg)) (\tau(m_P(kg))e(k_P(kg)),e(k)) dk$$
for all $\tau\in \Pi_2(M)$, $e\in \pi_K$ and $g\in G(F)$ where when we write $D(H_P(kg))$ we consider $D$ as a polynomial function on $\cA_M$. Clearly $\lvert D(H_P(kg))\rvert \ll \sigma(g)^{\deg(D)}$ for all $g\in G(F)$ and $k\in K$ and therefore
\begin{align}\label{eq 0huit Planch}
\displaystyle \left\lvert D\left(\lambda\mapsto (\pi_\lambda(g)e,e)\right)_{\lambda=0}\right\rvert \ll \sigma(g)^{\deg(D)} \int_K \delta_P(m_P(kg))^{1/2} \lvert(\tau(m_P(kg))e(k_P(kg)),e(k))\rvert dk
\end{align}
for all $\tau\in \Pi_2(M)$, $e\in \pi_K$ and $g\in G(F)$. By \cite[Theorem 2]{CHH}, we have
$$\displaystyle \lvert (\tau(m)v,v')\rvert\ll \dim(\tau(K_M)v)^{1/2} \dim(\tau(K_M)v')^{1/2} \Xi^M(m) \lVert v\rVert \lVert v' \rVert$$
for all $\tau\in \Pi_2(M)$, all $v,v'\in \tau$ and all $m\in M(F)$. Note that $\dim(\tau(K_M)e(k))\leqslant \dim(\pi(K)e)$ for all $\tau\in \Pi_2(M)$ and $e\in \pi_K$ and that by a new application of \cite[Theorem 2]{CHH} $\dim(\pi(K)e)\leqslant \dim(\rho)^2$ if $e\in \pi_K[\rho]$. Combining this with \eqref{eq 0huit Planch}, we obtain
\[\begin{aligned}
\displaystyle \left\lvert D\left(\lambda\mapsto (\pi_\lambda(g)e,e)\right)_{\lambda=0}\right\rvert & \ll \sigma(g)^{\deg(D)} \dim(\rho)^2 \int_K \delta_P(m_P(kg))^{1/2} \Xi^M(m_P(kg))\lVert e(k_P(kg))\rVert \lVert e(k)\rVert dk \\
 & \leqslant \sigma(g)^{\deg(D)} \dim(\rho)^2 \sup_{k\in K} \lVert e(k)\rVert^2 \int_K \delta_P(m_P(kg))^{1/2} \Xi^M(m_P(kg)) dk \\
 & \leqslant \sigma(g)^{\deg(D)} \dim(\rho)^3 \Xi^G(g)\lVert e\rVert
\end{aligned}\]
for all $\tau\in \Pi_2(M)$, $\rho\in \widehat{K}$, $e\in \pi_K[\rho]$ and $g\in G(F)$ where in the last inequality we have used \cite[Proposition 16(iv) p.329]{Var} and Lemma \ref{Sobolev}. Since $\lVert e\rVert=1$ and there exists $n\geqslant 1$ such that $\dim(\rho)\leqslant N(\rho)^n$ (\cite[p.291]{Wall}), this gives \eqref{eq 0sept Planch} and ends the proof of the Proposition \ref{prop 1 Planch}. $\blacksquare$

\flushright Rapha\"el Beuzart-Plessis \\
Aix Marseille University \\
CNRS, Centrale Marseille, I2M\\
Marseille, France\\ 
email: rbeuzart@gmail.com

\end{document}